\input amstex
\documentstyle{amsppt}

\centerline{\bf ACTION DE DERIVATIONS IRREDUCTIBLES SUR LES }

\vskip 3mm

\centerline{\bf  ALGEBRES QUASI-REGULIERES D'HILBERT}

\vskip 5mm

\centerline{ par Abderaouf Mourtada}

\vskip 5mm

\centerline{A la mémoire de Jean Martinet}

\document

\hcorrection{2cm}
\vcorrection{2cm}

\def \NN{\Bbb N}
\def \ZZ{\Bbb Z}
\def \RR{\Bbb R}
\def \CC{\Bbb C}

\vskip 1cm

\heading Université de Bourgogne, I.M.B.\\
        U.M.R. 5584 du C.N.R.S., U.F.R. des Sciences et Techniques\\
         9, avenue Alain Savary, B.P. 47 870, 21078 Dijon Cedex.\\
      \rm E-mail: mourtada\@u-bourgogne.fr  \endheading

\vskip 1cm

\noindent {\bf Mots-clés}: 16ème problème d'Hilbert, cycle limite, cyclicité, polycycle, singularité de champ de vecteurs, application de Dulac, application de retour de Poincaré, structure asymptotique quasi-analytique, action différentielle, idéal différentiel, faisceau différentiel, projection intégrale, multiplicité algébrique.\par

\vskip 5mm

\noindent {\bf Classification AMS}: 34C07 (Primary) 34C08, 37G15 (Secondary).\par

\vskip 1cm

\noindent{\bf Abstract.} {\sl We study the action of irreducible derivations $\chi$ on some Hilbert's quasi-regular algebras $QR{\Cal H}$ of germs at 0, of real analytic functions on $(U,0)$, where $U$ is some semi-algebraic open set. We show that these algebras are $\chi$-finite or locally $\chi$-finite: the degree of the projection $\pi_\chi$ restricted to fibers of $QR{\Cal H}$, is finite and the differential ideals are noetherian or locally noetherian. Moreover, these algebras satisfy to the double inclusion: for every germ $f$, there exist an algebraic multiplicity $ma_\chi(f)$ such that, modulo some algebraic factor which vanishs only on the germ of the boundary $(\partial U,0)$, and depends only on $ma$, the differential ideal of $f$ coincides with the saturation of his transverse ideal. In last, we give an application to the generalised Hilbert's 16th problem about limit cycles: there is no accumulation of limit cycles on hyperbolic polycycles in compact analytic families of vector fields on the sphere $S^2$. This is a highly non trivial result, it includes the case of a polycycle that is an accumulation of cycles.}\par

\vskip 1cm

\centerline{\bf Introduction} 

\vskip 5mm

 Le {\bf problème de Dulac} ([E], [I]), dit que tout champ de vecteurs analytique sur la sphère réelle $S^2$ a un nombre fini de cycles limites, autrement dit il n'y a pas accumulation de cycles limites sur les polycycles. Le {\bf 16ème problème d'Hilbert algébrique} ([Hi]), demande d'établir une majoration de ce nombre en fonction du degré pour les champs de vecteurs algébriques du plan réel. Plus généralement ([Ro1]), il s'agit de montrer qu'il n'y a pas accumulation de cycles limites sur les ensembles limites périodiques, dans les familles compactes de champs de vecteurs analytiques sur $S^2$: c'est le {\bf problème d'Hilbert analytique}. Soit $\Gamma_k$ un polycycle {\bf hyperbolique et monodromique réel}, à $k$ singularités, tangent \`a un champ de vecteurs analytique r\'eel $X_0$ d\'efini sur un voisinage $U_0$ de $\Gamma_k$. On suppose, uniquement pour simplifier la présentation, que le rapport des valeurs propres en chaque singularité de $\Gamma_k$ est égal à $-1$. Dans la section IVC, on montre le théorème fondamental\par
 
\vskip 3mm

\proclaim{Théorème 0} Soit $X_\nu$ un déploiement analytique de $X_0$ à $q$ paramètres. Alors il existe des entiers $N$ et $L$ et des voisinages $\Gamma_k\subset U\subset U_0$ et $V\in (\RR^q,0)$ tels que

\roster

\vskip 2mm

\item "$(i)$"  pour tout $\nu\in V$, le nombre de cycles limites de $X_\nu$ dans $U$ est majoré par $N$,

\vskip 2mm

\item "$(ii)$" la multiplicité de chacun de ces cycles limites est majorée par $L$.

\endroster
\endproclaim

\vskip 3mm

  Plusieurs travaux, sur des cas génériques, ont été établis par moi même ([M]), ou par Il'yashenko, Yakovenko et Kaloshin ([I-Y2], [Ka]), où les singularités semi-hyperboliques sont aussi considérées. La démarche adoptée dans ces travaux consiste, d'abord en une préparation générique des singularités (dans une classe de différentiabilité suffisament grande), suivie d'une application de la procédure d'élimi-\par\noindent nation de Khovanski ([K1]). Dans [M], les conditions génériques sont algébriquement controlées tout au long de cette procédure. Cependant, dans [I-Y2], il semble très difficile de relier les conditions génériques de cette procédure à celles, géométriques provenant du polycycle.\par
  
\vskip 3mm

 Dans le cas général du théorème 0, et en vue de couvrir les cas les plus dégénérés, l'approche mise en oeuvre peut se résumer ainsi: on prépare localement (dans une subdivision finie), l'application de retour $p_\nu$ du polycycle perturbé. Dans cette préparation, les {\bf propriétés de finitude} de $p_\nu$ (le nombre de ses points fixes et leur multiplicité), sont données par celles d'un certain jet fini qui est un {\bf fewnomial} [K1]. La théorie de Khovanski s'applique aisément à ces jets.\par
 
\vskip 3mm

  Voici les idées de base de cette approche: soient $(x,\alpha)$ des coordonnées analytiques locales sur $((\RR^{+*})^k\times\RR^q,0)$, et soit $B_k=\{\prod_{j=1}^k x_j=0\}$. Soit ${\Cal B}\supset\RR\{x,\alpha\}$ un anneau local de germes de fonctions analytiques sur $((\RR^{+*})^k\times\RR^q,0)$, continues sur $(B_k,0)$. Soit $\chi$ un germe en 0 de champs de vecteurs à composantes dans ${\Cal B}$. On s'intéresse aux dérivations qui induisent une action infinitésimale sur ${\Cal B}$, et plus particulièrement à celles qui satisfont aux conditions suivantes: $\chi ({\Cal B})\subset {\Cal B}$, $Sing(\chi)\subset (B_k,0)$ et $(B_k,0)$ est invariant par le flot $\varphi_\chi$ (cf. partie IB). Le but est d'étudier les propriétés de finitude topologiques et algébriques des éléments de l'algèbre ${\Cal B}$ relativement à la distribution induite par la dérivation $\chi$ dans $((\RR^{+*})^k\times\RR^q,0)$.\par

\vskip 3mm

 Les concepts suivants sont développés dans la partie IB. Soit $U\in ((\RR^{+*})^k\times\RR^q,0)$ un ouvert sur lequel est réalisée la dérivation $\chi$. Soit $\varphi_{\chi,U}$ le flot de $\chi$ dans $U$, et soit $\pi_{\chi,U}:\ U\to\ \widetilde{U}=U/\varphi_{\chi,U}$ la projection intégrale le long des orbites de $\chi$ dans $U$. Un germe $f\in{\Cal B}$ est dit $\chi$-régulier s'il existe un tel ouvert $U$ tel que le degré de $\pi_{\chi,U}$ restreinte aux fibres de $f$ soit fini. Il est dit $\chi$-fini s'il est $\chi$-régulier et si  son idéal différentiel $I_{\chi,f}$ est noethérien dans une extension étoilée de ${\Cal B}$. Il est dit localement $\chi$-fini s'il est $\chi$-régulier et s'il existe une subdivision finie $(U_i)$ de $U$, invariante par $\chi$, telle que chaque idéal restriction $I_{\chi,f|U_i}$ soit noethérien dans une extension étoilée de l'anneau restriction ${\Cal B}_{|U_i}$. Une sous-algèbre ou une sous-classe de ${\Cal B}$ est dite $\chi$-finie (resp. localement $\chi$-finie) si chacun de ses éléments est $\chi$-fini (resp. localement $\chi$-fini). Un résultat majeur de la partie IB est\par
 
\vskip 3mm

\proclaim{Lemme de finitude IB1} Soit ${\Cal M}$ l'idéal maximal de ${\Cal B}$ et soit ${\Cal M}_0\subset {\Cal M}$ un idéal stable par $\chi$. Soit ${\Cal B}'\subset{\Cal B}$ une sous-algèbre $\chi$-finie et stable par $\chi$. Alors la classe ${\Cal C}_{{\Cal M}_0,{\Cal B}'}=\{f\in g+{\Cal M}_0 I_{\chi,g};\ g\in{\Cal B}'\}$ est $\chi$-finie.
\endproclaim

\vskip 3mm

\noindent Dans ce cas, on dit que $f$ est $\chi$-équivalente à $g$ et plus généralement, on parle d'algèbres ou de classes $\chi$-équivalentes. Ce lemme donne une idée du type de préparation que l'on souhaite établir. Décrivons brièvement les outils pour atteindre cet objectif. Si $\gamma\subset U$ est une orbite de $\chi$, le lemme d'isomorphie IB4 dit que les fibres du faisceau différentiel ${\Cal I}_{\chi,f}[\gamma]$ sont isomorphes dans les anneaux analytiques locaux $\RR\{.\}$ correspondants. D'où l'existence d'un unique idéal transverse $J_{\chi,f,\gamma}$ dans un anneau analytique $\RR\{\beta\}$ dont les coordonnées $\beta$ sont des intégrales premières de $\chi$ le long de $\gamma$. De plus, le lemme de saturation IB5 permet de reconstruire chaque fibre du faisceau différentiel à partir de cet idéal transverse: pour tout $m\in\gamma$, $I_{\chi,f}(m)=\pi_{\chi_m}^*(J_{\chi,f,\gamma})$. La question naturelle qui se pose alors est: si $\gamma$ adhère à 0, quel est le lien entre la fibre différentiel en 0 et le saturé de cet idéal transverse? Il est clair que ce lien est d'autant plus fort que l'orbite $\gamma$ est principale dans $U$, i.e le saturé de toute transversale analytique à $\gamma$ est un voisinage de 0 dans $U$. L'idéal $I((B_k,0))$ est principal de générateur $\theta=\prod_{j=1}^k x_j$. Alors, ce lien s'exprime en général par une double inclusion qui relaxe l'égalité le long de $\gamma$, et qui s'inspire du Nullstellensatz d'Hilbert

$$
 (\theta^n)\pi_\chi^*(J_{\chi,f,\gamma})\subset I_{\chi,f}\subset \pi_\chi^*(J_{\chi,f,\gamma})
\tag $*$
$$

\noindent Le plus petit de ces entiers $n$ est la multiplicité $m_\chi(f)$ de $f$ relativement à $\chi$. Si l'anneau ${\Cal B}$ possède une structure asymptotique, il se trouve que cette multiplicité est intimement liée à la multiplicité algébrique $ma_\chi(f)$ qui est l'indice de stationnarité d'une suite croissante d'idéaux transverses, qui converge vers $J_{\chi,f,\gamma}$. Ainsi, en étudiant l'action de $\chi$ sur les jets finis de $f$, et en utilisant la double inclusion $(*)$, on montre que $f$ est $\chi$-équivalente à son jet d'ordre $ma_\chi(f)$ (cf. sections II, III, IV). Le lemme de finitude ci-dessus donnera les propriétés de finitude voulues.\par

\vskip 3mm

  La projection $\pi_\chi$ est unidimensionnelle. On peut généraliser cette approche du bord $(B_k,0)$, en considérant des projections $p$-dimensionnelles le long, par exemple des feuilles d'un feuilletage de dimension $p$.\par

\vskip 3mm

  Le germe de l'application de Dulac, et de ses déploiements, appartiennent à certaines algèbres $QR{\Cal H}^{1,.}$ décrites dans la partie IA. L'algèbre $QR{\Cal H}^{k,.}(x,.)\subset {\Cal B}$ est constituée des germes qui sont quasi-analytiques dans les coordonnées $x$\par
  
$$
QR{\Cal H}\cap (\cap_{n\in\NN}{\Cal M}_x^n)=\{0\}
$$

\noindent (${\Cal M}_x=\langle x_1,\ldots,x_k\rangle\subset{\Cal B}$); et qui possèdent une structure asymptotique élémentaire dans les coordonnées $x$: ces ingrédients sont des fonctions élémentaires d'Ecalle-Khovanski. Soit $\chi\in\Xi{\Cal H}_k$ une dérivation d'Hilbert réalisée sur un ouvert $U_k\in((\RR^{+*})^k\times\RR^q,0)$, de dimension de non trivialité $k-1$ (cf. partie IVA). Elle admet une action sur l'algèbre $QR{\Cal H}^{k,.}$ et elle possède une orbite principale $\gamma$ incluse dans $U_k$. Dans la partie IVC, on montre le théorème général suivant, dont le théorème 0 est une conséquence immédiate

\vskip 3mm

\proclaim{Théorème IVC1} L'algèbre $QR{\Cal H}^{k,.}$ est localement $\chi$-finie et satisfait localement à la double inclusion $(*)$.
\endproclaim

\vskip 3mm

 La dérivation $\chi$ n'est pas réduite. Il existe une désingularisation $(\pi_k,{\Cal N}_k)$ entièrement décrite pas les algèbres $QR{\Cal H}^{k,.}$ (cf. partie IVA), et dans laquelle les singularités réduites de $\chi$ sont de la forme 

$$
 \chi_\ell=\rho\frac{\partial}{\partial \rho}-\sum_{j=1}^\ell s_ju_j\frac{\partial}{\partial u_j}
$$ 

\noindent pour $\ell=0,\ldots,k-1$. Le théorème IVC1 est donc une conséquence de l'étude de l'action des dérivations réduites $\chi_\ell$ sur les algèbres $QR{\Cal H}^{p,.}(\rho,\rho',.)$, avec $p\leq k$. Or, une dérivation $\chi_\ell$, réalisée sur un ouvert $U_p$, admet une orbite principale incluse dans $U_p$ si et seulement si $p=1$. Dans ce cas, on montre les résultats principaux suivants dans les sections II et III\par

\vskip 3mm

\proclaim{Théorème principal II1} L'algèbre $QR{\Cal H}^{1,.}(\rho,.)$ est $\chi_0$-finie et satisfait à la double inclusion.
\endproclaim

\vskip 3mm

\noindent Ce théorème, de démonstration simple et basique est une introduction aux autres théorèmes. Notons $QR{\Cal H}^{1,.}_{cvg}$ la restriction d'un anneau analytique $\RR\{.\}$ au graphe des fonctions élémentaires de l'algèbre $QR{\Cal H}^{1,.}$ correspondante. Sa $\chi_\ell$-finitude est une conséquence simple de résultats de géométrie analytique classique et de la théorie de Khovanski-Tougeron ([K1], [T]).\par

\vskip 3mm

\proclaim{Théorème principal IIIA1} Pour tout $\ell$, l'algèbre $QR{\Cal H}^{1,.}_{cvg}$ satisfait à la double inclusion relativement à $\chi_\ell$.
\endproclaim

\vskip 3mm

\proclaim{Théorème principal IIIB1} Pour tout $\ell$, l'algèbre $QR{\Cal H}^{1,.}$ est localement $\chi_\ell$-finie et satisfait localement à la double inclusion.
\endproclaim

\vskip 3mm

 Si $p>1$, la dérivation $\chi_\ell$ admet une orbite principale $\gamma_p$ incluse dans le bord $B_p$: le saturé dans $U_p$ de toute semi-transversale analytique à $\gamma_p$ est un voisinage de 0 dans $U_p$. Les fibres différentielles le long de $\gamma_p$ ne sont pas forcément isomorphes et il n'existe pas forcément d'idéal transverse. Plus généralement, deux questions se posent concernant les anneaux $QR{\Cal H}^{p,.}$: sont-ils noethériens? et leurs semi-analytiques sont-ils induits par une structure $o$-minimale au dessus de $\RR$?\par
 
\vskip 3mm

  La preuve du théorème IVC1 s'appuie sur les trois théorèmes principaux ci-dessus, et sur 6 lemmes de base (cf. partie IB). Soit $\overline{{\Cal D}}_k$ le diviseur exceptionnel du morphisme $(\pi_k,{\Cal N}_k)$ de désingularisation de $\chi$ (cf. IVA). Soit $f\in QR{\Cal H}^{k,.}$ et soient $\widetilde{f}$ et $\widetilde{\chi}$ les relevés de $f$ et $\chi$ par $\pi_k$. Il s'agit de montrer que le faisceau ${\Cal I}_{\widetilde{\chi},\widetilde{f}}[\overline{{\Cal D}}_k]$ est localement $\widetilde{\chi}$-fini. La dérivation $\widetilde{\chi}$ admet une unique singularité $a_0$ sur ${\Cal D}_k$. Soit $\gamma_1\subset{\Cal D}_k$ une orbite de $\widetilde{\chi}$ et soit $a_1=\overline{\gamma}_1\cap\partial{\Cal D}_k$. Par compacité de $\overline{{\Cal D}}_k$, il suffit de montrer que le faisceau ${\Cal I}_{\widetilde{\chi},\widetilde{f}}[a_0\gamma_1 a_1]$ est localement $\widetilde{\chi}$-fini.\par

\vskip 3mm

 L'orbite $\gamma_0=\pi_k^{-1}(\gamma)$ est principale dans un voisinage $U_{1,a_0}$ de $a_0$. Le résultat en $a_0$ est donc une conséquence du théorème principal IIIB1. En tout point $a\in\gamma_1$, un représentant du germe $(\gamma_1,a)$ est principal dans un voisinage $U_{1,a}$ de $a$; cependant il est inclus dans le bord de $U_{1,a}$. Grâce au lemme de cohérence IB3, les résultats du théorème principal IIIB1 en $a_0$ se germifient en tout point $a\in\gamma_1$ suffisament proche de $a_0$: le germe en $a$ de $\widetilde{f}$ est $\widetilde{\chi}$-équivalent à un élément $g_a$ d'une algèbre convergente $QR{\Cal H}_{cvg}^{1,.}$, qui elle, satisfait au théorème principal IIIA1. Or, comme $\widetilde{f}$, ce germe se prolonge au dessus de $\gamma_1$ en une fonction $g$ dont tous les germes appartiennent à une algèbre convergente. Le lemme d'isomorphie s'applique aux faisceaux de cette algèbre le long d'orbites incluses dans le bord. Un recollement des idéaux de $g$ et de $\widetilde{f}$ donne le résultat au dessus de $\gamma_1$. En $a_1$, on utilise un argument de récurrence sur la dimension de non trivialité de la dérivation d'Hilbert (cf. lemmes de récurrence, parties IVB et IVC). Sa preuve est elle même construite autour des 6 lemmes de base et des 3 théorèmes principaux\par

\vskip 3mm

  L'application de Dulac de chaque singularité de $X_\nu$ est induite par un élément d'une algèbre $QR{\Cal H}^{1,.}$ (cf. appendice VA). Les cycles limites de $X_\nu$ correspondent aux intesections isolées des orbites d'une dérivation d'Hilbert $\chi\in\Xi{\Cal H}_k$ et des fibres d'un germe $f\in QR{\Cal H}^{k,.}$. Le théorème 0 est alors une conséquence simple du théorème IVC1: la propriété $(i)$ est équivalente à la $\chi$-régularité de $f$, et la propriété $(ii)$ est une conséquence de la noethérianité ou la locale noethérianité de l'idéal différentiel $I_{\chi,f}$. Cette approche algébrique et géomètrique est appliquable à tout ensemble limite périodique. Comme dans le problème de Dulac, la seule difficulté réside dans la {\bf complexité des structures asymptotiques} des meilleures algèbres et dérivations d'Hilbert correspondantes.\par
  
\vskip 3mm
  
   L'article est composé de 4 sections I,...,IV et un appendice V. Chaque section est subdivisée en parties A, B,... et chaque partie est subdivisée en paragraphes 1, 2,...\par

\vskip 5mm

\noindent{\bf Remerciements.} Je remercie vivement mes collègues A. Jebrane, P. Mardesic, R. Moussu, M. Pelletier, C. Rousseau et D. Schlomiuk, pour leur soutien durant ce pénible travail. Je tiens à remercier particulièrement R. Roussarie, qui a apporté d'énormes améliorations à certains résultats de ce travail.\par

\vskip 5mm

\centerline{\bf  I. Définitions et éléments de base.}\par

\vskip 5mm

\noindent {\bf A. Définitions des algèbres.}\par

\vskip 3mm

\noindent {\bf \S 1. L'algèbre} ${\Cal A}^{p,q}$.\par

\vskip 3mm

  Soit $(x,\alpha)=(x_1,\cdots,x_p,\alpha_1,\cdots,\alpha_q)$ des coordonn\'ees sur ${{\bold{\RR}}}^p \times {{\bold{\RR}}}^q$.\par

\proclaim{Définition IA1} On note ${\Cal A}^{p,q}(x,\alpha)$ (ou simplement ${\Cal A}^{p,q}$) l'alg\`ebre r\'eelle locale des germes analytiques r\'eels de $(({{{\bold{\RR}}}}^{+*})^p,0)\times ({{{\bold{\RR}}}}^q,0)$ qui sont continus sur le germe en 0 du bord $B_p=\{x_1\times \cdots \times x_p=0\}$. 
\endproclaim

 Sauf mention contraire, toutes les alg\`ebres de référence consid\'er\'ees dans la suite sont des sous-alg\`ebres locales de ${\Cal A}^{p,q}$ qui contiennent l'algèbre analytique $\RR \{ x,\alpha \}$. Ces algèbres ${\Cal A}^{p,q}$ ne sont pas stables par les dérivations les plus élémentaires: $f=x\sin (1/x)\in {\Cal A}^{1,0}$ mais $x\partial f/\partial x\not\in{\Cal A}^{1,0}$. Dans le problème d'Hilbert, ces algèbres serviront uniquement d'espaces d'intégrales premières pour les dérivations considérées (voir sections III et IV).\par

\vskip 3mm

\noindent {\bf \S 2. L'algèbre} $SB^{p,q}$ {\bf des germes sectoriellement bornés.}\par

\vskip 3mm

 Soit $\theta =(\theta_1,\cdots ,\theta_p)\in ]0,\pi/2[^p$ et $S_\theta$ le polysecteur 

$$
     S_{\theta}=\{w=(w_1,\cdots,w_p)\in ({{{\bold{\CC}}}}^*)^p;\quad 
                 |\arg (w_j)|< \theta_j \}      
\tag 1
$$
\noindent Soit $(S_\theta ,\infty )$ le germe de $S_\theta$ \`a l'infini.\par

\vskip 2mm

 Pour simplifier la présentation dans toute la suite, nous noterons souvent pareillement les germes, leurs représentants et les relevés de germes de fonctions dans la carte $w=-\log(x)$.\par

\vskip 2mm

\proclaim{Définition IA2} Les \'el\'ements de l'alg\`ebre $SB^{p,q}(x,\alpha)\subset {\Cal A}^{p,q}(x,\alpha)$ sont les germes $f$ qui admettent, pour tout $\theta \in ]0,\pi /2 [^p$, un prolongement holomorphe et born\'e sur $(S_\theta ,\infty)\times ({{{\bold{\CC}}}}^q,0)$ dans la carte $w=(w_j=-\log (x_j))_{j=1,\cdots ,p}$.\endproclaim

\vskip 2mm

 Ces alg\`ebres $SB^{p,q}$ sont le lieu naturel o\`u vivent les germes d'applications de Dulac des d\'eploiements holomorphes d'\'equations diff\'erentielles du $1^{\text{er}}$ ordre, y compris dans le domaine de Poincar\'e (voir appendice A). Ce sont les anneaux de références dans le problème d'Hilbert (voir sections II, III et IV). Leur intérêt premier réside dans leur structure holomorphe produit. Notons \'egalement $SB^{p,q}_0$ la sous-alg\`ebre de $SB^{p,q}$ des germes qui, pour tout $\theta$, tendent vers 0 quand $w\to \infty$ dans $(S_\theta,\infty)$, uniform\'ement en $\alpha$. Contrairement \`a l'algèbre ${\Cal A}^{p,q}$, les alg\`ebres $SB^{p,q}$ et $SB^{p,q}_0$ sont stables par les d\'erivations naturelles $\chi_j=x_j\partial/\partial x_j \sim -\partial /\partial w_j$ (conséquence immédiate des formules de Cauchy dans les secteurs $S_\theta$). Cependant, elles ne sont pas quasi-analytiques (voir ci-dessous).\par

\vskip 3mm

\noindent {\bf \S 3. Algèbres quasi-analytiques} $QA^{p,q}$.\par

\vskip 3mm

 Soit $P^+=\{w=u+iv\in{\bold{\CC}};\ u\geq 0\}$. Soient $u_0\geq 0$, $C>0$ et $K>1$. Les domaines de $P^+$ {\bf de type puissance} sont les domaines de la forme
 
$$
\Omega_{puis}(u_0,C,K)=\{w\in P^+;\ u>u_0,\ |v|<Cu^K\}
$$

\noindent Les domaines de $P^+$ {\bf de type exponentiel} sont les domaines de la forme

$$
\Omega_{exp}(u_0,C,K)=\{w\in P^+;\ u>u_0,\ |v|<C(\exp(u/K)-1)\}
$$

\comment
le -1 est pour la stabilité par addition
\endcomment

\noindent Ces domaines sont stables par addition (opération qui correspond à une multiplication dans la coordonnée $x=-\log w$).\par

\vskip 3mm

\proclaim{Définition IA3} Les domaines standards d'Ecalle-Il'yashenko sont les ouverts $\Omega$ de $P^+$ qui contiennent un domaine de type puissance. On note ${{\Cal E}}{{\Cal I}}$ l'ensemble de tels domaines.\par
\endproclaim

\vskip 3mm

 En particulier les domaines de type puissance et les domaines de type exponentiel, sont des domaines standards. Une intersection finie et une union finie de domaines standards est encore un domaine standard. Soit $\Omega\in{\Cal E}{\Cal I}$, si $t\in{\bold{\CC}}$ est tel que $t+\Omega$ soit inclus dans l'intérieur de $P^+$, alors $t+\Omega\in{\Cal E}{\Cal I}$, et si $t>0$, alors $t\Omega\in{\Cal E}{\Cal I}$ (la translation $w\mapsto t+w$ correspond à une homothétie $x\mapsto ax$, et l'homothétie $w\mapsto tw$ correspond à une ramification $x\mapsto x^s$).\par
 
 \vskip 3mm
 
  Les domaines de type puissance et les domaines de type exponeniel sont biholomorphiquement conjugués à un ouvert contenant $P^+$, par un difféomorphisme $\phi_.$ strictement réel et qui est équivalent à l'identité à l'infini
  
$$
\phi_{puis}(w)=w-\frac{1}{C^{1/K}\cos(\pi/2K)}w^{1/K}-U_0
$$

$$
\phi_{exp}(w)=w-K\log(w)-U_0\qquad U_0>0
$$

 Ecalle [E] et Il'yashenko [I] ont exhib\'e, pour l'application de Dulac d'un col hyperbolique réel, un type de tels
domaines: exponentiel pour le premier (et ceci est optimal pour les cols hyperboliques analytiquement normalisables: voir appendice VA), et polynomial pour le deuxi\`eme $(K=2)$.\par

\vskip 3mm

 L'application de Dulac d'une équation différentielle dans le domaine de Poincaré n'est pas bornée sur un domaine de ${\Cal E}{\Cal I}$ (voir appendice VA). Ceci motive la\par

\vskip 3mm

\proclaim{Définition IA4} L'alg\`ebre $QA^{p,q}(x,\alpha)\subset SB^{p,q}(x,\alpha)$ est l'ensemble des germes $f=\sum f_n\alpha^n$ dont les coefficients $f_n$ admettent un prolongement holomorphe et borné sur un même domaine $\Omega\in{\Cal E}{\Cal I}^p$ dans les coordonnées $w=(w_j=-\log (x_j))_{j=1,\ldots,p}$.\par
\endproclaim

\vskip 3mm

 Ces algèbres sont quasi-analytiques au sens suivant: soit ${\Cal M}_{x}=\langle x_1,\ldots,x_p\rangle$ l'idéal de $SB^{p,q}$ engendré par les fonctions coordonnées $x_j$, alors

$$ 
 QA^{p,q}\cap(\cap_n {\Cal M}_x^n)=\{0\}
\tag 2
$$ 

\noindent Pour $p=1$ et $q=0$, ce résultat a été démontré par Il'yashenko dans [I], par une double application du principe de Phragmen-Lindelof dans $P^+$ ([Ru, p.244]), en utilisant le difféomorphisme $\phi_{puis}$. Dans le cas général, notons ${\Cal M}_{x,0}$ l'idéal de $SB^{p,0}(x)$ engendré par les fonctions coordonnées $x_j$; si $f=\sum_k f_k\alpha^k\in QA^{p,q}\cap(\cap_n{\Cal M}_x^n)$, alors pour tout multi-indice $k$: $f_k\in QA^{p,0}\cap(\cap_n {\Cal M}_{x,0}^n)$ (par une simple identification des coefficients des séries en $\alpha$). Donc pour $p=1$, on obtient encore l'égalité (2). Supposons $p>1$ et considérons la restriction de $f_k$ à un voisinage de la diagonale de $({\bold{\RR}}^{+*})^p$: soit l'application $g:(y,\beta)\in {\bold{\RR}}^{+*}\times{\bold{\RR}}^{p-1}\mapsto x=g(y,\beta)=(y,y(1+\beta_1),\ldots,y(1+\beta_{p-1}))$. Soit $U\in(({\bold{\RR}}^{+*})^p,0)$ sur lequel est réalisée $f_k$ (il est indépendant de $k$), et soit $V\in({\bold{\RR}}^{+*}\times{\bold{\RR}}^{p-1},0)$ tel que $g(V)\subset U$. Soit $F_k$ le germe en 0 de $f_k\circ g_{|V}$, c'est un élément de l'anneau $SB^{1,p-1}(y,\beta)$ (car le germe à l'infini du translaté complexe de tout secteur $S_\theta$ est inclus dans le germe à l'infini d'un secteur $S_{\theta'}$). Soit ${\Cal M}_y$ l'idéal de $SB^{1,p-1}$ engendré par la coordonnée $y$, par la stabilité des domaines standards par intersection finie et par translation, on a $F_k\in QA^{1,p-1}\cap(\cap_n{\Cal M}_y^n)$ (on a même que $F_k\in QA^{1,0}(y)\{\beta\}$, la série étant convergente sur un produit $\Omega_0\times W$ où $\Omega_0$ est un domaine standard et $W$ est un voisinage de 0 dans ${\bold{\CC}}^{p-1}$). Le germe $F_k$ est donc identiquement nul, et il en est de même pour $f_k$ et pour $f$.\par

\vskip 3mm

  Ces alg\`ebres sont stables par les d\'erivations $\chi_j$ (par les formules de Cauchy dans les coordonnées $w_j$ dans les translatés $1+\Omega_j$). Leur localité est un problème ouvert. Elles sont strictement incluses dans les algèbres $SB^{p,q}$: en effet, le germe $f(x)=x^{\log(-\log(x)}\in SB^{1,0}(x)\cap(\cap_n{\Cal M}_x^n)$ (pour cela, il suffit de voir que pour tout $n\in{\NN}$, le relevé $f_n(w)=\exp(-w(\log(w)-n))$ est borné sur tout germe $(S_\theta,\infty)$ avec $\theta\in]0,\pi/2[$). Cependant, si $p>0$, la topologie de Krull de l'anneau $QA^{p,q}$ (induite par celle de l'anneau $SB^{p,q}$) n'est pas séparée (pour tout $s>0$, $x_1^s\in{\Cal M}$, où ${\Cal M}$ est l'idéal maximal de $SB^{p,q}$). Et les idéaux les plus simples des anneaux $QA^{p,q}$ ne sont pas noethériens (ni dans l'anneau $QA^{p,q}$ ni même dans son extension ${\Cal A}^{p,q}$): tel est le cas par exemple des idéaux engendrés par les {\bf fonctions élémentaires} $f_n=x(\log x)^n$ ou $g_n=x^{1/n}$. Ainsi, une structure asymptotique dans un {\bf nombre fini de fonctions élémentaires} est souhaitable.\par

\vskip 3mm

\noindent {\bf \S 4. Algèbres quasi-régulières d'Hilbert} $QR{\Cal H}^{p,q}$.\par

\vskip 3mm

  La derni\`ere condition de r\'egularit\'e qu'on impose sur les germes \'etudi\'es est l'existence d'une structure asymptotique élémentaire dans les variables quasi-ana-\par\noindent lytiques $x_j$. Pour $(y,\beta)\in\RR^{+*}\times\RR$, notons Ld (pour {\bf L}ogarithme {\bf d}\'eploy\'e) la fonction 

$$
\text{Ld} (y,\beta)=\int_1^y {t^{-1+\beta}} dt =
\left\lbrace
\aligned 
      \frac{y^\beta-1}{\beta}\quad & \text{ pour }\quad \beta\neq 0 \\
        \log y \quad & \text{ pour }\quad \beta=0
\endaligned
\right. 
\tag 3
$$

\noindent Ceci est simplement le compensateur \'el\'ementaire d'Ecalle-Roussarie [E], [Ro2]. La fonction $f(y,\beta)=yLd(y,\beta)\in QA^{1,1}(y,\beta)$: soit $F(w,\beta)=f(\exp(-w),\beta)$ et soit $\theta\in[0,\pi/2[$, en faisant le changement de coordonnées $t=\exp(-z)$ dans l'intégrale (3), on vérifie facilement que pour $|\beta|$ suffisament petit, on a $|F(w,\beta)|\leq 1/\cos(\theta)$ sur le secteur $S_\theta$. De plus, $F(w,\beta)=\sum F_n(w)\beta^n$ avec 

$$
F_n(w)=\frac{1}{(n+1)!}(-w)^{n+1}\exp(-w)
$$

\noindent chaque fonction $F_n$ est bornée sur le domaine exponentiel $\Omega_{exp}(0,1,n+2)$ et sur tout compact de $P^+$. Les fonctions $F_n$ sont donc bornées sur n'importe quel domaine de type puissance (voir appendice VA pour une preuve générale).\par

\vskip 3mm

  Soit $q=(q_1,q_2)\in\NN^2$ et $\alpha=(\mu,\nu)$ des coordonnées sur $\RR^{q_1}\times\RR^{q_2}$. Soient les fonctions élémentaires $z_{i,0}(x_i)=x_i\log x_i$ et $z_{i,j}$ leurs d\'eploiements

$$  
z_{i,j}(x_i,\mu_j)=x_i\text{Ld}(x_i,\mu_j)  
\tag 4
$$

\noindent Ces fonctions appartiennent \`a l'algèbre $QA^{1,1}$. Dans la suite, certaines notations (de sens clair dans le texte) désignent aussi bien des fonctions que les coordonnées correspondantes. Soient

$$
 X_i=(x_i,z_{i,0},z_{i,1},\cdots ,z_{i,q_1})\quad \text{ et }\quad X=(X_1,\ldots,X_p)
\tag 5
$$

\noindent Notons $\widehat x^i=(x_1,\cdots,x_{i-1},x_{i+1},\cdots,x_p)$ et $c_i$ et $c$ les immersions

$$ 
c_i(x,\alpha)=(X_i,\widehat x^i,\alpha)\quad\quad c(x,\alpha)=(X,\alpha)
\tag 6
$$

\vskip 3mm

\proclaim{Définition IA5} Convenons que $QR{\Cal H}^{0,q}(\alpha)=\RR\{\alpha \}$. Alors, l'alg\`ebre quasi-r\'egu\-li\`ere d'Hilbert $QR{\Cal H}^{p,q}(x,\alpha)\subset QA^{p,|q|}(x,\alpha)$ est l'ensemble des germes $f$ ayant un d\'eveloppement asymptotique de "type Hilbert": pour tout $i=1,\cdots ,p$, il existe une suite $(G_{i,m})_m$ dans $QR{\Cal H}^{p-1,q}(\widehat{x}_i,\alpha)[X_i]$ qui sont des polyn\^omes homog\`enes de degr\'e $m$ dans la variable $X_i$ telle que pour tout  $n\in{\bold{\NN}}$

$$  
f(x,\alpha)={\sum}_{m=0}^n G_{i,m}\circ c_i (x,\alpha) +x_i^n h_n\qquad\text{avec}\qquad h_n \in SB^{p,|q|}_0
\tag 7
$$

\endproclaim

\vskip 3mm

 Les variables analytiques $\nu$ n'interviennent pas dans la construction des fonctions élémentaires (4), d'où la distinction faite dans les variables analytiques $\alpha$. Dans le problème d'Hilbert, les variables $\mu$ sont les paramètres qui déploient les valeurs propres des singularités, et les variables $\nu$ sont tout autres paramètres.\par

\vskip 3mm

 L'unicité des séries formelles (7) ainsi que l'injectivité des morphismes série formelle associés sont démontrés dans la section II. On y démontre aussi l'existence et l'injectivité d'un morphisme série formelle $f\in QR{\Cal H}^{p,q}\mapsto \widehat{f}\in c^*(\RR\{\alpha\}[[X]])$; ceci implique en particulier que la topologie de Krull des algèbres $QR{\Cal H}^{p,q}$ est séparée. Les germes quasi-analytiques ($\in QA^{p,|q|}$) qui possèdent une telle série formelle forment une sur-algèbre de $QR{\Cal H}^{p,q}$ qui ne sera pas étudié dans ce travail. La sous-alg\`ebre $QR{\Cal H}_{\text{cvg}}^{p,q}$ des \'el\'ements "convergents" de l'alg\`ebre $QR{\Cal H}^{p,q}$ est d\'efinie par 

\proclaim{Définition IA6} On note   $QR{\Cal H}_{\text{cvg}}^{p,q}=c^* (\RR \{ X,\alpha \})$.
\endproclaim

\noindent Une conséquence algébrique de la transcendance du graphe de $c$ est que le morphisme $c^*$ est un isomorphisme sur son image. Ceci est démontré aussi dans le section II.\par

\vskip 5mm

\noindent {\bf B. Quelques généralités et six lemmes de base.}\par

\vskip 5mm

 Les {\bf anneaux de référence} sont les anneaux locaux ${\Cal B}$ telles que $\RR \{ x,\alpha\}\subset {\Cal B} \subset {\Cal A}^{p,q}(x,\alpha)$ et qui sont stables par les dérivations

$$
 \prod_{j=1}^p x_j \frac{\partial}{\partial y_i}\qquad\text{avec}\quad y=(x,\alpha)
$$

\noindent Ce sont des $\RR$-algèbres. Dans la suite, on parlera indifférement d'anneau ou d'algè-\par\noindent bre. L'idéal maximal de ${\Cal B}$ est l'idéal des germes nuls en 0. En cas d'ambiguité, on note ${\Cal B}(y)$ pour préciser le choix des coordonnées.\par

\vskip 3mm

\noindent {\bf §1. Anneaux restriction et anneaux extension}.\par

\vskip 3mm

   Soit $U$ un repr\'esentant de $(({\bold{\RR}}^{+*})^p\times {\bold{\RR}}^q,0)$ et soit $B_p=\{\prod_{j=1}^p x_j=0\}$ le bord associé de germe $(B_p,0)$ en 0. Soit $U_O\subset U$ un sous-ensemble quelconque dont l'adhérence contient 0. On note $(U_0,0)$ le germe de $U_0$ en 0 . On note aussi $\partial_0 U_0=B_p\cap \overline{U_0}$ le bord associé et $\partial_0(U_0,0)$ son germe en 0. Soit ${\Cal B}\subset{\Cal A}^{p,q}(y)$ un anneau de référence et soit $i_{U_0}:(U_0,0)\to (U,0)$ le germe de l'injection canonique (qu'on notera aussi $i_{U_0,U}$ en cas d'ambiguité).  On lui associe un morphisme étoilé
   
$$
i^*_{U_0}:f\in{\Cal B}\mapsto i^*_{U_0}(f)=f\circ i_{U_0}
$$
   
\noindent On généralise ainsi les anneaux de référence ${\Cal B}$ aux {\bf anneaux restriction} notés ${\Cal B}_{|U_0}$ et définis comme suit:

$$
 {\Cal B}_{|U_0}=i^*_{U_0}({\Cal B})
$$

\noindent C'est un anneau local (qui n'est pas forcément intègre). Il ne dépend de $U_0$ que par son germe $(U_0,0)$. Il est isomorphe à ${\Cal B}$ (par $i^*_{U_0}$) si  $(U_0,0)$ est d'intérieur non vide (ie. l'adhérence de l'intérieur de $U_0$ contient 0): en effet, une fonction analytique nulle sur un ouvert connexe est nulle sur la composante connexe de son domaine d'analycité contenant cet ouvert. Le sous-ensemble $U_0$ est dit {\bf semi-analytique élémentaire} de ${\Cal B}$ (ou décrit par ${\Cal B}$) s'il existe $V\in(U,0)$ et contenant $U_0$ et des germes $f_1,\ldots,f_n,g_1,\ldots,g_m\in{\Cal B}$ et représentés sur $V$ tels que

$$
U_0=\{y\in V;\ f_1(y)>0,\ldots,f_n(y)>0,g_1(y)=0,\ldots,g_m(y)=0\}
$$

\noindent Il est dit {\bf semi-analytique} de ${\Cal B}$ (ou décrit par ${\Cal B}$) si c'est une union finie de semi-analytiques élémentaires de ${\Cal B}$. Dans ce cas, le morphisme $i^*_{U_0}$ est un isomorphisme si et seulement si $(U_0,0)$ est d'intérieur non vide: en effet, si $U_{0,1},\ldots,U_{0,\ell}$ sont les semi-analytiques élémentaires formant $U_0$, l'un des $U_{0,j}$ est un ouvert dont l'adhérence contient 0, sinon il existe $g_{i_1,1},\ldots,g_{i_\ell,\ell}\in{\Cal B}\setminus\{0\}$ tels que $g=g_{i_1,1}\times\cdots\times g_{i_\ell,\ell}$ est nulle sur $(U_0,0)$; mais par l'isomorphisme $i^*_{U_0}$, $g$ est alors  nulle, ce qui contredit l'intégrité de ${\Cal B}$.\par

\vskip 3mm

  Soit $I$ un idéal de ${\Cal B}$. Le morphisme $i^*_{U_0}$ étant surjectif, le sous-ensemble $i^*_{U_0}(I)$ est un idéal de ${\Cal B}_{|U_0}$ dit {\bf idéal restriction}. On le note simplement $I_{|U_0}$; il ne dépend de $U_0$ que par son germe $(U_0,0)$. Inversement, tout idéal $J$ de ${\Cal B}_{|U_0}$ est un idéal restriction: $I=(i^*_{U_0})^{-1}(J)$ est un idéal de ${\Cal B}$ et $I_{|U_0}=J$.\par

\vskip 3mm

  Soit ${\Cal B}'\subset {\Cal A}^{p',q'}(x',\alpha')$ un anneau de référence. Soit $U'\in ((\RR^{+*})^{p'}\times \RR^{q'},0)$ et soit $U'_0$ un sous-ensemble de $U'$ dont l'adhérence contient 0. L'anneau ${\Cal B}'_{|U'_0}$ est dit {\bf anneau extension} de l'anneau ${\Cal B}_{|U_0}$ s'il existe un homomorphisme d'anneaux injectif $\Psi:\ {\Cal B}_{|U_0}\hookrightarrow{\Cal B}'_{|U'_0}$. Dans ce cas, si $J$ est un id\'eal de ${\Cal B}_{|U_0}$, on appelle {\bf idéal prolongé} associé à $J$ l'idéal de ${\Cal B}'_{|U'_0}$ engendr\'e par le sous-ensemble $\Psi(J)$. On le notera simplement $\Psi(J)$ si aucune confusion n'est à craindre. Cet idéal prolongé est aussi un idéal restriction.\par

\vskip 3mm

 Soit $\psi: U'_0 \rightarrow U_0$ un morphisme {\bf surjectif, continu sur } $U'_0\cup\{0\}$ et tel que $\psi(0)=0$. On note de la même façon son germe  $\psi:(U'_0,0)\rightarrow (U_0,0)$. Ce germe induit un morphisme étoilé $\psi^*$ qui agit sur l'anneau ${\Cal B}_{|U_0}$ et qui est injectif. On suppose que $\psi^*({\Cal B}_{|U_0})\subset {\Cal B}'_{|U'_0}$. Dans ce cas, on dira que l'anneau ${\Cal B}'_{|U'_0}$ est une {\bf extension étoilée} de l'anneau ${\Cal B}_{|U_0}$ et on la note $({\Cal B}'_{|U'_0},\psi)$. Si $g$ est un élément de ${\Cal B}_{|U_0}$, on a la relation suivante entre les germes en 0 des ensembles de zéros

$$
 Z(\psi^*(g))=\psi^{-1}(Z(g))
\tag 0
$$

\noindent Dans toute la suite, on ne considérera que des extensions étoilées.\par

\vskip 3mm

\proclaim{Définition IB1} Soit $I$ un idéal de ${\Cal B}$.\par

\roster

\vskip 2mm

\item"$(i)$" On dit que $I$ est {\bf "noethérien"} sur $U_0$ (ou que $I_{|U_0}$ est {\bf "noethérien"}) s'il existe une extension étoilée $({\Cal B}'_{|U'_0},\psi)$ de ${\Cal B}_{|U_0}$ dans laquelle l'idéal prolongé $\psi^*(I_{|U_0})$ est {\bf noethérien} dans le sens classique: il existe $h_1,\ldots,h_n\in \psi^*(I_{|U_0})$ tels que pour tout $h\in \psi^*(I_{|U_0})$, il existe $H_1,\ldots,H_n\in {\Cal B}'_{|U'_0}$ tels que $h=\sum_{j=1}^n H_j h_j$.\par

\vskip 2mm

\item"$(ii)$" On dit que $I$ est {\bf localement "noethérien"} sur $U_0$ (ou que $I_{|U_0}$ est {\bf localement "noethérien"}) s'il existe une subdivision finie de $U_0$ en des sous-ensembles $U_{i,0}$ qui adhérent à $0$ et telle que $I$ soit {\bf "noethérien"} sur chaque $U_{i,0}$.\par

\endroster

\endproclaim

\vskip 3mm

\noindent Si $(U_0,0)=(U,0)$, on dit simplement que $I$ est "noethérien" ou localement "noethé-\par\noindent rien". Dorénavant, on enlève les guillemets au mot "noethérien". En cas d'ambigui-\par\noindent té, on précisera l'anneau de référence pour les idéaux prolongés (qui est aussi l'extension étoilée associée). Soit $J\subset {\Cal B}_{|U_0}$ un idéal noethérien et soit $({\Cal B}'_{|U'_0},\psi)$ l'extension associée. On sait définir le germe en 0 de l'ensemble des zéros de l'idéal prolongé $\psi^*(J)$: c'est celui de n'importe quel système fini de générateurs de cet idéal dans l'anneau ${\Cal B}'_{|U'_0}$. Maintenant, si $g_1,\ldots,g_n\in J$ sont tels que $\psi^*(g_1),\ldots,\psi^*(g_n)$ forment un système de générateurs de $\psi^*(J)$ (il en existe), la relation (0) montre que le germe $Z(g_1)\cap\cdots\cap Z(g_n)$ est indépendant du système ainsi choisi, et qu'on a donc une notion {\bf d'ensemble des zéros de l'idéal} $J$, qu'on note $Z(J)$ et qui est donné par la formule ($\psi$ étant surjective)

$$
Z(\psi^*(J))=\psi^{-1}(Z(J))
\tag 1
$$

\noindent En particulier, pour tout $g\in J$, on a $g_{|Z(J)}=0$.\par

\vskip 3mm

\noindent {\bf §2. Projection unidimensionnelle}.\par

\vskip 3mm

 Soit $\Xi{\Cal B}$ la classe des germes en 0 de champs de vecteurs $\chi=\sum_{j=1}^{p+q} a_j(y)\partial/\partial y_j$ dont les composantes $a_j$ sont des \'el\'ements de ${\Cal B}$, et qui satisfont aux conditions suivantes\par

\roster

\vskip 2mm

\item"$(i)$" $\text{Sing} (\chi) \subset (B_p,0)$: il existe un ouvert $U\in ((\RR^{+*})^p\times\RR^q,0)$ tel que le champ $\chi$ se prolonge par continuité à $\overline{U}$ et tel que l'ensemble singulier de $\chi$ sur $\overline{U}$ soit inclus dans $\partial_0 U$.\par

\vskip 2mm

\item"$(ii)$" $(B_p,0)$ est invariant par $\chi$: il existe $U$ comme dans $(i)$ tel que $\partial_0 U$ est une union d'orbites de $\chi$ dans $\overline{U}$.\par

\vskip 2mm

\item"$(iii)$" $\chi ({\Cal B})\subset {\Cal B}$.\par

\endroster

\vskip 3mm

 Un ouvert $U$ satisfaisant aux conditions $(i)$ et $(ii)$ est dit {\bf admissible}. Soit $\chi\in\Xi{\Cal B}$ et soit $U$ un ouvert admissible. Soit $\varphi_{\chi,U}$ le flot de $\chi$ dans $U$; on note $\pi_{\chi,U}: U \mapsto \widetilde{U}=U/\varphi_{\chi,U}$ {\bf la projection le long des orbites de} $\chi$ {\bf dans} $U$ (on l'appelle aussi {\bf le morphisme intégral de} $\chi$ {\bf dans} $U$). L'espace $\widetilde{U}$ étant muni de la topologie quotient. Cette projection est donc continue et ouverte. Soit $S\subset (\RR^{+*})^p\times\RR^q$ dont l'adhérence contient 0\par

\vskip 3mm

\proclaim{Définition IB2} {\bf Le degr\'e de la projection} $\pi_{\chi,U}$ restreinte au sous-ensemble $S$ est 

$$
d^°\pi_{\chi,U|S}=\sup_{\gamma \subset U} b_0(\gamma \cap S)
$$ 

\noindent où $\gamma$ est une orbite de $\chi$ dans $U$, et $b_0$ est le premier nombre de Betti. {\bf On note}

$$
d^°\pi_{\chi|(S,0)}=\inf_U d^°\pi_{\chi,U|S}
$$

\noindent où la borne inférieure est prise sur tous les ouverts $U$ admissibles.\par
\endproclaim

\vskip 3mm

 La notation $\pi_\chi$ dans cette définition, ne désigne pas un germe. En général, il n'existe pas de notion de germe en 0, pour la projection intégrale, qui soit indépendante des ouverts $U$ (ou du moins d'une base d'ouverts $U$). Quand il en existe une, on note ce germe $\pi_\chi$; c'est par exemple le cas quand le champ $\chi$ admet une orbite "principale" dans $U$ (cf. fin de cette section), ou plus généralement, quand il admet dans $U$, $(p+q-1)$ intégrales premières $F_j$, analytiques et indépendantes: la (p+q-1)-forme $dF_1\wedge\cdots\wedge dF_{p+q-1}$ ne s'annule pas sur $U$.\par

\vskip 3mm

\noindent {\bf 2.1} $\chi${\bf -régularité et} $\chi${\bf -finitude}.\par

\vskip 3mm

 Soit $U_0\subset (\RR^{+*})^p\times\RR^q$ dont l'adhérence contient 0. On dit que le germe $(U_0,0)$ (ou simplement $U_0$) est {\bf invariant par} $\chi$ s'il existe un ouvert $U$ admissible tel que $U_0\cap U$ soit une union d'orbites de $\chi$ dans $U$.\par
 
\vskip 3mm

\proclaim{Définition IB3} Soit $f\in{\Cal B}$ et $Z(f)$ le germe en 0 de son ensemble des zéros. On dit que $f$ est $\chi${\bf -r\'eguli\`ere sur} $U_0$ si
  
$$
d^°\pi_{\chi|Z(f)\cap (U_0,0)}<+\infty
$$

\endproclaim

\vskip 3mm

\noindent Si $(U_0,0)=(U,0)$ où $U$ est un ouvert admissible, on dit simplement que $f$ est $\chi$-régulière. Van Den Dries parle dans l'un de ces travaux ([Dr]) d'une certaine "propriété de finitude" (qui porte justement sur des projections unidimensionnelles mais linéaires!), qui est équivalente à cette notion de $\chi$-régularité. Soit $I_{\chi,f}=\langle \chi^n f;\ n\in \NN \rangle$ l'idéal différentiel de $f$ dans l'anneau ${\Cal B}$. Si $U'\subset U$ est un ouvert admissible sur lequel est réalisée $f$, alors tous les germes $\chi^n f$ sont aussi réalisés sur $U'$. On note $Z(I_{\chi,f}$ le germe en 0 de $\cap_{n\in\NN} S_n$, où $S_n$ est un représentant de $Z(\chi^n f)$ sur $U'$. Cette définition de l'ensemble des zéros d'un idéal différentiel coincide avec la définition classique en cas de noethérianité. Remarquer que le germe $Z(I_{\chi,f})$ est invariant par $\chi$.\par

\vskip 3mm

\proclaim{Définition IB4} On suppose que $f$ est $\chi$-régulière sur $U_0$,\par

\roster

\vskip 2mm

\item"$(i)$" on dit que $f$ est $\chi${\bf -finie sur} $U_0$ si l'id\'eal différentiel $I_{\chi,f}$ est noeth\'erien sur $U_0$.\par

\vskip 2mm

\item"$(ii)$" On dit que $f$ est {\bf localement} $\chi${\bf -finie sur} $U_0$ si $I_{\chi,f}$ est localement noethé-\par\noindent rien sur $U_0$, la subdivision associée étant invariante par $\chi$.\par

\vskip 2mm

\endroster

\endproclaim

\vskip 3mm

\proclaim{Définition IB5} Une classe ${\Cal C}\subset {\Cal B}$ est dite $\chi${\bf -finie sur} $U_0$ (resp. {\bf localement} $\chi${\bf -finie sur} $U_0$) si tout $f\in {\Cal C}$ est $\chi$-finie sur $U_0$ (resp. localement $\chi$-finie sur $U_0$). Elle est dite $\chi${\bf -stable} si $\chi({\Cal C})\subset {\Cal C}$.\endproclaim 

\vskip 3mm

 Si $(U_0,0)=(U,0)$ où $U$ est un ouvert admissible, on dit simplement que $f$ (ou la classe ${\Cal C}$) est $\chi$-finie, ou alors que $f$ (ou la classe ${\Cal C}$) est localement $\chi$-finie.\par

\vskip 3mm

\noindent{\bf 2.2 Six lemmes de base.}\par

\vskip 3mm

 Le premier de ces lemmes, qui est un lemme fondamental, donne les propriétés de finitude d'une certaine classe de germes qui sont "bien préparés" et dans une extension appropriée. On suppose donc dans ce lemme que les anneaux de référence pour la noethérianité, sont les anneaux restriction ${\Cal B}_{|U_0}$. Soit ${\Cal M}$ l'idéal maximal de ${\Cal B}$\par

\vskip 3mm

\proclaim{Lemme IB1 (lemme de $\chi$-finitude)} Soit ${\Cal B}_0\subset {\Cal B}$ une alg\`ebre $\chi$-finie sur $U_0$ (resp. localement $\chi$-finie sur $U_0$) et $\chi$-stable. Soit ${\Cal M}_0\subset {\Cal M}$ un id\'eal $\chi$-stable. Soit ${\Cal N}_{U_0}\subset{\Cal B}$ l'idéal des germes nuls sur $(U_0,0)$. Alors la classe ${\Cal C}_{{\Cal B}_0,{\Cal M}_0}=\{f\in  g+ {\Cal M}_0 I_{\chi,g}+{\Cal N}_{U_0};\ g\in{\Cal B}_0\}\subset{\Cal B}$ est $\chi$-finie sur $U_0$ (resp. localement $\chi$-finie sur $U_0$).\par
\endproclaim

\vskip 3mm

\noindent {\bf Preuve.} Il suffit de la faire dans le cas de $\chi$-finitude sur $U_0$. Soit $g\in{\Cal B}_0$ et soit 

$$
 f\in g+{\Cal M}_0 I_{\chi,g}+{\Cal N}_{U_0}
\tag 2
$$

\noindent Le germe $(U_0,0)$ étant invariant par la dérivation $\chi$, l'idéal ${\Cal N}_{U_0}$ est $\chi$-stable. La relation (2) implique donc que $I_{\chi,f|U_0}\subset I_{\chi,g|U_0}$. Montrons d'abord que l'id\'eal $I_{\chi,f|U_0}$ est noeth\'erien dans l'anneau ${\Cal B}_{|U_0}$. L'algèbre ${\Cal B}_0$ étant $\chi$-finie sur $U_0$, l'idéal $I_{\chi,g|U_0}$ est noethérien dans l'anneau ${\Cal B}_{|U_0}$. Soit $(g_{|U_0}),\ldots,((\chi^\ell g)_{|U_0}))$ un système de générateurs de l'idéal $I_{\chi,g|U_0}$. L'idéal maximal de ${\Cal B}_{|U_0}$ est ${\Cal M}'={\Cal M}_{|U_0}$ (car les éléments de ${\Cal B}\subset{\Cal A}^{p,q}$ sont continues en 0). Comme ${\Cal M}_0\subset{\Cal M}$, on a ${\Cal M}'_{0}=({\Cal M}_{0|U_0})\subset {\Cal M}'$. En appliquant $\ell$ fois la dérivation $\chi$ à la relation (2), puis en prenant la restriction à $(U_0,0)$, on obtient

$$
 I_{\chi,g|U_0}\subset I_{\chi,f|U_0}+{\Cal M}'_0 I_{\chi,g|U_0}
$$

\noindent Comme l'idéal $I_{\chi,g|U_0}$ est noethérien, le Lemme de Nakayama (cf. [L]) implique que $I_{\chi,g|U_0}\subset I{\chi,f|U_0}$. On obtient alors $I_{\chi,f|U_0}=I_{\chi,g|U_0}$. L'idéal $I_{\chi,f|U_0}$ est donc noethérien dans l'anneau ${\Cal B}_{|U_0}$.\par

\vskip 3mm

 Montrons maintenant que $f$ est $\chi$-r\'eguli\`ere sur $U_0$. Soit $F=Z(I_{\chi,f|U_0})\subset (U_0,0)$ l'ensemble des z\'eros de l'id\'eal $I_{\chi,f|U_0}$; il est invariant par $\chi$ et $d^°\pi_{\chi|F}\leq 1$. Or, on a aussi $F=Z(I_{\chi,g})$. Soit $U$ un ouvert admissible tel que $U_0\cap U$ soit une union d'orbites de $chi$ dans $U$. Soit $U'\subset U$ un ouvert admissible où sont réalisés le germe $g$ et la relation (2). le sous-ensemble $U_0\cap U'$ est encore une union d'orbite de $\chi$ dans $U'$. Notons $S$ le représentant du germe $F$ sur $U'$ et considérons alors les sous-ensembles
 
$$
 S_j=\{m\in U'\cap U_0;\  |\chi^j g(m)|=\max_{i=0}^\ell |\chi^i g(m)| \}\setminus S
$$

\noindent Si $U'$ est suffisament petit, une orbite $\gamma$ de $\chi$ dans $U'\cap U_0$ rencontre $S_j$ en un nombre fini d'intervalles $\sigma$ de $\gamma$ (dont certains peuvent être réduits à un point), et ce nombre est uniform\'ement major\'e par un entier $n_j$ qui ne d\'epend que de $S_j$. En effet, les extr\'emit\'es $m$ de ces intervalles $\sigma$ qui appartiennent à $\gamma$ satisfont des \'equations du type

$$
 g_{i,j}^{\pm}(m)=\chi^j g(m)\pm \chi^i g(m)=0 \quad i\neq j
$$

\noindent et ces germes $g_{i,j}^{\pm}$ (en nombre fini) sont des éléments de l'algèbre ${\Cal B}_0$ qui est $\chi$-stable et $\chi$-finie sur $U_0$. Si on note $S_{i,j}^{\pm}$ le représentant de $Z(g_{i,j}^{\pm})$ sur $U'$, on peut donc prendre

$$
 n_j=\sum_{i\neq j} d^°\pi_{\chi,U'|S_{i,j}^{\pm}\cap U_0}
$$

\noindent Maintenant, sur un de ces intervalles $\sigma$ (non réduit à un point!), considérons la fonction $f_\sigma(t)=f\circ \varphi_{\chi,U'}(t,m)$ où $m$ est un point de $\sigma$ et $t$ décrit un intervalle $\tau$ de $\RR$ tel que $\varphi_{\chi,U'}(\tau,m)=\sigma$. Un calcul direct montre que les dérivvées successives de $f_\sigma$ sont données par

$$
f_\sigma^{(i)}(t)=\chi^i f\circ \varphi_{\chi,U'}(t,m)
$$

\noindent Et la relation (2) dérivée $j$ fois par rapport à $\chi$ et restreinte à $U_0\cap\sigma$ donne

$$
(\chi^j f)_{|U_0\cap\sigma}=(\chi^j g)_{|U_0\sigma}(1+O(m))
$$

\noindent Ceci montre que la dérivée $j$-ème de $f_\sigma$ ne s'annule pas sur l'intervalle $\tau$. Par conséquent et par le Lemme de Rolle appliqué à la fonction $f_\sigma$ sur l'intervalle $\tau$, on a

$$
 d^°\pi_{\chi,U'|S'\cap U_0}\leq \sum_{j=0}^\ell jn_j\leq \ell\sum_{j=0}^\ell n_j
$$

\noindent où $S'$ est le représentant du germe $Z(f)$ sur $U'$. \qed\par

\vskip 3mm

  L'hypothèse ${\Cal M}_0\subset{\Cal M}$ du lemme ne peut être affaiblie, comme le montre l'exem-\par\noindent ple suivant: $\chi=x\partial/\partial x$, $g=x$, $h=-1+\sin(1/x)\exp(-1/x)$ et $f=g+hg$; ${\Cal B}$ étant un anneau $\chi$-stable contenant l'anneau $\RR\{x\}$ et les dérivées successives $\chi^n h$, et ${\Cal B}_0=\RR\{x\}$. La "bonne préparation" du lemme  est à rapprocher de celle que l'on rencontre dans l'étude des singularités d'applications différentiables, et qui utilise uniquement "l'idéal jacobien". Ce qui consiste en fait à regarder l'action simultanée de plusieurs dérivations. Il faut noter que les algèbres étudiées dans ce travail {\bf ne sont pas différentiables en 0}.\par
  
\vskip 3mm
  
  Notons ${\Cal C}_{\chi,{\Cal B}}=\{f\in{\Cal B};\ I_{\chi,f} \text{ noethérien  dans } {\Cal B}\}$. Soit la relation sur ${\Cal C}_{\chi,{\Cal B}}$:  $f{\Cal R} g \Leftrightarrow$ il existe un idéal ${\Cal M}_0\subset {\Cal M}$ $\chi$-stable tel que $f-g\in {\Cal M}_0 I_{\chi,g}$. D'après la preuve précédente, c'est une relation d'équivalence sur ${\Cal C}_{\chi,{\Cal B}}$ dite $\chi${\bf -équivalence dans} ${\Cal B}$. Soient ${\Cal C}_1$ et ${\Cal C}_2$ des algèbres ou des classes telles que ${\Cal C}_1\subset{\Cal C}_2\subset{\Cal C}_{\chi,{\Cal B}}$. On dit que ${\Cal C}_2$ est $\chi${\bf -équivalente à} ${\Cal C}_1$ {\bf dans} ${\Cal B}$ si tout élément de ${\Cal C}_2$ est $\chi$-équivalent à un élément de ${\Cal C}_1$ dans ${\Cal B}$. Ainsi, si ${\Cal C}_1$ est une algèbre $\chi$-stable et $\chi$-finie, le lemme de finitude dit que la classe ${\Cal C}_2\subset\cup_{{\Cal M}_0} {\Cal C}_{{\Cal M}_0,{\Cal C}_1}$ est $\chi$-finie (l'union étant prise sur tous les idéaux ${\Cal M}_0\subset{\Cal M}$ et qui sont $\chi$-stables).\par

\vskip 3mm

\noindent{\bf Construction d'algèbre} $\chi${\bf -finie.} Il en existe un grand nombre d'après la riche littérature sur les propriétés de finitude des sous-ensembles analytiques avec de fortes conditions au bord ([T]), les sous-ensembles pfaffiens et les constructions dérivées de la théorie de l'o-minimalité ([W], [D-M-M], [L-R], [L-S], [D-S], [S]...); malheureusement, la question de la noethérianité est rarement étudiée dans ces derniers travaux. Un exemple simple de construction, déduit des idées originales de Khovanski-Tougeron, est illustré dans la section II par l'étude des sous-algèbres convergentes $QR{\Cal H}_{cvg}^{p,q}$. Il s'appuie sur les idées suivantes\par

\vskip 3mm

\proclaim{Lemme IB2 (lemme d'extension de Tougeron)} Soit ${\Cal B}_0$ une sous-alg\`ebre de ${\Cal B}$ stable par $\chi$. On suppose que\par

\roster

\vskip 2mm

\item"$(i)$" les semi-analytiques décrits par ${\Cal B}_0$ ont un nombre fini de composantes connexes.\par

\vskip 2mm

\item"$(ii)$" l'anneau ${\Cal B}_0$ est noeth\'erien dans ${\Cal B}$.\par

\vskip 2mm

\item"$(iii)$" Il existe $(p+q-1)$ 1-formes $\omega_j=\sum_{i=1}^{p+q} a_{i,j}dy_i$ avec $a_{i,j}\in{\Cal B}_0$, et un ouvert admissible $U$  telles que toute orbite de $\chi$ dans $U$ soit une intersection transverse de solutions séparantes des $\omega_j$.\par

\endroster

Alors, l'alg\`ebre ${\Cal B}_0$ est $\chi$-finie.\endproclaim

\vskip 3mm

\noindent {\bf Preuve.} L'hypothèse $(i)$ implique en particulier que les idéaux différentiels de ${\Cal B}_0$ sont noethériens dans ${\Cal B}$. Les hypothèses $(i)$ et $(ii)$ impliquent que l'algèbre ${\Cal B}_0$ est topologiquement noethérienne (cf. [T]). Si $S$ est un semi-analytique de ${\Cal B}_0$ et $\gamma$ une orbite r\'eguli\`ere de $\chi$ dans $U$

$$
 b_0(\gamma\cap S)=b_0( (\cap_{j=1}^{p+q-1} \Gamma_j)\cap S)
$$

\noindent où chaque $\Gamma_j$ est une solution séparante de $\omega_j$. Par l'hypothèse $(iii)$, ces $1$-formes sont à coefficients dans ${\Cal B}_0$. Donc, par la th\'eorie de Khovanski -Tougeron ([K1], [K2], [T]), ce nombre est majoré par un entier qui ne dépend que de $S$ et des $\omega_j$. Et donc, tout $f\in{\Cal B}_0$ est $\chi$-régulier.\qed

\vskip 3mm

 Si de plus, ${\Cal B}_0$ est un anneau de référence, on peut simplifier l'hypothèse $(iii)$ comme ceci: la dérivation $\chi$ admet $(p+q-1)$ intégrales premières indépendantes $F_j\in{\Cal B}_0$. Cette hypothèse ne peut être affaiblie, commme le montre l'exemple suivant: prenons $p=1$ et $q=2$, soient ${\Cal B}_0=\RR\{x,\alpha_1,\alpha_2\}$ et 
 
 $$
 \chi=-x\frac{\partial}{\partial x}+(a\alpha_1-b\alpha_2)\frac{\partial}{\partial\alpha_1}+(b\alpha_1+a\alpha_2)\frac{\partial}{\partial\alpha_2}
 $$
 
 \noindent avec $a>0$ et $b\neq 0$. Il est clair que l'anneau ${\Cal B}_0\subset{\Cal A}^{1,2}(x,\alpha)$ est un anneau de référence qui est $\chi$-stable et qui satifait aux hypothèses $(i)$ et $(ii)$. Cependant, le germe $f=\alpha_1\in{\Cal B}_0$ n'est pas $\chi$-régulier. Et le lemme I2 montre que la dérivation $\chi$ n'admet pas de paire d'intégrales premières indépendantes, dans aucun anneau de référence qui satisfait aux hypothèses $(i)$ et $(ii)$. Un calcul direct montre qu'il existe une paire d'intégrales premières indépendantes qui sont linéaires en $\alpha$ et dont les coefficients sont les fonctions oscillantes $x^a\cos(b\log(x))$ et $x^a\sin(b\log(x))$, qui appartiennent à l'anneau ${\Cal A}^{1,0}(x)$ mais pas à l'anneau $SB^{1,0}(x)$.\par

\vskip 3mm

\noindent {\bf Faisceau différentiel.}\par

\vskip 3mm

 Soient $f\in{\Cal B}$ et $\chi\in\Xi{\Cal B}$. Soit $U$ un ouvert admissible pour $\chi$ sur lequel $f$ est réalisée. On note ${\Cal I}_{\chi,f}[U]$ {\bf le faisceau différentiel} de $f$ sur $U$: sa {\bf fibre} ${\Cal I}_{\chi,f,m}$ (qu'on note aussi $I_{\chi,f}(m)$), en un point $m\in U$, est {\bf l'idéal différentiel} $I_{\chi_m,f_m}\subset {\Cal B}_m={\bold \RR}\{y-y_m \}$; $\chi_m$ et $f_m$ étant les germes de $\chi$ et $f$ en $m$, et $y-y_m$ sont les coordonnées locales en $m$ (les représentants d'éléments de $\RR\{y-y_m\}$ étant considérés sur des ouverts connexes). On suppose qu'en tout point $m$ du bord $\partial_0 U$, il existe un anneau de référence ${\Cal B}_{m}\subset{\Cal A}^{p_m,q_m}(y-y_m)$, qui est stable par $\chi_m$ (ainsi $\chi_m\in\Xi{\Cal B}_m$) et qui contient les germes en $m$ d'éléments de ${\Cal B}$ qui sont représentés sur un voisinage de $m$ ( en général, on choisi ${\Cal B}_m$ comme étant l'anneau de ces germes). Dans ce cas, on prolonge le faisceau ${\Cal I}_{\chi,f}[U]$ par le faisceau ${\Cal I}_{\chi,f}[U\cup\partial_0 U]$, dont les {\bf fibres différentielles} sur le bord $\partial_0 U$ sont les idéaux différentiels $I_{\chi_m,f_m}$ des anneaux  ${\Cal B}_m$ correspondants. Ce faisceau est donc le faisceau associé au préfaisceau ${\Cal F}[U\cup\partial_0 U]$ muni des morphismes restriction naturels, et dont les sections ${\Cal F}(V)$ (où $V$ est un ouvert de $U\cup\partial_0 U$) sont les idéaux engendrés par la famille $(\chi^n f)_{n\in\NN}$ dans l'anneau des fonctions sur $V$ dont le germe en tout point $m\in V$, est un élément de l'anneau local ${\Cal B}_m$.\par

\vskip 3mm

 On montre deux résultats importants concernant ces faisceaux. Conjointement au lemme de finitude, ces résultats (et leurs conséquences) constituent le socle de ce travail.\par

\vskip 3mm

\proclaim{Lemme IB3 (lemme de cohérence)} Le faisceau ${\Cal I}_{\chi,f}[U]$ est {\bf cohérent}: plus précisément, pour tout $m\in U$, il existe un ouvert $V_m\subset U$ contenant $m$, et un entier $\ell_m$ tels que en tout point $m'\in V_m$, la fibre $I_{\chi,f}(m')$ est engendrée par les germes en $m'$ de $(f,\chi f,\ldots,\chi^{\ell_m} f)$. Si de plus, la fibre en $0$ $I_{\chi,f}(0)=I_{\chi,f}$ est noethérienne dans ${\Cal B}$, alors il existe un ouvert $V\in (U,0)$ et contenu dans $U$ tel que le faisceau ${\Cal I}_{\chi,f}[U\cup\partial_0 V]$ soit {\bf cohérent}.\endproclaim

\vskip 3mm

\noindent {\bf Preuve.} Si la fibre en 0 est noethérienne et si $(f$, $\chi f,\ldots$, $\chi^\ell f)$ est un système de générateurs de $I_{\chi,f}$ dans ${\Cal B}$, alors il existe un voisinage ouvert $V$ de $0$ dans $U$ tel que la division

$$
 \chi^{l+1} f=\sum_{i=0}^\ell h_i \chi^i f
$$

\noindent soit réalisée sur un voisinage de $\overline{V}$ dans $U\cup\partial_0 U$. En dérivant plusieurs fois cette égalité, on obtient que pour tout $k$, la division

$$
 \chi^k f=\sum_{i=0}^\ell h_{i,k} \chi^i f
$$

\noindent est réalisée sur ce même voisinage. On obtient le résultat en germifiant ces égalités en n'importe quel point de $V\cup\partial_0 V$. Le reste du lemme s'obtient de la même façon, car pour tout $m\in U$, l'anneau ${\Cal B}_m=\RR\{y-y_m\}$ est noethérien.\qed\par

\vskip 3mm

 Soit $(\chi_i,{\Cal B}_i)_{i=1,2}$ deux couples tels que ${\Cal B}_i$ est un anneau de référence et $\chi_i$ est une dérivation appartenant à $\Xi{\Cal B}_i$. Soit $U_i$ un ouvert admissible pour $\chi_i$ et soit $\varphi: U_1\to U_2$ un difféomorphisme tel que $\varphi_*(\chi_1)=\chi_2$. On suppose que $\varphi$ se prolonge en un homéomorphisme de $U_1\cup \partial_0 U_1$ sur $U_2\cup\partial_0 U_2$ tel que $\varphi(0)=0$. On note de la même façon son germe en 0 $\varphi: (U_1,0)\to (U_2,0)$ (ainsi que celui de $\varphi^{-1}$), et on suppose que $\varphi^*({\Cal B}_2)={\Cal B}_1$ (donc $\varphi^*$ est un isomorphisme et $(\varphi^*)^{-1}=(\varphi^{-1})^*$). Un calcul direct montre que le diagramme suivant est commutatif
 
$$
\CD
{\Cal B}_2 @>\varphi^*>> {\Cal B}_1 \\
@V\chi_2 VV             @VV\chi_1 V  \\
{\Cal B}_2 @>\varphi^*>> {\Cal B}_1
\endCD
$$
 
\noindent Donc, si $f_1\in{\Cal B}_1$ et $f_2=(\varphi^{-1})^*(f_1)$, alors pour tout entier $n$
 
$$
\chi_1^n f_1=\varphi^*(\chi_2^n f_2)
$$

\noindent et par conséquent

$$
I_{\chi_1,f_1}=\varphi^*(I_{\chi_2,f_2})
$$

\noindent Ceci est à fortiori vrai aux voisinages de points $m_1\in U_1$ et $m_2=\varphi(m_1)\in U_2$, les anneaux locaux correspondants étant les anneaux de germes de fonctions analytiques en ces points (on les appelle aussi "anneaux analytiques" pour simplifier). On verra dans la suite des situations plus générales de tels "transfert", dont la première est\par
 
\vskip 3mm

\proclaim{Lemme IB4 (lemme d'isomorphie de Roussarie)} Le faisceau ${\Cal I}_{\chi,f}[U]$ est {\bf compatible avec la projection} $\pi_{\chi,U}$. Autrement dit, si $\gamma$ est une orbite de $\chi$ dans $U$, et si $m_1,m_2\in \gamma$, alors les fibres ${\Cal I}_{\chi,f,m_1}$ et ${\Cal I}_{\chi,f,m_2}$ sont {\bf isomorphes}, l'isomorphisme étant le germe du flot $\varphi_{\chi,U}$ aux voisinages de $m_1$ et $m_2$.\endproclaim

\vskip 3mm

\noindent {\bf Preuve.} Soit $t_0$ tel que $m_2=\varphi_{\chi,U} (t_0,m_1)$. Soient $m'_1\in U$ un point voisin de $m_1$, $t$ voisin de 0 et $m'_2=\varphi_{\chi,U} (t,m'_1)\in U$. Alors, un calcul direct donne

$$
 f(m'_2)=f(\varphi_{\chi,U} (t,m'_1))=\sum_{n\geq 0} \chi^n f(m'_1)\frac{t^n}{n!}
$$

\noindent la série étant uniformément convergente sur le produit d'un voisinage de $m_1$ dans $U$ et d'un voisinage de 0 dans $\RR$ (le flot étant analytique sur ce produit et $f$ est analytique sur un voisinage de $m_1$). Les id\'eaux de germes de fonctions analytiques en un point, sont {\bf ferm\'es pour la topologie de la convergence uniforme} ([B-M], [H]). Par conséquent, ${\Cal I}_{\chi ,f,m_1}\supset (\varphi_{\chi,U} (t,.))^* ({\Cal I}_{\chi,f,\varphi_{\chi,U}(t,m_1)})$, où $\varphi_{\chi,U}(t,.)$ désigne aussi le germe de ce difféomorphisme aux points $m_1$ et $\varphi_{\chi,U}(t,m_1)$. En considérant le flot inverse ($(\varphi_{\chi,U}(t,.))^*$ est un isomorphisme entre les anneaux analytiques correspondants), on obtient  ${\Cal I}_{\chi ,f,m_1}=(\varphi_{\chi,U} (t,.))^* ({\Cal I}_{\chi,f,\varphi_{\chi,U}(t,m_1)})$. Et, en utilisant un recouvrement fini du segment d'orbite $[m_1,m_2]$, on obtient

$$
  {\Cal I}_{\chi ,f,m_1}=(\varphi_{\chi,U} (t_0,.))^* ({\Cal I}_{\chi,f,m_2})
$$

\noindent où $\varphi_{\chi,U}(t_0,.)$ désigne aussi le germe aux points $m_1$ et $m_2$ de ce difféomorphisme.\qed\par

\vskip 3mm

  Soit $U\in ((\RR^{+*})^p\times\RR^q,0)$ et soit $c:m\in U\mapsto m'\in \RR^n$ un morphisme continu sur $U\cup\partial_0 U$ tel que $c(0)=0$. On note de la même façon son germe $c: (U,0)\to (\RR^n,0)$. Soit $u$ une coordonnée locale sur $(\RR^n,0)$. L'anneau ${\Cal B}=c^*(\RR\{ u\})$ est dit {\bf anneau restriction analytique} (ou encore {\bf anneau convergent}). Il est isomorphe à l'anneau restriction de référence $\RR\{u\}_{|c(U)}$. Il est clair que c'est un anneau local et noethérien (l'image réciproque par $c^*$ de tout idéal est un idéal, et l'anneau $\RR\{u\}$ est noethérien). Soit $y=(x,\alpha)$ des coordonnées sur $U$, on suppose maintenant que le morphisme $c$ est analytique, et qu'il s'écrit dans les coordonnées $y$ et $u$
  
$$
u=c(y)=(y,\psi(y))
$$

\noindent On suppose aussi que pour tout $i=1,\ldots,n-(p+q),\ j=1,\ldots,p+q$, il existe $h_{i,j}\in\RR\{u\}$ tel que

$$
\prod_{k=1}^p x_k \frac{\partial \psi_i}{\partial y_j}=c^*(h_{i,j})
$$

\noindent Dans ce cas, on vérifie facilement que l'anneau convergent ${\Cal B}$ est un anneau de référence. Soit 

$$
\chi=\sum_{j=1}^{p+q} a_j\frac{\partial}{\partial y_j}\in\Xi{\Cal B}
$$

\noindent Pour tout $i=1,\ldots,n$, le germe $g_i(u)=u_i\in\RR\{u\}$, donc $c^*(g_i)\in{\Cal B}$ et par conséquent $\chi c^*(g_i)\in{\Cal B}$ (en particulier, pour $i=1,\ldots,p+q$, on a $\chi c^*(g_i)=a_i$), il existe alors $h_i\in\RR\{u\}$ (il n'est pas unique en général!) tel que $\chi c^*(g_i)=c^*(h_i)$. Soit $U'\subset U$ un ouvert admissible pour $\chi$ tel que les germes $h_i$ soient réalisés sur un ouvert ${\Cal U}\in(\RR^n,0)$ et contenant $c(U')$. Soit 

$$
{\Cal X}=\sum_{i=1}^n h_i\frac{\partial}{\partial u_i}
$$

\noindent c'est une dérivation analytique réalisée sur ${\Cal U}$. Le morphisme $c$ est un difféomorphisme de $U'$ sur la variété analytique $c(U')$ ($c$ est une immersion), et par la construction de ${\Cal X}$, on a 

$$
c_*\chi={\Cal X}_{|c(U')}
\tag *
$$

\noindent (on note encore de la même façon le difféomorphisme associé à $c$, et son germe en 0). Ainsi, le diagramme suivant est commutatif

$$
\CD
\RR\{u\} @>c^*>> {\Cal B} \\
@V{\Cal X} VV             @VV\chi V  \\
\RR\{u\} @>c^*>> {\Cal B}
\endCD
$$

 Le difféomorphisme $c$ se prolonge en un homéomorphisme de $U'\cup\partial_0 U'$ sur son image. Soit $m\in \partial_0 U'$ et $m'=c(m)$, on note $c_m: (U',m)\to (\RR^n,m')$ le germe de $c$ aux points $m$ et $m'$. Soit $\RR\{u-u_{m'}\}$ l'anneau analytique local au point $m'$, et soit ${\Cal B}_m=c_m^*(\RR\{u-u_{m'}\})$, c'est un anneau restriction analytique qui est un anneau de référence (par les mêmes raisonnements que ci-dessus, en supposant bien sûr que les germes $h_{i,j}$ sont réalisés sur ${\Cal U}$). Il est clair qu'il contient les germes en $m$ d'éléments de ${\Cal B}$ représentés sur un voisinage de $m$. De plus, il est stable par la dérivation $\chi_m$ (par la commutativité du digramme obtenu à partir de la relation $(*)$ germifiée aux points $m$ et $m'$).\par
 
\vskip 3mm
  
  Ainsi, pour tout $f\in{\Cal B}$ réalisée sur un certain ouvert $V\in (U,0)$ (et $V\subset U'$), le lemme de cohérence dit que le faisceau ${\Cal I}_{\chi,f}[V\cup\partial_0 V]$ est cohérent. Par la continuité de la dérivation $\chi$ et du morphisme $c$ sur le bord $\partial_0 U'$, l'image par $c$ d'une orbite régulière de $\chi$ incluse dans le bord, est une partie connexe d'une orbite régulière de ${\Cal X}$ dans ${\Cal U}$. Le lemme d'isomorphie s'applique aux faisceaux différentiels ${\Cal I}_{{\Cal X},g}[{\Cal U}_g]$ le long de toute orbite régulière (${\Cal U}_g\subset {\Cal U}$ étant un un voisinage ouvert de 0, sur lequel $g$ est représentée). Ainsi, en transportant par le morphisme $c^*$, on obtient que le lemme d'isomorphie s'apllique aussi aux faisceaux ${\Cal I}_{\chi,f}[V\cup\partial_0 V]$ le long {\bf d'orbites régulières incluses dans le bord}. En général dans ce travail, le morphisme $c$ est une immersion dont l'image est le graphe de fonctions élémentaires (et simples!) de Khovanski. Un exemple, déjà rencontré, de tels anneaux convergents, est l'anneau $QR{\Cal H}^{p,q}_{cvg}$ (pour lequel le morphisme $c^*$ est un isomorphisme (cf. Section II)).\par

\vskip 3mm

 Si l'une des dérivations ${\Cal X}$ satisfait à l'hypothèse $(iii)$ du lemme d'extension (on bien si elle admet $(n-1)$ intégrales premières analytiques et indépendantes), alors ce lemme dit que l'algèbre $\RR\{u\}$ est ${\Cal X}$-finie. Il est clair, dans ce cas, que l'anneau convergent ${\Cal B}$ est aussi $\chi$-fini: si $f\in{\Cal B}$ et si $g\in\RR\{u\}$ est tel que $f=c^*(g)$, alors par transport par l'injection $c$
 
$$
d^°\pi_{\chi |Z(f)}\leq d^°\pi_{{\Cal X}|Z(g)}
$$

\vskip 3mm

\noindent {\bf Idéal} $\chi${\bf -transverse.}\par

\vskip 3mm

 Soit ${\Cal B}$ un anneau de référence. Soient $f\in{\Cal B}$ et $\chi\in\Xi{\Cal B}$. Soit $U$ un ouvert admissible sur lequel $f$ est réalisée. Soit $\gamma$ une orbite régulière de $\chi$ dans $U$. Soient $m,m'\in\gamma$ et $t_0>0$ tel que $\varphi_{\chi,U}(t_0,m)=m'$. Soit $\sigma_m\subset U$ une transversale à $\gamma$ en $m$, analytique. Quitte à réduire $\sigma_m$, on suppose que pour $t$ dans un voisinage de $[0,t_0]$, l'image de $\sigma_m$ par le flot $\varphi_{\chi,U}(t,.)$ est incluse dans $U$. Ainsi, $\sigma_{m'}=\varphi_{\chi,U}(t_0,\sigma_m)$ est une transversale à $\gamma$ en $m'$, analytique. Soient $\beta$ des coordonnées analytiques locales sur $(\sigma_m,m)$, il est clair que l'anneau restriction de référence $\RR\{y-y_m\}_{|\sigma_m}$ est isomorphe à l'anneau analytique $\RR\{\beta\}$ (qui lui même est isomorphe à l'anneau des intégrales premières analytiques locales en $m$ (cf. ci-dessous)). Par le difféomorphisme analytique $\varphi_{\chi,U}(t_0,.)$, on peut prendre $\beta$ comme coordonnées analytiques locales sur $(\sigma_{m'},m')$. Par cette identification, le lemme d'isomorphie dit alors que
$i^*_{\sigma_m}({\Cal I}_{\chi,f,m})=i^*_{\sigma_{m'}}({\Cal I}_{\chi,f,m'})$, d'où\par

\vskip 3mm 

\proclaim{Définition IB6} On appelle idéal $\chi$-{\bf transverse} de $f$ le long de $\gamma$ l'idéal restriction

$$
 J_{\chi,f,\gamma}=i^*_{\sigma_m}({\Cal I}_{\chi,f,m})={\Cal I}_{\chi,f,m|\sigma_m}\subset\RR\{\beta\}\quad \forall m\in\gamma
$$

\noindent où les coordonnées analytiques $\beta$ sur $\sigma_m$ sont des intégrales premières de $\chi$ le long de $\gamma$.
\endproclaim

\vskip 3mm

\noindent Et inversement, on peut reconstruire toute fibre ${\Cal I}_{\chi,f,m}$ le long de $\gamma$ à partir de l'idéal transverse  $J_{\chi,f,\gamma}$: soit $U_m\subset U$ le saturé de $\sigma_m$ par le flot $\varphi_{\chi,U}(t,.)$ pour $t$ dans un voisinage de 0. Il est clair que l'espace quotient $\widetilde{U_m}=U_m/\varphi_{\chi_m,U_m}$ s'identifie à la transversale analytique $\sigma_m$. Ainsi, on peut parler du germe de la projection intégrale $\pi_{\chi_m,U_m}: U_m\to \sigma_m$, qu'on note simplement $\pi_{\chi_m}: (U_m,m)\to (\sigma_m,m)$. On lui associe le morphisme étoilé $\pi_{\chi_m}^*: \RR\{\beta\}\to \RR\{y-y_m\}$ qui est injectif; son image est l'anneau des intégrales premières analytiques locales en $m$. Avant de poursuivre, faisons une remarque simple qui sera souvent utilisée dans la suite: soient $\psi_1^*,\psi_2^*$ deux morphismes étoilés, et soit $I$ un idéal, alors l'idéal prolongé (par $\psi_1^*$) associé à l'idéal prolongé (par $\psi_2^*$) associé à $I$, coincide avec l'idéal prolongé (par $(\psi_2\circ\psi_1)^*$) associé à $I$, ce que l'on note $(\psi_2\circ\psi_1)^*(I)=\psi_1^*\circ\psi_2^*(I)$.\par

\vskip 3mm

\proclaim{Lemme IB5 (lemme de saturation)} Pour tout $m\in \gamma$, on a $\quad {\Cal I}_{\chi,f,m}=\pi^*_{\chi_m}(J_{\chi,f,\gamma})$.
\endproclaim

\vskip 3mm

\noindent {\bf Preuve.} On se place dans un flow-box au voisinage de $m$. Quitte à réduire $U_m$, soit $\psi: (t,\beta)\in \tau\times \Sigma_m \mapsto m'\in U_m$ le difféomorphisme normalisant tel que $\psi(\{0\}\times\Sigma_m)=\sigma_m$ et $(\psi^{-1})_*(\chi_m)={\Cal Y}=\partial/\partial t$. Notons de la même façon son germe $\psi:(\tau\times\Sigma,(0,0))\to (U_m,m)$. On a le diagramme commutatif

$$
\CD
\RR\{y-y_m\} @>\psi^*>> \RR\{t,\beta\} \\
@V\chi_m VV             @VV{\Cal Y} V  \\
\RR\{y-y_m\} @>\psi^*>> \RR\{t,\beta\}
\endCD
\tag {$*$}
$$

\noindent où $\psi^*$ est un isomorphisme. L'orbite $\gamma\cap U_m$ est envoyée sur l'orbite $\Gamma=\tau\times \{ 0\}$. Notons $\psi_1:\Sigma_m\to\sigma_m$ (et son germe aux points 0 et $m$) la restriction de $\psi$ à $\{0\}\times\Sigma_m$. C'est un difféomorphisme. La projection intégrale $\pi_{{\Cal Y},\tau\times\Sigma_m}$ s'identifie avec la projection canonique $\pi: \tau\times \Sigma_m \to \Sigma_m$; notons de la même façon son germe $\pi:(\tau\times\Sigma_m,(0,0))\to (\Sigma_m,0)$. Les diagrammes suivants est commutatifs

$$
\CD
(\tau\times\Sigma_m,(0,0)) @>\psi>> (U_m,m) \\
@V\pi VV             @VV\pi_{\chi_m} V  \\
(\Sigma_m,0) @>\psi_1>> (\sigma_m,m)
\endCD
\tag {$**1$}
$$

$$
\CD
(\Sigma_m,0) @>\psi_1>> (\sigma_m,m) \\
@V i_{\Sigma_m} VV             @VV i_{\sigma_m} V  \\
(\tau\times\Sigma_m,(0,0)) @>\psi>> U_m
\endCD
\tag {$**2$}
$$

\noindent Par le choix des coordonnées $\beta$ sur $\Sigma_m$ et $\sigma_m$, le morphisme étoilé $\psi_1^*$ est l'identité de l'anneau $\RR\{\beta\}$; on a donc les diagrammes commutatifs

$$
\CD
\RR\{\beta\} @>id>> \RR\{\beta\} \\
@V\pi_{\chi_m}^* VV             @VV\pi^* V  \\
\RR\{y-y_m\} @>\psi^*>> \RR\{t,\beta\}
\endCD
\tag {$***1$}
$$

$$
\CD
\RR\{y-y_m\} @>\psi^*>> \RR\{t,\beta\} \\
@V i_{\sigma_m}^* VV             @VV i_{\Sigma_m}^* V  \\
\RR\{\beta\} @>id>> \RR\{\beta\}
\endCD
\tag {$***2$}
$$

\noindent Ces deux derniers diagrammes sont aussi valables pour les idéaux prolongés. Soit $F=\psi^*(f_m)=\sum a_n(\beta) t^n$, d'après le diagramme $(*)$, on a $\psi^*({\Cal I}_{\chi,f,m})=I_{\Cal Y,F}\subset \RR\{t,\beta\}$. D'après la définition de l'idéal transverse

$$
 J_{{\Cal Y},F,\Gamma}=i_{\Sigma_m})^*(I_{{\Cal Y},F})=i_{\Sigma_m}^*\circ\psi^*({\Cal I}_{\chi,f,m})
$$

\noindent Et d'après le diagramme $(***2)$, on a $J_{{\Cal Y},F,\Gamma}=i_{\sigma_m}^*({\Cal I}_{\chi,f,m})=J_{\chi,f,\gamma}$. Donc, d'après le diagramme $(***1)$, il suffit de montrer que $I_{\Cal Y,F}=\pi^*(J_{{\Cal Y},F,\Gamma})$. Par sa définition, on vérifie facilement que $J_{{\Cal Y},F,\Gamma}=\langle a_n;\ n\in \NN \rangle$ (d'où le nom d'{\bf idéal des coefficients} attribué par Roussarie à cet idéal). Soit ${\Cal M}$ l'idéal maximal de $\RR \{ t,\beta \}$. Notons $F_n={\bold j}^n_t(F)$. On a $F_n\in \pi^*(J_{\Cal Y,F,\Gamma})$ et $F-F_n\in {\Cal M}^n$ pour tout $n$. Donc, par le théorème d'intersection de Krull ([L]), $F\in \pi^*(J_{\Cal Y,F,\Gamma})$. Comme l'idéal prolongé $\pi^*(J_{{\Cal Y},F,\Gamma})$ est stable par la dérivation ${\Cal Y}$, on obtient $I_{{\Cal Y},F}\subset \pi^*(J_{{\Cal Y},F,\Gamma}$. Inversement, pour tout $n\in\NN$, on a ${\Cal Y}^n F=n!a_n+tH_n$ avec $H_n\in\pi^*(J_{{\Cal Y},F,\Gamma})$ (par le même raisonnement que pour $F$, en utilisant la série de $F$). Donc, pour tout $n$, $a_n(=\pi^*(a_n))\in I_{{\Cal Y},F}+{\Cal M}\pi^*(J_{{\Cal Y},F,\Gamma})$. Par conséquent

$$
\pi^*(J_{{\Cal Y},F,\Gamma})\subset I_{{\Cal Y},F}+{\Cal M}\pi^*(J_{{\Cal Y},F,\Gamma})
$$

\noindent Et, par le lemme de Nakayama ([L]), on obtient $\pi^*(J_{{\Cal Y},F,\Gamma})\subset I_{{\Cal Y},F}$.\qed\par

\comment
Une autre démonstration, qu'on verra dans la suite...
notons $J_n=\langle a_0,\ldots,a_n \rangle\subset \RR\{\beta\}$. Il est clair que $I_{{\Cal Y},F_n}\subset \pi^*(J_n)$. De l'égalité ${\Cal Y}^n F_n=n! a_n$ on déduit deux choses: d'abord  que $I_{\Cal Y,F_{n-1}} \subset I_{\Cal Y,F_{n}}$ (car $a_n\in I_{{\Cal Y},F_n}\Rightarrow a_n t^n\in I_{{\Cal Y},F_n}\Rightarrow F_n-a_nt^n=F_{n-1}\in I_{{\Cal Y},F_n}$); puis par une récurrence sur $n$, que $\pi^*(J_n)\subset I_{{\Cal Y},F_n}$, et donc que $\pi^*(J_n)=I_{{\Cal Y},F_n}$. Ces deux suites d'idéaux sont croissantes, et la première converge vers l'idéal $\pi^*(J_{{\Cal Y},F,\Gamma})$
\noindent La croissance des suites d'idéaux $(J_n)$ et $(I_{\Cal Y,F_n})$ et une deuxième application du théorème d'intersection de Krull fournissent le résultat.\qed\par
\endcomment

\vskip 3mm

\noindent {\bf Orbite principale et la double inclusion}.\par

\vskip 3mm

 Si l'orbite $\gamma$ adhère à $0$, et s'il existe une notion de germe de la projection intégrale $\pi_{\chi,U}$ en 0 (qu'on notera $\pi_\chi$), il se pose alors la question de comparer la fibre en 0 $I_{\chi,f}$ et l'idéal $\pi^*_{\chi}(J_{\chi,f,\gamma})$ (dans une extension commune aux anneaux $\pi_\chi^*(\RR\{\beta\})$ et ${\Cal B}$). Intuitivement, il est tout à fait légitime d'espérer établir {\bf la double inclusion}

$$
 (\prod_{j=1}^p x_j)^N \pi^*_{\chi}(J_{\chi,f,\gamma})\subset I_{\chi,f}\subset \pi^*_{\chi}(J_{\chi,f,\gamma})
$$

\noindent qui relaxe l'égalité du lemme de saturation, et qui s'inspire du Nullstelenzats d'Hilbert ([L]). Dans ce cas, il est nécessaire d'avoir l'égalité 

$$
Z(I_{\chi,f})=\pi^{-1}_{\chi}(Z(J_{\chi,f,\gamma}))
$$

\noindent (cf. formule (1)). Et pour cela, il suffit que $\gamma$ soit "principale" dans $U$\par

\vskip 3mm

\proclaim{Définition IB7} Soit $\gamma$ une orbite de $\chi$ dans $U$. Elle est dite {\bf principale} dans $U$ si 

\roster

\vskip 2mm

\item"$(i)$" elle adhère à 0;\par

\vskip 2mm

\item"$(ii)$" elle admet une transversale analytique $\sigma_0\subset U$ qui rencontre en au plus un point chaque orbite de $\chi$ dans $U$;\par

\vskip 2mm

\item"$(iii)$" pour toute  transversale analytique $\sigma\subset\sigma_0$, le saturé 

$$
\varphi_{\chi,U}(.,\sigma)=\pi_{\chi,U}^{-1}(\pi_{\chi,U}(\sigma))
$$ 

\noindent est un voisinage de $0$ dans $U$.\par

\endroster

\endproclaim

\vskip 3mm

 Ainsi, l'orbite principale $\gamma$ est la seule orbite de $\chi$ dans $U$ qui adhère à 0. C'est donc la seule orbite principale de $\chi$ dans $U$. Pour $\sigma\subset\sigma_0$, notons $U_\sigma=\varphi_{\chi,U}(.,\sigma)$, c'est un voisinage ouvert de 0, qui est une union d'orbites de $\chi$ dans $U$ (c'est un ouvert admissible). L'espace quotient $\widetilde{U_\sigma}$ s'identifie donc à un sous-espace de $\widetilde{U}$. Par conséquent, le morphisme $\pi_{\chi,U_\sigma}: U_\sigma\to\widetilde{U_\sigma}$ est simplement la restriction à $U_\sigma$ du morphisme $\pi_{\chi,U}$. Par la définition de $U_\sigma$, et par la condition $(ii)$, l'espace $\widetilde{U_\sigma}$ s'identifie à la transversale analytique $\sigma$. Par la condition $(iii)$, le morphisme $\pi_{\chi,U_\sigma}$ se prolonge continûment à $U_\sigma\cup\{0\}$ en posant $\pi_{\chi,U_\sigma}(0)=\gamma\cap\sigma$. On peut donc parler du germe en 0 du morphisme $\pi_{\chi,U}$ comme étant celui du morphisme $\pi_{\chi,U_\sigma}: U_\sigma\to\sigma$ aux points 0 et $\gamma\cap\sigma$. Ce germe est indépendant de la transversale $\sigma\subset\sigma_0$; on le note $\pi_\chi$. L'existence d'une orbite principale $\gamma$ implique en particulier l'existence de $(p+q-1)$ intégrales premières analytiques et indépendantes: si $\beta$ sont des coordonnées analytiques sur $\sigma_0$, et si on pose $g_j(\beta)=\beta_j$, les germes $G_j=\pi_\chi^*(g_j)$ admettent des représentants $F_j$ sur $U_{\sigma_0}$, qui sont des intégrales premières de $\chi$, analytiques et indépendantes.\par
 
\vskip 3mm

\noindent{\bf Autre cas de transfert.}\par

\vskip 3mm

  On en a vu beaucoup, et on verra beaucoup. La première apparition de chaque nouveau cas sera traîtée en détail. Un cas assez général, et abondamment employé dans la suite est le suivant: pour $i=1,2$, soit ${\Cal B}_i\subset{\Cal A}^{p_i,q_i}$ un anneau de référence. Soit $\chi_i\in\Xi{\Cal B}_i$ une dérivation, et $U_i$ un ouvert admissible pour $\chi_i$. On suppose que chaque dérivation $\chi_i$ admet une orbite principale $\gamma_i$ dans $U_i$. Pour simplifier, supposons que $U_i$ est le saturé d'une transversale analytique $\sigma_i$, comme ci-dessus (quitte à réduire $U_i$). Soit $W_i\subset\sigma_i$ une sous-variété analytique contenant le point $m_i=\gamma_i\cap\sigma_i$, et soit $V_i=\pi_{\chi,U_i}^{-1}(W_i)$, c'est une sous-variété analytique de $U_i$, qui adhère à 0, et qui est invariante par $\chi_i$.\par

\vskip 3mm

  Notons ${\Cal Y}_i=\chi_{i|(V_i,0)}$. En utilisant le flot de ${\Cal Y}_i$, on montre facilement que ${\Cal Y}_i({\Cal B}_{i|V_i})\subset{\Cal B}_{i|V_i}$. Soit $\Psi: V_2\to V_1$ un difféomorphisme analytique, continue en 0 avec $\Psi(0)=0$. On note de la même façon son germe $\Psi: (V_2,0)\to (V_1,0)$. On suppose que $\Psi_*{\Cal Y}_2={\Cal Y}_1$, et que $\Psi^*({\Cal B}_{1|V_1})\hookrightarrow {\Cal B}_{2|V_2}$. Dans ce cas, le diagramme suivant est commutatif

$$
\CD
{\Cal B}_{1|V_1} @>\Psi^*>> {\Cal B}_{2|V_2} \\
@V{\Cal Y}_1 VV             @VV{\Cal Y}_2 V  \\
{\Cal B}_{1|V_1} @>\Psi^*>> {\Cal B}_{2|V_2}
\endCD
\tag {$*$}
$$

  Une orbite principale étant la seule qui adhère à 0, et $\Psi$ étant continue en 0, on a $\Psi(\gamma_2)=\gamma_1$. Comme $\Psi$ est un difféomorphisme qui préserve les orbites, on peut supposer, pour simplifier que $\Psi(W_2)=W_1$ (quitte à déplacer les transversales $\sigma_i$). Notons $\psi: W_2\to W_1$ la restriction de $\Psi$ à $W_2$. C'est un difféomorphisme analytique, et on note de la même façon son germe $\psi: (W_2,m_2)\to (W_1,m_1)$. Soit $V_{m_i}$ un saturé suffisament petit, de $W_i$ par le flot de ${\Cal Y}_i$, et soit $\pi_{{\Cal Y}_i,m_i}$ le germe en $m_i$ de la projection intégrale. Le diagramme suivant est commutatif
  
$$
\CD
(V_{m_2},m_2) @>\Psi_{m_2}>> (V_{m_1},m_1) \\
@V\pi_{{Cal Y}_2,m_2} VV             @VV\pi_{{\Cal Y}_1,m_1} V  \\
(W_2,m_2) @>\psi>> (W_1,m_1)
\endCD
\tag {$**$}
$$
  
  La variété $W_i$ est analytique en $m_i$. Donc, par un bon choix de coordonnées analytiques locales $(\beta_i,\beta'_i)$ sur $(\sigma_i,m_i)$ (dans lesquelles $(W_i,m_i)$ est un graphe), on obtient que l'anneau local $\RR\{\beta_i,\beta'_i\}_{|W_i}$ sur $(W_i,m_i)$, est isomorphe à l'anneau analytique $\RR\{\beta_i\}$. On obtient alors le diagramme commutatif à partir du diagramme ci-dessus
  
$$
\CD
\RR\{\beta_1\} @>\psi^*>> \RR\{\beta_2\} \\
@V\pi_{{\Cal Y}_1,m_1}^* VV             @VV\pi_{{\Cal Y}_2,m_2}^* V  \\
{\Cal B}_{1,m_1|V_{m_1}} @>\Psi_{m_2}^*>> {\Cal B}_{2,m_2|V_{m_2}}
\endCD
\tag {$***$}
$$
  
\noindent Les variétés $V_i$ étant analytiques en $m_i$, les anneaux locaux ${\Cal B}_{i,m_i|V_{m_i}}$ sont aussi isomorphes à des anneaux analytiques.\par

\vskip 3mm

\proclaim{Lemme IB6 (lemme de transfert)} Soit $f_1\in{\Cal B}_1$ et $f_2\in{\Cal B}_2$ tel que $f_{2|V_2}=\Psi^*(f_{1|V_1})$. On a les égalités suivantes\par

\roster

\vskip 2mm

\item"$(a)$" $J_{\chi_i,f_i,\gamma_i|W_i}=i_{W_i,V_{i,m_i}}^*({\Cal I}_{\chi_i,f_i,m_i|V_{i,m_i}})$.\par

\vskip 2mm

\item"$(b)$" ${\Cal I}_{\chi_i,f_i,m_i|V_{i,m_i}}=\pi_{{\Cal Y}_i,m_i}^*(J_{\chi_i,f_i,\gamma_i|W_i})$.\par

\vskip 2mm

\item"$(c)$" $I_{\chi_2,f_2|V_2}=\Psi^*(I_{\chi_1,f_1|V_1})$.\par

\vskip 2mm

\item"$(d)$" $J_{\chi_2,f_2,\gamma_2|W_2}=\psi^*(J_{\chi_1,f_1,\gamma_1|W_1})$.\par

\endroster

\endproclaim

\vskip 3mm

\noindent{\bf Preuve.} Enlevons l'indice $i$ pour un moment. Pour l'égalité $(a)$, par définition d'une restriction, on a $J_{\chi,f,\gamma|W}=i_{W,\sigma}^*(J_{\chi,f,\gamma})$. Et par définition de l'idéal transverse, on a $J_{\chi,f,\gamma}=i_{\sigma,U_m}^*({\Cal I}_{\chi,f,m})$, où $U_m$ est un saturé suffisament petit, de $\sigma$ par le flot de $\chi$. D'un autre côté, on a ${\Cal I}_{\chi,f,m|V_{m}}=i_{V_m,U_m}^*({\Cal I}_{\chi,f,m})$. Et, on vérifie facilement que $i_{W,\sigma}^*\circ i_{\sigma,U_m}^*=i_{W,V_m}^*\circ i_{V_m,U_m}^*$.\par

\vskip 3mm

  Pour l'égalité $(b)$, par le lemme de saturation, on a 

$$
{\Cal I}_{\chi,f,m|V_m}=i_{V_m,U_m}^*\circ\pi_{\chi_m}^*(J_{\chi,f,\gamma})
$$

\noindent D'un autre côté, on a $\pi_{{\Cal Y}_m}^*(J_{\chi,f,\gamma|W})=\pi_{{\Cal Y}_m}^*\circ i_{W,\sigma}^*(J_{\chi,f,\gamma})$. Et, on vérifie facilement que les morphismes $\pi_{\chi_m}\circ i_{V_m,U_m}: V_m\to \sigma$ et $i_{W,\sigma}\circ\pi_{{\Cal Y}_m}: V_m\to\sigma$, coincident.\par
  
\vskip 3mm

  Notons $g_i=f_{i,|V_i}$. Il est clair que $I_{\chi_i,f_i|V_i}=I_{{\Cal Y}_i,g_i}$, et ceci est encore vrai germifié en n'importe quel point de $V_i$. L'égalité $(c)$ est donc une conséquence immédiate de la commutativité du diagramme $(*)$. Pour l'égalité $(d)$, on utilise l'égalité $(a)$

$$
J_{\chi_2,f_2,\gamma_2|W_2}=i_{W_2,V_{2,m_2}}^*({\Cal I}_{\chi_2,f_2,m_2|V_{2,m_2}})
$$

\noindent puis la relation $(c)$ germifiée aux points $m_1$ et $m_2$

$$
({\Cal I}_{\chi_2,f_2,m_2|V_{2,m_2}}=\Psi_{m_2}^*({\Cal I}_{\chi_1,f_1,m_1|V_{1,m_1}})
$$

\noindent et enfin, l'égalité $(b)$ pour obtenir

$$
J_{\chi_2,f_2,\gamma_2|W_2}=i_{W_2,V_{2,m_2}}^*\circ \Psi_{m_2}^*\circ \pi_{{\Cal Y}_1,m_1}^*(J_{\chi_1,f_1,\gamma_1|W_1})
$$

\noindent La composition de morphismes étoilés passe aux idéaux prolongés. On conclut en utilisant la commutativité du diagramme $(***)$, et la relation $\pi_{{\Cal Y}_2,m_2}\circ i_{W_2,V_{2,m_2}}=id_{W_2}$.\qed

\vskip 5mm

\centerline{{\bf  II. Action de la dérivation} $\chi=x\partial/\partial x$ {\bf et théorème principal 1}}\par

\vskip 5mm

 Les anneaux de références dans toute la suite sont les anneaux $SB^{p,|q|}(x,\alpha)$. Leur topologie de Krull n'est pas séparée. Leurs sous-algèbres $QA^{p,|q|}$ {\bf ne sont pas} $\chi${\bf -finies}: l'idéal différentiel du germe

$$
 f=\sum_{n>0} \alpha^n x^{1/n}
$$

\noindent n'est pas noethérien dans l'anneau $SB^{1,1}(x,\alpha)$: en effet, si tel est le cas, il existe $\ell\in\NN$, et des germes $h_i\in SB^{1,1}$ tels que

$$
\sum_{i=0}^\ell h_i\chi^i f +\chi^{\ell+1} f=0
\tag 0
$$

\noindent Soit $h_i=\sum_n h_{n,i}(x)\alpha^n$ la série de $h_i$. Les coefficients $h_{n,i}\in SB^{1,0}(x)$; notons $a_i=h_{0,i}(0)$. Par une identification des coefficients des séries en $\alpha$ dans l'égalité (0), on obtient les relations suivantes pour tout $n>0$

$$
\sum_{i=0}^\ell a_i(\frac{1}{n})^i+(\frac{1}{n})^{\ell+1}=0
$$

\noindent ce qui est impossible (en utilisant un système de Vandermonde adéquat).\par

\vskip 3mm

  Grâce à la propriété de {\bf quasi-régulari\-té}, la topologie de Krull des algèbres \par\noindent $QR{\Cal H}^{p,q}$ est séparée et, pour $p=1$, ces algèbres sont $\chi${\bf -finies} (cf. théorème II1). La question de leur noethérianité est un problème ouvert.\par

\vskip 3mm

  Plaçons nous dans l'anneau de référence $SB^{1,|q|}(x,\alpha)$. Soit $U\in ({\bold{\RR}}^{+*}\times{\bold{\RR}}^{|q|},0)$ tel que le sous-ensemble $\gamma=\{ (x,\alpha)\in U;\ \alpha=0\}$ soit connexe. Alors $\gamma$ est une orbite principale de $\chi$ dans $U$. Le morphisme int\'egral de $\chi$ dans $U$ est simplement la projection canonique $\pi:(x,\alpha)\mapsto \alpha$ (on notera de la même façon son germe en 0). Les id\'eaux $\chi$-transverses le long de $\gamma$ sont donc des id\'eaux de l'anneau analytique $\RR\{\alpha\}$. L'anneau $(SB^{1,|q|}(x,\alpha),\pi)$ est une extension étoilée de l'anneau $\RR\{\alpha\}$.\par

\vskip 3mm

\proclaim{théorème II1 (théorème principal 1)} L'algèbre $QR{\Cal H}^{1,q}$ est $\chi$-finie et satisfait à la double inclusion: pour tout $f\in QR{\Cal H}^{1,q}$ d'idéal $\chi$-transverse $J_{\chi,f,\gamma}$, il existe $n(f)$ tel que pour tout $\varepsilon>0$

$$
 (x^{n(f)+\varepsilon})\pi^*(J_{\chi,f,\gamma})\subset I_{\chi,f}\subset \pi^*(J_{\chi,f,\gamma})
\tag 1
$$

\noindent De plus, elle est $\chi$-équivalente à la sous-algèbre $QR{\Cal H}_{cvg}^{1,q}$.\par
\endproclaim

\vskip 3mm

 L'argument principal de la preuve de la $\chi$-finitude est que l'algèbre $QR{\Cal H}^{1,q}$ satisfait à la {\bf double inclusion} (1). La deuxième inclusion est une conséquence du\par

\vskip 3mm

\proclaim{Lemme II1 (lemme de division)} Soit $f\in SB^{1,|q|}$ et $J_{\chi,f,\gamma}$ son id\'eal $\chi$-transverse. Alors $I_{\chi,f}\subset\pi^*(J_{\chi,f,\gamma})$.
\endproclaim

\vskip 3mm

\noindent {\bf Preuve.} On utilise le th\'eor\`eme de division VB1 de l'appendice B. Soit $\Delta$ le complémentaire du diagramme des exposants initiaux de $J_{\chi,f,\gamma}$. Effectuons la division de $f(x,.)$ dans $J_{\chi,f,\gamma}$ pour tout $x\neq 0$. Soit $Q\in \pi^*(J_{\chi,f,\gamma})$ telle que

$$
 \text{Supp}(f(x,.)-Q(x,.))\subset \Delta\quad\text{pour tout}\quad x\neq 0
$$

\noindent Or, d'apr\`es le lemme de saturation IB5, on a $f(x,.)\in J_{\chi,f,\gamma}$ pour tout $x\neq 0$, et donc $f-Q\equiv 0$; d'o\`u le r\'esultat.\qed\par

\vskip 3mm

 La première inclusion est basée sur la notion de multiplicité et sur la propriété de quasi-analycité qui se substitue au théorème d'intersection de Krull ([L]) pour le passage à la limite. La mulptiplicité $m_\chi(f)$ est le plus petit des entiers $n$ tel que $(x^{n+\varepsilon})\pi^*(J_{\chi,f,\gamma})\subset I_{\chi,f}$ pour tout $\varepsilon>0$. Elle coîncide avec {\bf la multiplicité algébrique} $ma_\chi(f)$ qui se lit sur la série asymptotique de $f\in QR{\Cal H}^{1,q}$,  en faisant agir la dérivation $\chi$ d'abord sur les composantes élémentaires de cette série qui sont des fewnomials [K2], solutions d'équations différentielles simples.\par

\vskip 3mm

\noindent {\bf §1. Operateurs différentiels d'Euler} $E_{\Cal F}$.\par

\vskip 3mm

  Pour tout multi-indice $m=(m_0,\ldots,m_{q_1}) \in {\bold{\NN}}^{1+q_1}$ et tout $q_1$-uplet de nombre caract\'eristiques $(r_1,\ldots,r_{q_1})$ avec $r_j=1+\mu_j$, on pose

$$ 
r=(1,r_1,\ldots,r_{q_1})\quad \text{  et  }\quad e_m(\mu) =\langle m,r \rangle 
$$

\noindent A toute famille finie de multi-indices ${\Cal F} \subset {\bold{\NN}}^{1+q_1}$, on associe l'op\'erateur 

$$ 
E_{\Cal F}=\prod_{m\in \Cal F} (\chi -e_m(\mu) Id)  
\tag 2
$$

\noindent et on note $P_{\Cal F}$ son polyn\^ome caract\'eristique. On notera aussi 

$$
{\Cal F}_{?n} =\{m\in {\NN}^{1+q_1};\quad |m|?n\}.
$$

\noindent où $?$ est un opérateur binaire de comparaison. Les solutions de l'\'equation $E_{\Cal F}\cdot=0$ sont les combinaisons sur $\RR\{\alpha\}$ des m\^onomes $x^{e_m}$ pour $m\in {\Cal F}$, et des fonctions \'el\'ementaires $x^{e_m}(\log x)^p$ si $e_m$ est racine multiple de $P_{\Cal F}$. Donc, {\bf g\'en\'eriquement en} $\mu$, les mon\^omes $x^{e_m(\mu)}$ (pour $m\in {\Cal F}$), forment {\bf une base} du noyau de l'op\'erateur $E_{\Cal F}$.\par

\vskip 3mm

 Les fonctions \'el\'ementaires $z_j=xLd(x,\mu_j)$ satisfont aux \'equations diff\'erentielles

$$
 \chi z_j=r_jz_j +x
$$

\noindent donc, un mon\^ome $X^m=x^{m_0}z_1^{m_1}\cdots z_{q_1}^{m_{q_1}}$ satisfait \`a l'\'equation diff\'erentielle

$$
 \chi X^m=e_m X^m +\text{mon\^omes} \succ
\tag 3
$$

\noindent o\`u $\prec$ est un ordre adéquat sur ces mon\^omes (associé à l'ordre lexicographique sur $\NN^{1+q_1}$). Par suite, tout mon\^ome $X^m$ est dans le noyau de l'op\'erateur $E_{{\Cal F}_{=|m|}}$. D'autre part, de la relation $x^{r_j}=x+\mu_j z_j$ on déduit que

$$
 x^{e_m(\mu)}=\sum_{|m'|=|m|} c_{m'}(\mu) X^{m'}
$$

\noindent donc, pour tout $n$, la famille des mon\^omes $X^m$ de longeur $n$ est, {\bf g\'en\'eriquement en} $\mu$, {\bf une base du noyau de l'op\'erateur} $E_{{\Cal F}_{=n}}$.\par

\vskip 3mm

\noindent  {\bf §2. Le wronskien de l'opérateur} $E_{{\Cal F}_{=n}}$ {\bf et ses mineurs.}\par

\vskip 3mm

 Soit $n\in \bold{\NN}$ et $N(n)=\sharp {\Cal F}_{=n}$. Notons $M_n$ la matrice $N(n)\times N(n)$ dont les colonnes sont les d\'eriv\'ees successives par $\chi$ des mon\^omes $X^{m}$ ($|m|=n$),  et soit 

$$ 
 \Delta_n (x,\mu)=\det M_n    
\tag 4
$$

\noindent Ces monômes forment une base du noyau de l'opérateur $E_{{\Cal F}_{=n}}$. Donc, en d\'erivant les colonnes de $M_n$, et en utilisant (3), on obtient\par

\vskip 3mm

\proclaim{Lemme II2} Il existe une fonction algébrique $b_n(\mu)$, non-identiquement nulle telle que 

$$ 
 \Delta_n (x,\mu)=b_n(\mu)x^{s_n(\mu)} \quad\text{avec}\quad  s_n(\mu)=\langle {\sum}_{|m|=n} m,r\rangle
\tag 5
$$

\endproclaim

\vskip 3mm

\noindent Cette fonction algébrique est donnée par $b_n(\mu)=\Delta_n(1,\mu)$. Soit $A_n(x,\mu)$ la matrice compl\'ementaire de la matrice $M_n(x,\mu)$\par

\vskip 3mm

\proclaim{Lemme II3} Pour tout $x_0\neq 0$ et pour tout $\varepsilon>0$, les \'el\'ements de la matrice $M_n(x_0,\mu)A_n(x,\mu)$ appartiennent \`a l'id\'eal principal de $SB^{1,q_1}$ engendr\'e par\par
\noindent $x^{nN(n)-n-\varepsilon}b_n$.
\endproclaim

\vskip 3mm

\noindent {\bf Preuve.}  Ces \'el\'ements sont de la forme

$$
 B_{i,l}(x,\mu)=L_i(x_0,\mu)C_l(x,\mu)
\tag 6
$$

\noindent o\^u $L_i(x,\mu)$ est la ligne d'indice $i$ dans la matrice $M_n$ et $C_l(x,\mu)$ est la colonne d'indice $l$ dans la matrice $A_n$. Il est clair, d'apr\`es la d\'efinition de la matrice $A_n$ que ces \'el\'ements sont divisibles par $x^{nN(n)-n-\varepsilon}$ dans l'anneau $SB^{1,q_1}$ (chaque monôme $X^m$ est divisible par $x^{n-\varepsilon/s}$ dans cet anneau, pour tout $s>0$). Il suffit donc de montrer que l'id\'eal $\chi$-transverse $J$, de $B_{i,l}$ le long de $\gamma$ est inclus dans l'id\'eal principal engendr\'e par $b_n$. Le lemme de division II1 permettra de conclure.\par

\vskip 3mm

\noindent Pour cela, montrons par une récurrence sur $k$, que pour tout $(i,l)\in \{ 1,\ldots,N(n) \}^2$ et pour tout $k$

$$
 L_i(x,\mu)\chi^k C_l(x,\mu)\equiv 0\quad [(b_n)]\qquad \text{ dans } SB^{1,q_1} 
\tag 7
$$

\noindent Pour $k=0$, la relation (7) est une cons\'equence de l'égalité $M_n A_n=\Delta_n Id$.\par

\vskip 2mm

\noindent L'id\'eal $(b_n)$ \'etant stable par $\chi$, une d\'erivation de (7) donne

$$
 L_i \chi^{k+1} C_l \equiv -(\chi L_i)( \chi^k C_l)\quad [(b_n)]
$$

\noindent or, pour $i\neq N(n)$ on a $\chi L_i=L_{i+1}$ d'apr\`es la d\'efinition de la matrice $M_n$. Et en utilisant l'op\'erateur $E_{{\Cal F}_{=n}}$, on obtient

$$
 \chi L_{N(n)}=\sum_{j=1}^{N(n)} c_j(\mu) L_j
$$
 
\noindent Et ceci prouve la relation (7). Maintenant, d'après la définition IB6 de l'idéal $\chi$-transverse, l'idéal $J\subset\RR\{\mu\}$ est engendré par les dérivées successives 

$$
(\chi^k B_{i,l}(x_1,\mu))_{k\in\NN}
$$

\noindent pour tout $x_1>0$. On prend donc $x_1=x_0$ et on utilise (7) pour finir la preuve du lemme.\qed\par

\vskip 3mm

\noindent {\bf §3. Sur les blocs} $\chi${\bf -homogènes.}\par

\vskip 3mm

  Soit $c$ l'immersion de la partie IA, et soit $H{\Cal H}_n$ l'image par $c^*$ du $\RR \{ \alpha \}$-module $\RR\{\alpha\}[X]_n$, engendr\'e par les mon\^omes $X^m$ d'une m\^eme longueur $n$ ($c^*(X)=(x,(z_j(x,\mu_j))_{j=1,\ldots,q_1})$), rappelons qu'on note de la même façon ces fonctions $X,z_j$ et les coordonnées correspondantes). D'après le paragraphe 1, la restriction de $c^*$ à ce module est un {\bf isomorphisme}. Les éléments du module $H{\Cal H}_n$ sont dits {\bf blocs } $\chi${\bf -homogènes de degré }$n$. Ces modules sont stables par les opérateurs d'Euler définis ci-dessus.\par

\vskip 3mm

\proclaim{Lemme II4} Tout $g\in H{\Cal H}_n$ satisfait à la {\bf double inclusion}. Plus précisément, pour tout $\varepsilon >0$

$$
 (x^{n+\varepsilon})\pi^*(J_{\chi,g,\gamma})\subset I_{\chi,g}\subset\pi^*(J_{\chi,g,\gamma})
$$

\endproclaim

\vskip 3mm

\noindent {\bf Preuve.} La deuxième inclusion est donnée par le lemme de division II1. Un tel $g$ s'\'ecrit

$$
 g=\sum_{|m|=n} a_m(\alpha) X^m  
\tag 8
$$

\noindent il est donc dans le noyau de l'op\'erateur $E_{{\Cal F}_{=n}}$. Par cons\'equent, son id\'eal diff\'erentiel est engendr\'e par la famille $\{ \chi^jg;\ j<N(n) \}$, et donc son id\'eal $\chi$-transverse le long de $\gamma$ est engendr\'e par la famille $\{ \chi^j g(x_0,.);\ j<N(n) \}$ pour tout $x_0\neq 0$. Plusieurs d\'erivations de (8) donnent le syst\`eme

$$
 M_n (a_m)_m =(\chi^jg)_j    
\tag 9
$$

\noindent ($(a_m)_m$ et $(\chi^j g)_j$ désignent les vecteurs colonnes associés). En particulier

$$
 M_n(x_0,.) (a_m)_m =(\chi^jg(x_0,.))_j   
\tag 10
$$

\noindent Mutiplions (9) par la matrice $M_n(x_0,.)A_n$ et utilisons (10)

$$
 \Delta_n (\chi^j g(x_0,.))_j =M_n(x_0,.) A_n (\chi^j g)_j
$$

\noindent Par le lemme II2, on a $\Delta_n=x^{s_n}b_n$, et par le lemme II3 

$$
M_n(x_0,.)A_n\in (x^{nN(n)-n-\varepsilon}b_n)
$$

\noindent Or, $s_n(0)=nN(n)$. Ce qui finit la preuve du lemme.\qed\par

\vskip 3mm

\noindent{\bf Remarque II1.} Un cas simple, où cette première inclusion peut effectivement être donnée par le Nullstellensatz d'Hilbert, est le suivant: prenons $q_1=1$ et supposons $J=\RR\{\alpha\}$. L'immersion $c(x,\alpha)=(x,z(x,\mu),\alpha)$ est un difféomorphisme sur la variété image $V=c(U)$, qui se prolonge en un homéomorhisme sur le bord $\partial_0 U=\{x=0\}$. La dérivation sur $V$ $c_*(\chi)$, se prolonge sur un voisinage $W$ de 0, en une dérivation analytique

$$
{\Cal X}=x\frac{\partial}{\partial x}+(rz+x)\frac{\partial}{\partial z}
$$

\noindent Soit $G=(c^*)^{-1}(g)$. D'après le lemme de transfert IB6 et l'isomorphie de $c^*:\RR\{\alpha\}[x,z]_n\to H{\Cal H}_n$, on a $I_{\chi,g}=c^*(I_{{\Cal X},G}$ (car $\RR\{\alpha\}[x,z]_{n|V}\cong\RR\{\alpha\}[x,z]_n$). Plaçons nous dans le complexifié de $W$. L'ensemble des zéros de l'idéal différentiel $I_{{\Cal X},G}$ est un sous-ensemble analytique invariant par le flot de ${\Cal X}$. Or, pour $\alpha$ fixé générique, les seules feuilles analytiques de ${\Cal X}$ sont $\{x=0\}$ et $\{x+\mu z=0\}$. Par le Nullstellensatz d'Hilbert, il existe $M\in\NN$ tel que 

$$
(x(x+\mu z))^M\in I_{{\Cal X},G}
$$

\noindent et en appliquant $c^*$, on obtient $x^{(1+r)M}\in I_{\chi,g}$. Ceci est encore valable pour $g\in QR{\Cal H}^{1,(1,q_2)}_{cvg}$, par la définition de cet anneau. On peut généraliser cette idée à plusieurs fonctions élémentaires ($q_1>1$), seulement l'estimation de l'entier $M$ optimal semble difficile. Par la démarche adoptée dans ce travail, on obtient une estimation précise et optimal de cet entier, en se plaçant dans des anneaux plus larges.\par

\vskip 3mm

\noindent{\bf §4. Sur les fewnomials.}\par

\vskip 3mm

 Le passage des blocs $\chi$-homogènes aux fewnomials est basée sur l'inversion des opérateurs $E_{\Cal F}$\par
 
\vskip 3mm

\proclaim{Lemme II5} Soit $e(\alpha)$ un germe analytique en 0 tel que $e(0)\neq n$. Soit l'opérateur $E=\chi-eId$. Alors, pour tout $g\in H{\Cal H}_n$, $I_{\chi,E g}=I_{\chi,g}$.
\endproclaim

\vskip 3mm

\noindent {\bf Preuve.} Il est clair que $I_{\chi,Eg}\subset I_{\chi,g}$. De la relation $\chi g= eg +Eg$, on déduit

$$
 E_{{\Cal F}_{=n}}g\equiv P_{{\Cal F}_{=n}}(e)g\qquad [I_{\chi,Eg}]
$$

\noindent Or $E_{{\Cal F}_{=n}}g=0$ et le germe $e$ n'est pas valeur propre de cet op\'erateur, on obtient donc $g\in I_{\chi,Eg}$, par conséquent $I_{\chi,g}\subset I_{\chi,Eg}$.\qed\par

\vskip 3mm

 Soit $n\in\NN$ quelconque et soit $f=\sum_{p=0}^n g_p$, o\`u $g_p$ est un bloc $\chi$-homog\`ene de degr\'e $p$, $f$ est un fewnomial d'après la terminologie de Khovanski [K1].  Appliquons l'op\'erateur $E_{{\Cal F}_{<n}}$ \`a $f$
 
$$
E_{{\Cal F}_{<n}}f=E_{{\Cal F}_{<n}}g_n
\tag $10'$
$$

\noindent car deux opérateurs $\chi-eId$ et $\chi-e'Id$ commutent, et pour $p<n$, le bloc $g_p$ est dans le noyau de l'opérateur $E_{{\Cal F}_{<n}}$. En appliquant plusieurs fois le lemme II5 à $g_n$, on obtient que $I_{\chi,E_{{\Cal F}_{<n}}}=I_{\chi,g_n}$, et d'après (10'), $I_{\chi,g_n}\subset I_{\chi,f}$. De même, $I_{\chi,g_{n-1}}\subset I_{\chi,f-g_n}\subset I_{\chi,f}$, ...etc. On obtient donc\par
 
\vskip 3mm

\proclaim{Lemme II6} \qquad $I_{\chi,f}=\sum_{p=0}^n I_{\chi,g_p}$.
\endproclaim

\noindent Ces idéaux étant des idéaux différentiels, cette égalité est encore vraie germifiée en tout point de $\gamma$ voisin de 0. Prenons la restriction à une transversale à $\gamma$, cette égalité donne 

$$
J_{\chi,f,\gamma}\subset\sum_{p=0}^n J_{\chi,g_p,\gamma}
$$

\noindent et l'inclusion $I_{\chi,g_p}\subset I_{\chi,f}$ donne $J_{\chi,g_p,\gamma}\subset J_{\chi,f,\gamma}$ pour tout $p=0,\ldots,n$. D'où l'égalité

$$
J_{\chi,f,\gamma}=\sum_{p=0}^n J_{\chi,g_p,\gamma}
\tag $10"$
$$

\vskip 3mm

\proclaim{Lemme II7} Le germe $f$ satisfait à la double inclusion. Plus précisément, pour tout $\varepsilon>0$

$$
  (x^{n+\varepsilon})\pi^*(J_{\chi,f,\gamma})\subset I_{\chi,f}\subset\pi^*(J_{\chi,f,\gamma})
$$

\endproclaim

\vskip 3mm

\noindent{\bf Preuve.} Par les lemmes II4 et II6, on a

$$
\sum_{p=0}^n  (x^{p+\varepsilon})\pi^*(J_{\chi,g_p,\gamma})\subset I_{\chi,f}\subset\sum_{p=0}^n \pi^*(J_{\chi,g_p,\gamma})
$$

\noindent ou plus simplement

$$
(x^{n+\varepsilon})\sum_{p=0}^n \pi^*(J_{\chi,g_p,\gamma})\subset I_{\chi,f}\subset\sum_{p=0}^n \pi^*(J_{\chi,g_p,\gamma})
$$

\noindent L'égalité (10") permet de conclure.\qed\par

\vskip 3mm

\noindent {\bf §5. Sur les convergents.}\par

\vskip 3mm

\proclaim{Lemme II8} L'alg\`ebre $QR{\Cal H}^{1,q}_{\text{cvg}}$ est $\chi$-finie et satisfait à la {\bf double inclusion}. De plus, elle est $\chi$-équivalente à la sous-algèbre des fewnomials.\endproclaim

\vskip 3mm

\noindent {\bf Preuve.} Cette algèbre est clairement stable par $\chi$. Pour la $\chi$-finitude, on utilise le lemme d'extension IB2. L'anneau $QR{\Cal H}^{1,q}_{\text{cvg}}$ est la restriction d'un anneau analytique, c'st donc un anneau noethérien, d'où l'hypothèse $(ii)$ du lemme d'extension. L'hypothè\-se $(iii)$ est claire: prendre $\omega_j=d\alpha_j$. La preuve de l'hypothèse $(i)$ est une double application de la th\'eorie de Khovanski et d'un théorème de [T] (elle est aussi une conséquence d'un travail récent et général de Speissegger [S]). Soit $X=(x,z_1,\ldots,z_{q_1})$ et $z_0(x)=x\log x$. Soit $c_0$ l'immersion $c_0(X,\alpha)=(X,z_0(x),\alpha)$ et ${\Cal B}_1$ l'algèbre ${\Cal B}_1=c^*_0(\RR\{ X,z_0,\alpha \} )$. Le graphe de la fonction $z_0$ est une solution s\'eparante de l'\'equation diff\'erentielle

$$  
 \omega_0=xdz_0-(z_0+x)dx=0     
\tag 11
$$ 

\noindent considérée sur un ouvert connexe $U_0$, voisinage de 0 dans $\{ x>0,\ z_0<0\}$. La 1-forme $\omega_0$ est à coefficients dans l'anneau analytique $\RR\{ X,z_0,\alpha \}$. L'algèbre ${\Cal B}_1$ est donc topologiquement noethérienne [T]. Soit $c_1$ l'immersion $c_1(x,\alpha)=(X(x,\alpha),\alpha)$, alors $QR{\Cal H}^{1,q}_{cvg}=c_1^*({\Cal B}_1)$.  les graphes des fonctions $z_j$ sont des solutions s\'eparantes des \'equations diff\'erentielles

$$
 \omega_j=xdz_j-((1+\mu_j)z_j+x)dx-(x\log x(\mu_jz_j+x)-xz_j)\frac{d\mu_j}{\mu_j}=0
\tag 12
$$

\noindent considérées sur des ouverts connexes $U_{(?_j)}$, voisinages de 0 dans $\{ x>0,\ z_j<0,\ \mu_j ?_j 0\ j=1,\ldots,q_1\}$, où $?_j=>,<$ (si l'un des $\mu_j$ est nul, alors $z_j=z_0$). Ces 1-formes sont \`a coefficients dans l'alg\`ebre ${\Cal B}_1$. L'algèbre $QR{\Cal H}^{1,q}_{\text{cvg}}$ est donc topologiquement noethérienne. D'où l'hypothèse $(i)$.\par

\vskip 3mm

   Soit $f=\sum_{p=0}^{\infty} g_p \in QR{\Cal H}^{1,q}_{\text{cvg}}$ et soit $f_n=\sum_{p\leq n} g_p$, où $g_p$ est un bloc $\chi$-homogène de degré $p$. D'après le lemme II6, la suite des idéaux $(I_{\chi,f_n})$ est croissante dans l'anneau noethérien $QR{\Cal H}^{1,q}_{cvg}$. Elle est donc stationnaire. Soit $I$ sa limite et soit $n_0$ son indice de stationnarité: c'est le plus petit entier $n$ tel que $I_{\chi,f_n}=I_{\chi,f_{n'}}$ pour tout $n'\geq n$. Soit ${\Cal M}_X=\langle x,z_0,z_1,\ldots,z_{q_1}\rangle\subset QR{\Cal H}^{1,q}_{cvg}$, il est inclus dans l'idéal maximal de $QR{\Cal H}^{1,q}_{cvg}$, et il est stable par $\chi$. Pour tout $n\in\NN$, on a $f-f_n\in{\Cal M}_X^n$ d'après la série de $f$. Donc, pour tout $n\geq n_0$, on a
   
$$
I_{\chi,f}\subset I+{\Cal M}_X^n \quad\text{et}\quad I\subset I_{\chi,f}+{\Cal M}_X^n
$$

\noindent Par le théorème d'intersection de Krull, on a $I_{\chi,f}=I$. Or $I=I_{\chi,f_{n_0}}$, donc en germifiant ces égalités le long de $\gamma$, et en prenant les restrictions à une transversale, on obtient $J_{\chi,f,\gamma}=J_{\chi,f_{n_0},\gamma}$. Par conséquent, le lemme II7 appliqué au fewnomial $f_{n_0}$, montre que $f$ vérifie la double inclusion. Soit $k$ l'exposant d'Artin-Ress de l'idéal ${\Cal M}_X$ dans l'idéal $I$ ([L]): pour tout entier $n$, on a

$$
{\Cal M}_X^{n+k}\cap I={\Cal M}_X^k\cap ({\Cal M}_X^n I)
$$

\noindent Prenons $n>\max\{n_0,k\}$. On a $f-f_n\in{\Cal M}_X^n$ et $I_{\chi,f-f_n}\subset I$, donc $f-f_n\in{\Cal M}_X I$. Or, par le choix de $n$, $I=I_{\chi,f_n}$, le germe $f$ est donc $\chi$-équivalent au fewnomial $f_n$, et ceci finit la preuve du lemme.\qed\par

\vskip 3mm

\noindent {\bf §6. Enfin, sur l'algèbre} $QR{\Cal H}^{1,q}$.\par

\vskip 3mm

  Soit $f\in QR{\Cal H}^{1,q}$ et soit $\widehat{f}=\sum_{p=0}^{\infty} g_p$ {\bf une} s\'erie formelle associ\'ee \`a $f$. Les sommes finies $f_n=\sum_{p\leq n} g_p$ sont des fewnomials. La suite des idéaux différentiels $(I_{\chi,f_n})$ est croissante (dans chacun des anneaux $QR{\Cal H}^{1,q}_{cvg}$, $QR{\Cal H}^{1,q}$ ou $SB^{1,|q|}$). Il en est de même de la suite des idéaux $\chi$-transverses $(J_{\chi,f_n,\gamma})$ dans l'anneau $\RR\{\alpha\}$. Soit $I$ la limite "différentielle" de la suite $(I_{\chi,f_n})$ dans  l'anneau $QR{\Cal H}^{1,q}_{cvg}$, et soit $J$ la limite "$\chi$-transverse" de la suite $(J_{\chi,f_n,\gamma})$ dans l'anneau $\RR\{\alpha\}$. On va montrer que la limite différentielle $I$, prolongée dans l'anneau $SB^{1,|q|}$, coincide avec l'idéal différentiel de $f$ dans cet anneau (mais pas forcément dans l'anneau $QR{\Cal H}^{1,q}$). Par contre, on montre que la limite $\chi$-transverse $J$ coincide avec l'idéal $\chi$-transverse de $f$ dans l'anneau $\RR\{\alpha\}$; et cela grâce à deux arguments principaux: la propriété de quasi-analycité et le théorème de division VB1.\par

\vskip 3mm

\proclaim{Lemme II9} Pour tout $n$, on a $\qquad J_{\chi,g_n,\gamma}\subset J_{\chi,f,\gamma}$.
\endproclaim

\vskip 3mm

\noindent {\bf Preuve.} Supposons d'abord que $f=g_n+o(x^n)$. D'après la définition des séries asymptotiques de $f$ (cf. Déf. IA5), on peut trouver $\varepsilon>0$ suffisament petit, et $h\in SB^{1,|q|}$ tels que

$$
 f=g_n+x^{n+2\varepsilon}h   
\tag 13
$$

\noindent En germifiant l'égalité $x^{n+2\varepsilon}h=f-g_n$ en un point de $\gamma$ ($x>0$), et en prenant la retriction à une transversale des dérivées successives, on obtient $J_{\chi,h,\gamma}\subset J_{\chi,f,\gamma}+J_{\chi,g_n,\gamma}$. Donc, d'apr\`es le lemme de division II1, appliqué à $h$, il existe $h_1\in \pi^*(J_{\chi,f,\gamma})$ et $h_2\in \pi^*(J_{\chi,g_n,\gamma})$ telles que $h=h_1+h_2$. Maintenant, en consid\'erant les id\'eaux diff\'erentiels et en utilisant l'égalité (13), on obtient

$$
I_{\chi,g_n}\subset I_{\chi,f}+(x^{n+2\varepsilon})(I_{\chi,h_1}+I_{\chi,h_2})
$$

\noindent soit

$$
I_{\chi,g_n}\subset I_{\chi,f}+(x^{n+2\varepsilon})(\pi^*(J_{\chi,f,\gamma})+\pi^*(J_{\chi,g_n,\gamma}))
$$

\noindent Le lemme de division II1, appliqué à $f$, donne

$$
I_{\chi,g_n}\subset \pi^*(J_{\chi,f,\gamma})+(x^{n+2\varepsilon})\pi^*(J_{\chi,g_n,\gamma})
$$

\noindent et le lemme II4 appliqué à $g_n$ fournit

$$
I_{\chi,g_n}\subset \pi^*(J_{\chi,f,\gamma})+(x^{\varepsilon}) I_{\chi,g_n}
$$

\noindent L'idéal $I_{\chi,g_n}$ étant noethérien, et l'idéal principal $(x^{\varepsilon})$ étant inclus dans l'idéal maximal de l'anneau $SB^{1,q}$, le lemme de Nakayama ([L]) donne

$$
I_{\chi,g_n}\subset \pi^*(J_{\chi,f,\gamma})
$$

\noindent En prenant la restriction à une transversale à $\gamma$, on obtient $J_{\chi,g_n,\gamma}\subset J_{\chi,f,\gamma}$. Maintenant, si $\widehat{f}=\sum_{p\geq 0} g_p$, on montre comme ci-dessus que $J_{\chi,g_0,\gamma}\subset J_{\chi,f,\gamma}$, puis en supposant que pour tout $p<n$, on a $J_{\chi,g_p,\gamma}\subset J_{\chi,f,\gamma}$, on montre comme ci-dessus que 

$$
J_{\chi,g_n,\gamma}\subset J_{\chi,f-\sum_{p<n} g_p,\gamma}\subset J_{\chi,f,\gamma}
$$

\noindent et ceci finit la preuve du lemme.\qed\par

\vskip 3mm

\noindent Ce lemme et la propriété de quasi-analycité impliquent le\par

\vskip 3mm

\proclaim{Lemme II10} Il existe une application s\'erie formelle injective

$$
 f\in QR{\Cal H}^{1,q} \mapsto \widehat{f}=\sum_{p\geq 0} g_p \in c^*(\RR\{\alpha\}[[X]])
$$

\endproclaim

\vskip 3mm

\noindent{\bf Preuve.} Pour prouver l'existence, il suffit de montrer que le germe nul admet une unique série qui est la série nulle. En effet, si $\widehat{0}=\sum_{p\geq 0} g_p$, le lemme II9 dit que pour tout $n$, $J_{\chi,g_n,\gamma}\subset\{0\}$, et le lemme de division II1 dit alors que $g_n\equiv 0$. Pour prouver l'injection, prenons $f\in QR{\Cal H}^{1,q}$ tel que $g_p(f)=0$ pour tout $p$. On a donc $f=o(x^n)$ pour tout $n$; comme $QR{\Cal H}^{1,q}\subset QA^{1,|q|}$, on obtient $f\equiv 0$.\qed\par

\vskip 3mm

\noindent Une conséquence de cela est que le morphisme $c^*:\RR\{X,\alpha\}\to QR{\Cal H}^{1,q}_{cvg}$ est un isomorphisme. En effet, si $F=\sum_{p\geq 0} G_p\in\RR\{X,\alpha\}$ (les $G_p$ étant ses parties homogènes en $X$ de degré $p$), son image $f=c^*(F)$ admet comme série $\sum_{p\geq 0} c^*(G_p)$, et on a vu que la restriction de $c^*$ aux modules $\RR\{\alpha\}[X]_p$ est un isomorphisme sur les modules $H{\Cal H}_p$.\par

\vskip 3mm

 Le résultat clé de cette section est le suivant\par
 
\vskip 3mm

\proclaim{Lemme II11} $\quad J=J_{\chi,f,\gamma}$.\endproclaim

\vskip 3mm

\noindent {\bf Preuve.} D'apr\`es le lemme II9 et l'égalité (10"), on a 

$$
J_{\chi,f_n,\gamma}\subset J_{\chi,f,\gamma}\quad\text{pour tout}\ n
\tag $13'$
$$

\noindent donc $J\subset J_{\chi,f,\gamma}$. Pour montrer l'autre inclusion, on utilise le th\'eor\`eme de division VB1 (appendice VB), sur les algèbres $QA^{1,|q|}[.]$. Les fonctions $f_n$ sont algébriques dans les fonctions élémentaires $z_j$; il existe donc $\Omega\in {\Cal E}{\Cal I}$ tel que $f,\ f_n\in QA[\Omega]$ pour tout $n$. Soit $\varphi_1,\ldots,\varphi_l$ une base de $J$. Pour tout $F\in QA[\Omega]$, on note

$$
 Q(F)=\sum_i Q_i(F) \varphi_i \quad\text{ et }\quad R(F)=F-Q(F)
$$

\noindent où les fonctions $Q_i(F)$ et $R(F)$ sont données par le théorème de division VB1. D'après la définition de l'idéal limite transverse $J$, on a $J_{\chi,f_n,\gamma}\subset J$ pour tout $n$. Donc, d'apr\`es le lemme de division II1, on a $R(f_n)\equiv 0$ pour tout $n$. Et, par unicit\'e de la division dans le théorème VB1, on a

$$
 R(f)=R(f-f_n)\quad\text{pour tout}\ n
$$

\noindent Or, $f-f_n=o(x^n)$. Le lemme VB1 implique que $R(f)=o(x^n)$ pour tout $n$. Comme $R(f)\in QA[\Omega]$, on obtient $R(f)\equiv 0$. Ceci prouve que $f\in\pi^(J)$, et donc (par un raisonnement maintenant classique), que $J_{\chi,f,\gamma}\subset J$.\qed\par

\vskip 3mm

\proclaim{Définition II1} L'indice de stationnarité de la suite d'idéaux $(J_{\chi,f_n,\gamma})$ est dit {\bf multiplicité algébrique} de $f$ relativement à $\chi$ le long de $\gamma$. On la note $ma_{\chi}(f)$.\endproclaim

\vskip 3mm

\noindent {\bf §7. Preuve du théorème principal II1.}\par

\vskip 3mm

Soit $n\geq ma_{\chi}(f)$ et soit $h_n=f-f_n=o(x^{n+2\varepsilon})$ (d'après la série asymptotique de $f$). On a $J_{\chi,h_n,\gamma}\subset J$ d'apr\`es le lemme II11. Et, par le lemme de division II1, on a

$$
 h_n\in (x^{n+2\varepsilon})\pi^*(J)
$$

\noindent D'apr\`es le lemme II7, on a $(x^{n+\varepsilon})\pi^*(J)\subset I_{\chi,f_n}$. Donc, $h_n\in (x^{\varepsilon}) I_{\chi,f_n}$. L'idéal $(x^{\varepsilon})$ est inclus dans l'idéal maximal de $SB$, et il est satble par $\chi$. Comme $f_n$ est un fewnomial, le lemme de finitude IB1 et le lemme II8 impliquent que l'algèbre $QR{\Cal H}^{1,q}$ est $\chi$-finie, et qu'elle est $\chi$-équivalente à la sous-algèbre des fewnomials.\par

\vskip 3mm

 Ceci implique en particulier que $I_{\chi,f}=I_{\chi,f_n}$ pour tout $n\geq ma_\chi (f)$. En prenant $n=ma_\chi (f)$ et en utilisant la double inclusion du lemme II7, on obtient la double inclusion pour l'algèbre $QR{\Cal H}^{1,q}$
 
$$
  (x^{ma_\chi(f)+\varepsilon})\pi^*(J_{\chi,f,\gamma})\subset I_{\chi,f}\subset\pi^*(J_{\chi,f,\gamma})
$$

\qed\par

\vskip 3mm

\noindent{\bf Remarque II2.} En fait, cette multiplicité algébrique $ma_\chi(f)$ est le plus petit entier $n$ tel que $I_{\chi,f}\supset (x^{n+\varepsilon})\pi^*(J_{\chi,f,\gamma})$ pour tout $\varepsilon>0$ (dans l'anneau $SB^{1,|q|}$). En effet, supposons que $I_{\chi,f}\supset (x^{ma_\chi(f)-1+\varepsilon})\pi^*(J_{\chi,f,\gamma})$ pour tout $\varepsilon>0$. On a 

$$
f-f_{ma_\chi(f)-1}\in (x^{ma_\chi(f)-1+2\varepsilon})\pi^*(J_{\chi,f,\gamma})
$$

\noindent pour $\varepsilon>0$ suffisament petit (par les mêmes raisonnements que ci-dessus). Donc

$$
I_{\chi,f}\subset I_{\chi,f_{ma_\chi(f)-1}}+(x^{ma_\chi(f)-1+2\varepsilon})\pi^*(J_{\chi,f,\gamma})
$$

\noindent et en utilisant l'hypothèse ci-dessus

$$
(x^{ma_\chi(f)-1+\varepsilon})\pi^*(J_{\chi,f,\gamma})\subset I_{\chi,f_{ma_\chi(f)-1}}+ (x^{ma_\chi(f)-1+2\varepsilon})\pi^*(J_{\chi,f,\gamma})
$$

\noindent Par le lemme de Nakayama, on obtient

$$
(x^{ma_\chi(f)-1+\varepsilon})\pi^*(J_{\chi,f,\gamma})\subset I_{\chi,f_{ma_\chi(f)-1}}
$$

\noindent Mais ceci implique (par un raisonnement classique), que $J_{\chi,f,\gamma}\subset J_{\chi,f_{ma_\chi(f)-1},\gamma}$, ce qui contredit la définition de la multiplicité algébrique $ma_\chi(f)$.\par

\vskip 3mm

\noindent{\bf §8. Morphismes série formelle des algèbres} $QR{\Cal H}^{p,q}$.\par

\vskip 3mm

  Soit $f\in QR{\Cal H}^{p,q}(x,\alpha)$ et soit $\widehat{f}^i=\sum_{n\geq 0} g_{i,n}$ une série formelle de $f$ relativement à la variable $x_i$ (cf. Déf. IA5). Montrons qu'elle est unique. Notons $f_{i,n}=\sum_{\ell\leq n} g_{i,n}$. Soient $h_n\in SB_0^{p,|q|}$ tels que
  
$$
f-f_{i,n}=x_i^n h_n
$$

\noindent Soient $U_{i,n}\in ((\RR^{+*})^p\times\RR^{|q|},0)$ une suite décroissante d'ouverts telle que $f_{i,n}$ et $h_n$ soient réalisées sur $U_{i,n}$. Soient $x^{(n)}=(x_{1,n},\ldots,x_{i-1,n},0,x_{i+1,n},\ldots,x_{p,n})$ tel que $x_{j,n}>0$ et $(x^{(n)},0)$ appartient à l'intérieur de $\overline{U}_{i,n}\cap\{x_i=0\}$ dans $\RR^{p+|q|-1}$. Le germe de $f,f_{i,n}$ et $h_n$ en $(x^{(n)},0)$ est un élément d'une algèbre $QR{\Cal H}^{1,q'}(x_i,\alpha')$ (qui dépend de $n$ par les coordonnées analytiques $\alpha'_j=x_j-x_{j,n}$ pour $j\neq i$). Le germe en $(x^{(n)},0)$ de $g_{i,\ell}$ (pour $\ell\leq n$), est un élément du $\RR\{\alpha'\}$-module correspondant $H{\Cal H}_\ell\subset QR{\Cal H}^{1,q'}(x_i,\alpha')$. Donc, par le lemme II9, si $f\equiv 0$ alors $f_{i,n}\equiv 0$, et ceci pour tout $n$. Par conséquent la série $\widehat{f}^i$ est identiquement nulle. L'injectivité de ce morphisme série formelle est une conséquence facile de la propriété de quasi-analycité.\par

\vskip 3mm

 En partant de chaque série formelle $\widehat{f}^i$, et en utilisant une récurrence sur $p$, on construit une série formelle unique

$$
 \widehat{f}\in c^*(\RR\{\alpha\}[[X]])
$$ 

\noindent au sens suivant: pour tout $n$

$$
f-{\bold j}^n_X(\widehat{f})\in{\Cal M}^n_x SB_0^{p,|q|}
$$

\noindent où ${\Cal M}_x$ est l'idéal de $SB^{p,|q|}$ engendré par les coordonnées $x_1,\ldots,x_p$. Pour prouver l'unicité de cette série, il suffit de se placer dans n'importe quel carte projective dans la coordonnée $x$, par exemple: $x_1=y_1$ et $x_i=y_iy_1$ pour $i\neq 1$. On utilise alors la méthode de la première partie de ce paragraphe, aux voisinages de points tels que $y_1=0,\ \alpha=0$ et $y_{i,n}>0$ pour $i\neq 1$. Les formules suivantes (obtenues par un calcul direct)

$$
z_{i,j}(x_i,\mu_j)=z_{i,j}(y_i,\mu_j)y_1+(y_i+\mu_jz_{i,j}(y_i,\mu_j))z_{1,j}(y_1,\mu_j)
$$

\noindent montrent que les jets dans les fonctions élémentaires $X_1(y_1,\mu)$, sont préservés (les fonctions $z_{i,j}(y_i,\mu_j)$ sont analytiques dans les coordonnées $y_i-y_{i,n}$). L'injectivité de ce morphisme est encore une conséquence facile de la propriété de quasi-analycité.\par

\vskip 5mm

\centerline{{\bf III. Action de la dérivation} 
$\chi=x\partial/\partial x-\sum_{j=1}^\ell s_j(\alpha)u_j\partial/\partial u_j$}\par

\vskip 5mm

 Soient $\alpha=(\mu,\nu)$ des coordonn\'ees sur ${\bold{\RR}}^{q_1}\times {\bold{\RR}}^{q_2}$ et soit $u$ une coordonn\'ee sur ${\bold{\RR}}^\ell$. Posons $\alpha'=(\alpha,u)$ et $q=(q_1,q_2+\ell)$. On veut étudier l'action de la dérivation $\chi$ sur l'algèbre $QR{\Cal H}^{1,q}(x,\alpha')$. Les germes $s_j$ sont analytiques et on suppose pour simplifier la présentation que $s_j(\alpha)=1+\mu_j$. Soit $U\in ({\bold{\RR}}^{+*}\times {\bold{\RR}}^{|q|},0)$ un ouvert tel que le sous-ensemble $\gamma=\{(x,\alpha')\in U;\ \alpha'=0\}$ soit connexe. C'est une orbite principale de $\chi$ dans $U$. Soit $\pi_\chi$ (le germe de) la projection intégrale de $\chi$ dans $U$. Quitte à réduire $U$, l'espace $\pi_{\chi}(U)$ s'identifie \`a une transversale analytique \`a $\gamma$: $W=\{ x=x_0>0 \}$. Des coordonn\'ees analytiques naturelles sur $W$ sont les int\'egrales premi\`eres de $\chi$, c'est à dire la coordonnée $\alpha$ et les (coordonnées) germes  

$$
 \lambda_j(x,\alpha')=x^{s_j(\alpha)}u_j\in QR{\Cal H}^{1,q}
$$

\noindent Les idéaux $\chi$-transverses le long de $\gamma$ sont donc des idéaux de l'anneau $\RR\{ \alpha,\lambda \}$, et on a le morphisme étoilé $\pi_\chi^*:\RR\{\alpha,\lambda\}\to SB^{1,|q|}(x,\alpha,u)$.\par

\vskip 3mm

  La question qui se pose est: l'algèbre $QR{\Cal H}^{1,q}$ est-elle $\chi$-finie? Le problème est ouvert. Cependant, on montre dans cette section, qu'elle est {\bf localement} $\chi${\bf -finie}.\par

\vskip 3mm

\noindent{\bf A. Etude globale et théorème principal 2.}\par

\vskip 3mm

  Le germe en 0 de la projection $\pi_\chi$ {\bf n'est pas isomorphe à celui d'une projection linéaire}, d'où des difficultés nouvelles dans l'étude de cette action. La première difficulté apparaît dans l'étude de l'action formelle de $\chi$.\par

\vskip 3mm

\noindent{\bf §1. Action formelle.}\par

\vskip 3mm

   Notons $r=(1,(1+\mu_j))$ et $s=(s_j)$. Les valeurs propres de l'opérateur $\chi$ sont

$$
 e_{m,n}(\alpha)=\langle m,r\rangle -\langle n,s\rangle
$$

\noindent Les mon\^omes $X^m u^n$ satisfont aux \'equations diff\'erentielles

 $$
 \chi(X^m u^n)=e_{m,n}X^m u^n +\text{monômes}\succ
$$

\noindent pour un ordre $\prec$ adéquat sur ces monômes (induit par l'ordre sur les monômes $X^m$). Le degré d'un monôme $X^mu^n$ relativement à $\chi$ est $|m|-|n|$. Soit $c$ l'immersion

$$
 c(x,\alpha,u)=(X(x,\mu),\alpha,u)
$$

\vskip 2mm

\proclaim{Définition IIIA1} On note $\widehat{H{\Cal H}}_p\subset c^*(\RR\{\alpha\}[[X,u]])$ le $\RR\{\alpha\}$-module des séries de monômes de degré $p\in\ZZ$. Ses éléments sont dits {\bf blocs} $\chi${\bf -homogènes formels} de degré $p$.
\endproclaim

\vskip 2mm

\noindent Dans l'action de $\chi_0=x\partial/\partial x$ (cf. section II), on a montré l'existence d'une application série formelle injective. Par une {\bf construction explicite}, on en déduit facilement le

\vskip 3mm

\proclaim{Lemme IIIA1} Il existe une application série {\bf doublement formelle} et injective

$$
 f\in QR{\Cal H}^{1,q}\to \widehat{f}=\sum_{p\in\ZZ} g_p(f)\in \sum_{p\in\ZZ}\widehat{H{\Cal H}}_p
$$

\endproclaim

\vskip 3mm

\noindent {\bf Preuve.} Soit $f\in QR{\Cal H}^{1,q}(x,\alpha')$. Sous l'action du champ $\chi_0$, on lui associe une unique s\'erie formelle

$$
 \widehat{f}^0=\sum_{M\geq 0} g_{0,M}(f)
\tag 1
$$

\noindent o\`u $g_{0,M}$ est un bloc $\chi_0$-homog\`ene de degr\'e $M$, i.e

$$
 g_{0,M}=\sum_{|m|=M} a_m(\alpha') X^m
$$

\noindent Or $g_{0,M}(f)\in QR{\Cal H}^{1,q}_{cvg}(x,\alpha')$, on lui associe donc une unique s\'erie convergente

$$
 g_{0,M}(f)=\sum_{p\leq M} g_{p}(g_{0,M}(f))
\tag 2
$$

\noindent o\`u $g_{p}(g_{0,M}(f))$ est un bloc $\chi$-homog\`ene convergent de dgré $p$. Les blocs formels $g_p(f)$ sont donnés par

$$
 g_p(f)=\sum_{M\geq p} g_{p}(g_{0,M}(f))
$$

\noindent Cette application est injective comme les applications (1) et (2).\qed

\vskip 3mm

 L'anneau $c^*(\RR[[X,\alpha']])\hookleftarrow c^*(\RR\{\alpha\}[[X,u]])$ est noethérien. Grâce au théorème d'intersection de Krull et au lemme d'Artin-Ress ([L]), on va montrer que toute série $\widehat{f}$ est $\chi$-équivalente à une somme finie de blocs $\chi$-homogènes (par des raisonnements abondamment employés dans la section II). Là apparaissent les premières difficultés de cette action: ces blocs $\chi$-homogènes $g_p(f)$ peuvent être formels, et même en cas de convergence d'un bloc $g_p$, les valeurs propres correspondantes $e_{m,n}$ peuvent avoir des points d'accumulation dans ${\bold{\ZZ}}\setminus\{ p \}$. \par

\vskip 3mm

\proclaim{Lemme IIIA2} Soit $e(\alpha)$ un germe analytique en 0, tel   que $e(0)\neq p$, et soit l'opérateur $E=\chi-e\text{Id}$. Pour tout $g\in \widehat{H{\Cal H}}_p$, on a $I_{\chi,Eg}=I_{\chi,g}$ dans l'anneau $c^*(\RR[[X,\alpha']])$.
\endproclaim

\vskip 3mm

\noindent {\bf Preuve.} Il suffit de montrer que $I_{\chi,g}\subset I_{\chi,Eg}$. Soit $G_k={\bold j}^k_u(g)$. Comme la dérivation $\chi$ est diagonale dans la coordonnée $u$, on a ${\bold j}^k_u(Eg)=EG_k$. Or le germe $G_k$ est dans le noyau de l'opérateur

$$
 \prod_{|m|-|n|=p,\quad |n|\leq k} (\chi -e_{m,n}Id) 
$$

\noindent comme $e(0)\neq p$, on montre par les mêmes arguments que dans le lemme II5 que $I_{\chi,EG_k}=I_{\chi,G_k}$. Soit ${\Cal M}$ l'idéal maximal de $c^*(\RR[[X,\alpha']])$. On a donc pour tout k

$$
 I_{\chi,g}\subset I_{\chi,Eg}+{\Cal M}^k
$$

\noindent on obtient le résultat par le théorème d'intersection de Krull.\qed\par

\vskip 3mm

\proclaim{Lemme IIIA3} Notons $f_{p_1,p_2}=\sum_{p_1\leq p\leq p_2} g_p(f)$, alors $I_{\chi,f_{p_1,p_2}}\subset I_{\chi,f_{p_1,p_2+1}}$.
\endproclaim

\vskip 3mm

\noindent {\bf Preuve.} Le germe ${\bold j}^k_u(f_{p_1,p_2})$ est dans le noyau de l'opérateur

$$
 \prod_{p_1\leq |m|-|n|\leq p_2;\quad |n|\leq k} (\chi -e_{m,n}Id) 
$$

\noindent et par le lemme IIIA2, on a 

$$
 I_{\chi,g_{p_2+1}}\subset I_{\chi,f_{p_1,p_2+1}}+{\Cal M}^k
$$

\noindent on conclut par le théorème d'intersection de Krull.\qed\par

\vskip 3mm

\proclaim{Lemme IIIA4} Il existe $p_1,p_2\in \ZZ$ tels que $\widehat{f}\equiv f_{p_1,p_2}\quad [c^*({\Cal M}_{X,u)})I_{\chi,f_{p_1,p_2}}]$.
\endproclaim

\vskip 3mm

\noindent {\bf Preuve.} D'après le lemme IIIA3, la suite d'idéaux $(I_{\chi,f_{-p,p}})$ est croissante et par le théorème d'intersection de Krull, sa limite est l'idéal différentiel de $\widehat{f}$. Donc pour tout $p\in\ZZ$, on a $I_{\chi,g_p(f)}\subset I_{\chi,\widehat{f}}$. Soit $k_1$ et $k_2$ les exposants d'Artin-Ress des idéaux ${\Cal M}_u$ et $c^*({\Cal M}_X)$ dans l'idéal $I_{\chi,\widehat{f}}$. Il suffit de prendre $p_1=-k_1$ et $p_2=k_2$.\qed\par

\vskip 3mm

\noindent{\bf §2. Multiplicité algébrique et la double inclusion.}\par

\vskip 3mm

 Soit $f\in QR{\Cal H}^{1,q}(x,\alpha,u)$. Dans l'action de $\chi_0$ sur $QR{\Cal H}^{1,.}$, la première inclusion est établie grâce à la structure asymptotique de $QR{\Cal H}^{1,.}$ (bien adaptée à $\chi_0$), et à la notion de multiplicité algébrique: c'est le plus petit entier $m$ tel que pour tout $\varepsilon>0$, $(x^{m+\varepsilon})\pi_{\chi_0}^*(J_{\chi_0,f,\gamma_0})\subset I_{\chi_0,f}$ dans l'anneau $SB^{1,.}$ correspondant. C'est aussi l'indice de stationnarité de la suite des idéaux $\chi_0$-transverses des éléments de la série de $f$ relativement à $\chi_0$. D'après l'étude de l'action formelle de $\chi$, cette structure asymptotique n'est plus adaptée à $\chi$ et les les sous-espaces stables $\widehat{H{\Cal H}}_p$ sont formels. Ceci motive la

\vskip 3mm

\proclaim{Définition IIIA2} Soit $\widehat{f}=\sum_{p\in\ZZ} g_p(f)$ la série de $f$ relativement à $\chi$. Le germe $f$ est dit {\bf quasi-convergent} si $g_p(f)\in QR{\Cal H}^{1,q}_{cvg}$ pour tout $p$. On note $QR{\Cal H}^{1,q}_{{\bold qcvg}}$ le $\RR\{\alpha\}$-module correspondant.
\endproclaim

\vskip 3mm

 Soit $f\in QR{\Cal H}^{1,q}_{qcvg}$. Les sommes finies correspondantes $f_{p_1,p_2}$ sont des éléments de l'anneau restriction analytique $QR{\Cal H}^{1,q}_{cvg}$. Les lemmes IIIA2 et IIIA3 s'appliquent sur cet anneau. Soit $I_f$ l'idéal limite différentiel de la suite $(I_{\chi,f_{-p,p}})_{p\in\ZZ}$, et soit $J_f$ l'idéal limite transverse de la suite $(J_{\chi,f_{-p,p},\gamma})_{p\in\ZZ}$. On obtient (en utilisant les lemmes IIIA2 et IIIA3, et les raisonnements de la section II)

$$
 I_f=\sum_{p\in\ZZ} I_{\chi,g_p(f)}\quad\text{et}\quad J_f=\sum_{p\in\ZZ} J_{\chi,g_p(f),\gamma}
$$

\noindent Il s'agit de comparer les idéaux $I_f$ et $I_{\chi,f}$ dans l'anneau $SB^{1,|q|}$, et les idéaux $J_f$ et $J_{\chi,f,\gamma}$ dans l'anneau $\RR\{\alpha,\lambda\}$. Soit ${\Cal M}'$ l'id\'eal maximal de $\RR\{\alpha,\lambda\}$ et soit $K({\Cal M}'_\lambda,J_f)$ l'exposant d'Artin-Ress de l'id\'eal ${\Cal M}'_{\lambda}=\langle\lambda_1,\ldots,\lambda_\ell\rangle\subset\RR\{\alpha,\lambda\}$, dans l'idéal $J_f$. Soit $p_0$ le plus petit des entiers $p\geq -K$ tel que $J_{\chi,g_p,\gamma}\not\subset {\Cal M}' J_f$.

\vskip 3mm

\proclaim{Définition IIIA3} {\bf La multiplicit\'e alg\'ebrique} $ma_{\chi}(f)$ est l'indice $p$ de stationnarit\'e  de la suite croissante des id\'eaux $(J_{\chi,f_{p_0,p},\gamma})_{p\geq p_0}$. La multiplicit\'e alg\'ebrique positive est $ma^+_{\chi}(f)=\max \{ ma_{\chi}(f),0 \}$.
\endproclaim

\vskip 3mm

\noindent Cette multiplicité algébrique est invariante dans les perturbations de $f$ dans l'idéal $\pi_\chi^*({\Cal M}' J_f)\cap QR{\Cal H}^{1,q}_{qcvg}$. Ceci permet d'étendre $QR{\Cal H}^{1,q}_{qcvg}$ en une classe d'éléments possédant la notion de multiplicité algébrique

\vskip 3mm

\proclaim{Définition IIIA4} Un germe $f\in QR{\Cal H}^{1,q}$ est dit {\bf presque quasi-convergent} s'il est limite dans la ${\Cal M}_{(\alpha,u)}$-topologie de $SB^{1,|q|}$, d'une suite $(f_n)_n$ de germes quasi-convergents, dont la suite des idéaux limites transverses $(J_{f_n})_n$ est croissante.
\endproclaim

\vskip 3mm

\proclaim{Lemme et définition IIIA5} Tout germe $f$ presque quasi-convergent possède une {\bf multiplicité algébrique} $ma_\chi(f)$ qui est la limite de la suite $(ma_\chi(f_n))_n$, et un idéal limite transverse $J_f$ qui est la limite de la suite $(J_{f_n})$.
\endproclaim

\vskip 3mm

\noindent{\bf Preuve.} Montrons que les suites $(ma_\chi(f_n))$ et $(J_{f_n})$ admettent des limites qui ne dépendent que de $f$. Soit $J$ l'idéal limite de la suite $(J_{f_n})$ d'indice de stationnarité $n_0$ et soit $K$ l'exposant d'Artin-Ress de l'idéal ${\Cal M}'_{(\alpha,\lambda)}$ dans $J$. Soit $k>K$, pour $n$ assez grand et $n'>n$, on a $f_{n'}-f_n\in {\Cal M}_{(\alpha,u)}^k$, et donc $g_p(f_{n'})-g_p(f_n)\in {\Cal M}_{(\alpha,u)}^k$ pour tout $p\in\ZZ$. Par conséquent, si $n\geq n_0$, on a $J_{\chi,g_p(f_{n'}),\gamma}\subset J_{\chi,g_p(f_n),\gamma}+{\Cal M}' J$ et 
$J_{\chi,g_p(f_n),\gamma}\subset J_{\chi,g_p(f_{n'}),\gamma}+{\Cal M}' J$. En sommant sur $p$ la première inclusion par exemple, on obtient 

$$
 J=\cup_{p\leq ma_\chi(f_{n'})} J_{\chi,g_p(f_{n'}),\gamma}\subset \cup_{p\leq ma_\chi(f_{n'})} J_{\chi,g_p(f_n),\gamma}+{\Cal M}' J
$$

\noindent le lemme de Nakayama et la définition de $J=J_{f_n}$ donnent $ma_\chi(f_n)\leq ma_\chi(f_{n'})$. La deuxième inclusion donne l'inégalité inverse. La suite $(ma_\chi(f_n))$ est donc stationnaire. Notons $ma$  sa limite.\par

\vskip 3mm

 Soit $(f'_n)$ une autre suite de germes quasi-convergents qui tend vers $f$ dans la ${\Cal M}_{(\alpha,u)}$-topologie et dont la suite des idéaux limites $(J_{f'_n})$ est croissante. Soit $J'$ son idéal limite d'indice de stationnarité $n'_0$ et soit $ma'$ la limite de la suite $(ma_\chi(f'_n))$. Montrons que $J'=J$ et $ma'=ma$. Pour tout $k$ et $n\geq\max(n_0,n'_0)$ assez grand, on a $f'_n-f_n\in {\Cal M}_{(\alpha,u)}^k$ et donc $g_p(f'_n)-g_p(f_n)\in{\Cal M}_{(\alpha,u)}^k$ pour tout $p$. Par conséquent, $J'\subset J+({\Cal M}')^k$ et $J\subset J'+({\Cal M}')^k$ pour tout $k$. Par le théorème d'intersection de Krull, on obtient $J'=J$. Notons $J_f$ cet idéal. En appliquant le raisonnement ci-dessus à $f'_n-f_n$ et à $J_f$, on obtient $ma'=ma$.\qed

\vskip 3mm

 Soit $W_0\subset W$ un semi-analytique de l'anneau $\RR\{\alpha,\lambda\}$ qui adhère à 0. Soit $U_0=\pi_\chi^{-1}(W_0)$. On définit de la même façon, une {\bf multiplicité algébrique restreinte} $ma_\chi(f_{|U_0})$ en considérant les idéaux $\chi$-transverses et leurs idéaux limites dans l'anneau restriction $\RR\{\alpha,\lambda\}_{|W_0}$. On verra dans la partie B que tout $f\in QR{\Cal H}^{1,q}$ admet une {\bf localisation finie} dans laquelle il est (à extension étoilée près), {\bf presque quasi-convergent}, et donc possède une {\bf multiplicité algébrique restreinte}.\par

\vskip 3mm

\noindent{\bf §2.1. Exemples de germes quasi-convergents.}\par

\vskip 3mm

  Omettons les exposants $1,q$ et $1,|q|$ pour un moment.\par
  
\vskip 3mm

\noindent{\bf (a) La sous-algèbre intégrale-linéaire } $\Lambda^*(QR{\Cal H}^{1,q}(x,\alpha,\lambda))$.\par

\vskip 3mm

 Soit $U_0\in (\RR^{+*}\times\RR^{|q|},0)$ et soit $\Lambda$ le difféomorphisme sur $U$ définie par

$$
 \Lambda(x,\alpha,u)=(x,\alpha,\lambda(x,\alpha,u))\in U_0
$$

\noindent L'image $\Lambda(U)$ étant ouverte, le morphisme $\Lambda^*:SB^{1,|q|}(x,\alpha,\lambda)\to SB^{1,|q|}(x,\alpha,u)$ est injectif, la sous-algèbre $QR{\Cal H}^{1,q}(x,\alpha,\lambda)$ est donc isomorphe à la sous-algèbre 

$$
 \Lambda^*(QR{\Cal H}^{1,q}(x,\alpha,\lambda)\subset QR{\Cal H}^{1,q}(x,\alpha,u)
$$

\noindent (l'inclusion s'obtient en remarquant que $x^{s_j}=x+\mu_jz(x,\mu_j)$). Le champ $\Lambda_*\chi$ coincide avec la restriction à $\Lambda(U)$ du champ $\chi_0$ et le diagramme suivant est commutatif\par

$$
\CD
SB(x,\alpha,\lambda)  @>\Lambda^*>>   SB(x,\alpha,u) \\
@V\chi_0 VV                      @VV\chi V \\
SB(x,\alpha,\lambda)  @>\Lambda^*>>   SB(x,\alpha,u) 
\endCD \tag {$*$}
$$

\vskip 3mm

\proclaim{Lemme IIIA6}\hskip 1cm $\Lambda^*(QR{\Cal H}(x,\alpha,\lambda))\subset QR{\Cal H}_{qcvg}$. 
\endproclaim

\vskip 3mm

\noindent{\bf Preuve.} Soit $F\in QR{\Cal H}(x,\alpha,\lambda)$ et soit $f=\Lambda^*(F)$. Soit
$\widehat{f}=\sum_{p\in\ZZ} g_p(f)$ la s\'erie formelle de $f$ relativement \`a $\chi$, et soit $\widehat{F}=\sum_{p\geq 0} g_{0,p}(F)$ la s\'erie formelle de $F$ relativement \`a $\chi_0$. Prolongeons le morphisme $\Lambda^*$ aux séries formelles: $\Lambda^*(\widehat{F})=\sum_p \Lambda^*(g_{0,p}(F))$. D'après le diagramme ci-dessus, l'image par $\Lambda^*$ d'un sous-espace stable par $\chi_0$, est incluse dans un sous-espace stable par $\chi$: plus précisément, pour $p\geq 0$, $\Lambda^*(H{\Cal H}_p)\subset\widehat{H{\Cal H}}_p\cap QR{\Cal H}_{cvg}$ . Donc, en utilisant les formules (1) et (2) du lemme IIIA1, qui définissent la série de $f$, on obtient $g_p(f)=\Lambda^*(g_{0,p}(F))\in QR{\Cal H}_{cvg}$, pour $p\geq 0$ et $g_p(f)=0$ pour $p<0$. En effet, écrivons

$$
g_{0,p}(F)=\sum_{|n|=p,|m|\geq 0} b_{n,m}(\alpha)\lambda^m X^n
$$

\noindent les séries $\sum_{|m|\geq 0} b_{n,m}(\alpha)\lambda^m$ étant convergentes sur un voisinage de 0. Soit

$$
\sum_{M\geq 0} F_M \quad\text{avec}\quad F_M=\sum_{|n|+|m|=M}b_{n,m}(\alpha)X^n\lambda^m
$$

\noindent la série de $F$ dans les variables $(X,\lambda)$. Elle est bien définie d'après la série $\widehat{F}$ et l'analycité dans les coordonnées $\lambda$. Le germe $\Lambda^*(F_M)$ est un bloc $\chi_0$-homogène de degré $M$. De plus, on vérifie facilement que pour tout $M$

$$
\Lambda^*(F-\sum_{M'\leq M}F_{M'})=o(x^M) \quad\text{dans l'anneau}\ SB
$$

\noindent Ceci prouve que $\Lambda^*(F_M)=g_{0,M}(f)$. Or, il est facile de voir que $g_p(g_{0,M}(f))=0$ si $p<0$ ou $p>M$, et que $g_p(g_{0,M}(f))=\sum_{|n|=p,|m|=M-p} \Lambda^*(b_{n,m}X^n\lambda^m)$ pour $0\leq p\leq M$. D'après la formule (2), on a

$$
g_p(f)=\sum_{M\geq p} g_p(g_{0,M}(f))=\sum_{M\geq p}(\sum_{|n|=p,|m|=M-p}\Lambda^*( b_{n,m}X^n\lambda^m))
$$

\noindent donc

$$
g_p(f)=\Lambda^*(g_{0,p}(F))
$$
 
\noindent et ceci finit la preuve du lemme.\qed\par

\vskip 3mm

  Soit $\gamma_0=\Lambda(\gamma)$, elle est principale dans $U_0$, et le morphisme $\Lambda$ est un difféomorphisme de $U$ sur un voisinage de $\gamma_0$. En identifiant les intégrales premières de $\chi$ et $\chi_0$, on obtient le diagramme commutatif

$$
\CD
 U     @>\pi_\chi >>    W \\
@V\Lambda VV           @VVid V \\
U_0    @>\pi_{\chi_0}>> W
\endCD \tag {$**$}
$$

\noindent Les idéaux $\chi$-transverse et $\chi_0$-transverse coincident (lemme de transfert IB6). Donc si $F\in QR{\Cal H}(x,\alpha,\lambda)$ et si $f=\Lambda^*(F)$, on a $J_{\chi,f,\gamma}=J_{\chi_0,F,\gamma_0}$ ($=J_F$ d'après le théorème principal II1). Or, on a $J_F=J_f$ d'après la définition de l'idéal limite transverse, et les séries de $f$ et $F$. Et donc $ma_\chi(f)=ma_{\chi_0}(F)\geq 0$. On a aussi $I_f=\Lambda^*(I_F)=\Lambda^*(I_{\chi_0,F})=I_{\chi,f}$ en utilisant le diagramme (*). Par conséquent, en appliquant le morphisme $\Lambda^*$ à la double inclusion 

$$
 (x^{ma_{\chi_0}(F)+\varepsilon})\pi_{\chi_0}^*(J_F)\subset I_{\chi_0,F}\subset\pi_{\chi_0}^*(J_F)
$$

\noindent et en utilisant le diagramme (**), on obtient la généralisation facile du théorème principal II1

\vskip 3mm

\proclaim{Lemme IIIA7} L'algèbre $\Lambda^*(QR{\Cal H}(x,\alpha,\lambda))$ est $\chi$-finie et satisfait à la double inclusion. Elle est $\chi$-équivalente à la sous-algèbre $\Lambda^*(QR{\Cal H}_{cvg}(x,\alpha,\lambda))$.
\endproclaim

\vskip 3mm

\noindent{\bf (b) La sous-algèbre algébrique }$QR{\Cal H}(x,\alpha)[u]$. Elle est clairement incluse dans $QR{\Cal H}_{qcvg}$ (les blocs $\chi$-homogènes sont algébriques dans les variales $(X,u)$). Soit $f\in QR{\Cal H}(x,\alpha)[u]$ de degré $N(f)$ en $u$, et soit $S(\mu)=N(f)\sum_{j=1}^\ell s_j$. Alors, $\widehat{f}=\sum_{p\geq -N(f)} g_p(f)$, et le germe 

$$
h=x^{S(\mu)}f
$$

\noindent s'identifie à un élément de l'algèbre $\Lambda^*(QR{\Cal H}(x,\alpha)[\lambda])$ (en utilisant les formules $x^{s_j}u_j=\lambda_j$). De plus, $g_p(h)=x^{S(\mu)}g_{p-S(0)}(f)$ pour tout $p\geq (\ell-1)N(f)$, et $g_p(h)=0$ si $p<(\ell-1)N(f)$.  La multiplicité algébrique de $h$ est (par sa définition), $ma_\chi(h)=ma_\chi(f)+\ell N(f)$. Donc en appliquant les résultats du {\bf (a)}, on obtient

\vskip 3mm

\proclaim{Lemme IIIA8} L'algèbre $QR{\Cal H}(x,\alpha)[u]$ est $\chi$-finie et satisfait à la double inclusion

$$
 (x^{ma^+_\chi(f)+\varepsilon})\pi_\chi^*(J_f)\subset I_{\chi,f}\quad\text{et}\quad (x^{\ell N(f)+\varepsilon})I_{\chi,f}\subset\pi_\chi^*(J_f)
$$

\noindent Elle est $\chi$-équivalente à la sous-algèbre $QR{\Cal H}_{cvg}(x,\alpha)[u]$.
\endproclaim

\vskip 3mm

\noindent Le facteur $x^{\ell N(f)}$ dans la deuxième inclusion, est optimal comme le montre l'exem\-ple suivant: $f=u^N x\log x$ ($\ell=1$ et $s=1$), l'idéal $\chi$-transverse de $f$ est $(\lambda^N)$. l'idéal saturé $\pi_\chi^*(J_{\chi,f,\gamma})$ est l'idéal $(x^N u^N)$. Le plus petit entier $n$ tel que $x^n f\in (x^N u^N)$ (dans l'anneau $SB$), est $N$.\par

\vskip 3mm

\noindent{\bf (c) La sous-algèbre convergente} $QR{\Cal H}_{cvg}$. Elle est $\chi$-finie d'après le lemme d'extension IB2 (même démarche que dans la section II). Le lemme IIIA4 s'applique à cette algèbre noethérienne: $I_f=I_{\chi,f}$ et $J_f=J_{\chi,f,\gamma}$. Soit $H{\Cal H}_p\subset QR{\Cal H}_{cvg}$ le $\RR\{\alpha\}$-module des blocs $\chi$-homogènes de degré $p\in\ZZ$. Si $\ell>0$, il est de dimension infinie et la méthode de réduction de la section II est inopérante. Une première approche à la première inclusion est

\vskip 3mm

\proclaim{Lemme IIIA9} Soit $g\in H{\Cal H}_p$ d'idéal $\chi$-transverse $J_g$ et soit $p^+=\max(p,0)$. Pour tout $\varepsilon>0$ et pour tout $N$ 

$$
 (x^{p^++\varepsilon})\pi_\chi^*(J_g)\subset I_{\chi,g}+\pi_\chi^*(({\Cal M}'_\lambda)^N)
$$
\endproclaim

\vskip 3mm

\noindent{\bf Preuve.} Par récurrence sur $\ell$. Le cas $\ell=0$ est donné par le lemme II4. Supposons $\ell>0$. L'idéal $I_{\chi,g}$ étant noethérien dans $QR{\Cal H}_{cvg}$, $g$ satisfait à une équation différentielle résolue

$$
 \chi^{M_0+1}g-\sum_{i=0}^{M_0}h_i\chi^i g=0
$$

\noindent Son idéal $\chi$-transverse est donc \ engendré \ par la \ restriction \ de \ la \ famille $(g,\ldots,\chi^{M_0}g)$ à toute transversale $\{x=x_0>0\}$. Quitte à effectuer une homothètie en $x$, on peut supposer que le polydisque de convergence de la série de $g$ (dans les variables $X$), contient la transversale $\{x=1\}$. Et quitte à effectuer une ramification $x'=x^{s_\ell}$, et à augmenter le nombre des variables $\mu$, on peut supposer que $s_\ell\equiv 1$. Ordonnons $g$ en puissances croissantes des fonctions élémentaires $z$: $g=\sum_{|m|\geq 0} a_m z^m$ et posons 

$$
g_M=\sum_{|m|\leq M} a_m z^m
$$

$$
u'=(u_1,\ldots,u_{\ell-1})\quad\text{et}\quad \lambda'=(\lambda_1,\ldots,\lambda_{\ell-1})
$$

\noindent On vérifie facilement (sur les monômes de $g_M$), que le germe $x^{M-p}g_M$ s'identifie à un élément $g'_M$ de l'algèbre $(\Lambda')^*(QR{\Cal H}(x,\alpha,\lambda_\ell,u'))$ avec 
$$
\Lambda'(x,\alpha,u)=(x,\alpha,\lambda_\ell,u')\quad\text{et}\quad \lambda_\ell=xu_\ell
$$

\noindent Comme on l'a vu au {\bf (a)}, l'action de $\chi$ sur cette algèbre est équivalente à celle de la dérivation

$$
 \chi'=x\frac{\partial}{\partial x}-\sum_{j=1}^{\ell-1} s_j u_j\frac{\partial}{\partial u_j}
$$

\noindent sur l'algèbre $QR{\Cal H}(x,\alpha,\lambda_\ell,u')$. Soit $G_M=((\Lambda')^*)^{-1}(g'_M)$. Il est $\chi'$-homogène de degré $M$. Par \ l'hypothèse \ de \ \ récurrence

$$
(x^{M+\varepsilon})\pi_{\chi'}^*(J_{G_M})\subset I_{\chi',G_M}+\pi_{\chi'}^*(({\Cal M}'_{\lambda'})^N)
$$

\noindent Or l'idéal $\chi'$-transverse de $G_M$ coincide avec l'idéal $\chi$-transverse de $g_M$ et le relevé de son idéal différentiel est $(x^{M-p})I_{\chi,g_M}$. Donc en relevant cette inclusion, on obtient 
$$
(x^{p^++\varepsilon})\pi_\chi^*(J_{g_M})\subset I_{\chi,g_M}+\pi_\chi^*(({\Cal M}'_\lambda)^N)
$$

  Maintenant, si $M\geq M_0$ et si on se restreint à la transversale $\{x=1\}$, on voit que l'idéal $\chi$-transverse de $g_M$ contient celui de $g$. Et si on prend $M\geq N+p^+$, on a $g-g_M\in\pi_\chi^*(({\Cal M}'_\lambda)^N)$. D'où le résultat.\qed

\vskip 3mm

\proclaim{Lemme IIIA10} Soit $f\in QR{\Cal H}_{cvg}$. Pour tout $\varepsilon>0$ et pour tout $N$ 

$$
 (x^{ma^+_\chi(f)+\varepsilon})\pi_\chi^*(J_f)\subset I_f+\pi_\chi^*(({\Cal M}'_\lambda)^N)
$$

\endproclaim

\vskip 3mm

\noindent{\bf Preuve.} Soit $p_0$ comme dans la définition IIIA3, et soient $p_1\leq p_0$ et $p_2\geq ma^+_\chi(f)$ tels que $I_f=I_{\chi,f_{p_1,p_2}}$. On a 

$$
 I_f=\sum_{p=p_1}^{p_2} I_{\chi,g_p}\supset \sum_{p=p_1}^{ma^+_{\chi}(f)} I_{\chi,g_p}\quad\text{et}\quad J_f=\sum_{p=p_1}^{p_2} J_{g_p}=\sum_{p=p_1}^{ma^+_{\chi}(f)} J_{g_p}
$$

\noindent on obtient le résultat en appliquant le lemme IIIA9 aux fonctions $g_p$.\qed

\vskip 3mm

  l'anneau $SB$ n'étant pas noethérien, on ne peut passer à la limite dans ces inclusions. Comme dans la section II, on contourne cette difficulté grâce au théorème de division VB2 (appendice VB), et grâce au lemme de Nakayama. On obtient le très important résultat suivant\par

\vskip 3mm

\proclaim{Théorème IIIA1 (théorème principal 2)} L'algèbre $QR{\Cal H}_{cvg}$ satisfait globalement à la double inclusion: pour tout $f\in QR{\Cal H}_{cvg}$ de multiplicité algébrique $ma^+$, il existe un entier $n(f)$ tel que pour tout $\varepsilon>0$

$$
(x^{ma^++\varepsilon})\pi_\chi^*(J_f)\subset I_{\chi,f}\quad\text{et}\quad (x^{n(f)})I_{\chi,f}\subset \pi_\chi^*(J_f)
$$

\endproclaim

\vskip 3mm

\noindent{\bf Preuve.} La difficulté réside dans la transcendance des intégrales premières non triviales de $\chi$. Grâce à la convergence des séries d'éléments de $QR{\Cal H}_{cvg}$, on montre d'abord que l'action de $\chi$ sur cette algèbre est équivalente à celle d'une dérivation d'intégrales premières non triviales algébriques, sur une autre algèbre convergente. Le théorème de division VB2 s'applique à cette nouvelle dérivation.\par

\vskip 3mm

  Généralisons le difféomorphisme $\Lambda$ en les difféomorphismes $\Lambda_{s}$ suivants: soit $s\geq 0$ et $\Lambda_{s}:(x,\alpha,u)\in U\mapsto (y,\alpha,v)\in U_s$ avec

$$
 y=x^{\frac{1}{s+1}}\quad\text{et}\quad v_j=x^{\frac{1}{s+1}+\mu_j}u_j
$$

\noindent Le champ $(s+1)(\Lambda_s)_*\chi$ coincide avec la restriction à $\Lambda_s(U)$ du champ 

$$
 {\Cal Y}_s=y\frac{\partial}{\partial y}-s\sum_{j=1}^{\ell}v_j\frac{\partial}{\partial v_j}
$$

\noindent Le morphisme $\Lambda_s^*:SB(y,\alpha,v)\to SB(x,\alpha,u)$ est un isomorphisme sur son image et le diagramme suivant est commutatif

$$
\CD
SB(y,\alpha,v)  @>\Lambda_s^*>>   SB(x,\alpha,u) \\
@V{\frac{1}{s+1}{\Cal Y}_s} VV                      @VV\chi V \\
SB(y,\alpha,v)  @>\Lambda_s^*>>   SB(x,\alpha,u) 
\endCD \tag {$*$}
$$

 Soit $\gamma_s=\Lambda_s(\gamma)$, elle est principale dans $U_s$ et le morphisme $\Lambda_s$ est un difféomorphisme de $U$ sur $\Lambda_s(U)$ qui est un voisinage de $\gamma_s$. En identifiant les intégrales premières de $\chi$ et ${\Cal Y}_s$, on obtient le diagramme commutatif

$$
\CD
 U     @>\pi_\chi >>    W \\
@V\Lambda_s VV           @VVid V \\
U_s    @>\pi_{{\Cal Y}_s}>> W
\endCD \tag {$**$}
$$

\noindent Les idéaux $\chi$-transverses et ${\Cal Y}_s$-transverses coincident.\par 

\vskip 3mm

  Si $s$ est un entier$\geq 2$, un calcul simple sur les monômes montre que l'image $\Lambda_s^*(SB(y,\alpha,v))$ contient les $\RR\{\alpha\}$-modules $H{\Cal H}_p(x,\alpha,u)$ pour $p\geq 0$. Et il existe un morphisme linéaire: $\alpha\mapsto\beta(\alpha)$ tel que l'image r\'eciproque $(\Lambda_s^*)^{-1}(H{\Cal H}_p(x,\alpha,u))$ est  incluse dans le $\RR\{\beta\}$-module $H{\Cal H}_{(s+1)p}(y,\beta,v)$ des blocs ${\Cal Y}_s$-homog\`enes de degr\'e $(s+1)p$, restreint à l'image $\beta(\alpha)$. En effet, si $z_j$ est une fonction élémentaire de l'algèbre $QR{\Cal H}(x,\alpha,u)$ qui satisfait à l'équation $\chi_0 z_j=r_jz_j+x$, la fonction $Z_j=y^{-s}(\Lambda_s^*)^{-1}(z_j)$ satisfait à l'équation ${\Cal Y}_0 Z_j=(r_j+s(r_j-1))Z_j+(s+1)y$. Si $m$ et $n$ sont des multi-indices tels que $|m|-|n|=p\geq 0$

$$
 (\Lambda_s^*)^{-1}(X^m u^n)=y^{m_0+s|m|-2|n|}\prod Z_j^{m_j}\prod_{j=1}^\ell y^{n_j(1-(s+1)\mu_j)} v_j^{n_j}
\tag $2'$
$$

\noindent Maintenant, il suffit de remarquer que $y^{1-(s+1)\mu_j}=y-(s+1)\mu_j y\text{Ld}(y,-(s+1)\mu_j)$ et d'utiliser la convergence des séries d'éléments de $H{\Cal H}_p$.\par

\vskip 3mm

  Soit $f\in QR{\Cal H}_{cvg}(x,\alpha,u)$ et soient $p_1$ et $p_2$ comme dans le lemme IIIA10. Quitte \`a remplacer $f$ par $x^{-p_1}f$, on peut supposer que $p_1\geq 0$. Soit $s$ un entier$\geq 2$ et $F=(\Lambda_s^*)^{-1}(f_{p_1,p_2})$. C'est un élément de l'algèbre $QR{\Cal H}(y,\beta,v)_{cvg}$ retreinte à l'image $\beta(\alpha)$. La multiplicité algébrique de $F$ relativement à ${\Cal Y}_s$ est $(s+1)ma_\chi(f)$. Le lemme IIIA10 s'applique à l'action de ${\Cal Y}_s$ sur cette algèbre: pour tout $\varepsilon>0$ et pour tout $N$

$$
 (y^{(s+1)ma^++\varepsilon})\pi_{{\Cal Y}_s}^*(J_F)\subset I_{{\Cal Y}_s,F}+\pi_{{\Cal Y}_s}^*(({\Cal M}'_\lambda)^{N})
$$

 Le théorème de division VB2 s'applique à l'action de ${\Cal Y}_s$ sur l'anneau $SB(y,\alpha,v)$. Soit $n(F)$ l'entier donné par ce théorème, il ne dépend que de $J_F=J_f$ et de $s$. Le théorème de division VB2 appliqué à $F$ fournit la deuxième inclusion

$$
 (y^{n(F)}) I_{{\Cal Y}_s,F}\subset \pi_{{\Cal Y}_s}^*(J_F)
\tag 3
$$

  Soit $(\varphi_1,\ldots,\varphi_L)$ un système de générateurs de $J_F$ dans l'anneau $\RR\{\alpha,\lambda\}$ et soit $N>(s+1)ma^++n(F)+1$. Pour tout $i=1,\ldots,L$ l'inclusion ci-dessus montre qu'il existe $h_i\in SB(y,\alpha,v)$ tel que

$$
 y^{(s+1)ma^++\varepsilon}\pi_{{\Cal Y}_s}^*(\varphi_i)-y^N h_i\in I_{{\Cal Y}_s,F}
\tag 4
$$

\noindent Ceci implique que l'idéal ${\Cal Y}_s$-transverse de $h_i$ est inclus dans $J_F$. Donc par le théorème de division VB2, les inclusions (4) donnent

$$
 (y^{(s+1)ma^++\varepsilon})\pi_{{\Cal Y}_s}^*(J_F)\subset I_{{\Cal Y}_s,F}+(y^{(s+1)ma^++1})\pi_{{\Cal Y}_s}^*(J_F)
$$

\noindent et par le lemme de Nakayama, on obtient

$$
 (y^{(s+1)ma^++\varepsilon})\pi_{{\Cal Y}_s}^*(J_F)\subset I_{{\Cal Y}_s,F}
\tag 5
$$

\noindent On applique alors le morphisme $\Lambda_s^*$ aux inclusions (3) et (5) et la commutativité du diagramme $(*)$ pour obtenir le résultat.\qed

\vskip 3mm

\noindent{\bf Remarque IIIA1.} De cette preuve, et du théorème de division VB2, on déduit le résultat suivant: si $h\in QR{\Cal H}^{1,q}_{cvg}$, et si son idéal $\chi$-transverse est inclus dans celui de $f$, alors

$$
(x^{n(f)})I_{\chi,h}\subset \pi_\chi^*(J_{\chi,f,\gamma})
$$

\vskip 3mm

\noindent{\bf Remarque IIIA2.} Soit $W_0\subset W$ un semi-analytique de l'anneau $\RR\{\alpha,\lambda\}$, qui adhère à 0. Soit $U_0=\pi_\chi^{-1}(W_0)$. Il est clair que le lemme IIIA9 est vrai sur les restrictions $W_0$ et $U_0$ (la multiplicité algébrique restreinte coincide avec la multiplicité algébrique d'un bloc $g_p$). Par conséquent, le lemme IIIA10 et le théorème principal IIIA1, sont encore vrais sur les restrictions $W_0$, $U_0$, en utilisant la multiplicité algébrique restreinte $ma_\chi(f_{|U_0})$.\par

\vskip 3mm

\noindent{\bf §2.2. Exemple de germe presque quasi-convergent.}\par

\vskip 3mm

 Soit $f\in QR{\Cal H}^{1,q}(x,\alpha)$ d'idéal $\chi$-transverse $J_{\chi,f,\gamma}$. On dit que $f$ (ou $J_{\chi,f,\gamma}$) satisfait à l'hypoth\`ese $(H\lambda)$ s'il existe un entier $N(f)$ tel que  
 
$$
J_{\chi,f,\gamma}\supset ({\Cal M}'_{\lambda})^{N(f)}
$$

  On note ${\Cal C}^1_\lambda$ la classe des germes qui satisfont à l'hypothèse $(H\lambda)$. Tout $f\in {\Cal C}^1_\lambda$ est  presque quasi-convergent: soit $f_n={\bold{j}}_u^{N(f)+n}(f)$, il est quasi-convergent et par le lemme IIIA8, $J_{f_n}=J_{\chi,f_n,\gamma}$. Par l'hypothèse $(H\lambda)$, on a $J_{\chi,f,\gamma}\subset J_{f_n}+{\Cal M}' J_{\chi,f,\gamma}$. Donc, par le lemme de Nakayama, on a $J_{\chi,f,\gamma}=J_{f_n}$. De plus, la suite $(f_n)$ converge vers $f$ dans la ${\Cal M}_{(\alpha,u)}$-topologie. Par conséquent, $f$ possède une multiplicité algébrique $ma_\chi(f)$ et $J_f=J_{\chi,f,\gamma}$. 

\vskip 3mm

\proclaim{Lemme IIIA11} La classe ${\Cal C}^1_\lambda$ est $\chi$-finie et satisfait à la double inclusion

$$
 (x^{ma^+_\chi(f)+\varepsilon})\pi_\chi^*(J_f)\subset I_{\chi,f}\quad\text{et}\quad (x^{\ell N(f)+\varepsilon})I_{\chi,f}\subset\pi_\chi^*(J_f)
$$

\noindent De plus, elle est $\chi$-équivalente à la sous-classe $QR{\Cal H}^{1,q}_{cvg}\cap {\Cal C}^1_\lambda$.
\endproclaim

\vskip 3mm

\noindent{\bf Preuve.} Par l'hypothèse $(H\lambda)$, $x^{\ell N(f)+\varepsilon}(f-f_0)\in\pi_\chi^*(J_f)$. Et par le lemme IIIA8, $(x^{\ell N(f)+\varepsilon})I_{\chi,f_0}\subset \pi_\chi^*(J_f)$, d'où la deuxième inclusion pour $f$. Soit $F={\bold j}_X^{\ell N(f)+ma^+}(f-f_0)$, on a $f-f_0-F\in (x^{ma^++\varepsilon})\pi_\chi^*(J_f)$, et l'idéal $\chi$-transverse de $F$ est inclus dans ${\Cal M}'J_f$. Le germe $f_0$ a la même multiplicité algébrique que $f$, et d'après le lemme IIIA8, il est $\chi$-équivalent au germe

$$
 h=\sum_{-N(f)\leq p\leq ma^+} g_p(f_0)
$$

\noindent De plus on a $f_0-h\in (x^{ma^++\varepsilon})\pi_\chi^*(J_f)$. Maintenant, le germe $h+F\in QR{\Cal H}_{cvg}$ a la même multiplicité algébrique que $f$, et il est algébrique dans la variable  $X$, il est donc $\chi$-équivalent à un germe $H$ algébrique dans les variables $(X,u)$ (cf. {\bf (c)}). D'après le lemme IIIA8, le germe $H$ satisfait à la première inclusion. Comme $f-(h+F)\in (x^{ma^++\varepsilon})\pi_\chi^*(J_f)$, le germe $f$ satisfait à la première inclusion et il est $\chi$-équivalent à $H$.\qed

\vskip 3mm

 Ce lemme est encore vrai pour la classe ${\Cal C}^1_{\lambda,loc}$ des germes $f$ qui satisfont à l'hypothèse $(H\lambda)$ sur une restriction $W_0(f)\subset W$.\par

\vskip 5mm

\noindent{\bf B. Etude localisée et théorème principal 3.}\par

\vskip 3mm

 Pour $\ell=0$, on a montré dans la section II, que l'algèbre $QR{\Cal H}^{1,q}$ est $\chi_0$-finie, et que tout élément $f$ se divise dans son idéal $\chi$-transverse $J_{\chi_0,f,\gamma_0}$, dans l'anneau $QA^{1,|q|}[\Omega_f]$. La question qui se pose alors est la suivante: $f$ se divise-t-elle dans $J_{\chi_0,f,\gamma_0}$ dans un sous-anneau de $QA^{1,|q|}$, qui possède une {\bf structure asymptotique élémentaire}? La réponse est oui modulo une désingularisation de $J_{\chi_0,f,\gamma_0}$ (voir ci-dessous).\par
 
\vskip 3mm

 Pour $\ell>0$, on a vu dans la partie A que si le germe $f$ est presque quasi-convergent, il existe un entier $n$ tel $x^n f$ se divise dans son idéal $\chi$-transverse (qu'on note simplement $J_f$), disons dans l'anneau $SB$ (on peut montrer que la division se fait dans l'anneau $QA$, mais ceci ne sera pas utilisé). Je pense qu'il est possible de généraliser le théorème de division VB2 à l'action de la dérivation $\chi$ (ie. avec valeurs propres $s_j(\mu)$ transcendantes), en utilisant encore les techniques de la référence [B-M]. Cependant, si le germe $f$ est général, sa série formelle relativement à $\chi$ (cf. lemme IIIA1), est {\bf doublement formelle}. D'où la difficulté de pouvoir approcher son idéal $\chi$-transverse $J_f$, par l'intermédiaire d'idéaux $\chi$-transverses de {\bf germes simples}.\par
 
\vskip 3mm

 Si cet idéal $J_f$ était {\bf principal et monomial}, sa reconnaissance dans la série de $f$, serait une simple étude sur les séries Tayloriennes. D'où l'idée du lemme IIIB1 ci-dessous, basé sur une {\bf désingularisation d'Hironaka}. Grâce à cette désingularisation, on montre le théorème principal suivant\par

\vskip 3mm

\proclaim{Théorème IIIB1 (théorème principal 3)} L'algèbre $QR{\Cal H}^{1,q}$ est localement $\chi$-finie.
\endproclaim

\vskip 3mm

 Néanmoins, cette désingularisation donne un caractère fortement technique à la preuve du théorème. Il est donc très souhaitable d'avoir un résultat global (de $\chi$-finitude), par les méthodes algébrico-géométriques des sections II et IIIA.\par
 
\vskip 3mm

\proclaim{Lemme IIIB1} Soit $J$ un idéal de l'anneau $\RR\{\alpha,\lambda\}$. Il existe une désingularisation d'Hironaka $(\psi,{\Cal N})$ telle que dans chaque carte $(a,V_a)$, l'idéal $\psi_a^*(J)$ est principal et monomial. Autrement dit, il existe une coordonnée $v$ sur $V_a$, et un multi-indice $n$, tels que $\psi_a^*(J)=(v^n)$.
\endproclaim

\vskip 3mm

\noindent{\bf Preuve.} Elle est basée sur une récurrence sur le nombre de générateurs (indépendants) de $J$. Soit $(\varphi_1,\ldots,\varphi_L)$ un système de générateurs de $J$. Si $L=1$, par une désingularisation d'Hironaka ([Hir1,2]), le génératuer $\varphi_1$ se relève localement en un germe qui est équivalent à un monôme. Si $L>1$, on applique une désingularisation d'Hironaka au germe

$$
 \varphi_1\times\cdots\times\varphi_L\times(\varphi_L-\varphi_{L-1})
$$

\noindent s'il est non identiquement nul. Localement dans cette désingularisation, les trois germes $\varphi_{L-1}$, $\varphi_L$ et $\varphi_L-\varphi_{L-1}$ sont équivalents à des monômes. Donc, en utilisant l'ordre lexicographique sur les monômes, on a ou bien $(\varphi_{L-1})\subset (\varphi_L)$, ou bien $(\varphi_L)\subset (\varphi_{L-1})$ et on applique l'hypothèse de récurrence.\qed

\vskip 3mm

\noindent{\bf Preuve du théorème IIIB1.} Pour bien illustrer les étapes cruciales de cette preuve, commençons par le cas $\ell=0$\par

\vskip 3mm

\noindent {\bf §1. Cas} $\ell=0$.\par

\vskip 3mm

 Soit $f\in QR{\Cal H}^{1,q}(x,\alpha)$ ($q=(q_1,q_2)$), et soit $J_f$ son idéal $\chi$-transverse le long de $\gamma$ ($\chi=\chi_0$ et $\gamma=\gamma_0$). Montrons que $f$ {\bf se divise localement dans} $J_f$ {\bf dans un sous-anneau de} $QA$ {\bf qui admet une structure asymptotique élémentaire}. Soit $\widehat{f}=\sum_{n\geq 0} g_n$ sa série asymptotique formelle relativement à la dérivation $\chi$ (section II). Les fonctions $g_n$ sont des blocs $\chi$-homogènes de degré $n$, autrement dit si 

$$
\text{Ld}(x,\mu_j)=\frac{x^{\mu_j}-1}{\mu_j}
$$

\noindent $z_j(x,\mu)=x\text{Ld}(x,\mu_j)$, $z=(z_1,\ldots,z_{q_1})$ et $X=(x,z)$, alors 

$$
g_n=\sum_{|m|=n} a_m(\alpha) X^m
$$

\noindent avec $a_m\in \RR\{\alpha\}$. On a vu dans la section II que si on note $f_n=\sum_{n'\leq n} g_{n'}$, alors $f-f_n\in (x^n)SB_0$ pour tout $n$. Soient $J_{g_n}$, $J_{f_n}$ et $J_f$ les idéaux $\chi$-transverses le long de $\gamma$. Ce sont des idéaux de l'anneau $\RR\{\alpha\}$. La suite $(J_{f_n})$ est croissante et elle converge vers $J_f$ grâce à la propriété de quasi-analycité dans la coordonnée $x$. De plus, pour tout $n$, on a $f-f_n\in (x^n)\pi_\chi^*(J_f)$ dans l'anneau $SB_0$. Et, $J_{f_n}=\sum_{n'\leq n} J_{g_{n'}}$. Donc, par unicité de la division (voir ci-dessous), il suffit de montrer que pour tout $n$, le bloc $g_n$ se divise localement dans $J_f$, dans un sous-anneau de $QA$ qui possède une structure asymptotique élémentaire (ce sous-anneau étant indépendant de $n$!).\par

\vskip 3mm

 Soit $(\psi,{\Cal N})$ une désingularisation dans laquelle l'idéal $J_f$ est principal et monomial (lemme IIIB1), et soit $(a,V_a)$ une carte de coordonnée $v$ telle que 
 
$$
\psi_a^*(J_f)=(\varphi(v))\quad\text{avec}\quad \varphi(v)=\prod_{j=1}^p v_j^{n_{a,j}}
$$

\noindent On peut supposer, sans perte de généralité, que $p=|q|$. Les variables analytiques de grande complexité dans l'algèbre $QR{\Cal H}^{1,q}(x,\mu,\nu)$, sont les variables $\mu$. Par une suite d'éclatements sphériques dans la coordonnée $v$, on suppose que les relevés $\mu_{j,a}$ des fonctions $\mu_j$ sont préparées de la façon suivante: pour tout $i=1,\ldots,p$

$$
 \mu_{j,a}(v)=\mu_{j,p+1-i}+\mu'_{j,p+1-i}
\tag 1
$$

\noindent les fonctions $\mu_{j,p+1-i}$ étant {\bf indépendantes des coordonnées}  $(v_i,\ldots,v_p)$ et 

$$
 \mu'_{j,p+1-i}\ \ \in (v_1\times\cdots\times v_i)
$$

\noindent En effet, effectuons un premier éclatement sphérique dans la coordonnée $v$: $v=t_1 w$ où $t_1$ est une coordonnée locale sur $(\RR,0)$, et $w\in S^{p-1}$ la sphère de $\RR^p$. Soit $v^1$ une coordonnée analytique locale sur $(S^{p-1},w)\cong(\RR^{p-1},0)$. Le relevé du germe $\mu_{j,a}$ s'écrit

$$
\rho_1+t_1\theta_1(t_1,v^1)
$$

\noindent (avec $\rho_1\equiv 0$!). Supposons que, après $i$ éclatements sphériques des coordonnées $v^j$ ($j=0,\ldots,i-1$, $v^0=v$), le relevé du germe $\mu_{j,a}$ s'écrive

$$
\rho_i(t_1,\ldots,t_{i-1})+t_1\times\cdots\times t_i\theta_i(t^i,v^i)
$$

\noindent où $t^i=(t_1,\ldots,t_i)$ est une coordonnée analytique sur $(\RR^i,0)$, et $v^i$ est une coordonnée analytique sur $(\RR^{p-i},0)$. On applique alors la première étape au germe $\theta_i(0,v^i)$ (si $p-i>1$), et on pose $\rho_{i+1}=\rho_i+t_1\times\cdots\times t_i\theta_i(t^i,0)$. On obtient les formules (1) en posant 

$$
\mu_{j,p+1-i}=\rho_i\quad\text{et}\quad \mu'_{j,p+1-i}=t_1\times\cdots\times t_i\theta_i
$$

 Notons toujours  $(\psi,{\Cal N})$ la composée de ces deux désingularisations (l'idéal $(\varphi)$ est toujours monomial dans cette désingularisation).  Notons ${\Cal B}$ l'une des algèbres $SB$, $QA$ ou $QR{\Cal H}$ (sans préciser les exposants). Dans la carte $(a,V_a)$, une extension naturelle de l'anneau ${\Cal B}(x,\alpha)$, appropriée à $f$, est l'anneau $({\Cal B}(x,\alpha,v),\pi)$, où $\pi$ est (le germe de) la projection canonique: $(x,\alpha,v)\in U\times V_a\mapsto (x,\alpha)\in U$ ($U\in(\RR^{+*}\times\RR^{|q|},0)$ est un ouvert produit). Précisons un peu, dans ce cas simple, les propriétés du transfert par le morphisme $\pi$. Notons toujours $\chi=x\partial/\partial x$ la dérivation définit sur $U\times V_a$. Le morphisme $\pi$ est une submersion, et on a $\pi_*\chi=\chi$ (la coordonnée $v$ est une intégrale première de $\chi$, et la restriction de $\pi$ aux fibres $\{v=const\}$ est un difféomorphisme sur son image $U$). L'orbite $\Gamma=\{(\alpha,v)=0\}$ est principale dans $U\times V_a$, et $\pi(\Gamma)=\gamma$. Notons $\pi_0$ la restriction de $\pi$ à une transversale $\{x=x_0>0\}$. L'idéal $\chi$-transverse de $\pi^*(f)=f$, est l'idéal prolongé $\pi_0^*(J_f)\subset\RR\{\alpha,v\}$, il est engendré par un système de générateurs de $J_f$ dans $\RR\{\alpha\}$. Le morphisme $\Psi_a:\ v\in V_a\mapsto (\psi_a(v),v)\in W_a$ est une immersion analytique, qui est un difféomorphisme sur son image $W_a\subset\{x=x_0\}$. L'anneau restriction $\RR\{\alpha,v\}_{|W_a}$ est donc isomorphe à l'anneau $\RR\{v\}$. Par conséquent, on a $\pi_0^*(J_f)_{|W_a}=(\varphi)$.\par
 
\vskip 3mm
 
 De même, posons $U_a=\pi_\chi^{-1}(W_a)$, le morphisme $\Phi_a:\ (x,v)\in(\RR^{+*},0)\times V_a\mapsto (x,\Psi_a(v))\in U_a$ est une immersion analytique, qui est un difféomorphisme sur son image $U_a$. Donc, si ${\Cal B}$ est l'une des algèbres $SB$ ou $QA$, l'anneau restriction ${\Cal B}(x,\alpha,v)_{|U_a}$ est isomorphe à l'anneau ${\Cal B}(x,v)$. La dérivation $\chi$ est préservée par ce difféomorphisme. Soit $g$ un bloc $\chi$-homogène de la série de $f$. Notons $f_a=\Phi_a^*(f_{|U_a})$ et $g_a=\Phi_a^*(g_{|U_a})$. D'après le lemme de transfert IB6, l'idéal $\chi$-transverse de $f_a$ le long de l'orbite $\{v=0\}$, est $(\varphi)$. On a $g\in\pi_\chi^*(J_f)$, donc $g_a\in\pi_\chi^*((\varphi))$. Il s'agit de montrer que le quotient de la division de $g_a$ par $\varphi$, est dans un sous-anneau de $QA^{1,p}(x,v)$ qui possède une structure asymptotique élémentaire. {\bf La complexité de cette structure ne dépendra que du germe} $f$ par l'intermédiaire de son idéal $\chi$-transverse $J_f$. Pour cela, on va construit une suite de $p$ extensions $\widetilde{QR{\Cal H}}_i^{1,p}(x,v)$ de l'algèbre $QR{\Cal H}^{1,q+(0,p)}(x,\alpha,v)_{|U_a}$ (rappelons que $p=|q|$). La récurrence de cette construction est basée sur la récurrence de construction de nouvelles fonctions élémentaires d'Ecalle-Khovanski $z^{...}$, $y^{\ldots}$ à partir des fonctions élémentaires $z$ de l'algèbre $QR{\Cal H}^{1,q}$, et des formules $(1)$.\par

\vskip 3mm

\noindent{\bf (a) Construction des extensions }$\widetilde{QR{\Cal H}}_i^{1,p}$.\par

\vskip 3mm

 Pour $j=1,\ldots,q_1$, posons $r_{j,0}(v)=1+\mu_{j,a}(v)$ et $r_{j,1}(v_1,\ldots,v_{p-1})=1+\mu_{j,1}(v_1,\ldots,\ \ v_{p-1})$. D'après (1), le germe $\tau_{j,1}=\mu_{j,a}-\mu_{j,1}=\mu'_{j,1}\in rad(\varphi)=(v_1\times\cdots\times v_p)$. Pour tout $m_1\in \NN$, développons les (relevés des) fonctions $z_j$ à l'ordre $m_1$ dans la variable 
$\tau_{j,1}=\mu_{j,a}-\mu_{j,1}$

$$
 \Phi_a^*(z_j)=\sum_{n=0}^{m_1} \tau_{j,1}^n z_{j,n} +\tau_{j,1}^{m_1+1} y_{j,m_1}
\tag $2_1$
$$

\noindent En utilisant l'équation $\chi \Phi_a^*(z_j)=r_{j,0}\Phi_a^*(z_j)+x$ et en identifiant les coefficients des puissances de $\tau_{j,1}$, on obtient des équations différentielles simples pour les fonctions $z_{j,n}$ et $y_{j,m_1}$

$$
 \chi z_{j,0}=r_{j,1}z_{j,0} + x
$$

$$
 \chi z_{j,n}=r_{j,1}z_{j,n} + z_{j,n-1}
\tag $3_1$
$$

$$
 \chi y_{j,m_1}=r_{j,0} y_{j,m_1} + z_{j,m_1}
$$

\vskip 2mm

\noindent avec les conditions initiales $z_{j,n|x=1}\equiv 0$ et $y_{j,m_1|x=1}\equiv 0$. Comme pour les fonctions élémentaires $z_j$, on montre grâce à l'opérateur intégral de Dulac (cf. appendice VA), que $z_{j,n}\in QA^{1,p-1}(x,v_1,\ldots,v_{p-1})$ et $y_{j,m_1}\in QA^{1,p}(x,v)$. Notons $z^n=(z_{1,n},\ldots,z_{q_1,n})$, $y^{m_1}=(y_{1,m_1},\ldots,y_{q_1,m_1})$ et $X_{1,m_1}=(x,z^0,\ldots,z^{m_1},y^{m_1})$. Un germe $G_{1,n}\in SB^{1,p}(x,v)$ est un {\bf 1-bloc }$\chi${\bf -homogène de degré} $n$ {\bf et de complexité} $m_1$, si

$$
 G_{1,n}=\sum_{|m'|=n} a_{m'}(v) X_{1,m_1}^{m'}
$$

\noindent avec $a_{m'}\in\RR\{v\}$. Le germe $G_{1,n}$ est donc un élément de l'algèbre $QA^{1,p}(x,v)$. La première extension $\widetilde{QR{\Cal H}}_{1,m_1}^{1,p}(x,v)\subset QA^{1,p}(x,v)$ est {\bf l'algèbre des germes qui possèdent une série asymptotique formelle dans les 1-blocs} $\chi${\bf -homo\-gènes}, de complexité $m_1$: $F\in QA^{1,p}$ est un élément de $\widetilde{QR{\Cal H}}_{1,m_1}^{1,p}$ s'il existe une suite $(G_{1,n})$ de 1-blocs $\chi$-homogènes de degré $n$ et de complexité $m_1$, telle que pour tout $n$

$$
F-\sum_{n'\leq n} G_{1,n'}\in (x^n)SB_0^{1,p}
$$
 
\noindent On note $\widetilde{QR{\Cal H}}_1^{1,p}$ la {\bf limite inductive} de ces algèbres, quand $m_1$ décrit $\NN$.\par

\vskip 3mm

 Soit $r_{j,2}=1+\mu_{j,2}$. D'après les formules (1), le germe $\tau_{j,2}=\mu_{j,1}-\mu_{j,2}\in (v_1\times\cdots\times v_{p-1})$.  Soit $m_2\in\NN$ et soit $n_1\leq m_1$. En utilisant les équations $(3_1)$, on vérifie facilement que les coefficients $z_{j,n_1,n_2}$ et $y_{j,n_1,m_2}$ du développement $(2_2)$ à l'ordre $m_{2}$ de la fonction $z_{j,n_1}$ dans la variable $\tau_{j,2}$, satisfont des équations différentielles $(3_2)$ du même type que les équations $(3_1)$, de valeurs propres $r_{j,1}$ ou $r_{j,2}$. De plus, ces coefficients vérifient des relations linéaires déduites des équations

$$
 \frac{\partial^n z_{j,0}}{\partial \tau_{j,1}^n}=n!z_{j,n}
$$

\noindent Comme ci-dessus, on montre que $z_{j,n_1,n_2}\in QA^{1,p-2}(x,v_1,\ldots,v_{p-2})$ et $y_{j,n_1,m_{2}}\in QA^{1,p-1}(x,v_1,\ldots,v_{p-1})$. Pour $n_1=0,\ldots,m_1$ et $n_2=0,\ldots,m_{2}$, notons $z^{n_1n_2}=(z_{1,n_1,n_2},\ldots,z_{q_1,n_1,n_2})$,   $y^{n_1,m_{2}}=(y_{1,n_1,m_{2}},\ldots,y_{q_1,n_1,m_{2}})$ et $y^{m_1m_2}=(y^{n_1,m_{2}})$. Soit 

$$
 X_{2,m_1,m_2}=(x,(z^{n_1n_2}),y^{m_1},y^{m_1m_2})
$$

\noindent Un {\bf 2-bloc} $\chi${\bf -homogène de degré} $n$ {\bf et de complexité} $(m_1,m_2)$, est un germe de la forme

$$
 G_{2,n}=\sum_{|m'|=n} a_{m'}(v) X_{2,m_1,m_2}^{m'}
$$

\noindent avec $a_{m'}\in\RR\{v\}$. La deuxième extension $\widetilde{QR{\Cal H}}^{1,p}_{2,m_1,m_2}(x,v)\subset QA^{1,p}(x,v)$ est l'alg\-èbre des germes qui possèdent un développement asymptotique dans les 2-blocs $\chi$-homogènes de complexité $(m_1,m_2)$. On note $\widetilde{QR{\Cal H}}^{1,p}_2$ la limite inductive de ces algèbres.\par

\vskip 3mm

 En répétant ce procédé $p$ fois, on obtient l'extension $\widetilde{QR{\Cal H}}^{1,p}_{p}\subset QA^{1,p}(x,v)$ qui est la limite inductive d'algèbres $\widetilde{QR{\Cal H}}^{1,p}_{p,m_1,\ldots,m_p}(x,v)$ dont les éléments ont une {\bf structure asymptotique élémentaire dans les} $p${\bf -blocs} $\chi${\bf -homogènes dans la variable} 

$$
X_{p,m_1,\ldots,m_p}=(x,(z^{n_1\ldots n_p}),y^{m_1},y^{m_1m_2},\ldots,y^{m_1\ldots m_p})
$$ 

\vskip 2mm
\noindent Le uplet $(m_1,\ldots,m_p)$ étant {\bf la complexité} de ces $p$-blocs. Les fonctions $z^{...}$ et $y^{...}$ satisfont des équations différentielles $(3_p)$ du même type que $(3_1)$, de valeurs propres $r_{j,0}$, $r_{j,1}$, ... et $r_{j,p}\equiv 1$. Ainsi construites, les fonctions élémentaires $z^{n_1\ldots n_i}\in QA^{1,p-i}(x,v_1,\ldots,v_{p-i}))$ sont indépendantes des coordonnées $v_{p-i+1},\ldots,v_p$. Un résultat qui généralise le théorème principal II1 est le suivant\par

\vskip 3mm

\proclaim{Lemme IIIB2} Pour tout $i=1,\ldots,p$ et pour tout choix d'entiers $m_1,\ldots,m_i$; l'algè\-bre $\widetilde{QR{\Cal H}}^{1,p}_{i,m_1,\ldots,m_i}(x,v)$ est $\chi$-finie et satisfait à la double inclusion. De plus, elle admet un morphisme série formelle $f\mapsto \widehat{f}^i$ qui est injectif, et on a les inclusions

$$
 \Phi_a^*(QR{\Cal H}^{1,q+(0,p)}_{|U_a})\hookrightarrow\widetilde{QR{\Cal H}}^{1,p}_{1,m_1}\hookrightarrow\widetilde{QR{\Cal H}}^{1,p}_{2,m_1,m_2}\hookrightarrow\ldots\hookrightarrow\widetilde{QR{\Cal H}}^{1,p}_{p,m_1,\ldots,m_p}
$$
\endproclaim

\vskip 3mm

\noindent{\bf Preuve.} D'après la preuve du théorème principal II1, la $\chi$-finitude, la double inclusion et l'existence d'un morphisme série formelle injectif sont conséquences du théorème de division VB1, de la quasi-analycité: $\widetilde{QR{\Cal H}}^{1,p}_{i,...}\subset QA^{1,p}(x,v)$, et de l'étude de l'action de $\chi$ sur les $i$-blocs $\chi$-homogènes.\par

\vskip 3mm 

 Pour cela, on se place dans la situation générale: on suppose que les valeurs propres $r_{j,0}$ sont indépendantes et on généralise les valeurs propres $r_{j,1}$ en $(r_{j,1,n})_n$, $r_{j,2}$ en $(r_{j,2,n_1,n_2})_{n_1,n_2}$, ...etc. Ainsi le nombre des valeurs propres indépendantes coincide avec le nombre des fonctions élémentaires $z^{...}$ et $y^{...}$ qui satisfont aux équations différentielles $(3_i)$ et à la condition initiale $z^{...}_{|x=1}=y^{...}_{|x=1}\equiv 0$. On généralise aussi les germes $\tau_{j,i}(v)$ en des variables indépendantes $\tau=(\tau_{j,i})$. Notons $\mu_{...}=r_{...}-1$, $\mu^i=(\mu_{...})$ et $\alpha^i=(\mu^i,\tau,v)$. Soit $H{\Cal H}_{i,n}$ le $\RR\{\alpha^i\}$-module des $i$-blocs $\chi$-homogènes de degré $n$ (et de complexité $(m_1,\ldots,m_i)$). Il est stable par $\chi$ d'après la linéarité du système $(3_i)$. Il s'agit donc de montrer que tout élément de $H{\Cal H}_{i,n}$ satisfait à la double inclusion, et que sa multiplicité algébrique est $n$. Mettons un ordre sur les monômes de $H{\Cal H}_{i,n}$ compatible avec la triangularité du système $(3_i)$. Pour tout multi-indice $m'$ de longueur $n$, posons $e_{m'}=\sum m'_{...}r_{...}$. Alors, d'après le système $(3_i)$, le $\RR\{\alpha^i\}$-module $H{\Cal H}_{i,n}$ est inclus dans le noyau de l'opérateur

$$
 E_n=\prod_{|m'|=n} (\chi-e_{m'}Id)
$$

 Pour pouvoir appliquer les méthodes de la section II à $H{\Cal H}_{i,n}$, il suffit de montrer que, génériquement en $\mu^i$, ce module coincide avec le noyau de l'opérateur $E_n$. Or, génériquement en $\mu^i$, la famille des monômes $x^{e_{m'}}$ forme une base de ce noyau. Or, par une récurrence sur les fonctions $z^{...}$ et $y^{...}$, et par la résolution triangulaire du système $(3_i)$, on montre que ces fonctions sont des combinaisons linéaires des fonctions $x^{r_{...}}$, par un système triangulaire inversible: en effet, si $\chi f=r_{\ell+1}f+g$ avec $f_{|x=1}\equiv 0$, et si $g=\sum_{j=0}^\ell a_j(\tau,r_1,\ldots,r_\ell)x^{r_j}$ (convention $r_0=1$), où les $a_j$ sont des polynômes en $\tau$ dont les coefficients sont des fonctions rationnelles non identiquement nulles, alors

$$
 f=\sum_{j=0}^\ell \frac{a_j}{r_j-r_{\ell+1}}x^{r_j}-(\sum_{j=0}^\ell\frac{a_j}{r_j-r_{\ell+1}})x^{r_{\ell+1}}
$$

\noindent Ainsi, $f=\sum_{j=0}^{\ell+1}b_j(\tau,r_1,\ldots,r_{\ell+1})x^{r_j}$ et les fonctions $b_j$ sont des polynômes en $\tau$, dont les coefficients sont rationnelles et non identiquement nulles. Par conséquent, les monômes $X^{m'}_{i,m_1,...,m_i}$ forment une base du noyau de l'opérateur $E_n$ pour des valeur génériques de $\mu^i$.\par

\vskip 3mm

 La double inclusion est donc réalisée sur le $\RR\{\alpha^i\}$-module $H{\Cal H}_{i,n}$, et en prenant la restriction au graphe $v\mapsto\alpha^i(v)=(\mu^i(v),\tau(v),v)$, elle est réalisée sur le $\RR\{v\}$-module des $i$-blocs $\chi$-homogènes de degré $n$. Dés lors, on applique exactement la démarche de la section II à l'algèbre $\widetilde{QR{\Cal H}}^{1,p}_{i,m_1,\ldots,m_i}(x,v)$. Si on note $c$ l'immersion: $(x,v)\mapsto (X_{i,m_1,\ldots,m_i}(x,v),v)$, cette algèbre est $\chi$-équivalente à la sous-algèbre 

$$
\widetilde{QR{\Cal H}}^{1,p}_{i,m_1,\ldots,m_i,(cvg)}=c^*(\RR\{X_{i,m_1,\ldots,m_i},v\}
$$

\noindent à laquelle s'applique le lemme d'extension IB2. De plus, on obtient (comme dans la section II) l'existence et l'injectivité du morphisme série formelle $f\mapsto \widehat{f}^i$.\par

\vskip 3mm

 En substituant les formules $(2_i)$ (linéaires dans les fonctions élémentaires), dans un $(i-1)$-bloc, on obtient un $i$-bloc de même degré (en convenant qu'un 0-bloc est un bloc $\chi$-homogène dans les variables $(x,z)$) .  Par l'existence et l'injectivité des morphismes séries formelles $f\mapsto\widehat{f}^{i-1}$ et $f\mapsto\widehat{f}^i$, cette application se prolonge injectivement de $\widetilde{QR{\Cal H}}^{1,p}_{i-1,m_1,...,m_{i-1}}$ vers $\widetilde{QR{\Cal H}}^{1,p}_{i,m_1,...,m_i}$. Ceci donne les inclusions du lemme. Ces extensions sont donc des extensions étoilées, les morphismes associés étant simplement l'identité.\qed\par

\vskip 3mm

\noindent{\bf (b) Division de} $g_a$ {\bf dans l'idéal} $(\varphi)$ {\bf dans l'extension} $\widetilde{QR{\Cal H}}^{1,p}_p$.\par

\vskip 3mm

\noindent Choisissons $m_1\geq \max\{n_{a,j};\ j=1,\ldots,p\}$. Soit $G_1$ l'image de $g_a$ dans l'extension $\widetilde{QR{\Cal H}}^{1,p}_{1,m_1}$. Il est de même degré $n$ que $g_a$, et son idéal $\chi$-transverse est inclus dans $(\varphi)$. D'après les égalités (1), les germes $\tau_{j,1}(v)$ appartiennent à l'idéal $(v_1\times\cdots\times v_p)=rad((\varphi))$. Donc d'après les formules $(2_1)$ et le choix de $m_1$, on a la division

$$
 G_1=H_1+\varphi H_2
$$

\noindent les germes $H_1$ et $H_2$ sont des 1-blocs de degré $n$, et les monômes de $H_1$ sont indépendants des fonctions élémentaires $y^{m_1}$. D'après les égalités (1) et les équations $(3_1)$, les fonctions élémentaires $z^.$ sont indépendantes de la coordonnée $v_p$. Comme l'idéal $\chi$-transverse de $H_1$ est inclus dans $(\varphi)$, un développement taylorien de ses coefficients donne la division

$$
 H_1=v_p^{n_{a,p}}F_{1}
$$

\noindent le germe $F_{1}$ étant un 1-bloc de degré $n$, dont l'idéal $\chi$-transverse est inclus dans $(\varphi_1)=(v_1^{n_{a,1}}\times\cdots\times v_{p-1}^{n_{a,p-1}})$. Ainsi, par une récurrence sur $i$, et en choisissant $m_2=m_1$, ..., $m_p=m_1$, on obtient que l'image $G_p$ de $g_a$ dans l'algèbre $\widetilde{QR{\Cal H}}^{1,p}_{p,m_1,...,m_p}$ s'écrit

$$
 G_p=\varphi H_{p}
$$

\noindent le germe $H_{p}$ étant un $p$-bloc de degré $n$.\par

\vskip 3mm

\noindent{\bf (c) Division de} $f_a$ {\bf dans l'idéal} $(\varphi)$ {\bf dans l'extension} $\widetilde{QR{\Cal H}}^{1,p}_p$.\par

\vskip 3mm

 Soient $(g_{n,a}(f_a)$ les 0-blocs de la série de $f_a$. D'après le théorème de division VB1, il existe un unique $Q\in QA^{1,p}(x,v)$ tel que $f=\varphi Q$. D'autre part, pour tout $n\in\NN$, il existe un $p$-bloc $G_{p,n}$ tel que $g_{n,a}(f_a)=\varphi G_{p,n}$. Soit $N\in\NN$, en utilisant le lemme de division II1, appliqué au germe $f_a-\sum_{n\leq N} g_{n,a}$, on obtient 
 
$$
Q-\sum_{n\leq N} G_{p,n}\in (x^N)SB^{1,p}_0
$$

\noindent Et ceci prouve que $Q\in\widetilde{QR{\Cal H}}^{1,p}_{p,m_1,\ldots,m_1}$.\par

\vskip 3mm

 Il est clair que l'idéal $\chi$-transverse de $Q$ n'est pas propre. Soit $ma_a$ sa multilicité algébrique (lemme IIIB2). On a $I_{\chi,Q}\supset (x^{ma_a+\varepsilon})$ pour tout $\varepsilon>0$. Donc, par la définition de la multiplicité algébrique, on a $ma_a\leq ma_\chi(f)$ pour tout $a$.\par

\vskip 3mm

\noindent{\bf Remarque IIIB1.} Notons $\varphi_i(v)=\prod_{j=1}^i v_j^{n_{a,j}}$. Si $h\in \Phi_a^*(QR{\Cal H}^{1,q+(0,p)}_{|U_a})$ et si on note $h_i$ son image dans l'extension $\widetilde{QR{\Cal H}}^{1,p}_{i,m_1,\ldots,m_1}$, l'étude précédente  montre que la série formelle de $h_i$ est la somme d'une série formelle $\widehat{h}_{i,1}$ indépendante des fonctions élémentaires $y^{...}$, et d'une série formelle $\widehat{h}_{i,2}$ qui se divise par $\varphi_i$, dans l'anneau des séries formelles associé à $\widetilde{QR{\Cal H}}^{1,p}_{i,m_1,\ldots,m_1}$.\par

\vskip 3mm

\noindent{\bf Remarque IIIB2.} Une conséquence facile de cette étude est que tout élément d'une algèbre $\widetilde{QR{\Cal H}}_{...}$ se divise dans son idéal $\chi$-transverse dans une extension $\widetilde{\widetilde{QR{\Cal H}}}_{...}$ obtenue en prenant aussi les développements $(2_{...})$ pour les fonctions élémentaires $y^{...}$. Ces nouvelles algèbres satisfont aussi au lemme IIIB2,...etc.\par
 
\vskip 3mm

\noindent{\bf §2. Cas } $\ell>0$.\par

\vskip 3mm

 Reprenons la dérivation $\chi=x\partial/\partial x-\sum_{j=1}^\ell s_j(\mu) u_j\partial/\partial u_j$. Soit 

$$
f\in QR{\Cal H}^{1,q}(x,\alpha,u)
$$

\noindent et soit $J_f\subset\RR\{\alpha,\lambda\}$ son idéal $\chi$-transverse le long de $\gamma=\{(\alpha,u)0\}$. Les coordonnées $\lambda=(\lambda_1,\ldots,\lambda_\ell)$ sont les intégrales premières non triviales de $\chi$\ : $\lambda_j=x^{s_j} u_j$. Il s'agit d'abord de diviser localement $f$ dans l'idéal $\pi_\chi^*(J_f)$, dans un anneau qui possède une structure asymptotique élémentaire. L'idée générale est la suivante: en général, les blocs formels $\chi$-homogènes de la série de $f$ sont divergents, et leurs idéaux $\chi$-transverses sont dans un anneau formel. Mais, on a vu dans la partie IIIA, que si le germe $f$ est presque quasi-convergent (bien approché par les germes quasi-convergents dans la ${\Cal M}_{(\alpha,u)}$-topologie), alors il possède une multiplicité algébrique et un idéal limite transverse. On applique donc à $J_f$ une désingularisation dans laquelle (ie. localement) le germe $f$ est {\bf presque quasi-convergent}, et mieux encore, le quotient de cette division locale par $\pi_chi^*(J_f)$, satisfait à lhypothèse $(H\lambda)$ (ou mieux encore, son idéal $\chi$-transverse n'est pas propre!). Ceci nécéssite une préparation des intégrales premières non triviales $\lambda_j$, simlaire à celle des intégrales premières $\mu_j$.\par

\vskip 3mm

  Soit donc $(\psi,{\Cal N})$ une désingularisation d'Hironaka dans laquelle l'idéal $J_f$ est principal et monomial (lemme IIIB1). Soit $(a,V_{a})$ une carte de cette désingularisation de coordonnée $v$. Soit $(\varphi)=\psi_a^*(J_f)$ avec $\varphi=\prod_{j=1}^p v_j^{n_j}$ et $\mu_{j,a}$ et soient $\lambda_{j,a}$ les relevés des coordonnées $\mu_j$ et $\lambda_j$. Notons $s_{j,0}=1+\mu_{j,a}$ ($=r_{j,0}$, voir §1). On suppose (pour simplifier les notations!) que $p=|q|$ et que les germes $\mu_{j,a}$ et $\lambda_{j,a}$ sont  préparés sphériquement comme dans les formules (1)

$$
 \mu_{j,a}(v)=\mu_{j,p+1-i}+\mu'_{j,p+1-i}
\tag $4$
$$

$$
\lambda_{j,a}=\lambda_{j,p+1-i}+\lambda'_{j,p+1-i}
\tag $4'$
$$

\noindent les fonctions $\mu_{j,p+1-i}$ et $\lambda_{j,p+1-i}$ étant indépendantes des coordonnées  $(v_i,\ldots,v_p)$ et les fonctions $\mu'_{j,p+1-i}$ et $\lambda'_{j,p+1-i}$ appartiennent à l'idéal $(v_1\times\cdots\times v_i)$ dans l'anneau $\RR\{v\}$. Notons toujours  $(\psi,{\Cal N})$ la composée de ces deux désingularisations.\par

\vskip 3mm

  Comme dans le paragraphe 1, on commence par relever les germes de $f$ et $\chi$ (qu'on note de la même façon), dans l'extension naturelle dans les coordonnées $(x,\alpha,v,u)$: $f\in QR{\Cal H}^{1,q+(0,p)}(x,\alpha,v,u)$ et $v$ sont des intégrales premières de $\chi$. On note aussi de la même façon l'idéal $\chi$-transverse $J_f$, dans ce relevé. Soient les morphismes analytiques
  
$$
\psi_a:v\in V_a\mapsto (\alpha_a(v),\lambda_a(v))\  \text{et}\  \Psi_a:\ v\in V_a\mapsto (\psi_a(v),v)\in W_a=\Psi_a(V_a)
$$

\noindent Le morphisme $\Psi_a$ étant un difféomorphisme sur son image, les anneaux $\RR\{v\}$ et $\RR\{\alpha,\lambda,v\}_{|W_a}$ sont isomorphes; donc $\Psi_a^*(J_{f|W_a})=(\varphi)$. Soit $U_a=\pi_\chi^{-1}(W_a)$, il contient l'orbite principale $\gamma_a=\{(\alpha,v,u)=0\}$. Par le difféomorphisme $\Psi_a$, identifions les variétés analytiques $V_a$ et $W_a$, et les idéaux $J_{f|W_a}$ et $(\varphi)$. Notons $\pi_{\chi|U_a}:(U_a,0)\to(V_a,0)$ le germe en 0 de la restriction à la variété analytique $U_a$, de la projection intégrale $\pi_\chi$. Il s'agit de diviser $f_a=f_{|U_a}$ dans l'idéal $\pi_{\chi|U_a}^*(\varphi)$. La démarche suit et généralise celle du cas $\ell=0$: on construit $p$ extensions $\widetilde{QR{\Cal H}}_i\subset QA$, de l'algèbre $QR{\Cal H}^{1,q+(0,p)}_{|U_a}$, telles que l'image de $f_{a}$ dans la $p$-ème extension, se divise par $\varphi$, dans cette $p$-ème extension. Ces extensions possèdent bien sûr une structure asymptotique élémentaire.\par

\vskip 3mm

\noindent{\bf (a) Division par l'idéal} $(v_p^{n_p})$ {\bf dans une extension} $\widetilde{QR{\Cal H}}_1$.\par

\vskip 3mm

  Reprenons les notations du §1 relatives aux formules (4) (qui remplacent les formules (1)). Soient $u^{(1)}=(u_{1,1},\ldots,u_{\ell,1})$ des coordonnées analytiques locales sur $(\RR^\ell,0)$. La dérivation

$$
 {\Cal X}_1=\chi -\sum_{j=1}^\ell r_{j,1}(v)u_{j,1}\frac{\partial}{\partial u_{j,1}}
$$

\noindent agit sur l'algèbre $QR{\Cal H}^{1,q+(0,p+\ell)}(x,\alpha,v,u,u^{(1)})$. Il admet $\gamma_1=\{(\alpha,v,u,u^{(1)})=0\}$ comme orbite principale. Notons $\lambda^{(1)}=(\lambda_{1,1},\ldots,\lambda_{\ell,1})$ des coordonnées sur $(\RR^\ell,0)$. Des coordonnées analytiques transverses à $\gamma_1$ sont $(\alpha,\lambda,v,\lambda^{(1)})$. Soit ${\Cal W}_1=\{  (\Psi_a(v),\lambda^{(1)}(v));\ v\in V_a\}$ et ${\Cal U}_1=\pi_{{\Cal X}_1}^{-1}({\Cal W}_1)$. Les germes $\lambda_{j,1}(v)$ sont donnés par les formules (4'). Sur la variété analytique ${\Cal U}_1$ (de dimension $p+1$), on a donc les relations supplémentaires $x^{r_{j,1}(v)}u_{j,1}=\lambda_{j,1}(v)$. Comme précédement, on identifie les variétés analytiques ${\Cal W}_1$ et $V_a$ (de dimension $p$), et on note $\pi_{{\Cal X}_1|{\Cal U}_1}:{\Cal U}_1\to V_a$ la restriction associée. On a l'injection canonique 

$$
QR{\Cal H}(x,\alpha,v,u)\hookrightarrow 
QR{\Cal H}(x,\alpha,v,u,u^{(1)})
$$

\noindent par la projection canonique $\pi':(x,\alpha,v,u,u^{(1)})\mapsto (x,\alpha,v,u)$. Notons toujours $f$ son image dans cette extension. Par la définition de ${\Cal U}_1$, la restriction $\pi'_{|{\Cal U}_1}:{\Cal U}_1\to U_a$ est un difféomorphisme sur son image $U_a$. On a donc aussi l'injection canonique

$$
QR{\Cal H}(x,\alpha,v,u)_{|U_a}\hookrightarrow 
QR{\Cal H}(x,\alpha,v,u,u^{(1)})_{|{\Cal U}_1}
$$

\noindent L'image de la dérivation restreinte ${\Cal X}_{1|{\Cal U}_1}$ est la dérivation $\chi$, et les idéaux transverses (le long de $\gamma_a$ et $\gamma_1$, dans la variété $V_a$), sont préservés (lemme de transfert IB6). Notons toujours $f_a$ la restriction de $f$ à ${\Cal U}_1$.\par

\vskip 3mm

 Soient $w^{(1)}=(w_{1,1},\ldots,w_{\ell,1})$ des coordonnées sur $(\RR^\ell,0)$. Soit le difféomorphisme sur $(\RR^{+*}\times\RR^{2(p+\ell)},0)$

$$
 \Phi_1(x,\alpha,v,u,u^{(1)})=(x,\alpha,v,u^{(1)},w^{(1)})=(x,\alpha,v,u^{(1)},u-u^{(1)})
$$

\noindent et soient ${\Cal X}'_1=(\Phi_1)_*{\Cal X}_1$, ${\Cal U}'_1=\Phi_1({\Cal U}_1)$ et $\gamma'_1=\Phi_1(\gamma_1)$. On a trivialement 

$$
\Phi_1^*(QR{\Cal H}(x,\alpha,v,u^{(1)},w^{(1)}))=QR{\Cal H}(x,\alpha,v,u,u^{(1)})
$$ 

\noindent et 

$$
\Phi_1*(QR{\Cal H}(x,\alpha,v,u^{(1)},w^{(1)})_{|{\Cal U}'_1})=QR{\Cal H}(x,\alpha,v,u,u^{(1)})_{|{\Cal U}_1}
$$

\noindent Soit $f'=(\Phi_1^{-1})^*(f)$ et $f'_a$ sa restriction à ${\Cal U}'_1$. Pour poursuivre, on a besoin de la définition suivante

\vskip 3mm

\proclaim{Définition IIIB1} Soit $t$ une coordonnée locale sur $(\RR^k,0)$ et soit $\chi_0=x\partial/\partial x$. L'extension 

$$
 \widetilde{QR{\Cal H}}^{1,(p,k)}_{1,m_1}(x,v,t)\subset QA^{1,p+k}(x,v,t)
$$

\noindent est l'algèbre des germes qui possèdent une série asymptotique formelle dans les 1-blocs $\chi_0$-homogènes de complexité $m_1$, construits (comme au §1) grâce aux formules $(4)$, et dont les coefficients appartiennent à l'anneau $\RR\{v,t\}$.\par
\endproclaim

\vskip 3mm

\noindent{\bf Remarque IIIB3.} Pour toute partition $t=(t',t")\in\RR^{k_1}\times\RR^{k_2}$

$$
 \widetilde{QR{\Cal H}}^{1,(p,k)}_{1,m_1}(x,v,t)
\subset \widetilde{QR{\Cal H}}^{1,(p,k_1)}_{1,m_1}(x,v,t')\{t"\}
$$

\noindent les séries étant convergentes sur un produit, au sens suivant 

$$
\widetilde{QR{\Cal H}}^{1,(p,k_1)}_{1,m_1}(x,v,t')\{t"\}\subset SB^{1,p+k}(x,v,t)
$$

\noindent et pour tout $h\in\widetilde{QR{\Cal H}}^{1,(p,k_1)}_{1,m_1}(x,v,t')\{t"\}$, de série 

$$
h=\sum_{|m|\geq 0} h_m (t")^m
$$

\noindent il existe un domaine standard $\Omega$ tel que $h_m\in QA^{1,p+k_1}[\Omega]$ pour tout $m$.\par  

\vskip 3mm

  Soient $U_{1,a}$ et $U_{2,a}$ les variétés analytiques, images de ${\Cal U}'_1$ par les projections canoniques 

$$
 (x,\alpha,v,u^{(1)},w^{(1)})\mapsto (x,v,u^{(1)})\quad\text{et}\quad (x,\alpha,v,u^{(1)},w^{(1)})\mapsto (x,v,u^{(1)},w^{(1)})
$$

\noindent et soient $\gamma_{1,a}$, $\gamma_{2,a}$ les images de $\gamma'_1$. Les restrictions de ces projections à ces variétés sont des difféomorphismes. Par une généralisation facile des résultats du §1, on a l'injection suivante pour tout $m_1$

$$
 QR{\Cal H}^{1,q+(0,p+\ell)}(x,\alpha,v,u^{(1)},w^{(1)})_{|{\Cal U}'_1}\hookrightarrow 
 \widetilde{QR{\Cal H}}^{1,(p,2\ell)}_{1,m_1}(x,v,u^{(1)},w^{(1)})_{|U_{2,a}}
 \tag 5
$$

\noindent La dérivation ${\Cal X}'_1$ (restreinte à ${\Cal U}'_1$) est préservée (par restriction à $U_{2,a}$). Choisissons $m_1>\max\{n_j;\ j=1,\ldots,p\}$. Soit $h\in \widetilde{QR{\Cal H}}^{1,(p,2\ell)}_{1,m_1}$ tel que l'image de $f'_a$ dans cette extension soit égale à $h_a=h_{|U_{2,a}}$. La motivation de ces extensions est la suivante: par la définition de ${\Cal U}_1$, on a

$$
x^{s_{j,0}}u_{j|{\Cal U}_1}=\lambda_{j,a}(v)\quad\text{et}\quad x^{r_{j,1}}u_{j,1|{\Cal U}_1}=\lambda_{j,1}(v)
$$

\noindent pour tout $j=1,\ldots,\ell$. Par conséquent, par le difféomorphisme $\Phi_1$, un calcul direct donne
  
$$
 x^{s_{j,0}+r_{j,1}}w_{j,1|{\Cal U}'_1}=\lambda'_{j,1}x^{r_{j,1}}+\lambda_{j,1}(x^{r_{j,1}}-x^{s_{j,0}})
\tag 6
$$

\noindent (voir formules (4')). Notons $\delta_j(x,v)=x^{s_{j,0}}-x^{r_{j,1}}$. En utilisant les développements $(2_1)$, on voit que l'image de $\delta_j$ dans l'extension (5), est un élément de l'idéal $(\tau_{j,1}(v))\subset rad((\varphi))$ (plus précisément, dans l'anneau $\widetilde{QR{\Cal H}}^{1,p}_{1,m_1}(x,v)$). 
D'après les formules (4'), on a $\lambda'_{j,1}\in rad((\varphi))$, donc l'image du germe (6) dans l'extension (5), est un élément de l'idéal $rad((\varphi))$ dans l'anneau $\widetilde{QR{\Cal H}}^{1,p}_{1,m_1}(x,v)$).\par

\vskip 3mm

 Posons donc $S_1=m_1\sum_{j=1}^\ell (s_{j,0}+r_{J,1})$ et notons

$$
 h_{1}={\bold j}_{w^{(1)}}^{m_1}(x^{S_1} h)\quad\text{et}\quad  h_2=x^{S_1}h-h_1
\tag 7
$$

\noindent (en prenant bien sûr l'image du monôme $x^{S_1}$ dans l'extension (5)). Notons aussi $h_{1,a}$ et $h_{2,a}$ leurs restrictions à $U_{2,a}$, par la remarque IIIB3, ce sont des éléments de l'extension (5). Et, par la remarque faite sur le germe (6), et par le choix de $m_1$, le germe $h_{2,a}$ se divise dans l'idéal $(\varphi){\Cal M}'$ dans l'extension (5) (${\Cal M}'$ étant l'idéal maximal de $\RR\{v\}$). Comme l'idéal ${\Cal X}'_1$-transverse de $h_a$ est $(\varphi)$, il en est de même de celui de $h_{1,a}$.\par

\vskip 3mm

  Le germe $h_1$ est polynomial dans la coordonnée $w^{(1)}$. L'égalité (6) montre que sa restriction $h_{1,a}$ s'identifie à la restriction à $U_{1,a}$ d'un élément $F_1$ de l'algèbre $\widetilde{QR{\Cal H}}^{1,(p,\ell)}_{1,m_1}(x,v,u^{(1)})$. L'image dans cette extension, de la dérivation ${\Cal X}'_{1}$ (restreinte à $U_{2,a}$), est la dérivation 

$$
 \chi_1=\chi_0 -\sum_{j=1}^\ell r_{j,1}u_{j,1}\frac{\partial}{\partial u_{j,1}}
$$

\noindent (restreinte à $U_{1,a}$). Les idéaux transverses sont préservés dans cette extension. Notons $F_{1,a}$ la restriction de $F_1$ à $U_{1,a}$. L'idéal $\chi_1$-transverse de $F_{1,a}$, le long de $\gamma_{1,a}$, est $(\varphi)$.\par

\vskip 3mm

 Il s'agit maintenant de montrer que $F_{1,a}$ appartient à l'idéal $(v_p^{n_p})$ dans l'anneau  $\widetilde{QR{\Cal H}}^{1,(p,\ell)}_{1,m_1|U_{1,a}}(x,v,u^{(1)})$. Pour cela, on commence par effectuer une division de $F_1$ par $v_p^{n_p}$ dans l'anneau  $\widetilde{QR{\Cal H}}^{1,(p,\ell)}_{1,m_1}(x,v,u^{(1)})$. En effet, d'après le théorème de division VB1 dans l'anneau $QA^{1,p+\ell}(x,v,u^{(1)})$, il existe $f_1\in QA^{1,p+\ell}(x,v,u^{(1)})$ et $R_1,\ldots,\ R_{n_p}\in QA^{1,p+\ell-1}(x,v_1,\ldots,v_{p-1},u^{(1)})$ tels que
 
$$
F_1=v_p^{n_p} f_1+\sum_{j=0}^{n_p-1} v_p^{j} R_j=v_p^{n_p} f_1+R
$$

\noindent Par la remarque IIIB3, les opérations de prise de jet fini en $w^{(1)}$, et d'extension (5) commutent. Donc, d'après la remarque IIIB1, la série formelle $\widehat{F}_1$ relativement à la dérivation $\chi_0$, est la somme d'une série formelle $\widehat{F}_{1,1}$ indépendante des fonctions élémentaires $y^{m_1}$, et d'une série formelle $\widehat{F}_{1,2}$ qui se divise dans l'idéal $(\varphi)$ dans l'anneau formel associé à $\widetilde{QR{\Cal H}}^{1,(p,\ell)}_{1,m_1}(x,v,u^{(1)})$. En divisant dans l'idéal $(v_p^{n_p})$, les 1-blocs $\chi_0$-homogènes et les restes de la série $\widehat{F}_1$, et en utilisant l'unicité de la division de $F_1$, on voit que $f_1\in  \widetilde{QR{\Cal H}}^{1,(p,\ell)}_{1,m_1}(x,v,u^{(1)})$ et que pour tout $j=0,\ldots,n_p-1$, $R_j\in \widetilde{QR{\Cal H}}^{1,(p-1,\ell)}_{1,m_1}(x,v_1,\ldots,v_{p-1},u^{(1)})$.\par

\vskip 3mm 
 
 Maintenant, pour montrer que $R_{|U_{1,a}}\equiv 0$, on utilise le lemme de saturation IB5. Soit $x_0>0$ suffisament petit, et soit $m_0=(x_0,0)\in\gamma_{1,a}$. Redressons le champ $\chi_{1|U_{1,a}}$ dans un voisinage de $m_0$ inclus dans $U_{1,a}$
 
$$
x=x_0\exp(t),\ u_{j,1}=\frac{\lambda_{j,1}(v)}{x_0}\exp(-r_{j,1}(v)t)
$$

\noindent D'après le lemme de saturation IB5, le germe de $F_{1,a}$ en $m_0$, appartient au saturé de l'idéal $(\varphi)$. Donc, dans les coordonnées locales $(t,v)$ sur $(\RR^{p+1},0)$, ce germe est divisible par $v_p^{n_p}$. Par conséquent, le germe de $R_{|U_{1,a}}$ en $m_0$, est divisible par $v_p^{n_p}$. Or, les germes $r_{j,1}(v)$ et $\lambda_{j,1}(v)$ sont indépendants de la coordonnée $v_p$ (formules (4) et (4')). Comme $R$ est polynomial dans la coordonnée $v_p$, de degré$\leq n_p-1$, il s'ensuit que $R_{|U_{1,a}}\equiv 0$.\par

\vskip 3mm

 On a donc construit $f_1\in \widetilde{QR{\Cal H}}^{1,(p,\ell)}_{1,m_1}(x,v,u^{(1)})$ tel que 

$$
 F_{1,a}=v_p^{n_p} f_{1,a}
$$

\noindent (où $f_{1,a}$ est la restriction de $f_1$ à $U_{1,a}$). L'idéal $\chi_1$-transverse de $f_{1,a}$ est donc égal à l'idéal $(\varphi_1)=(v_1^{n_1}\times\cdots\times v_{p-1}^{n_{p-1}})$ et, toujours par la remarque IIIB1, la série formelle de $f_1$ relativement à $\chi_0$, est la somme d'une série formelle $\widehat{f}_{1,1}$ indépendante des fonctions élémentaires $y^{m_1}$, et d'une série formelle $\widehat{f}_{1,2}$ divisible par l'idéal $(\varphi_1)$, dans l'anneau formel associé à $\widetilde{QR{\Cal H}}^{1,(p,\ell)}_{1,m_1}(x,v,u^{(1)})$.\par

\vskip 3mm

 {\bf En résumé}, notons ${\Cal U}_{1,a}$ l'image de ${\Cal U}_1$ par la projection canonique 

$$
 (x,\alpha,v,u,u^{(1)})\mapsto (x,v,u,u^{(1)})
$$

\noindent et $\Phi_{1,0}$ la restriction de $\Phi_1$ à $\{\alpha=0\}$. Soit $H\in\widetilde{QR{\Cal H}}^{1,(p,2\ell)}_{1,m_1}(x,v,u^{(1)},w^{(1)})$ tel que $h_{2,a}=\varphi H_a$ ($H_a$ étant sa restriction à $U_{2,a}$). Notons $H_1=\Phi_{1,0}^*(H)$ et $H_{1,a}$ sa restriction à ${\Cal U}_{1,a}$. L'image de la variété ${\Cal U}_{1,a}$ par la projection canonique $(x,v,u,u^{(1)})\mapsto (x,v,u^{(1)})$ est la variété $U_{1,a}$, et par une identification triviale, $\Phi_{1,0}^*(f_1)=f_1$. Donc, en composant toutes ces extensions, on obtient une extension

$$
(\widetilde{QR{\Cal H}}^{1,(p,2\ell)}_{1,m_1}(x,v,u,u^{(1)})_{|{\Cal U}_{1,a}},\pi_1)\hookleftarrow QR{\Cal H}^{1,q+(0,p)}(x,\alpha,v,u)_{|U_a}
$$

\noindent telle que

$$
 \pi_1^*(x^{S_1}f_a)=v_p^{n_p} f_{1,a}+\varphi H_{1,a}
$$

\noindent avec $f_1\in \widetilde{QR{\Cal H}}^{1,(p,\ell)}_{1,m_1}(x,v,u^{(1)})$, dont la série formelle relativement à $\chi_0$ est la somme des deux séries formelles $\widehat{f}_{1,1}$ et $\widehat{f}_{1,2}$, et telle que l'idéal $\chi_1$-transverse de $f_{1,a}$ est $(\varphi_1)$. L'image dans cette extension, de la dérivation $\chi$ (restreinte à $U_a$), est la dérivation ${\Cal X}_1$ (restreinte à ${\Cal U}_{1,a}$).\par

\vskip 3mm 

\noindent{\bf (b) Division dans l'idéal} $(\varphi)$ {\bf dans une extension} $\widetilde{QR{\Cal H}}_{p}$.\par

\vskip 3mm

\proclaim{Définition IIIB2} Soient $m_1,\ldots,m_i\in\NN$ ($i\leq p$), et soit $t$ une coordonnée locale sur $(\RR^k,0)$. L'extension 

$$
 \widetilde{QR{\Cal H}}^{1,(p,k)}_{i,m_1,\ldots,m_i}(x,v,t)\subset QA^{1,p+k}(x,v,t)
$$

\noindent est l'algèbre des germes qui possèdent une série asymptotique formelle dans les $i$-blocs $\chi_0$-homogènes de complexité $(m_1,\ldots,m_i)$, construits grâce aux formules $(4)$, et dont les coefficients appartiennent à l'anneau $\RR\{v,t\}$.\par
\endproclaim

\vskip 3mm

\noindent La remarque IIIB3 s'appliquent à ces algèbres.\par

\vskip 3mm

 Choisissons $m_2=\cdots=m_p=m_1$. On répéte le procédé du §2a, $p-1$ fois, appliqué au germe $f_1$. On construit une suite de germes $f_i$, $H_i$ (avec $f_0=f$ et $H_0=0$), et une suite d'extensions 
 
$$
(\widetilde{QR{\Cal H}}^{1,(p,(i+2)\ell)}_{i+1,m_1,...,m_{i+1}}(x,v,u,u^{(1)},\ldots,u^{(i+1)})_{|{\Cal U}_{i+1,a}},\pi_{i,i+1})
$$

$$
\hookleftarrow 
\widetilde{QR{\Cal H}}^{1,(p,(i+1)\ell)}_{i,m_1,\ldots,m_i}(x,v,u,u^{(1)},\ldots,u^{(i)})_{|{\Cal U}_{i,a}}
$$

\noindent telles que

$$
\pi_{i,i+1}^*(x^{S_{i+1}}f_{i,a})=\frac{\varphi}{\varphi_{i+1}} f_{i+1,a} + \varphi H_{i+1,a}
$$

\noindent avec $S_i=m_i\sum_{j=1}^\ell (r_{j,i-1}+r_{j,i})$ pour $i=1,\ldots,p$. Le germe 

$$
H_i\in\widetilde{QR{\Cal H}}^{1,(p,(i+1)\ell)}_{i,m_1,\ldots,m_i}(x,v,u,u^{(1)},\ldots,u^{(i)})
$$

\noindent et le germe $f_i\in\widetilde{QR{\Cal H}}^{1,(p,\ell)}_{i,m_1,\ldots,m_i}(x,v,u^{(i)})$ (en posant $u^{(0)}=u$). les germes $H_{i,a}$ et $f_{i,a}$ dénotent leurs restrictions à ${\Cal U}_{i,a}$. La série formelle de $f_i$ relativement à $\chi_0$, est la somme d'une série formelle $\widehat{f}_{i,1}$ indépendante des fonctions élémentaires $y^{m_1},\ y^{m_1m_2},\ldots,\ y^{m_1\cdots m_i}$; et d'une série formelle $\widehat{f}_{i,2}$ divisible par $\varphi_i$ dans l'anneau formel associé à $\widetilde{QR{\Cal H}}^{1,(p,\ell)}_{i,m_1,\ldots,m_i}(x,v,u^{(i)})$. La variété analytique ${\Cal U}_{i,a}$, de dimension $p+1$, est donnée par les conditions

$$
{\Cal U}_{i,a}=\{ (x,v,u,\ldots,u^{(i)})\in {\Cal U}_i\in (\RR^{+*}\times\RR^{p+(i+1)\ell},0);
$$

$$
 x^{r_{j,k}(v)}u_{j,k}=\lambda_{j,k}(v),\ j=1,\ldots,\ell;\ k=0,\ldots,i\}
$$

\noindent L'image de la dérivation

$$
 {\Cal X}_i=\chi -\sum_{j=1,\ldots,\ell,\ k=1,\ldots,i} r_{j,k}u_{j,k}\frac{\partial}{\partial u_{j,k}}
$$

\noindent (restreinte à ${\Cal U}_{i,a}$), est la dérivation ${\Cal X}_{i+1}$ (restreinte à ${\Cal U}_{i+1,a}$).

\vskip 3mm

  Posons $S_a=S_1+\cdots +S_p$. On a donc construit une extension

$$
(\widetilde{QR{\Cal H}}^{1,(p,(p+1)\ell)}_{p,m_1,...,m_p}(x,v,u,u^{(1)},\ldots,u^{(p)})_{|{\Cal U}_{p,a}},\pi_p)\hookleftarrow  QR{\Cal H}^{1,q+(0,p)}(x,\alpha,v,u)_{|U_a}
$$ 

\noindent telle que

$$
 \pi_p^*(x^{S_a} f_a)=\varphi Q_a
$$

\noindent avec $Q\in(\widetilde{QR{\Cal H}}^{1,(p,(p+1)\ell)}_{p,m_1,...,m_p}(x,v,u,u^{(1)},\ldots,u^{(p)})$ ($Q_a$ étant sa retsriction à ${\Cal U}_{p,a}$). L'idéal ${\Cal X}_p$-transverse de $Q_a$ le long de $\gamma_p=\{(v,u,u^{(p)})=0\}$, {\bf n'est pas propre}. Le germe $Q_a$ est donc {\bf presque quasi-convergent}, car il satisfait à l'hypothèse $(H\lambda)$ (cf. partie IIIA §2.2). Pour finir la preuve du théorème principal IIIB1, il suffit de montrer que le germe $Q_a$ est ${\Cal X}_p$-fini.\par

\vskip 3mm

  Plus généralement, indiquons briévement comment on adapte les principaux résultats de la partie A, à l'action de la dérivation ${\Cal X}_p$ (qu'on notera ${\Cal X}$ pour simplifier), sur l'algèbre étendue $\widetilde{QR{\Cal H}}^{1,(p,(p+1)\ell)}_{p,m_1,...,m_p}(x,v,u,u^{(1)},\ldots,u^{(p)})$ (qu'on notera simplement $\widetilde{QR{\Cal H}}^{1,(p,\ell')}_{p,m}(x,v,u')$). Soit
  
$$
X_{p,m}(x,v)=(x,z^{...}(x,v),y^{...}(x,v))
$$

\noindent les fonctions élémentaires de cette algèbre (cf. §1 pour les notations $z^{...}$, $y^{...}$). Les $p$-blocs formels ${\Cal X}$-homogènes de degré $k\in\ZZ$ (et de complexité $m$), sont les séries formelles de la forme

$$
G_k=\sum_{|n|-|n'|=k} a_{n,n'}(v) X_{p,m}^n (u')^{n'}
$$

\noindent avec $a_{n,n'}\in\RR\{v\}$. On définit alors, de la même façon que dans la partie A, les germes quasi-convergents, la multiplicité algébrique relativement à ${\Cal X}$, et les germes presque quasi-convergents. Soit $c_{p,m}$ l'immesion

$$
c_{p,m}(x,v)=(X_{p,m}(x,v),v,u')
$$

\noindent et soit la sous-algèbre {\bf convergente} $\widetilde{QR{\Cal H}}^{1,(p,\ell'}_{p,m,(cvg)}=c^*(\RR\{X_{p,m},v,u'\}$ (ici, $X_{p,m}$ désignent des variables). Si on généralise les valeurs propres $r_{j,i}$ (comme dans le lemme IIIB2), on voit que cette algèbre convergente est simplement la restriction d'une algèbre convergente $QR{\Cal H}^{1,q'}(x,\alpha',u')_{cvg}$ au graphe du germe analytique

$$
v\mapsto \mu'(v)=(r_{j,i}(v)-1)_{j=1,\ldots,q_1;\ i=0,\ldots,p}
$$

\noindent (Ceci n'est pas le cas pour l'algèbre $\widetilde{QR{\Cal H}}^{1,(p,\ell')}_{p,m}(x,v,u')$: un élément de cette algèbre n'est pas forcément une restriction d'un élément d'une algèbre $QR{\Cal H}$). Donc, par une vérification aisée des formules (2') (partie A), appliquées aux fonctions élémentaires $z^{...}$, $y^{...}$, et utilisant uniquement les équations différentielles $(3_p)$ (partie A), et par le lemme IIIB1, on obtient\par

\vskip 3mm

\proclaim{Lemme IIIB3} L'algèbre $\widetilde{QR{\Cal H}}^{1,(p,\ell')}_{p,m,(cvg)}$ est ${\Cal X}$-finie, et elle satisfait globalement à la double inclusion.\par
\endproclaim

\vskip 3mm

 On définit de la même façon que dans la partie A, la classe $\widetilde{{\Cal C}}^1_\lambda$ des germes qui satisfont à l'hypothèse $(H\lambda)$. Donc, en utilisant le lemme IIIB1, qui remplace le théorème principal II1, et le lemme IIIB3, qui remplace le théorème principal IIIA1, on obtient\par
 
\vskip 3mm

\proclaim{Lemme IIIB4} La classe $\widetilde{{\Cal C}}^1_\lambda$ est ${\Cal X}$-finie, et satisfait à la double inclusion.\par
\endproclaim
  
\vskip 3mm

  Ce lemme est encore vrai sur la restriction ${\Cal U}_{p,a}\subset {\Cal U}_p$, ce qui finit la preuve du théorème.\qed
  
\vskip 3mm

 Cette preuve suggère la définition suivante de la multiplicité algébrique relativement à $\chi$, pour tout $f\in QR{\Cal H}^{1,q}$: soit ${\Cal D}$ le diviseur exeptionnel de la désingularisation $(\psi,{\Cal N})$, alors on pose
 
$$
ma_{\chi}(f)=\inf_{(\psi,{\Cal N})}\sup_{a\in{\Cal D}}(ma_{{\Cal X}}(Q_a)-S_a(0))
$$
  
\noindent et ce nombre est fini, car ${\Cal D}$ est compact, et la multiplicité algébrique des germes d'idéal transverse non propre (et même des germes presque quasi-convergents), est semi-continue supérieurement, comme fonction des variables $v$ (ceci est étudié en détail dans l'article [M']).\par

\vskip 5mm

\centerline{{\bf IV. Démonstration du théorème 0.}}\par

\vskip 5mm

\noindent {\bf A. Désingularisation de la dérivation d'Hilbert.}\par

\vskip 3mm

 Soit $\Xi[QR{\Cal H}^{p,q}]$ le $QR{\Cal H}^{p,q}$-module de germes en $0$ de champs de vecteurs \`a composantes dans $QR{\Cal H}^{p,q}$ et qui laissent invariant l'alg\`ebre $QR{\Cal H}^{p,q}$. Ce module contient le sous-module engendr\' e par les d\' erivations \' el\' ementaires $x_j\partial/\partial x_j$ pour $j=1,\ldots,p$. On s'intéresse plus particuli\`erement \`a une sous-classe de $\Xi[QR{\Cal H}^{p,q}]$ qui apparait dans le probl\`eme d'Hilbert hyperbolique, et qui a la propri\' et\' e d'être stable dans une certaine d\' esingularisation. Cette désingularisation est inspirée de la géométrie du polycycle déployé. On note cette classe $\Xi{\Cal H}[QR{\Cal H}^{p,q}]$ et elle est d\' efinie de la façon suivante: soit $k\leq p$; posons $r_j=1+\mu_j$ pour $j=1,\ldots,k-1$, $x=(x_1,\cdots,x_k)$ et $x'=(x_{k+1},\cdots,x_p)$. Les \' el\' ements de $\Xi{\Cal H}_k[QR{\Cal H}^{p,q}]$ sont les germes en 0, de champs de vecteurs $\chi$ qui satisfont aux conditions suivantes\par

\vskip 3mm

\roster

\item"(a)" si $k=1$, $\Xi{\Cal H}_1[QR{\Cal H}^{p,q}]=\{ x_1\partial/\partial x_1,\cdots,
x_p\partial/\partial x_p \}$.\par

\vskip 2mm

\item"(b)" Si $k>1$, on a $\chi x_1=\prod_{j=1}^k x_j$ et $\chi$ admet comme int\' egrales premi\`eres les fonctions coordonn\' ees $\alpha'=(x',\alpha)$ et $(k-1)$ germes $g_j(x_j,x_{j+1},\alpha')=d_j(x_j,\alpha')-x_{j+1}$ avec 
$$
d_j=x_j^{r_j}(1+D_j) \quad \text{et}\quad D_j=O(x_j)\in QR{\Cal H}^{p-k+1,q}(x_j,x',\alpha)
$$

\endroster

\vskip 3mm

\noindent L'entier $k-1$ est appelé {\bf dimension de non trivialité} de $\chi$, et les germes $g_j$ sont dits ses{\bf intégrales premières non triviales}. On pose 

$$
 \Xi{\Cal H}[QR{\Cal H}^{p,q}]=\cup_{k=1}^p \Xi{\Cal H}_k[QR{\Cal H}^{p,q}]
$$

\vskip 3mm

\noindent {\bf 1§. Réduction de certains éléments de} $\Xi[QR{\Cal H}^{p,q}]$ {\bf à des éléments de}\par
\noindent $\Xi{\Cal H}[QR{\Cal H}^{p,q}]$.\par

\vskip 3mm

 Soit $\chi \in \Xi[QR{\Cal H}^{p,q}]$ tel que $\chi x_1=\prod_{j=1}^k x_j$ et qui admet comme int\' egrales premi\`eres la coordonn\' ee $\alpha'$ et $(k-1)$ germes $g_j(x_j,x_{j+1},\alpha')=a_j(\alpha')d_j(x_j,\alpha')-f_{j+1}(x_{j+1},\alpha')$ tels que $d_j$ a la même structure que dans le cas (b) ci-dessus et  $f_j=b_j(\alpha')x_j(1+O(x_j))\in QR{\Cal H}^{p-k,q+(0,1)}(\alpha',x_j)$. Les fonctions  $a_j$ et  $b_j$ appartiennent à l'algèbre $\in QR{\Cal H}^{p-k,q}$ avec $a_j(0)>0$ et  $b_j(0)>0 $. Par le théorème d'inversion VB4 (appendice VB), il existe un diff\' eomorphisme 

$$
H(x,\alpha')=(h_1(x_1,\alpha'),\ldots,h_k(x_k,\alpha'),\alpha')
$$

\noindent o\`u les germes $h_j$ ont la m\^ eme structure que les germes $f_j$ et tel que le champ $H_*\chi$ soit équivalent \`a un \' el\' ement de $\Xi{\Cal H}_k[QR{\Cal H}^{p,q}]$.\par

\vskip 3mm

\noindent {\bf §2. Désingularisation d'éléments de} $\Xi{\Cal H}$.\par

\vskip 3mm

 Soit $1<k\leq p$ et $\chi \in \Xi{\Cal H}_k[QR{\Cal H}^{p,q}]$ repr\' esent\' e sur $U\in ( (\RR^{+*})^p \times \RR^{|q|},0 )$ qu'on choisit de la forme $U=U_1\times U_2$ de coordonn\' ees $(x,\alpha')$.\par

\vskip 3mm

\noindent {\bf 2.1. Premier éclatement de} $\chi$.\par

\vskip 3mm

  Posons $r_k=1$ et pour $i\leq j$, $r_{i,j}=r_i\times \cdots \times r_j$. Soit ${\chi}_{pr}\in \Xi{\Cal H}_k$ le champ dont $(k-1)$ int\' egrales premi\`eres non triviales, sont les parties principales des $g_j$: $g_{j,pr}=x_j^{r_j}-x_{j+1}$. Ces int\' egrales premi\`eres  sont invariantes sous l'action du sous-groupe du groupe lin\' eaire $GL(\RR,k)$, constitu\' e des transformations de matrice $T_{k,(\rho,\alpha')}=Diag(\rho,\rho^{r_{1,1}},\ldots,\rho^{r_{1,k-1}})$ avec $\rho>0$. En utilisant l'action de $\chi_{pr}$ sur $x_1$, on obtient $(T_{k,(\rho,\alpha')})^{-1}_*\chi_{pr}=\rho^{s_k}{\Cal Y}_{pr}$, avec $s_k=\sum_{j=1}^{k-1} r_{1,j}$ et le champ ${\Cal Y}_{pr}$ a la m\^ eme expression que $\chi_{pr}$ dans la coordonn\' ee $y=T_{k,(\rho^{-1},\alpha')}(x)$. La d\' esingularisation adapt\' ee au probl\`eme est donc quasi-sph\' erique dans la coordonn\' ee $x$ et est fibr\' ee dans la coordonn\' ee $\alpha'$. Soit la fonction 
  
$$
Q(y,\alpha')=\sum_{j=1}^{k} y_j^{r_{j,k}}
$$

\noindent et pour $\varepsilon>0$, les quasi-sphères

$$
 S_k^{+*}(\varepsilon)=\{ (y,\alpha')\in ({\bold{\RR}}^{+*})^k\times U_2;\quad
                                          Q(y,\alpha')=\varepsilon \}
$$

\noindent on note $S_k^{+*}(\alpha',\varepsilon)$ leurs sections par les fibres $\alpha'=\text{constante}$. Le morphisme d'\' eclatement est

$$
 T_k:(\rho,y,\alpha')\in {\bold{\RR}}^{+*}\times S_k^{+*}(1)\mapsto (x,\alpha')=(T_{k,(\rho,\alpha')}(y),\alpha')
$$

\noindent Soit ${\Cal Y}$ le champ de vecteurs tel que $\rho^{s_k}{\Cal Y}=(T_k)^{-1}_*\chi$. Il admet comme int\' egrales premi\`eres, outre la coordonn\' ee $\alpha'$, les $(k-1)$ germes $G_j=g_j\circ T_k$ et on v\' erifie facilement que $G_j=\rho^{r_{1,j}}L_j$ avec $L_j(\rho,y,\alpha')=y_j^{r_j}(1+O(\rho))-y_{j+1}$, et elle est localement induite par un élément de l'algèbre $QR{\Cal H}^{p'(0),q'(0)}$ avec $p'(0)=p-k+1$ et $q'(0)=(kq_1+k-1,q_2+k-1)$. Soit l'application $L=(L_1,\ldots,L_{k-1})$ et $u$ une coordonn\' ee sur ${\RR}^{k-1}$. Pour simplifier l'expression du champ ${\Cal Y}$, il est naturel d'introduire le morphisme suivant

$$
 {\Cal L}_k(\rho,y,\alpha')=(\rho,L(\rho,y,\alpha'),\alpha')=(\rho,u,\alpha')                                                    $$

\noindent Soit ${\Cal L}_{k,(\rho,\alpha')}$ ses fibres par $(\rho,\alpha')=\text{constante}$
et ${\Cal D}_k(\alpha')={\Cal L}_{k,(0,\alpha')}( S_k^{+*}(\alpha',1))$. Notons simplement ${\Cal D}_k$ au lieu de ${\Cal D}_k(0)$.\par

\vskip 3mm

\proclaim{Proposition IVA1}

\vskip 2mm

\roster

\item"$(i)$" Le morphisme ${\Cal L}_{k,0}$ est un difféomorphisme de ${S_k^{+*}(0,1)}$ sur 
${{\Cal D}_k}$ qui se prolonge continument en une bijection de $\overline{S_k^{+*}(0,1)}$ sur $\overline{{\Cal D}_k}$.\par

\vskip 2mm

\item"$(ii)$" Soit $K$ un ouvert relativement compact de $S_k^{+*}(0,1)$. Alors, quitte \`a r\' eduire le voisinage $U_2$, le morphisme ${\Cal L}_k$ est un difféomorphisme de $V_k=(\RR^{+*},0)\times ( (K\times U_2)\cap S_k^{+*}(1) )$ sur son image $\widetilde{U}_k={\Cal L}_k(V_k)$. Les composantes de son inverse sont localement induites par des éléments de l'algèbre $QR{\Cal H}^{p'(0),q'(0)}$.\par

\vskip 2mm

\item"$(iii)$" Sur $\widetilde{U}_k$, on a 

$$
 ({\Cal L}_k)_*{\Cal Y}\sim \widetilde{\chi}=
                    \rho\frac{\partial}{\partial \rho}-
                    \sum_{j=1}^{k-1}r_{1,j}u_j\frac{\partial}{\partial u_j}
$$

\endroster

\endproclaim

\vskip 3mm

\noindent {\bf Preuve.} $(i)$ Les $(k-1)$ germes $L_j(0,.,.)$ sont des int\' egrales premi\`eres du champ ${\Cal Y}_{pr}$ qui est transverse aux sph\`eres $S_k^{+*}(\varepsilon)$ car ${\Cal Y}_{pr}Q>0$. Donc, ${\Cal L}_{k,0}$ est un difféomorphisme de $S_k^{+*}(0,1)$ sur ${\Cal D}_k$ qui se prolonge continument à $\overline{S_k^{+*}(0,1)}$. Comme $r_j(0)=1$ pour tout $j$, le système 

$$
 Q(y,0)=1+u_0,\quad y_j-y_{j+1}=u_j
$$

\noindent est linéaire inversible, et ceci prouve que ${\Cal L}_{k,0}$ est une bijection de $\overline{S_k^{+*}(0,1)}$ sur $\overline{{\Cal D}_k}$.\par

\comment
 Si certains $r_j(0)$ sont différents de 1, il est nécessaire d'utiliser la réduction de la proposition IVA2 ci-dessous qui décrit la topologie du champ ${\Cal Y}$ au voisinage du bord de la quasi-sphère $S_k^{+*}(0,1)$.\par
\endcomment

\vskip 3mm

\noindent $(ii)$ C'est une cons\' equence du $(i)$, par transversalit\' e et par la relative compacit\' e de $K$. La structure des composantes de l'inverse est conséquence du théorème des fonctions implicites dérivé du théorème de division VB3.\par

\vskip 3mm

\noindent $(iii)$ Etudions l'action de ${\Cal Y}$ sur la coordonnée $\rho$. On a $\rho^{r_{1,k}}=Q(x,\alpha')=Q\circ T_k(\rho,y,\alpha')$ et donc $\rho^{s_k}({\Cal Y}\rho^{r_{1,k}})=(\chi Q)\circ T_k$. Les sections de cette d\' erni\`ere fonction par $(\rho,\alpha')=$constante, sont $(\chi Q_{\alpha'})\circ T_{k,(\rho,\alpha')}$. Or pour $y\in ({\bold{\RR}}^{+*})^k$ quelconque, on a $Q_{\alpha'}\circ T_{k,(\rho,\alpha')}(y)=\rho^{r_{1,k}}Q_{\alpha'}(y)$, et donc

$$
(\chi_{pr} Q_{\alpha'})\circ T_{k,(\rho,\alpha')}=
         \rho^{s_k+r_{1,k}}({\Cal Y}_{pr} Q_{\alpha'})
$$

\noindent La condition (b) sur la dérivation d'Hilbert $\chi$ donne

$$
((\chi-\chi_{pr}) Q)\circ T_k=\rho^{s_k+r_{1,k}}O(\rho)
$$

\noindent et donc

$$
 {\Cal Y}\rho=\rho F_k(\rho,y,\alpha')\quad\text{avec}\quad F_k=\frac{1}{r_{1,k}}({\Cal Y}_{pr} Q)+O(\rho)
$$

 Soit $\widehat{\chi}=({\Cal L}_k)_*{\Cal Y}$ défini sur $\widetilde{U}_k$. On a $\widehat{\chi}\rho=\rho(F_k\circ {\Cal L}_k^{-1})$. Or, les fonctions $G_j\circ {\Cal L}_k^{-1}=\rho^{r_{1,j}}u_j$ sont des int\' egrales premi\`eres de $\widehat{\chi}$. Par cons\' equent $\widehat{\chi}u_j=-r_{1,j}u_j(F_k\circ {\Cal L}_k^{-1})$, et le champ recherch\' e est

$$
 \widetilde{\chi}=\frac{1}{F_k\circ {\Cal L}_k^{-1}}\widehat{\chi}
$$

\qed

\vskip 3mm

\noindent {\bf 2.2. Sur le bord de} $S_k^{+*}(0,1)$.\par

\vskip 3mm

  Ce bord est une union finie de sous-ensembles $B_{k',i}$ isomorphes chacun \`a l'un des sous-ensembles

$$
 \{ 0 \} \times S_{k-k',i}^{+*}(0,1)\subset ({\bold{\RR}}^+)^{k'}\times                                                            ({\bold{\RR}}^{+*})^{k-k'}
$$

\noindent avec $k'<k$, l'indice $i$ étant \'enum\' eratif

\vskip 3mm

\proclaim{Proposition IVA2} Soit $y^0\in B_{k',i}$. Il existe un voisinage $V_{k',i}$ de $(0,y^0,0)$ dans ${\bold{\RR}}^{+*}\times S_k^{+*}(1)$ et un difféomorphisme ${\Cal L}_{k',i}$ de $V_{k',i}$ sur son image $\widetilde{U}_{k',i}$ tels que le champ $({\Cal L}_{k',i})_*{\Cal Y}$ soit \' equivalent \`a un \' el\' ement de $\Xi{\Cal H}_{k'}[QR{\Cal H}^{p'(k'),q'(k')}]$ avec $p'(k')=p-k+k'+1$ et $q'(k')=(kq_1+k-1,q_2+k-k'-1)$. Ce morphisme ${\Cal L}_{k',i}$ se prolonge continument et bijectivement sur
$\overline{V}_{k',i}\cap \{0\}\times\overline{S_k^{+*}(0,1)}$. Ces composantes et celles de son inverse sont induits par des éléments de l'algèbre $QR{\Cal H}^{p',q'}$.
\endproclaim

\vskip 3mm

\noindent {\bf Preuve.} Posons $y^0=(y_{1,0},\ldots,y_{k,0})$ et supposons d'abord que $y_{0,1}=0$. Soit $j_0$ le plus petit des entiers $j$ tels que $y_{j,0}=0$ et $y_{j+1,0}>0$. La preuve consiste à trivialiser le champ ${\Cal Y}$ dans la coordonn\' ee $\rho$ en utilisant l'int\' egrale premi\`ere $G_{j_0}$ puis, on le trivialise de la même façon dans les coordonnées $y_j$ telles que $y_{j,0}>0$ et on utilise une récurrence sur le nombre d'intégrales premières non triviales. On a $G_{j_0}=\rho^{r_{1,j_0}}L_{j_0}$ et $L_{j_0}(\rho,y^0,\alpha')=-y_{j_0+1,0}<0$, par conséquent

$$
 G_{j_0}(\rho,y,\alpha')=-y_{j_0+1,0}\rho^{r_{1,j_0}}( 1+O(||y-y^0||) )
$$

\noindent Soit $G=(-G_{j_0})^{1/r_{1,j_0}}=c\rho ( 1+O(||y-y^0||) )$ avec $c>0$ et soit le morphisme

$$
 E_{k',i}(\rho,y,\alpha')=(G(\rho,y,\alpha'),T_{k,(\rho/G,\alpha')}(y),\alpha')=(\rho',y',\alpha')
$$

\noindent par le théorème d'inversion VB4, généralisé à plusieurs variables, il existe un voisinage $V_{k',i}$ tel que le morphisme $E_{k',i}$ soit un difféomorphisme de $V_{k',i}$ sur son image $V'$ qui est une sous-vari\' et\' e de $({{\RR}}^{+*})^{p'}\times {{\RR}}^{|q'|}$ de codimension $1$. De plus, il  se prolonge continument en une bijection de $\overline{V}_{k',i}\cap \{0\}\times\overline{S}_k^{+*}(0,1)$ sur $\overline{V'}\cap \{(\rho',\alpha')=0\}$. Soit ${\Cal Y}'=(E_{k',i})_*{\Cal Y}$, $G$ étant une intégrale première de ${\Cal Y}$, on a ${\Cal Y}'\rho'=0$. Soit le morphisme $T'_k(\rho',y',\alpha')=(T_{k,(\rho',\alpha')}(y'),\alpha')$ défini sur $V'$. D'après la définition du morphisme $E_{k',i}$, on a le diagramme commutatif

$$
\CD
 V_{k',i}    @>E_{k',i}>>      V'  \\
@VT_kVV                    @VVT'_kV \\
U_{k',i}     @>id>>         U_{k',i}
\endCD \tag {$*$}
$$

%$T_k(\rho,y,\alpha')=T'_k(\rho',y',\alpha')=(x,\alpha')$ 

\noindent et donc les fonctions $G'_j=g_j\circ T'_k$ sont des intégrales premières de ${\Cal Y}'$. Or $G'_j=(\rho')^{r_{1,j}}L'_j$, donc les fonctions $L'_j$ sont des intégrales premières de ${\Cal Y}'$ et la sous-vari\' et\' e $V'$ est donn\' ee par l'equation
$L'_{j_0}(\rho',y',\alpha')=-1$ qui est un graphe: $y'_{j_0+1}=f(\rho',y'_{j_0},\alpha')=1+(y'_{j_0})^{r_{j_0}}(1+O(y'_{j_0}))$. Soit la projection canonique $\pi':(\rho',y',\alpha')\in V'\mapsto (\rho',\widehat{y'}^{j_0+1},\alpha')\in V'_1$, le champ ${\Cal Y}'_1=\pi'_*{\Cal Y}'$ a pour intégrales premières les fonctions $L'_j$ pour $j\not\in\{j_0,j_0+1\}$ et la fonction 

$$
 L'_{j_0+1}(\rho',f,y'_{j_0+2},\alpha')-L'_{j_0+1}(\rho',1,0,\alpha')=a(\rho',\alpha')(y'_{j_0})^{r_{j_0}}(1+O(y'_{j_0}))-y'_{j_0+2}
$$

\noindent avec $a(0)>0$. D'après le diagramme $(*)$, on a $(T'_k)^{-1}_*\chi=(\rho')^{s_k}{\Cal Y}'$ et en utilisant l'action de $\chi$ sur $x_1$, on obtient ${\Cal Y}' y'_1=\prod_{j=1}^k y'_j$. Donc si $k-k'=1$, le champ ${\Cal Y}'_1$ est équivalent à un élément de $\Xi{\Cal H}_{k'}$ par la réduction du paragraphe 1. Supposons $k-k'>1$ et réindéxons les coordonnées $y'_j$ et les intégrales premières $L'_j$. Il existe $j_1\geq j_0$ tel que $y'_{j_1,0}=0$ et $y'_{j_1+1,0}>0$. Soit le morphisme 

$$
 {\Cal L}'_1:(\rho',y',\alpha')\in V'_1\mapsto (\rho',\widehat{y'}^{j_1+1},\alpha',L'_{j_1})=(\rho',\widehat{y'}^{j_1+1},\alpha',u'_{j_1})\in V'_2
$$

\noindent L'équation $L'_{j_1}=u'_{j_1}$ est un graphe: $y'_{j_1+1}=f_1(\rho',y'_{j_1},\alpha',u'_{j_1})$. Par la méthode ci-dessus, on montre que le morphisme ${\Cal L}'_1$ est un difféomorphisme sur son image et que le champ ${\Cal Y}'_2=({\Cal L}'_1)_*{\Cal Y}'_1$ admet pour intégrales premières la coordonnée $u'_{j_1}$, les fonctions $L'_j$ pour $j\not\in\{j_1,j_1+1\}$ et la fonction $L'_{j_1+1}(\rho',f_1,y'_{j_1+2},\alpha')-L'_{j_1+1}(\rho',y'_{j_1+1,0},0,\alpha')=a_1(\rho',\alpha',u'_{j_1})(y'_{j_1})^{r_{j_1}}(1+O(y'_{j_1}))-y'_{j_1+2}$. Par l'hypothèse de récurrence, il est équivalent à un élément de $\Xi{\Cal H}_{k'}$.\par

\vskip 3mm

  Supposons $y_{1,0}>0$ et soit $j_0$ le plus petit des entiers $j$ tels que $y_{j,0}>0$ et $y_{j+1,0}=0$. On trivialise le champ ${\Cal Y}$ dans la coordonnée $\rho$ en utilisant l'intégrale première $G_{j_0}$ et le morphisme $E_{k',i}$ associé. Puis, on trivialise le champ ${\Cal Y}'$ dans les coordonnées $y'_1,\ldots,y'_{j_0}$ en utilisant les intégrales premières $L'_1,\ldots,L'_{j_0}$ par l'intérmédiaire du morphisme

$$
 {\Cal L}'_0(\rho',y',\alpha')=(\rho',y'_{j_0+1},\ldots,y'_k,L'_1,\ldots,L'_{j_0-1},\alpha')
$$

\noindent le syst\`eme d'\' equations

$$
 L'_1=u'_1,\ldots,\quad L'_{j_0-1}=u'_{j_0-1},\quad L'_{j_0}=1
$$

\noindent dont les inconnues sont $y'_1,\ldots,y'_{j_0}$, s'inverse ligne par ligne dans l' alg\`ebre $QR{\Cal H}^{p',q'}$. Le champ ${\Cal Y}'_0=({\Cal L}'_0)_*{\Cal Y}'$ admet pour intégrales premières la coordonnée $u'$ et les fonctions $L'_{j_0+1},\ldots,L'_k$. D'après  l'action de $\chi$ sur $x_{j_0+1}$, on a 

$$
 {\Cal Y}' y'_{j_0+1}=r_{1,j_0}(y'_1)^{r_1}\times\cdots\times (y'_{j_0})^{r_{j_0}}\prod_{j>j_0}y'_j(1+O(\rho'))
$$

\noindent et on conclut en utilisant la première partie de la preuve.\qed

\vskip 3mm

 On note $(\pi_k,{\Cal N}_k)$ cette première étape de la désingularisation de $\chi$ de diviseur exceptionnel $\overline{{\Cal D}_k}$. Soit $f\in QR{\Cal H}^{p,q}$ et $\widetilde{f}$ son relevé par $\pi_k$. Si $a\in {\Cal D}_k$, $\widetilde{f}_a$ est induit par un élément de $QR{\Cal H}^{p'(0),q'(0)}$ et si $a\in\partial {\Cal D}_k$ tel que $\widetilde{\chi}_a$ soit de dimension de non trivialité $k'-1$, $\widetilde{f}_a$ est induit par un élément de $QR{\Cal H}^{p'(k'),q'(k')}$.\par

\vskip 5mm

\noindent {\bf B. L'hypothèse } $(H\lambda)$ {\bf et le lemme de récurrence 1.}\par

\vskip 3mm

 Considérons une dérivation $\chi\in\Xi{\Cal H}_k[QR{\Cal H}^{k,q}]$ réalisée sur un voisinage $U$ de 0. Soient $(g_j)$ ses $(k-1)$ intégrales premières non triviales. Son morphisme intégral est $\pi_\chi (x,\alpha)=(\alpha,(g_j(x,\alpha)))=(\alpha,\lambda)$. D'après les propositions IVA1 et IVA2, l'orbite $\gamma=\{\alpha=0,g_1=\cdots=g_{k-1}=0\}$ est principale dans $U$. Soit $W$ une transversale à $\gamma$, analytique de coordonnées $(\alpha,\lambda)$. On définit de la même façon que dans la partie IIIA, les classes ${\Cal C}^k_{\lambda}$ et ${\Cal C}^k_{\lambda,loc}$ des germes $f\in QR{\Cal H}^{k,q}$ qui satisfont globalement ou localement à l'hypothèse $(H\lambda)$.\par

\vskip 3mm

\proclaim{Théorème IVB1} La classe ${\Cal C}^k_{\lambda,loc}$ est localement $\chi$-finie.
\endproclaim

 La preuve est basée sur les théorèmes principaux II1 et IIIA2, et sur un argument de récurrence sur la dimension de non trivialité de la dérivation d'Hilbert. Définissons d'abord les anneaux des intégrales premières qui apparaîtront dans cette partie

\vskip 3mm

\proclaim{Définition IVB1} Soit $t=(t_1,\ldots,t_n)$ des coordonnées sur $\RR^n$. Notons ${\Cal A}_0^n(t)=\RR\{t\}$. Soit $V_0$ un semi-analytique de l'anneau ${\Cal A}_0^n(t)$ qui est ouvert et qui adhère à 0. On note ${\Cal A}_1^n(V_0)$ l'anneau des germes analytiques et bornés sur (un germe en 0 de) $V_0$. Supposons défini l'anneau ${\Cal A}_i^n(.)$ et soit $V_i$ un semi-analytique de cet anneau qui est ouvert et qui adhère à 0. On note ${\Cal A}_{i+1}^n(V_i)$ l'anneau des germes analytiques et bornés sur (un germe en 0 de) $V_i$.
\endproclaim

\vskip 3mm

\noindent{\bf L'argument de récurrence.}\par

\vskip 3mm

 Soit $x=(x_1,\ldots,x_k)$, $\rho=(\rho_1,\ldots,\rho_p)$ avec $p>0$ et $\alpha=(\mu,\nu,\nu')$ des coordonnées sur $\RR^{q_1}\times\RR^{q_2}\times\RR^{q_3}$. Soit $q=(q_1,q_2+q_3)$ et soit une dérivation $\chi \in \Xi{\Cal H}_k[QR{\Cal H}^{p+k,q}]$ réalisée sur un ouvert $U$ et d'int\' egrales premi\`eres non-triviales $g_j=x_j^{r_j}(1+D_j)-x_{j+1}$. Toujours d'après les propositions IVA1 et IVA2, la dérivation $\chi$ admet une orbite principale dans $U$ incluse dans le bord de $U$: $\gamma=\{(\rho,\alpha,(g_j)_j)=0 \}$. Son morphisme intégral est $\pi_\chi:(x,\rho,\alpha)\in U\mapsto (\rho,\alpha,\lambda)=(\rho,\alpha,(g_j(x,\rho,\alpha))_j)\in W$ où $W$ est isomorphe à une semi-transversale à $\gamma$.\par

\vskip 3mm

\noindent $(t_k)$ Soit $V$ un semi analytique d'un anneau ${\Cal A}^{q_2}_i(.)$ qui est ouvert et qui adhère à 0 et soit 

$$
 im:\nu\in V\mapsto (\rho(\nu),\mu(\nu),\nu,\nu'(\nu),\lambda(\nu)))\in W
$$

\noindent une immersion dont les composantes appartiennent à l'anneau ${\Cal A}_{i+1}^{q_2}(V)$.\par

\vskip 3mm

 Soit $W_0=im(V)\subset W$ et soit $U_0=\pi_\chi^{-1}(W_0)$. Notons $U_{0,m}$ et $\pi_{\chi,m}$ les germes de $U_0$ et $\pi_\chi$ en tout point $m\in\gamma$. Soit les morphismes $\pi_\chi^*:SB^{p,|q|+k-1}\to SB^{p+k,|q|}$ et $\pi_{\chi,m}^*:SB^{p,|q|+k-1}\to SB^{p,|q|+k}$

\vskip 3mm

\proclaim{Lemme IVB1 (lemme de récurrence 1)} Soit $f\in QR{\Cal H}^{k+p,q}$, on suppose que
\roster

\item"$(i_k)$" Il existe $N_0=(n_{i,0})_i$ tel que pour tout $m\in\gamma$

$$
 \pi_{\chi,m|U_{0,m}}^*((\rho^{N_0}))\subset {\Cal I}_{\chi,f,m|U_{0,m}}
$$

\item"$(ii_k)$" Il existe $N_1=(n_{i,1})$ avec $n_{i,1}>n_{i,0}$ tel que les germes $f$ et $D_j$ admettent des représentants $\in QR{\Cal H}_{cvg}$ dans les quotients $SB/(\rho^{N_1})$.\par

\endroster

\noindent Alors $f$ est localement $\chi$-finie sur $U_0$.
\endproclaim

\vskip 3mm

\noindent Remarque. Ce lemme est encore vrai dans des situations plus générales pour le sous-ensemble $W_0$. On garde cette formulation simple pour plus de cohérence avec le deuxième lemme de récurrence et parce qu'elle est suffisante pour la preuve du théorème.\par

\vskip 3mm

\noindent {\bf Preuve du lemme.} Par récurrence sur $k$. Le cas $k=1$ est une conséquence de la première partie de la preuve avec $\widetilde{\chi}=\chi$. Soit $k>1$ et $(\pi_k,{\Cal N}_k)$ le premier éclatement de $\chi$ donné par les propositions IVA1 et IVA2, et soit $\overline{{\Cal D}}_k$ son diviseur exceptionnel. Soient $\widetilde{\chi}$ et $\widetilde{f}$ les relevés de $\chi$ et $f$ par $\pi_k$. Soit $\gamma_0=\pi_k^{-1}(\gamma)$. D'après la proposition IVA1, la semi-transversale $W$ est isomorphe à une semi-transversale à $\gamma_0$ de même coordonnée. Soit $\widetilde{U}_0=\pi_k^{-1}(U_0)$, c'est aussi le saturé de $W_0$ par le flot de $\widetilde{\chi}$ dans ${\Cal N}_k$. Il s'agit de montrer que le faisceau ${\Cal I}_{\widetilde{\chi},\widetilde{f}}[\overline{{\Cal D}}_k\cap\overline{\widetilde{U}}_0]$ est localement $\widetilde{\chi}$-fini.\par

\vskip 3mm

\noindent {\bf (a) Au dessus de } ${\Cal D}_k$.\par

\vskip 3mm

 Le désingularisé de $\chi$ est

$$
 \widetilde{\chi}=\rho_0\frac{\partial}{\partial \rho_0}-
                                   \sum_{j=1}^{k-1}r_{1,j} u_j\frac{\partial}{\partial u_j}
$$

\noindent soit $a_0$ son unique singularité sur ${\Cal D}_k$. D'après l'hypothèse $(ii_k)$, il existe une fonction $g$ définie au dessus de tout relativement compact de ${\Cal D}_k$ telle que pour tout $a\in{\Cal D}_k$, le germe $g_a$ est induit par un élément de $QR{\Cal H}^{1,.}_{cvg}(\rho_0,.)$ sur un semi-analytique de $QR{\Cal H}_{cvg}^{p,.}(\rho,.)$ et 

$$
\widetilde{f}_{a}-g_a\in (\rho^{N_1})
\tag 1
$$

\noindent Soit $X(\rho,\mu)$ les fonctions élémentaires de l'algèbre $QR{\Cal H}^p(\rho,.)$ \ et \ soit \ l'immersion $c(\rho_0,\rho,\alpha,u)=(\rho_0,\alpha,X,u)$. Soit $G\in QR{\Cal H}^{1,.}_{cvg}(\rho_0,.)$ tel que $g_{a_0}=c^*(G)$ et soit $ma^+(G)$ sa multiplicité algébrique positive relativement à $\widetilde{\chi}_{a_0}$ le long de $\gamma_0$. Soit ${\Cal U}_0=c(\widetilde{U}_{0,a_0})$ , Le point clé est que la multiplicité algébrique restreinte $ma^+(G_{|{\Cal U}_0})\leq ma^+(G)$ est indépendante des représentants convergents de $f$ et $D_j$ dans les quotients $SB/(\rho^{N_1})$. En effet, d'après la proposition IVA1, on a $\rho_0^{s_k}\widetilde{\chi}\sim (\pi_k^{-1})_*(\chi)$, par conséquent les fibres différentielles le long de $\gamma$ et $\gamma_0$ sont isomorphes et d'après l'hypothèse $(i_k)$, les fibres ${\Cal I}_{\widetilde{\chi}_{a_0},\widetilde{f}_{a_0},m_0|\widetilde{U}_{0,m_0}}$ contiennent l'idéal $\pi_{\widetilde{\chi}_{a_0}|\widetilde{U}_{0,m_0}}^*(\rho^{N_0})$ pour tout $m_0\in\gamma_0$. D'après (1), il en est de même des fibres de $g_{a_0}$ le long de $\gamma_0$ restreintes à $\widetilde{U}_{0,a_0}$. Or le germe $g_{a_0}$ et les intégrales premières de $\widetilde{\chi}$ sont des éléments d'un anneau restriction analytique; le lemme d'isomorphie I4 s'applique: les fibres différentielles de $g_{a_0}$ le long de $\gamma_0$ sont isomorphes et il en est de même de celles de $\widetilde{f}_{a_0}$ restreintes à $\widetilde{U}_{0,a_0}$.\par

\vskip 3mm

 Il s'ensuit deux choses: d'une part, $\widetilde{f}_{a_0|\widetilde{U}_{0,a_0}}$ possède un idéal $\widetilde{\chi}_{a_0}$-transverse le long de $\gamma_0$, qui coincide avec celui de $g_{a_0|\widetilde{U}_{0,a_0}}$ et qui contient l'idéal $(\rho^{N_0})$ (Il en est donc de même pour $f_{|U_0}$ le long de $\gamma$, par l'isomorphisme $\pi_k$). D'autre part, la relation (1) et le lemme IIIA5 montrent que $\widetilde{f}_{a_0|\widetilde{U}_{0,a_0}}$ est {\bf presque quasi-convergente}, et possède donc une {\bf multiplicité algébrique restreinte}, qui conicide avec celle de $g_{a_0|\widetilde{U}_{0,a_0}}$ et qui ne dépend donc pas du choix de $g$. Notons simplement $ma^+$ cette multiplicité.\par

\vskip 3mm

 Soit $J_G$ l'idéal $\widetilde{\chi}_{a_0}$-transverse de $G$ le long de $c(\gamma_0)$ et ${\Cal W}_0=c(W_0)$. Le théorème principal IIIA1 s'applique à l'action de $\widetilde{\chi}_{a_0}$ sur $G$ par restriction à ${\Cal U}_0$

$$
I_{\widetilde{\chi}_{a_0},G|{\Cal U}_0}\supset (\rho_0^{ma^++1})\pi_{\widetilde{\chi}_{a_0}|{\Cal U}_0}^*(J_{G|{\Cal W}_0})
$$ 

\noindent et en appliquant le morphisme $c^*$, on obtient

$$
 I_{\widetilde{\chi}_{a_0},g_{a_0}|\widetilde{U}_{0,a_0}}
\supset (\rho_0^{ma^++1})\pi_{\widetilde{\chi}_{a_0}|\widetilde{U}_{0,a_0}
}^*(J_{g_{a_0}|W_0})
\tag 2
$$ 

\noindent Un plongement de $f$ et $\chi$ est alors nécessaire. Soit $f'$ et $D'_j$ des représentants convergents de $f$ et $D_j$ dans les quotients $SB/(\rho^{N_1})$. Soit $X(x,\mu)$ les fonctions élémentaires de l'algèbre $QR{\Cal H}^{k,.}(x,.)$ et soient  $f"={\bold j}_X^{ma^++3}(f-f')$ et $D"_j={\bold j}_X^{ma^++3}(D_j-D'_j)$. Les germes $f"$ et $D"_j$ appartiennent à l'idéal $(\rho^{N_1})$ et sont induits par des éléments de l'algèbre $QR{\Cal H}^{k,.}_{cvg}(x,.)$ sur un semi-analytique de $QR{\Cal H}^{p,.}(\rho,.)$. En effet, soit $(a_m(\rho,\alpha))$ la famille des coefficients de $f"$ et des $D"_j$ dans leur développement en série de $X$. Notons $v=(a_m-a_m(0))$ ces nouvelles coordonnées, les fonctions $v(\nu)$ appartiennent aussi à l'anneau ${\Cal A}^{q_2}_{i+1}(V)$. Remplaçons les coordonnées $\alpha$ par $\alpha'=(\alpha,v)$ et gardons les mêmes notations pour $im$, $W_0$, $U_0$ ....\par

\vskip 3mm

 Ainsi, quitte à remplacer $f'$ par $f'+f"$ et $D'_j$ par $D'_j+D"_j$, on peut supposer que $f$ et $D_j$ admettent des représentants convergents dans les quotients $SB/(\rho^{N_1}){\Cal M}_x^{ma^++3}$. La relation (1) est donc remplacée par

$$
\widetilde{f}_{a}-g_a\in (\rho_0^{ma^++2}\rho^{N_1})
\tag 3
$$

\noindent et la relation (2) est encore satisfaite par l'invariance de la multiplicité algébrique restreinte $ma^+$.\par

\vskip 3mm

 Soit $D_0={\Cal D}_k\cap\overline{\widetilde{U}}_0$. Montrons que pour tout $a\in D_0$, le germe $\widetilde{f}_a$ est $\widetilde{\chi}_a$-équivalent à $g_a$ sur $\widetilde{U}_{0,a}$. Le théorème de finitude IB1 permettera de conclure. Soit $\gamma_1$ une trajectoire incluse dans $D_0$. Etudions le faisceau ${\Cal I}_{\widetilde{\chi},\widetilde{f}}(\{a_0\}\cup\gamma_1)_{|\widetilde{U_0}}$. En $a_0$, on a $J_{g_{a_0}|W_0}\supset (\rho^{N_0})$ et d'après (2) 

$$
 I_{\widetilde{\chi}_{a_0},g_{a_0}|\widetilde{U}_{0,a_0}}
\supset (\rho_0^{ma^++1}\rho^{N_0})
\tag 4
$$ 

\noindent donc d'après (3), on a $(\widetilde{f}_{a_0}-g_{a_0})_{|\widetilde{U}_{0,a_0}}\in {\Cal M} I_{\widetilde{\chi}_{a_0},g_{a_0}|\widetilde{U}_{0,a_0}}$. Par conséquent, $\widetilde{f}_{a_0}$ est $\widetilde{\chi}_{a_0}$-finie sur $\widetilde{U}_{0,a_0}$.\par

\vskip 3mm

 Soit $b\in\gamma_1$ suffisament proche de $a_0$. D'après (4) et le lemme de cohérence IB3, on a 

$$
 {\Cal I}_{\widetilde{\chi},g,b|\widetilde{U}_{0,b}}
\supset (\rho_0^{ma^++1}\rho^{N_0})
$$

\noindent et donc d'après (3), $\widetilde{f}_{b}$ est $\widetilde{\chi}_{b}$-finie sur $\widetilde{U}_{0,b}$. Soit $a$ un point quelconque de $\gamma_1$, d'après le lemme d'isomorphie IB4, les fibres ${\Cal I}_{\widetilde{\chi},g,b|\widetilde{U}_{0,b}}$ et 
${\Cal I}_{\widetilde{\chi},g,a|\widetilde{U}_{0,a}}$ 
sont isomorphes par le flot de $\widetilde{\chi}$ qui préserve l'idéal $(\rho_0^{ma^++1}\rho^{N_0})$. Par conséquent, $\widetilde{f}_{a}$ est $\widetilde{\chi}_{a}$-finie sur $\widetilde{U}_{0,a}$ et on a 

$$
{\Cal I}_{\widetilde{\chi},\widetilde{f},a|\widetilde{U}_{0,a}}\supset (\rho_0^{ma^++1}\rho^{N_0})
\tag 5
$$

\vskip 3mm

\noindent {\bf (b) Sur le bord de } ${\Cal D}_k$.\par

\vskip 3mm

 Soit $a_1=\overline{\gamma}_1\cap\partial {\Cal D}_k$. Soit $\widetilde{\chi}_{a_1}$ le désingularisé de $\chi$ en $a_1$. D'après la proposition IVA2, c'est une dérivation d'Hilbert de dimension de non trivialité $k'-1<k-1$. Montrons d'abord la propriété $(t_{k'})$. Dans la réduction de $\chi$ au voisinage de $a_1$, on trivialise dans la coordonnée $\rho_0$ en posant par exemple $\rho'_0=\rho_0^{r_{1,1}}|u_1|=|\lambda_1|$ si $u_1(a_1)\neq 0$. Les intégrales premières non triviales $\lambda'_j$ au voisinage de $a_1$ sont des fonctions régulières des rapports

$$
 u'_j=\frac{\lambda_j}{(\rho'_0)^{r_{1,j}}}
$$

\noindent des intégrales premières non triviales en $a_0$. Notons $w$ ces nouvelles coordonnées germifiées autour de $\gamma_1$ et $V'$ le semi-analytique de ${\Cal A}^{q_2}_{i+1}(V)$ correspondant. On obtient ainsi la propriété $(t_{k'})$.\par

\vskip 3mm

 Montrons que $\widetilde{f}_{a_1}$ et $\widetilde{\chi}_{a_1}$ satisfont aux hypothèses $(i_{k'})$ et $(ii_{k'})$ du lemme. Soit $a\in\gamma_1$ suffisament proche de $a_1$. D'après la proposition IVA2, les germes en $a\in\gamma_1$ des champs $\widetilde{\chi}$ et $\widetilde{\chi}_{a_1}$ sont équivalents. Donc, d'après (5) et l'expression de $\pi_k$ en $a_1$, on a

$$
{\Cal I}_{\widetilde{\chi}_{a_1},\widetilde{f}_{a_1},a|\widetilde{U}_{0,a}}\supset ((\rho'_0)^{ma^++1}\rho^{N_0})
$$

\noindent ceci prouve l'hypothèse $(i_{k'})$ avec $N'_0=(ma^++1,N_0)$. Soient $g'_j=(y')_j^{r_j}(1+D'_j)-y'_{j+1}$ les intégrales premières non triviales de $\widetilde{\chi}_{a_1}$. Les germes $f$ et $D_j$ admettent des représentants convergents dans les quotients $SB/(\rho^{N_1}){\Cal M}_x^{ma^++3}$. Donc, d'après l'expression de $\pi_k$, les germes $\widetilde{f}_{a_1}$ et $D'_j$ admettent des représentants convergents dans les quotients $SB/((\rho'_0)^{ma^++2}\rho^{N_1})$. Ceci prouve l'hypothèse $(ii_{k'})$ avec $N'_1=(ma^++2,N_1)$.\qed

\vskip 3mm

\noindent {\bf Preuve du théorème.}\par

\vskip 3mm

\comment
 Soit $(X_\nu,\Gamma_k)$ un déploiement analytique de $X_0$ à $q_0$ paramètres. D'après le théorème de l'appendice B de [M-R], il existe des transversales analytiques $\sigma_j$ et $\tau_j$ au voisinage de chaque sommet $P_j$ tels que l'application de Dulac du coin $P_j$ soit induite par un élément de $QR{\Cal H}^{1,(1,q_0)}(x_j,\mu_j,\nu)$. Posons 
$\alpha=(\mu_1,\ldots,\mu_k,\nu)=(\mu,\nu)$ et $q=(q_1,q_2)=(k,q_0)$. Soit $d_j\in QR{\Cal H}^{1,q}(x_j,\alpha)$ qui induit l'application de dulac de $P_j$ et $\lambda_j(\nu)$ les germes qui déploient les connexions $\gamma_{j,j+1}$. Les cycles du champ $X_{\nu}$ proches de $\Gamma_k$ rencontrent les transversales $\sigma_j$ aux points dont les abscisses $x_j$ sont solutions du système
$$
  d_1(x_1,\alpha(\nu))-f_{2}(x_{2},\alpha(\nu))=\lambda_1(\nu),\ldots,  d_k(x_k,\alpha(\nu))-f_{1}(x_{1},\alpha(\nu))=\lambda_k(\nu)
\tag 1
$$
\noindent les germes $f_j$ \' etant des diff\' eomorphismes analytiques dans la variable $x_j$ qui préservent l'origine et l'orientation. Par la réduction de la partie C1, le système (1) est équivalent au système
$$
  d_1(x_1,\alpha(\nu))-x_{2}=\lambda_1(\nu),\ldots,d_k(x_k,\alpha(\nu))-x_{1}=\lambda_k(\nu)
$$
\noindent Posons $x=(x_1,\ldots,x_k)$ et $\lambda=(\lambda_1,\ldots,\lambda_{k-1})$. Soit $\chi\in\Xi{\Cal H}_k[QR{\Cal H}^{k,q}]$ dont $(k-1)$ intégrales non triviales sont $g_j(x_j,x_{j+1},\alpha) =d_j(x_j,\alpha)-x_{j+1}$. Son morphisme intégral est $\pi_\chi(x,\alpha)=(\alpha,(g_j)_j)=(\alpha,\lambda)$. Soit $W_0$ l'image d'un voisinage de $0$ dans $\RR^{q_0}$ par l'immersion $\nu \mapsto (\alpha(\nu),\lambda(\nu))$ et soit $U_0=\pi_\chi^{-1}(W_0)$, c'est un semi-analytique de $QR{\Cal H}$. Soit $f(x,\alpha)=d_k(x_k,\alpha)-x_1-\lambda_k(\nu)$. Alors $\text{cyc}(X_\nu,\Gamma_k)<\infty$ si et seulement si $f$ est $\chi$-régulière sur $U_0$ et la multiplicité des cycles limites proches de $\Gamma_k$ est majorée si l'idéal différentiel $I_{\chi,f}$ est localement noethérien sur $U_0$.\par
\vskip 2mm
 La trajectoire $\gamma=\{ \alpha=0, g_1=\cdots=g_{k-1}=0 \}$ est principale dans un voisinage de $0$. 
\endcomment

 Soit $f\in{\Cal C}^k_{\lambda,loc}$ et soit $W_0\subset W$ un semi-analytique de $\RR\{\alpha,\lambda\}$ et $N_0$ tels que

$$
 J_{\chi,f,\gamma|W_0}\supset {\Cal M}_{\lambda|W_0}^{N_0} 
\tag 6
$$

\noindent En utilisant la stratification analytique de $W_0$, on peut supposer que c'est un graphe analytique, d'où la propriété $(t_k)$. Notons que la preuve ci-dessous s'applique encore à des sous-ensembles $W_0$ plus généraux.\par

\vskip 3mm

  Soit $U_0=\pi_\chi^{-1}(W_0)$. La preuve reprend certains arguments de la preuve du lemme de récurrence IVB1 dont on reprend les notations. Soit $(\pi_k,{\Cal N}_k)$ le premier éclatement de $\chi$ de diviseur exceptionnel $\overline{{\Cal D}_k}$ et soit $\widetilde{U}_0=\pi_k^{-1}(U_0)$.\par

\vskip 3mm

\noindent {\bf (a) Au dessus de } ${\Cal D}_k$. Soit 

$$
 \widetilde{\chi}=\rho\frac{\partial}{\partial \rho}-
                                   \sum_{j=1}^{k-1}r_{1,j} u_j\frac{\partial}{\partial u_j}
$$

\noindent le désingularisé de $\chi$ et $\widetilde{f}$ le relevé de $f$ par $\pi_k$. Soit 
$\gamma_0=\pi_k^{-1}(\gamma)$ et $a_0$ l'unique singularit\' e de $\widetilde{\chi}$ sur ${\Cal D}_k$. La relation (6) est équivalente à  

$$
 J_{\widetilde{\chi}_{a_0},\widetilde{f}_{a_0},\gamma_0|W_0}\supset {\Cal M}_{\lambda|W_0}^{N_0}
\tag 7
$$

\noindent le germe $\widetilde{f}_{a_0}$ est donc presque quasi-convergent sur $\widetilde{U}_{0,a_0}$ et possède une multiplicité algébrique restreinte $ma^+$. D'après le lemme IIIA11, il est $\widetilde{\chi}_{a_0}$-fini sur $\widetilde{U}_{0,a_0}$ et

$$
 I_{\widetilde{\chi}_{a_0},\widetilde{f}_{a_0}|\widetilde{U}_{0,a_0}}\supset (\rho^{ma^++1})\pi_{\widetilde{\chi}_{a_0}|\widetilde{U}_{0,a_0}}^*(J_{\widetilde{\chi}_{a_0},\widetilde{f}_{a_0},\gamma_0|W_0})
$$

\noindent Donc d'après (7), on a pour tout $j$

$$
 I_{\widetilde{\chi}_{a_0},\widetilde{f}_{a_0}|\widetilde{U}_{0,a_0}}\supset (\rho^{ma^++1}(\rho^{r_{1,j}}u_j)^{N_0})_{|\widetilde{U}_{0,a_0}}
\tag 8
$$

\noindent Soit $b\in \gamma_1$ proche de $a_0$. L'une des coordonnées $u_j(b)$ est non nulle. Le lemme de cohérence IB3 appliqué à (8) donne

$$
{\Cal I}_{\widetilde{\chi},\widetilde{f},b|\widetilde{U}_{0,b}}\supset
(\rho^{N_1})
\tag 9
$$

\noindent Soit $X(\rho,\mu)$ les fonctions élémentaires de l'algèbre $QR{\Cal H}^{1,.}(\rho,.)$ et soit $g={\bold j}_X^{N_1+1}(\widetilde{f})$. Son germe en tout $a\in{\Cal D}_k$ est induit par un élément de $QR{\Cal H}^{1,.}_{cvg}(\rho,.)$ et $\widetilde{f}_a-g_a\in (\rho^{N_1+1})$. D'après (9) et le lemme de Nakayama, les fibres en $b$ de $\widetilde{f}$ et $g$ coincident, et par le lemme d'isomorphie IB4, leurs faisceaux le long de $\gamma_1$ coincident, et on a

$$
{\Cal I}_{\widetilde{\chi},\widetilde{f},a|\widetilde{U}_{0,a}}\supset
(\rho^{N_1})
\tag 10
$$

\vskip 2mm

\noindent {\bf (b) Sur le bord de } ${\Cal D}_k$.\par

\vskip 3mm

  On utilise le lemme de récurrence IVB1. D'après (10), les fibres de $\widetilde{f}$ le long de $\gamma_1$ satisfont à l'hypothèse $(i_{k'})$. L'hypothèse $(ii_{k'})$ est une conséquence de la structure du morphisme $\pi_k$.\qed

\vskip 5mm

\noindent{\bf C. Cas général et lemme de récurrence 2.}\par

\vskip 3mm

  Soient $x=(x_1,\ldots,x_k)$, $\alpha=(\mu,\nu)$ et $\chi\in\Xi{\Cal H}_k$. Grâce au théorème principal IIIB1, on généralise le théorème IVB1 dans le

\vskip 3mm

\proclaim{Théorème IVC1} L'algèbre $QR{\Cal H}^{k,.}(x,\alpha)$ est localement $\chi$-finie.
\endproclaim

\vskip 3mm

\noindent{\bf Preuve.} Soit $f\in QR{\Cal H}^{k,.}$. Nous allons montrer que $f$ est localement $\chi$-finie sur un voisinage $U$ de 0. Soit $(\pi_k,{\Cal N}_k)$ le premier éclatement de la désingularisation de $\chi$ de diviseur exceptionnel $\overline{{\Cal D}}_k$ et soit $\widetilde{\chi}$ et $\widetilde{f}$ les relevés par $\pi_k$ de $\chi$ et $f$. Montrons que $\widetilde{f}$ est localement $\widetilde{\chi}$-finie sur un voisinage dans ${\Cal N}_k$ de tout point de $\overline{{\Cal D}}_k$.\par

\vskip 3mm

\noindent{\bf §1. Au dessus de} ${\Cal D}_k$.\par

\vskip 3mm

 Soit $u=(u_1,\ldots,u_{k-1})$ la coordonnée globale sur ${\Cal D}_k$ et $(\rho,\alpha,u)$ les coordonnées sur ${\Cal N}_k$ au dessus de ${\Cal D}_k$. D'après la partie A, on a 

$$
 \widetilde{\chi}=\rho\frac{\partial}{\partial\rho}-\sum_{j=1}^{k-1} s_ju_j\frac{\partial}{\partial u_j}
$$ 

\noindent et $\widetilde{f}_a\in QR{\Cal H}^{1,.}(\rho,\alpha,u-u_a)$ pour tout $a\in{\Cal D}_k$. Soit $a_0$ l'unique singularité de $\widetilde{\chi}$ sur ${\Cal D}_k$ et $\gamma_0=\{\alpha=0,\ \ u=0\}$ la trajectoire principale de $\widetilde{\chi}_{a_0}$ dans un voisinage $U_{a_0}$ de $a_0$. Soit $\pi_{\widetilde{\chi}_{a_0}}:\ (\rho,\alpha,u)\in U_{a_0}\mapsto (\alpha,\lambda)\in W$ le morphisme intégral de $\widetilde{\chi}_{a_0}$.\par

\vskip 3mm

\noindent{\bf 1.1. En} $a_0$. Le résultat est une conséquence imm\'ediate du théorème principal IIIB1.\par

\vskip 3mm

\noindent{\bf 1.2. En dehors de} $a_0$.\par

\vskip 3mm

 Soit $J\subset\RR\{\alpha,\lambda\}$ l'idéal $\widetilde{\chi}_{a_0}$-transverse de $\widetilde{f}_{a_0}$ le long de $\gamma_0$. L'idée générale est la suivante: on veut préparer $\widetilde{f}$ dans $J$ globalement au dessus de ${\Cal D}_k$. Ceci repose sur une préparation de l'idéal $J+{\Cal M}_\lambda$ en vue de la détermination de la perte d'analycité dans $J$ à la traversée de l'ensemble singulier $\{\rho=0,\ \ u=0\}$. En effet, l'ensemble limite du saturé du sous-ensemble $Z(J)\cap \{\lambda=0\}$ est inclus dans le bord $\{\rho=0\}$. Soit donc $(\psi,{\Cal N})$ une désingularisation dans laquelle l'idéal $J+{\Cal M}_\lambda$ est principal, monomial et ordonné (lemme IIIB1). Soit $(c,V_c)$ une carte de cette désingularisation de coordonnée $v=(v_1,\ldots,v_p)$. Notons $J_c=(\varphi)=\psi_c^*(J)$ avec $\varphi=\prod_{j=1}^p v_j^{n_j}$. Soient $\mu_{j,c}=\psi_c^*(\mu_j)$, $s_{j,c}=1+\mu_{j,c}$ et $\lambda_{j,c}=\psi_c^*(\lambda_j)$. Quitte à réindéxer, on suppose que

$$
 (\lambda_{k-1,c})\subset\cdots\subset (\lambda_{1,c})
\tag 1
$$

\noindent La perte d'analycité dans cette carte est alors l'image dans $W$ du sous-ensemble $\{\lambda_{1,c}=0\}$. Deux cas se présentent:\par

\vskip 2mm

\noindent {\bf (a)}  $(\lambda_{1,c})\subset \text{rad}(\varphi)$. Dans ce cas, la perte d'analycité est totale: l'idéal $J_c$ satisfait à l'hypothèse $(H\lambda)$ relativement à l'idéal $\psi_c^*({\Cal M}_\lambda)$. Le résultat est donc une conséquence du théorème IVB1.\par

\vskip 3mm

\noindent {\bf (b)}  Dans le cas contraire, soit $p'<p$ tel que $\lambda_{1,c}=v_1^{n'_1}\times\cdots\times v_{p'}^{n'_{p'}}(1+O(v))$ avec $n'_j\leq n_j$. Soit $\varphi=\varphi'\varphi"$ l'unique factorisation de $\varphi$ telle que $\text{rad}(\varphi')=\text{rad}(\lambda_{1,c})$ et $\varphi'\wedge\varphi"=1$. Pour obtenir une division globale (au dessus de ${\Cal D}_k$) de $\widetilde{f}$ par $\varphi"=v_{p'+1}^{n_{p'+1}}\times\cdots\times v_p^{n_p}$, il faut plutôt étudier les intégrales premières ramifiées $|\lambda_{j,c}|^{1/s_{j,c}}$ (cf. (3) ci-dessous). Pour cela, une deuxième préparation des intégrales premières $\lambda_{j,c}$ et $\mu_{j,c}$ est nécessaire.\par

\vskip 3mm

 Pour simplifier la présentation, on notera toujours $s_j$ les relevés des fonctions $s_{j,c}$ dans cette préparation. Soient $p"=p-p'$, $v'_0=(v'_{j,0})=(v_1,\ldots,v_{p'})$, et $v"_0=(v"_{j,0})=(v_{p'+1},\ldots,v_p)$. Plaçons nous dans un quadrant dans les coordonnées $v'_0$, par exemple $v'_{j,0}>0$ pour tout $j$. Effectuons une désingularisation dans les coordonnées $v"_0$, et prenons par exemple la carte

$$
 v"_{j,0}=v"_{1,0}\widehat{v}_j
\tag 2
$$

\noindent puis faisons un éclatement dans le couple $(v"_{1,0},\lambda_{1,c})$. Deux cas se présentent:\par

\vskip 3mm

\noindent {\bf (b.1)}  $|v"_{1,0}|<\epsilon \lambda_{1,c}$. Dans ce cas, on pose 

$$
v"_{1,0}=v"_{1,1}\lambda_{1,c},\  v"_{j,1}=\widehat{v}_j \ \text{pour}\  j>1,\ v"_1=(v"_{j,1})\ \text{et}\  v'_1=v'_0
$$

\vskip 3mm

\noindent {\bf (b.2)} $\lambda_{1,c}<(2/\epsilon)|v"_{1,0}|$. Cette situation est couverte par un nombre fini de cartes qui sont de deux types: pour l'un, il existe $j_0$ tel que $(v'_{j_0,0})^{n'_{j_0}}<(2/\epsilon)|v"_{1,0}|$. Dans ce cas, on pose 

$$
v'_{j_0,0}=v'_{j_0,1}((2/\epsilon)|v"_{0,1}|)^{1/n'_{j_0}},\  v'_{j,1}=v'_{j,0}\ \text{pour}\  j\neq j_0,
$$

$$
 v'_{p'+1,1}=|v"_{0,1}|^{1/n'_{j_0}},\ v'_1=(v'_{j,1})\ \text{et} v"_1=(\widehat{v}_2,\ldots,\widehat{v}_{p"})=(v"_{j,1})
$$

\noindent Pour l'autre type, on a $(v'_{j,0})^{n'_{j}}\geq (2/\epsilon)|v"_{1,0}|$ pour tout $j$. Dans ce cas, on pose 

$$
v'_{j,0}=v'_{p'+1,1}v'_{j,1}\ \text{avec}\ \prod_{j=1}^{p'} (v'_{j,1})^{n'_j}=(2/\epsilon)|v"_{0,1}|,
$$

$$
 v'_1=(v'_{j,1})\ \text{et}\  v"_1=(\widehat{v}_2,\ldots,\widehat{v}_{p"})=(v"_{j,1})
$$

\vskip 3mm

 Dans les deux cas, notons $\widetilde{v}_1$ les coordonnées locales au voisinage des coordonnées $v'_.$ ou $v"_.$ qui ne sont pas voisines de 0. Soient $\mu_{j,1}$, $\lambda_{j,1}$ et $\varphi_1$ les relevés des fonctions $\mu_{j,c}$, $\lambda_{j,c}$ et $\varphi$ dans les coordonnées 

$$
v'_1=(v'_{1,1},\ldots,v'_{p'_1,1}),\ v"_1=(v"_{1,1},\ldots,v"_{p"_1,1})\ \text{et}\  \widetilde{v}_1=(\widetilde{v}_{1,1},\ldots,\widetilde{v}_{\widetilde{p}_1,1})
$$

\noindent On a $(\varphi_1)=(\varphi'_1(v'_1))(\varphi"_1(v"_1))$ avec $\text{rad}(\varphi'_1)=\text{rad}(\lambda_{1,1})$. On répète alors ce procédé au plus $p"$ fois appliqué aux cartes (b.2). A une certaine étape $i$ de ce procédé, on obtient $p"_i=0$, auquel cas la perte d'analycité est totale et on est dans la situation de l'hypothèse $(H\lambda)$. Le résultat est alors une conséquence du théorème IVB1.\par

\vskip 3mm

 Plaçons nous maintenant dans une carte (b.1) à une certaine étape $i$. Soient $\mu_{j,i}$, $\lambda_{j,i}$ et $\varphi_i$ les relevés des fonctions $\mu_{j,c}$, $\lambda_{j,c}$ et $\varphi$ dans les coordonnées 

$$
v'_i=(v'_{1,i},\ldots,v'_{p'_i,i}),\ v"_i=(v"_{1,i},\ldots,v"_{p"_i,i})\ \text{et}\  \widetilde{v}_i=(\widetilde{v}_{1,i},\ldots,\widetilde{v}_{\widetilde{p}_i,i})
$$

\noindent On a $(\varphi_i)=(\varphi'_i(v'_i))(\varphi"_i(v"_i))$ avec $\text{rad}(\varphi'_i)=\text{rad}(\lambda_{1,i})$. Soit $(c_i,V_i)$ cette carte de coordonnées $v_i=(v'_i,\widetilde{v}_i,v"_i)$ avec $V_i=V'_i\times\widetilde{V}_i\times V"_i$ et $V'_i$ est un voisinage de 0 dans un quadrant, par exemple $v'_{j,i}>0$ pour tout $j$. Soit $W_i\subset W$ l'image de cette carte et soit $U_i$ le saturé de $W_i$ par le flot de $\widetilde{\chi}$. Il s'agit de diviser $\widetilde{f}$ par $\varphi"_i$ globalement au dessus de $D_i=\overline{U_i}\cap {\Cal D}_k$.\par

\vskip 3mm

 Supposons que les fonctions $\mu_{j,i}$ et $\lambda_{j,i}$ sont préparées sphériquement dans les coordonnées $v"_i$ comme dans le théorème principal\  IIIB1, dont on reprend les notations. Divisons $\rho^S\widetilde{f}_{a_0|U_{i,a_0}}$ par $\varphi"_i$ dans une extension adaptée $(\widetilde{QR{\Cal H}}_{p"_i|{\Cal U}_{0,p"_i}},\pi_{p"_i})$ sur laquelle agit la dérivation ${\Cal X}_{p"_i}$. Soit ${\Cal U}_{p"_i}$ le saturé de ${\Cal U}_{0,p"_i}$ par le flot de ${\Cal X}_{p"_i}$ au dessus de ${\Cal D}_k^{p"_i+1}$ et soit $\Delta_i=\overline{{\Cal U}}_{p"_i}\cap {\Cal D}_k^{p"_i+1}$. Les séries de $\rho^S\widetilde{f}$ dans les variables $w^{(n)}=u^{(n-1)}-u^{(n)}$ sont convergentes sur un voisinage de 0 uniformément au dessus de tout compact de la diagonale de ${\Cal D}_k^{p"_i+1}$. Donc, la division en $a_0$ est globale au dessus de $D_i$ si $\Delta_i$ est inclus dans cette diagonale.\par

\vskip 3mm

 Or, en tout point $a\in D_i\setminus \{a_0\}$, il existe $\ell$ tel que $u_{\ell}(a)\neq 0$. Sur un voisinage de $a$, les fonctions 

$$
 \frac{|u_j|^{1/s_j}}{|u_{\ell}|^{1/s_{\ell}}}=
 \frac{|\lambda_j|^{1/s_j}}{|\lambda_{\ell}|^{1/s_{\ell}}}
\tag 3
$$

\noindent sont des intégrales premières de $\widetilde{\chi}$. Notons $\tau_{j,\ell}(v_i)$ ces derniers rapports écrits dans la carte $V_i$ pour des indices $j,\ell<k$. Soit $V_{\ell,i}=\{ v_i\in V_i;\ \tau_{j,\ell}(v_i)<3\ \text{pour tout}\ j\}$. D'après (1) et l'éclatement (b.1)

$$
 \lambda_{j,i}\in {\Cal M}_{v"_i} \Longrightarrow \lambda_{j,i}\in (\lambda_{1,i}^2){\Cal M}_{v"_i}
\tag 4
$$

\noindent Soit $k_i$ le plus grand indice $j$ tel que $\lambda_{j,i}\not\in{\Cal M}_{v"_i}$. D'après (1) et (4), les sous-ensembles $V_{\ell,i}$ sont vides pour tout $\ell=k_i+1,\ldots,k-1$.\par

\vskip 3mm

 Plaçons nous dans l'un des sous-ensembles $V_{\ell,i}$, par exemple $V_{1,i}$ et soient $W_{1,i}$, $U_{1,i}$, ${\Cal U}_{1,p"_i}$, $D_{1,i}$ et $\Delta_{1,i}$ les sous-ensembles correspondants. D'après l'éclatement (b.1)

$$
 \mu_{j,i}-\mu_{j,i}(v'_i,\widetilde{v}_i,0)\in (\lambda_{1,i}){\Cal M}_{v"_i}
\tag 5
$$

\noindent Donc les rapports $\tau_{j,1}$ tendent vers $\tau_{j,1}(v'_i,\widetilde{v}_i,0)$ quand $v"_i$ tend vers 0 uniformément en $(v'_i,\widetilde{v}_i)$. De plus, le sous-ensemble $V_{1,i}$ contient un produit $V'_{1,i}\times V"_{1,i}$ où $V'_{1,i}=\{(v'_i,\widetilde{v}_i)\in V'_i\times \widetilde{V}_i;\ \tau_{j,1}(v'_i,\widetilde{v}_i)<2\ \text{pour tout}\ j\}$ et $V"_{1,i}$ est un voisinage de 0.\par

\vskip 3mm

 Soit $n\leq p"_i$. Sur $D_{1,i}$ privé d'un voisinage de 0, la coordonnée $u_1$ ne s'annule pas. Et d'après (5), $|u_{1,n}|^{1/s_{1,n}}/|u_1|^{1/s_1}$ tend vers 1 quand $v_i$ tend vers 0 (et même uniformément en $(v'_i,\widetilde{v}_i)$). Donc, sur $\Delta_{1,i}$ on a $u_{1,n}=u_1$. De même, en utilisant les rapports

$$
 \tau'_{j,1,n}=\frac{|u_{j,n}|}{|u_{1,n}|^{s_{j,n}/s_{1,n}}}
$$

\noindent et la relation (5), on montre que sur $\Delta_{1,i}$ on a $u_{j,n}=u_j$ pour tout $j$ et pour tout $n$. Le sous-ensemble $\Delta_{1,i}$ est donc inclus dans la diagonale de ${\Cal D}_k^{p"_i+1}$.\par

\vskip 3mm

 Soit $h$ le quotient de la division de $\rho^S\widetilde{f}$ par $\varphi"_i$ au dessus de $\Delta_{1,i}$. Pour tout $A\in \Delta_{1,i}$, on a $h_A\in \widetilde{QR{\Cal H}}_{p"_i}(\rho,v_i,u-u(A),u^{(1)}-u^{(1)}(A),\cdots,u^{(p"_i)}-u^{(p"_i)}(A))$. 
Soit $A_0=\{\rho=0,\ v_i=0,\ u=u^{(1)}=\cdots=u^{(p"_i)}=0\}$ la singularité principale de ${\Cal X}_{p"_i}$. L'idéal transverse de $h_{A_0|{\Cal U}_{1,p"_i}}$ est $(\varphi'_i)$ qui satisfait à l'hypothèse $(H\lambda)$. On conclut donc au résultat par les méthodes de la partie B appliquées au faisceau ${\Cal I}_{{\Cal X}_{p"_i},h}[\Delta_{1,i}]$.\par

\vskip 3mm

\noindent{\bf §2. Sur le bord de }${\Cal D}_k$.\par 

\vskip 3mm

 On utilise le lemme de récurrence 2 ci-dessous, qu'on a choisi de présenter au §3 pour deux raisons: d'une part, ses idées généralisent simplement celles des §1 et 2 et celles du lemme de récurrence 1, et d'autre part sa présentation est beaucoup plus difficile essentiellement à cause des notations. Soit $\gamma\subset D_{1,i}$ une trajectoire de $\widetilde{\chi}$, $a=\partial {\Cal D}_k\cap\overline{\gamma}$ et $k'$ la dimension de non trivialité de la dérivation d'Hilbert $\widetilde{\chi}_a$. Il s'agit de prouver les hypothèses algébriques $(i_{k'})$ et $(ii_{k'})$ et la propriété topologique $(t_{k'})$ au voisinage de $a$.\par

\vskip 3mm

 Soit $C_{k'}$ le sous-ensemble de $\partial {\Cal D}_k$ constitué des points où la dimension de non trivialté de la dérivation d'Hilbert est $k'-1$. Quitte à réduire les rapports $\tau_{j,1}$ dans un voisinage de leurs valeurs sur $\gamma$, on suppose que l'adhérence de $D_{1,i}$ ne rencontre qu'une seule composante connexe de $C_{k'}$. Soit $\Gamma\subset \Delta_{1,i}$ l'orbite du champ ${\Cal X}_{p"_i}$ correspondante à $\gamma$ et soit $A$ le point de $\overline{\Gamma}$ correspondant à $a$. D'après le §1.2, pour tout point $B$ de $\Gamma$ voisin de $A$, les fibres ${\Cal I}_{{\Cal X}_{p"_i},h,B|{\Cal U}_{1,p"_i}}$ sont noethériennes et satisfont, par restriction à ${\Cal U}_{1,p"_i}$ à l'inclusion

$$
 (\rho^{n_0})\subset {\Cal I}_{{\Cal X}_{p"_i},h,B}
$$

\noindent Or si $b$ est le point de $\gamma$ correspondant à $B$, les germes ${\Cal X}_{p"_i,B}$ et $\widetilde{\chi}_b$ sont réguliers et donc "équivalents": plus précisément, les fonctions 

$$
  u'_{j,n}=\frac{u_{j,n}}{|u_1|^{s_{j,n}/s_1}}
$$

\noindent sont des intégrales premières analytiques de ${\Cal X}_{p"_i}$ le long de $\Gamma$. Soit $u'_{j,n}(\Gamma)$ leurs valeurs sur $\Gamma$. Dans le changement de coordonnées $\nu'_{j,n}=u'_{j,n}-u'_{j,n}(\Gamma)$ germifié au dessus de $\Gamma$, le champ ${\Cal X}_{p"_i}$ est transformé dans le champ $\widetilde{\chi}$. Donc, dans une extension évidente de $SB_{|U_{1,i,b}}$ obtenue par adjonction des coordonnées $\nu'_{j,n}$, les fibres ${\Cal I}_{\widetilde{\chi},\rho^S\widetilde{f},b}$ sont noethériennes et satisfont à la double inclusion
 
$$
(\rho^{n_0})\pi^*_{\widetilde{\chi},b}(\varphi"_i)\subset {\Cal I}_{\widetilde{\chi},\rho^S\widetilde{f},b}\subset\pi^*_{\widetilde{\chi},b}(\varphi"_i)
$$

\noindent et ceci prouve l'hypothèse $(i_{k'})$ car le germe en $b$ de la dérivation d'Hilbert $\widetilde{\chi}_a$ est équivalent à la dérivation $\widetilde{\chi}_b$. L'hypothèse $(ii_{k'})$ est une conséquence de la structure du morphisme de désingularisation $\pi_k$.\par

\vskip 3mm

 Dans la réduction de la dérivation d'Hilbert $\chi$ au voisinage de $a$, on trivialise dans la coordonnée $\rho$ en posant $(\rho')^{s_1}=\rho^{s_1}|u_1|=|\lambda_{1,i}|$. Les intégrales premières non triviales $\lambda'_j$ au voisinage de $a$ sont des fonctions régulières des rapports $\nu'_{j}=\lambda_{j,i}/(\rho')^{s_j}$ d'intégrales premières en $a_0$. On obtient donc la propriété $(t_{k'})$ en posant 

$$
\rho^{(1)}=v'_i,\ L'=L,\ \widetilde{\nu}=\widetilde{v}_i,\ \nu=v"_i,\  \nu'=((\nu'_j),(\nu'_{j,n})),\ V'=V'_{1,i}\ \text{et}\  V=V"_{1,i}
$$

\noindent \qed\par

\vskip 3mm 

\noindent{\bf §3. L'argument de récurrence.}\par

\vskip 3mm

 Définissons d'abord les espaces des intégrales premières correspondants à cette situation et qui généralisent ceux de la partie B

\vskip 3mm

\proclaim{Définition IVC1} Soit $(x,\alpha)$ des coordonnées sur $\RR^k\times \RR^n$. Notons 

$$
 {\Cal A}_0^{k,n}(x,\alpha)={\Cal A}^{k,n}(x,\alpha)
$$

\noindent Soit $V$ un semi-analytique ouvert de ${\Cal A}_0^{k,n}(x,\alpha)$ qui adhère à 0. On note ${\Cal A}_1^{k,n}(V)$ l'anneau des germes de fonctions analytiques et bornées sur (un germe en 0 de) $V$. Supposone définis les anneaux ${\Cal A}_i^{k,n}(.)$. Soit $V$ un semi-analytique ouvert de ${\Cal A}_i^{k,n}(.)$ qui adhère à 0. On note ${\Cal A}_{i+1}^{k,n}(V)$ l'anneau des germes de fonctions analytiques et bornées sur (un germe en 0 de) $V$.\par
\endproclaim

\vskip 3mm

\proclaim{Définition IVC2} Soient $(x,\alpha,\beta)$ des coordonnées sur $\RR^k\times\RR^n\times\RR^m$. Notons 

$$
 {\Cal A}_0^{k,(n,m)}(x,\alpha,\beta)={\Cal A}^{k,n+m}(x,\alpha,\beta)
$$

\noindent Soit $V'$ un semi-analytique ouvert d'un anneau ${\Cal A}_i^{k,n}(.)$ qui adhère à 0, et soit $V$ un voisinage de 0 dans $\RR^m$. On note ${\Cal A}_i^{k,(n,m)}(V'\times V)$ l'anneau des germes de fonctions analytiques et bornées sur (le germe en 0 de) $V'\times {\bold V}$ où ${\bold V}$ est le complexifié de $V$.\par  
\endproclaim

\vskip 3mm

\noindent {\bf Remarque IVC1.} Tout $f\in {\Cal A}_i^{k,(n,m)}(V'\times V)$ est la somme d'une série 

$$
 f=\sum_N f_N(x,\alpha) \beta^N
\tag 6
$$

\noindent convergente sur le produit d'un représentant $V'_f$ de $V'$ et d'un polydisque de $\CC^m$ et dont les coefficients $f_N$ appartiennent à l'anneau ${\Cal A}_{i+1}^{k,n}(V')$ et sont tous réalisés sur $V'_f$.\par

\vskip 3mm

 Soient $x=(x_1,\ldots,x_{k})$, $\rho=(\rho_1,\ldots,\rho_L)$ avec $L>0$ et $\alpha=(\mu,\nu,\widetilde{\nu},\nu')$ une coordonnée sur $\RR^{q_1}\times\RR^{q_2}\times\RR^{q_3}\times\RR^{q_4}$. Soient $q=(q_1,q_2+q_3+q_4)$ et $\chi\in\Xi{\Cal H}_{k}[QR{\Cal H}^{k+L,q}]$. Soit $U$ un voisinage de 0 sur lequel $\chi$ est réalisée et soient $g_j=x_j^{r_j}(1+D_j)-x_{j+1}$ ces intégrales premières non triviales. Posons $g=(g_1,\ldots,g_{k-1})$ et soit $\gamma=\{ (\rho,\alpha,g)=0\}$ l'orbite principale dans $U$. Soit $\pi_\chi:(x,\rho,\alpha)\in U\mapsto (\rho,\alpha,\lambda)=(\rho,\alpha,g)\in W$ le morphisme intégral de $\chi$. Soit $\rho^{(1)}=(\rho_1,\cdots,\rho_{L'})$, $\rho^{(2)}=(\rho_{L'+1},\cdots,\rho_L)$, $\beta=(\rho^{(1)},\widetilde{\nu})$ et soit

$$
 im:(\beta,\nu)\in V'\times V\mapsto (\rho^{(1)},\rho^{(2)}(\beta),\mu(\beta,\nu),\nu,\widetilde{\nu},\nu'(\beta,\nu),\lambda(\beta,\nu))\in W
$$ 

\noindent une immersion telle que\par

\vskip 2mm

\noindent $(t_k)$  $V'$ est un semi-analytique ouvert d'un anneau ${\Cal A}_.^{L',q_3}(.)$ et $V$ est un voisinage de 0 dans $\RR^{q_2}$. De plus, les fonctions composantes de $\rho^{(2)}$ appartiennent à un anneau ${\Cal A}_.^{L',q_3}(V')$ et les fonctions composantes de $\mu$, $\nu'$ et $\lambda$ appartiennent à un anneau ${\Cal A}_.^{L',(q_3,q_2)}(V'\times V)$.\par

\vskip 3mm

 Soit $W_0=im(V'\times V)$ et soit $U_0=\pi_\chi^{-1}(W_0)$. Si $m\in\gamma$, on note $U_{0,m}$ le germe de $U_0$ en $m$.\par

\vskip 3mm

\proclaim{Lemme IVC1 (lemme de récurrence 2)} Soit $f\in QR{\Cal H}^{k+L,q}$. On suppose que\par

\roster

\item"$(i_{k})$" il existe un idéal $J_0\subset\RR\{\nu\}$ et $N_0=(n_{i,0})$ tels que pour tout $m\in\gamma$, la fibre ${\Cal I}_{\chi,f,m|U_{0,m}}$ satisfait, par restriction à $U_{0,m}$ à la double inclusion

$$
 \pi_{\chi,m}^*((\rho^{N_0})J_0)\subset{\Cal I}_{\chi,f,m}\subset\pi_{\chi,m}^*(J_0)
$$

\item"$(ii_{k})$" il existe $N_1=(n_{i,1})$ avec $n_{i,1}>n_{i,0}$ tel que les germes $f$ et $D_j$ admettent des représentants convergents dans les quotients $SB/(\rho^{N_1})$.

\endroster

\noindent Alors $f$ est localement $\chi$-finie sur $U_0$.

\endproclaim

\vskip 3mm

\noindent{\bf Preuve.} Elle est basée sur une récurrence sur la dimension de non-trivialité $k-1$ et sur les arguments de l'algorithme de finitude utilisés dans le lemme de récurrence IVB1. L'étape $k=1$ est une conséquence de ce qui suit. L'idée générale est la suivante: après désingularisation de $\chi$ si nécessaire, on divise le relevé de $f$ dans l'idéal $J_0$. Par l'hypothèse $(i_k)$, le quotient de cette division satisfait alors aux hypothèses du lemme de récurrence IVB1.\par

\vskip 3mm

 Soit $(\pi_k,{\Cal N}_k)$ le premier éclatement de la désingularisation de $\chi$ de diviseur exceptionnel $\overline{{\Cal D}}_k$ et soit $\widetilde{\chi}$ et $\widetilde{f}$ les relevés par $\pi_k$ de $\chi$ et $f$. En identifiant $W_0$ à son image par $\pi_k^{-1}$, soit $\widetilde{U}_0$ le saturé de $W_0$ par le flot de $\widetilde{\chi}$. Montrons que $\widetilde{f}$ est localement $\widetilde{\chi}$-finie sur un voisinage dans $\widetilde{U}_0$ de tout point de $\overline{{\Cal D}}_k\cap\overline{\widetilde{U}}_0$.\par

\vskip 3mm

\noindent{\bf 3.1. Au dessus de} ${\Cal D}_k$.\par

\vskip 3mm

 Soit $u=(u_1,\ldots,u_{k-1})$ la coordonnée globale sur ${\Cal D}_k$ et $(\rho_0,\rho,\alpha,u)$ les coordonnées sur ${\Cal N}_k$ au dessus de ${\Cal D}_k$. On a 

$$
 \widetilde{\chi}=\rho_0\frac{\partial}{\partial \rho_0}-\sum_{j=1}^{k-1} s_ju_j\frac{\partial}{\partial u_j}
$$ 

\noindent et $\widetilde{f}_a\in QR{\Cal H}^{1+L,.}(\rho_0,\rho,\alpha,u-u_a)$ pour tout $a\in{\Cal D}_k$. Soit $a_0$ l'unique singularité de $\widetilde{\chi}$ sur ${\Cal D}_k$ et $\gamma_0=\{\rho=0,\ \alpha=0,\ \ u=0\}$ la trajectoire principale de $\widetilde{\chi}_{a_0}$ dans un voisinage de $a_0$.\par

\vskip 3mm

\noindent{\bf 3.1.1. En} $a_0$.\par

\vskip 3mm

 La preuve reprend, en la généralisant, la méthode du théorème principal IIIB1. Soit $(\psi,{\Cal N})$ une désingularisation (dans les coordonnées $\nu$) dans laquelle l'idéal $J_0$ est principal et monomial. En utilisant (6), on suppose que les fonctions $\mu_j$ et $\lambda_j$ sont préparées sphériquement dans cette désingularisation. Soit donc $(c,V_{c})$ une carte de cette désingularisation de coordonnée $v$ et soit $W_{c}\subset W_0$ et $U_{c}\subset\widetilde{U}_0$ les ensembles correspondants. Soit $J_{c}=(\varphi)=\psi_{c}^*(J_0)$, avec $\varphi=\prod_{j=1}^p v_j^{n_j}$. Soient $\mu_{j,c}=\psi_{c}^*(\mu_j)$ et $\lambda_{j,c}=\psi_{c}^*(\lambda_j)$ avec

$$
 \mu_{j,c}=\mu_{j,p+1-i}+\mu'_{j,p+1-i}\quad\text{et}\quad \lambda_{j,c}=\lambda_{j,p+1-i}+\lambda'_{j,p+1-i}
\tag 7
$$

\noindent les fonctions $\mu_{j,p+1-i}$ et $\lambda_{j,p+1-i}$ sont indépendantes des coordonnées $(v_i,\ldots,v_p)$. Les fonctions $\mu'_{j,p+1-i}$ et $\lambda'_{j,p+1-i}$ appartiennent à l'idéal $(v_1\times\cdots\times v_i)$ dans l'anneau ${\Cal A}(V'\times V_{c})$.\par

\vskip 3mm

 Dans la division de $\rho_0^{S_0}\widetilde{f}_{a_0}$ par $\varphi$, les nouvelles fonctions élémentaires des extensions $\widetilde{QR{\Cal H}}$ portent sur les coordonnées $(\rho_0,\rho)$ de l'algèbre $QR{\Cal H}^{1+L,.}$. De plus les (fonctions) coordonnées $(\mu,\nu')$ analytiques de cette algèbre et les intégrales premières $\lambda$ donnent lieu, dans cette division, à des (fonctions) coordonnées analytiques de l'extension qu'on notera $(\nu')^{(0)}$.\par

\vskip 3mm

  Soit donc $(\widetilde{QR{\Cal H}}^{1+L}_p(\rho_0,\rho,v,\widetilde{\nu},(\nu')^{(0)},(u_0^{(j)}))_{|{\Cal U}_p},\pi_p)$ l'extension dans laquelle la division de $\pi_p^*(\rho_0^{S_0}\widetilde{f}_{a_0})$ par $\varphi$ est réalisée. Cette extension est définie par les formules sphériques (7), et par une induction sur $L$ comme pour les algèbres $QR{\Cal H}^{1+L,.}$. Soit $h_0$ le quotient de cette division. Soit ${\Cal X}_p$ la dérivation qui relève $\widetilde{\chi}_{a_0}$ sur $U_{c,0}$ germe de $U_{c}$ en $a_0$. Soit $\Gamma_0$ son orbite principale, et $A_0$ sa singularité principale. On a $\rho_0\neq 0$ le long de $\gamma_0$. Donc en appliquant le morphisme $\pi_p^*$ ($\pi_p$ étant germifié le long de $\Gamma_0$) à la double inclusion $(i_k)$, on obtient que les fibres de $h_0$ le long de $\Gamma_0$ contiennent $\rho^{N_0}$ (par restriction à ${\Cal U}_p$).\par

\vskip 3mm

  Dés lors, on applique la méthode du lemme de récurrence IVB1. Soit $ma^+$ la multiplicité algébrique d'un représentant convergent de $h_0$ modulo $(\rho^{N_1})$, il en existe un d'après l'hypothèse $(ii_k)$. Ce représentant satisfait globalement à la double inclusion en $A_0$ d'après le théorème principal IIIA1. Donc, en prenant un plongement de $f$ et des fonctions $D_j$ à un ordre quelconque $n_{0,1}>ma^+$ (cf. lemme de récurrence IVB1), on obtient que $h_0$ est ${\Cal X}_p$-finie sur une restriction ${\Cal U}_{p,n_{0,1}}$ de ce plongement (et il en est de même pour $\widetilde{f}_{a_0}$ sur $U_{c,0,n_{0,1}}$ germe en $a_0$ d'une restriction $U_{c,n_{0,1}}$ de ce plongement). De plus on a la double inclusion par restriction à ${\Cal U}_{p,n_{0,1}}$

$$
 (\rho_0^{ma^++1}\rho^{N_0}\varphi)\subset I_{{\Cal X}_p,\pi_p^*(\rho_0^{S_0}\widetilde{f}_{a_0})}\subset (\varphi)
\tag 8
$$

\vskip 3mm

\noindent{\bf 3.1.2 En dehors de} $a_0$.\par

\vskip 3mm

 Comme dans le paragraphe 1.2, on veut déterminer la factorisation $J_c=(\varphi')(\varphi")$ telle que le sous-ensemble image dans $W_c$ de $\{\varphi'=0\}$ représente la perte d'analycité à la traversée des singularités $a_0(\alpha)$ et telle que la division de $\widetilde{f}$ par $\varphi"$ soit globale au dessus du sous-ensemble correspondant $D_c\subset{\Cal D}_k$. On généralise la préparation du paragraphe 1.2 aux anneaux ${\Cal A}(V'\times V_c)$ de la façon suivante: soit  $\lambda'_{j,c}(\beta)=\lambda_{j,c}(\beta,0)$ et $\lambda"_{j,c}=\lambda_{j,c}-\lambda'_{j,c}$. D'après (2), on a $\lambda"_{j,c}\in (v_1)$. Donc, l'ensemble limite du saturé par le flot de $\widetilde{\chi}$ de l'image dans $W_c$ du sous-ensemble $\{(\lambda'_{j,c})=0,\ v_1=0\}$ est inclus dans l'ensemble singulier  de $\widetilde{\chi}$ au dessus de $a_0$. Considérons donc les ensembles $V'_{\ell,c}=\{|\lambda'_{j,c}|<2|\lambda'_{\ell,c}|;\ \text{pour tout}\ j\}$ pour $\ell=1,\ldots,k-1$ et $V'_{k,c}=\{|\lambda'_{j,c}|<2|v_1|,\ \text{pour tout}\ j\}$. La pr\'eparation (pr) est la suivante\par

\vskip 3mm

\noindent {\bf (pr1)}  Sur $V'_{k,c}$, on pose 

$$
v'_1=v_1=(v'_{1,1}),\ v"_1=(v_2,\ldots,v_p)=(v"_{1,1},\ldots,v"_{p"_1,1})\ \text{et}\  \beta_1=(\beta,v'_1)
$$

\noindent Remarquons que $p"_1<p$. Les fonctions $\lambda_{j,c}$ se relèvent dans les fonctions   $\lambda_{j,1}=v'_{1,1}(\lambda'_{j,1}(\beta_1)+O(v"_1))$.\par

\vskip 3mm

\noindent {\bf (pr2)}  Sur $V'_{1,c}$ par exemple, deux cas se présentent:\par

\vskip 3mm

\noindent {\bf (pr2.1)} $|v_1|\leq \epsilon |\lambda'_{1,c}|$. On pose alors 

$$
v_1=v"_{1,1}\lambda'_{1,c},\ v'_1=\lambda'_{1,c}=(v'_{1,1}),\ \beta_1=(\beta,v'_1)\ \text{et}\ v"_1=(v"_{1,1},v_2,\ldots,v_p)=(v"_{j,1})
$$

\noindent Les fonctions $\lambda_{j,c}$ se relèvent dans les fonctions $\lambda_{j,1}=v'_{1,1}(\lambda'_{j,1}(\beta_1)+O(v"_1))$. Remarquons que $\lambda'_{1,1}\equiv 1$.\par

\vskip 3mm

\noindent {\bf (pr2.2)} $\epsilon|\lambda'_{1,c}|<|v_1|<2|\lambda'_{1,c}|$. Ce cas se traîte comme le cas (pr1).\par

\vskip 3mm

 On répète ce procédé au plus $p$ fois appliqué aux cartes (pr1). A une certaine étape, on obtient $p"_.=0$. La perte d'analycité est alors totale.\par
 
\vskip 3mm

 Plaçons nous dans une carte (pr2.1) à une étape $i$. Les relevés des fonctions $\lambda_{j,c}$ s'écrivent 

$$
\lambda_{j,i}=\prod_{\ell=1}^{p'_i} v'_{\ell,i}(\lambda'_{j,i}(\beta_i)+O(v"_i))\ \text{avec}\ \lambda'_{1,i}\equiv 1
$$

\noindent D'après (7), on a $O(v"_i)=O(v"_{1,i})$. Dans ce cas, on applique la méthode préparatoire (pr) ci-dessus aux (fonctions) coordonnées $\lambda'_{j,i}\not\equiv 1$ et à la coordonnée $v"_{1,i}$. Puis, on applique la méthode du paragraphe 1.2 (cas (b)) au couple constitué du facteur dominant de la  préparation (pr) et de la (fonction) coordonnée $\lambda_{1,i}$ qui représente la perte d'analycité. On répète ce procédé au plus $(k-2+p"_i)$ fois appliqué aux cartes (b.2). A une certaine étape de cette préparation, on obtient $p"_.=0$.\par

\vskip 3mm

 Dans tous les cas, notons $\widetilde{v}_i$ les coordonnées locales au voisinage des coordonnées $v"_.$ qui ne sont pas voisines de 0. Posons $v_i=(v'_i,\widetilde{v}_i,v"_i)$ et $\widetilde{\beta}_i=(\beta_i,\widetilde{v}_i)$. On obtient la factorisation du relevé $(\varphi_i)=(\varphi'_i(v'_i))(\varphi"_i(v"_i))$ et la préparation des relevés  

$$
\lambda_{1,i}=\prod_{j=1}^{p'_i} (v'_{j,i})^{n'_{1,j,i}}(1+O(v"_i))
$$

\noindent et pour $\ell\neq 1$

$$
\lambda_{\ell,i}=\prod_{j=1}^{p'_i} (v'_{j,i})^{n'_{\ell,j,i}}(\lambda'_{j,i}(\widetilde{\beta}_i)+O(v"_i))\ \text{avec}\  n'_{\ell,j,i}\geq n'_{1,j,i}
$$

\vskip 3mm

 D'après cette double préparation, il existe une partition $v'_i=(\rho'_i,\rho"_i)$ telle que les coordonnées $\rho"_i$ soient des fonctions des coordonnées $\beta'_i=(\beta,\rho'_i,\widetilde{v}_i)$ appartenant à un anneau ${\Cal A}_.(V'_i)$ où $V'_i$ est un semi-analytique ouvert d'un anneau ${\Cal A}_.(.)$ donné par la préparation ci-dessus. Comme dans le paragraphe 1, on utilise (7) pour montrer que les rapports $\tau_{j,\ell}(\beta'_i,v"_i)$ convergent vers $\tau_{j,\ell}(\beta'_i,0)$ quand $v"_i$ tend vers 0 uniformément en $\beta'_i$ sur les sous-ensembles non vides $V_{\ell,i}=\{(\beta'_i,v"_i)\in V'_i\times V_i;\ \tau_{j,\ell}<3\ \text{pour tout}\ j\}$ où $V_i$ est un voisinage de 0 dans $\RR^{p"_i}$ donné par la préparation ci-dessus. De plus ce sous-ensemble $V_{\ell,i}$ contient le produit de $V'_{\ell,i}=\{ \beta'_i\in V'_i;\ \tau_{i,\ell}<2\ \text{pour tout}\ j\}$ qui est un semi-analytique ouvert de ${\Cal A}_.(V'_i)$, et d'un voisinage $V"_{\ell,i}$ de 0.\par

\vskip 3mm

 Plaçons nous dans l'un des sous-ensembles non vides $V_{\ell,i}$, par exemple $V_{1,i}$ et 
supposons que les fonctions $\mu_{j,i}$ et $\lambda_{j,i}$ sont préparées sphériquement dans les coordonnées $v"_i$. Soit $h$ le quotient de la division de $\rho_0^{S_1}\widetilde{f}$ par $\varphi"_i$ en $a_0$ (si $p"_i=0$, on prend $h=\widetilde{f}$). Soit $n'_i$ le plus petit entier tel que $\lambda_{1,i}^{n'_i}\in (\varphi'_i)$. Choisissons un ordre $n_{0,1}$ de plongement de $f$ et des fonctions $D_j$ (cf. §3.1.1) tel que

$$
 n_{0,1}>n_{0,0}=ma^++1+S_1(0)+2n'_i
$$

\noindent Soit $(\widetilde{QR{\Cal H}}^{1+L}_{p"_i|{\Cal U}_{p"_i}},\pi_{p"_i})$ l'extension associée à cette division. Soient $W_{1,i}$, $U_{1,i}\subset U_{c,n_{0,1}}$, ${\Cal U}_{1,p"_i}$, $D_{1,i}$ et $\Delta_{1,i}$ les sous-ensembles correspondants à $V_{1,i}$. Comme dans le paragraphe 1, on montre que $\Delta_{1,i}$ est inclus dans la diagonale de ${\Cal D}_k^{p"_i+1}$. Donc cette division est globale au dessus de $\Delta_{1,i}$ et pour tout $A\in \Delta_{1,i}$, on a $h_A\in \widetilde{QR{\Cal H}}^{1+L}_{p"_i}(\rho_0,\rho,v_i,\widetilde{\nu},(\nu')^{(1)},u-u(A),u^{(1)}-u^{(1)}(A),\cdots,u^{(p"_i)}-u^{(p"_i)}(A))$. Les coordonnées $(v')^{(1)}$ sont définies comme les coordonnées $(v')^{(0)}$.\par

\vskip 3mm

 En tout point $b\in D_{1,i}\setminus\{a_0\}$, les champs $\widetilde{\chi}_b$ et ${\Cal X}_{p,B}$ (cf. §3.1.1) sont équivalents. Comme dans le paragraphe 1, dont on reprend les notations, on construit une extension de $SB_{|U_{1,i,b}}$ par adjonction de (fonctions) coordonnées $(\nu')^{(2)}$ qui appartiennent à un anneau ${\Cal A}_.(V'_{1,i}\times V"_{1,i})$ d'après (7) et la préparation ci-dessus. Et d'après (8) et le lemme de cohérence I3, si en plus $b\in\overline{U}_{c,0,n_{0,1}}$, on a la double inclusion dans cette extension

$$
 (\rho_0^{ma^++1}\rho^{N_0}\varphi_i)\subset {\Cal I}_{\widetilde{\chi},\rho_0^{S_0}\widetilde{f},b}\subset (\varphi_i)
\tag 9
$$

\noindent De même, en tout point $b\in D_{1,i}\setminus \{a_0\}$, les champs $\widetilde{\chi}_b$ et ${\Cal X}_{p"_i,B}$ sont équivalents. On construit donc une extension de cette dernière extension par adjonction de (fonctions) coordonnées $(\nu')^{(3)}\in{\Cal A}_.(V'_{1,i}\times V"_{1,i})$ et on note $H$ l'image de $h$ dans cette extension. On a donc sur une restriction adéquate au dessus de $D_{1,i}\setminus\{a_0\}$

$$
 \rho_0^{S_1}\widetilde{f}=\varphi"_i H
\tag 10
$$

\noindent et si $b$ est voisin de $a_0$, la double inclusion (9) donne

$$
 (\rho_0^{S_1+ma^++1}\rho^{N_0}\varphi'_i)\subset {\Cal I}_{\widetilde{\chi},H,b}
\tag 11
$$

\noindent Or sur $D_{1,i}\setminus\{a_0\}$, la coordonnée $u_1$ ne s'annule pas et on a $\lambda_{1,i}^{n'_i}\in (\varphi'_i)$. Donc l'inclusion (11) donne

$$
 (\rho_0^{n_{0,0}}\rho^{N_0})\subset {\Cal I}_{\widetilde{\chi},H,b}
\tag 12
$$

\noindent et par le plongement ci-dessus à l'ordre $n_{0,1}$, le germe de $H$ en tout point de $D_{1,i}$ est convergent modulo $(\rho_0^{n_{0,1}}\rho^{N_1})$. Pour conclure, on applique la méthode du lemme de récurrence IVB1 à $H$ au dessus de $D_{1,i}$ et on utilise (10).\par

\vskip 3mm

\noindent{\bf 3.2 Sur le bord de } ${\Cal D}_k$.\par

\vskip 3mm

 Soit $a\in \overline{D}_{1,i}\cap\partial {\Cal D}_k$ et $k'-1$ la dimension de non trivialitéde la dérivation d'Hilbert en $a$. Il reste simplement à nommer les données du lemme en fonction de $k'$. On pose $\rho_{k'}^{(1)}=(\rho^{(1)},\rho'_i)$, $\widetilde{\nu}_{k'}=(\widetilde{\nu},\widetilde{v}_i)$, $\nu_{k'}=v"_i$, $\rho_{k'}^{(2)}=(\rho^{(2)},\rho"_i)$, $\nu'_{k'}=(((\nu')^{(j)})_{j=0,\ldots 4},(\nu'_j))$ où les (fonctions) coordonnées $(\nu')^{(4)}$ proviennent du plongement à l'ordre $n_{0,1}$ et les $\nu'_j$ sont déterminés par les rapports $\tau_{j,1}$ comme dans le paragraphe 2. On pose aussi $V'_{k'}=V'_{1,i}$, $V_{k'}=V"_{1,i}$, $N_{0,k'}=(n_{0,0},\ldots,n_{0,0},N_0)$ et $N_{1,k'}=(n_{0,1},\ldots,n_{0,1},N_1)$, $n_{0,0}$ et $n_{0,1}$ étant répétés $p'_i$ fois. Le reste de la preuve est maintenant classique!\qed

\vskip 3mm

\noindent {\bf Preuve du théorème 0.}\par

\vskip 3mm

 Soit $(X_\nu,\Gamma_k)$ un déploiement analytique de $X_0$ à $q_0$ paramètres. D'après le théorème VA1 (appendice VA), il existe des transversales analytiques $\sigma_j$ et $\tau_j$ au voisinage de chaque sommet $P_j$ tels que l'application de Dulac du coin $P_j$ soit induite par un élément $d_j\in QR{\Cal H}^{1,(1,q_0)}(x_j,\mu_j,\nu)$. Posons 
$\alpha=(\mu_1,\ldots,\mu_k,\nu)=(\mu,\nu)$ et $q=(q_1,q_2)=(k,q_0)$. Soient  $\lambda_j(\nu)$ les germes qui déploient les connexions $\gamma_{j,j+1}$. Les cycles du champ $X_{\nu}$ proches de $\Gamma_k$ rencontrent les transversales $\sigma_j$ aux points dont les abscisses $x_j$ sont solutions du système

$$
  d_1(x_1,\alpha(\nu))-f_{2}(x_{2},\alpha(\nu))=\lambda_1(\nu),\ldots,  d_k(x_k,\alpha(\nu))-f_{1}(x_{1},\alpha(\nu))=\lambda_k(\nu)
\tag 13
$$

\noindent les germes $f_j$ \' etant des diff\' eomorphismes analytiques dans la variable $x_j$ qui préservent l'origine et l'orientation. Par la réduction de la partie A, le système (13) est équivalent au système

$$
  d_1(x_1,\alpha(\nu))-x_{2}=\lambda_1(\nu),\ldots,d_k(x_k,\alpha(\nu))-x_{1}=\lambda_k(\nu)
$$

 Posons $x=(x_1,\ldots,x_k)$ et $\lambda=(\lambda_1,\ldots,\lambda_{k-1})$. Soit $\chi\in\Xi{\Cal H}_k[QR{\Cal H}^{k,q}]$ dont $(k-1)$ intégrales non triviales sont $g_j(x_j,x_{j+1},\alpha) =d_j(x_j,\alpha)-x_{j+1}$. Soit $W_0$ l'image d'un voisinage de $0$ dans $\RR^{q_0}$ par l'immersion $\nu \mapsto (\alpha(\nu),\lambda(\nu))$ et soit $U_0=\pi_\chi^{-1}(W_0)$. Soit $f(x,\alpha)=d_k(x_k,\alpha)-x_1-\lambda_k(\nu)$. Alors la partie $(i)$ du théorème est équivalente à la $\chi$-régularité de $f$ sur $U_0$ et la partie $(ii)$ est vérifiée si l'idéal différentiel $I_{\chi,f}$ est localement noethérien sur $U_0$ (à extensions étoilées près). Le théorème 0 est donc une conséquence simple du théorème IVC1 ci-dessus.\qed\par

\vskip 5mm

\centerline{{\bf Appendice V.}}\par

\vskip 5mm

\noindent {\bf A. Déploiements d'applications de Dulac.}\par

\vskip 3mm

 Il est connu ( [E] et [I] ) que l'application de Dulac d'une singularit\'e hyperbolique r\'eelle est un \'el\'ement de $QA^{1,0}$. Un th\'eor\`eme de [M-M] dit que l'application de Dulac est "convergente" si et seulement si la singularit\'e est analytiquement normalisable. Ceci limite consid\'erablement le champ d'application des arguments classiques de la g\'eom\'etrie analytique et justifie l'approche quasi-analytique adopt\'ee qui prolonge celles d'Ecalle et d'Il'yashenko.\par

\vskip 3mm

 On se limitera au cas d'une \'equation diff\'erentielle r\'esonnante de nombre caract\'eristique $r=1$; le cas $r=n/m$ s'en d\'eduit par une double ramification et le cas quasi-r\'esonnant pr\'esente moins de difficult\'e. Soit $\omega_\nu =xdy+y(1+\mu(\nu)+a(x,y,\nu))dx$ un  d\'eploiement analytique \`a $q$ param\`etres $\nu=(\nu_1,\ldots,\nu_q)$ d'une 1-forme analytique r\'eelle r\'esonnante. Ce d\'eploiement est induit par le d\'eploiement

$$
  \omega_\alpha =xdy+y(1+\mu+a(x,y,\nu))dx
$$

\noindent avec $\alpha =(\mu,\nu)\in W_1\times W_q \in ({{\bold{\RR}}},0)\times ({{\bold{\RR}}}^q,0)$. On pose $r=1+\mu$. L'objet de cette partie est le\par

\vskip 3mm

\proclaim{Théorème VA1} Soit $d(.,\alpha)$ l'application de Dulac de $\omega_\alpha$, alors il existe $D\in QR{\Cal H}^{1,(1,q)}$ tel que $d=x^r(1+D)$ et $D(0)=0$.
\endproclaim

\vskip 3mm

 La propri\'et\'e de quasi-analycité a \'et\'e  d\'emontr\'ee dans [E-M-R] sur des domaines de ${\Cal E}{\Cal I}$ de type puissance 

$$
\{ w=u+iv;\quad |v|<Cu^n;\quad u\gg 1\},\quad C>0,\quad n\geq 2 
\tag 1
$$

\noindent en utilisant l'id\'ee g\'eom\`etrique d'Il'yashenko [I] bas\'ee sur la structure de l'holono\-mie de l'une des s\'eparatrices. Or les domaines de ${\Cal E}{\Cal I}$ qui s'imposent naturellement dans le probl\`eme sont les domaines qui sont optimaux pour les formes normales; i.e les domaines d'Ecalle [E] de type exponentiel 

$$
\{ w=u+iv;\quad |v|<C\exp(u/K)-1;\quad u\gg 1\},\quad C>0,\quad K>1
\tag 2
$$

\noindent On adopte une d\'emarche qui combine l'id\'ee g\'eom\`etrique d'Il'yashenko en modifiant les chemins d'int\'egration, et celle de Dulac [D] qui consiste \`a construire une int\'egrale premi\`ere analytique dans l'une des variables. On suppose que $\omega_\alpha$ est pr\'epar\'ee \`a l'ordre 1

$$
a=xy\sum_{n\geq 1} a_n(x,\nu)y^{n-1}
\tag 3
$$

\noindent et qu'on a la majoration suivante

$$
\sum_{n\geq 1} \| a_n\|_{\overline{D}(0,1)\times {{\bold{W}}}_q}\leq 1/4  
\tag 4
$$

\noindent o\`u ${{\bold{W}}}_q$ est le complexifi\'e de $W_q$. Soit

$$
 f(x,y,\alpha)=\sum_{n\geq 1} f_n(x,\alpha)y^n
\tag 5
$$

\noindent  l'int\'egrale premi\`ere de $\omega_\alpha$ telle que $f(1,y,\alpha)\equiv y$. L'application de Dulac de $\omega_\alpha$ est analytiquement conjugu\'ee \`a \par

$$
 d(x,\alpha)=f(x,1/2,\alpha)
\tag 6
$$

\noindent et le th\'eor\`eme est cons\'equence de la

\vskip 3mm

\proclaim{Proposition VA1} Il existe  $F\in QR{\Cal H}^{1,(1,q+1)}$ tel que $f=x^ry(1+F) \text{ et } F(0)=0$.
\endproclaim
 
\vskip 3mm

\noindent {\bf §1. Opérateur intégral de Dulac.}\par

\vskip 3mm

 Notons $P^+=\{w\in{{\bold{\CC}}};\quad \text{Re}(w)>0\}$. Soit $\Omega$ un ouvert simplement connexe de $P^+$ tel que $0\in \overline \Omega$ et $H_c (\Omega)$ l'espace des fonctions holomorphes sur $\Omega$ et continues sur $\overline \Omega$. Pour tout $s\in {{{\bold{\CC}}}}^*$, on d\'efinit sur $H_c (\Omega)$ l'op\'erateur ${\Cal L}_s$ par

$$ 
 {\Cal L}_s(f)(w)=
          s\exp(-sw)\int_{\gamma_w} \exp((s-1)z) f(z) dz 
\tag 7
$$

\noindent o\`u $\gamma_w\subset \Omega$ est un chemin ${\Cal C}^1$ régulier joignant $0$ et $w$. Pour $s$ fix\'e, l'ensemble $L_s=\{ w\in P^+;\quad \text{Re}(sw)=0 \}$ est dit direction singuli\`ere de l'operateur ${\Cal L}_s$ et $\Omega_s=\{ w\in P^+;\text{ Re}(sw)>0 \}$ est dit domaine non-singulier de ${\Cal L}_s$. On montre que, sous une certaine condition g\'eom\`etrique sur le chemin $\gamma_w$, et donc sur l'ouvert $\Omega$, l'op\'erateur ${\Cal L}_s$ est 2-lipschitzien.\par

\vskip 3mm

\proclaim{Lemme VA1} Si $\Omega \subset \Omega_s$ et si le long du chemin $\gamma_w$ la condition suivante est satisfaite 

$$ 
| \tan (\arg (s\frac{dz}{dt})) |\leq |\exp(z) |   
\tag 8
$$

\noindent alors pour tout $f\in H_c(\Omega)$

$$
 | {\Cal L}_s (f)(w) |\leq 2\| f \|_{\gamma_w}  
\tag 9
$$

\endproclaim

\vskip 3mm

\noindent {\bf Preuve.} Soit $t_0$ tel que $z(t_0)=w$. Gr\^ace \`a (8), on a la majoration 

$$ 
|{\Cal L}_s (f)(w) |\leq 2\|f\|_{\gamma_w} |\exp(-sw)|\int_0^{t_0}
                                    |\text{Re}(sz^{'})|\exp(\text{Re}(sz)) dt  
$$

\noindent le chemin $\gamma_w$ \'etant ${\Cal C}^1$-régulier, la condition (8) implique que la fonction $\text{Re}(sz')$ ne s'annule pas sur $\gamma_w$. Et comme $\Omega \subset \Omega_s$, le r\'esultat en d\'ecoule facilement.\qed\par

\vskip 3mm

\noindent {\bf §2. Preuve de la proposition.}\par

\vskip 3mm

  C'est une conséquence du lemme VA1 et des lemmes ci-dessous. Les coefficients $f_n$ de la s\'erie de $f$ v\'erifient les \'equations diff\'erentielles

$$  
x\frac{\partial f_n}{\partial x}-nrf_n=
           x{\sum}_{p=1}^{n-1} pa_{n-p}f_p  
\quad\text{ et }\quad f_1=x^r
$$

\noindent Il est clair que chaque fonction $f_n$ est holomorphe, mais n'est pas forc\'ement born\'ee sur $P^+\times {\bold{W}}_{q+1}$. Pour $n>1$, posons 

$$
 h_n=-\frac{1}{nr}{\sum}_{p=1}^{n-1} pa_{n-p}f_p
$$

\noindent alors, les coefficients $f_n$ sont donn\'es par

$$
 f_n ={\Cal L}_{nr} (h_n)  
$$

\noindent et d'apr\`es (4) on obtient 

$$
\| h_n \|_* \leq \frac{1}{2}\max_{1\leq p\leq n-1}\| f_p \|_* 
$$

\noindent o\`u $*$ est un domaine qui d\'ependra du contexte. De m\^eme, soit $F_n=f_n/f_1$ les coefficients de la s\'erie de  $1+F$ et

$$
 H_n=-\frac{1}{(n-1)r}\sum_{p=1}^{n-1} pa_{n-p} F_p
$$

\noindent alors, on a

$$
 F_n={\Cal L}_{(n-1)r} (H_n)
$$

\noindent En particulier, la fonction \'el\'ementaire $z(x,\mu)=x\text{Ld}(x,\mu)$ est donnée par

$$
 z={\Cal L}_r(-\frac{1}{r})
$$

\vskip 2mm

\proclaim{Lemme VA2} Les germes de $f$ et $F$ sont des \'el\'ements de $SB^{1,q+2}$.\endproclaim

\noindent {\bf Preuve.}  Le domaine non-singulier de l'opérateur ${\Cal L}_{nr}$ coincide avec $\Omega_r$. Soit $\theta \in ]0,\pi/2[$ et $S_\theta =\{ w\in P^{+};\quad |\arg(w)|\leq \theta \}$. Si $\arg(r)$ est suffisament petit, le secteur $S_\theta$ est inclus dans $\Omega_r$. Soit $u_0>0$ tel que

$$
 \exp(u_0)\gg \tan(\theta)  
\tag 10
$$

\noindent et $w_0\in S_{\theta}$ tel que

$$
 \arg(w_0)=\theta \quad\text{ et }\quad   \arg(w_0 -u_0)\sim \theta 
\tag 11
$$

\noindent Soit $S_{0,\theta}=S_\theta\cap (S_{\arg(w_0-u_0)}+u_0)$. On joint $w\in S_{0,\theta}$ \`a $0$ par un chemin ${\Cal C}^1$-régulier voisin du chemin $\gamma_w=[0,u_0]\cup [u_0,w]$ qui satisfait \`a la condition (8) du lemme VA1 gr\^ace \`a (10) et (11). D'où le résultat.\qed\par

\vskip 3mm

\noindent {\bf Chemins exponentiels.} Pour tout $u_0\geq 1$ et $K\geq 1$, notons 

$$
 V_{u_0,K}=\{ w=u_0+u+iv \in P^+;\quad u\geq 0,\quad |v|\leq \exp(u/K)-1 \}
$$

\noindent On joint un \'el\'ement $w\in V_{u_0,K}$ \`a $0$ par le chemin

$$
 \gamma_w=[0,u_0]\cup \{ z=u_0+u+iC(\exp(u/K)-1);\quad u\in [0,\text{ Re}(w)-u_0] \}
$$

\noindent o\`u $C\in[-1,1]$ est une constante qui ne d\'epend que de $w$. Par un calcul simple, il existe $M(K)>0$ tel que pour tout $u_0\geq 1$ et sur tout chemin exponentiel

$$
 |\frac{z'(u)}{z(u)}|\leq M(K)  
\tag 12
$$

\vskip 3mm

\proclaim{Lemme VA3} Le germe de $f$ est un \'el\'ement de $QA^{1,q+2)}$.\endproclaim

\vskip 3mm

\noindent {\bf Preuve.} Montrons que $f$ satisfait à la condition de quasi-analycité. Remarquons d'abord que si $r$ est r\'eel, le domaine non-singulier de l'op\'erateur ${\Cal L}_{nr}$ est $P^+$. Les chemins $\gamma_w$ de $V_{1,1}$ satisfont clairement \`a la condition (8) du lemme VA1.\par

\vskip 2mm
 
\noindent Soit $\varphi_0$ le morphisme d'\'eclatement du point $(\infty,0)$ de $P^+\times {{\bold{\CC}}}$ dans la direction r\'eelle $0$ de ${{\bold{\RR}}}P^1$

$$
\varphi_0 :\quad (w,\widetilde \mu) \mapsto \mu=\frac{\widetilde \mu}{ w+1}
\tag 13
$$

\noindent et soit $\Phi_0$ l'application 
$$
\Phi_0 (w,y,\widetilde \mu,\nu)=(w,y,\varphi_0 (w,\widetilde \mu ),\nu). 
\tag 14
$$

\noindent Soit $K>1$. Si $\widetilde {{{\bold{W}}}}_1 \in ({{\bold{\CC}}},0)$ est suffisament petit, les projections sur $[0,1]\cup V_{1,K}$ des fibres $\varphi^{-1}_0(\mu)$ de $([0,1]\cup V_{1,K}) \times \widetilde {{{\bold{W}}}}_1$ contiennent les sous-ensembles 

$$
 \Delta_\mu(w)=\{ z\in [0,1]\cup V_{1,K};\quad\text{ Re}(z)\leq \text{ Re}(w),\quad |z|\leq |w| \} 
\tag 15
$$

\noindent qui sont inclus dans le domaine non-singulier $\Omega_r$ des op\'erateurs ${\Cal L}_{nr}$. En effet, si $z\in ]0,1]$, on a $\text{Re}(rz)=z\text{Re}(r)>0$ et si $z\in V_{1,K}$

$$
 \text{Re}(rz)=\text{Re}(z)+\text{Re}(\frac{\widetilde{\mu}}{w+1}z)>0
$$

\noindent Soit $u_0>1$. Les chemins $\gamma_w$ de $[0,u_0]\cup V_{u_0,K}$ sont inclus dans $\Delta_{\mu}(w)$ et satisfont \`a la condition (8) du lemme VA1 si $u_0$ est suffisament grand. En effet, elle est clairement satisfaite sur le segment r\'eel. Sur le chemin exponentiel

$$
 rz'=z'+\mu z \frac{z'}{z}
$$

\noindent et d'apr\`es (12)

$$
 |\frac{\text{Im}(rz')}{\text{Re}(rz')}|\leq \exp(u_0+u)
$$

\noindent Soit $\overline{D(0,\rho)}$ un disque de $\widetilde{{\bold{W}}}_1$. Soit $\widetilde{f}$, $\widetilde{F}$ et $\widetilde{f}_1$ les relev\'es par $\Phi_0$ de $f$, $F$ et $f_1$. On a donc

$$
 ||1+\widetilde{F}||_{V_{u_0,K}\times \overline{D(0,1/2)}\times \overline{D(0,\rho)}\times {\bold{W}}_q}\leq C_0
$$

\noindent et

$$
 ||\widetilde{f}_1(w,.)||_{\overline{D(0,\rho)}}\leq C_1|\exp(-w)|
$$

\noindent par cons\'equent

$$
 ||\widetilde{f}(w,.)||_{\overline{D(0,1/2)}\times \overline{D(0,\rho)}\times {\bold{W}}_q}\leq C_2|\exp(-w)|
$$

\noindent Par les formules intégrales de Cauchy sur $\overline{D(0,\rho)}$, les coefficients de la série $\widetilde{f}={\sum}_{k\geq 0} \widetilde{c}_k \widetilde{\mu}^k$ admettent les majorations

$$
 ||\widetilde{c}_k(w,.)||_{\overline{D(0,1/2)}\times {\bold{W}}_q}\leq C_k|\exp(-w)|
$$

 Soit maintenant $f={\sum}_{k\geq 0} c_k \mu^k$ la s\'erie de $f$. On a $c_k=(w+1)^k\widetilde{c}_k$ et donc chaque coefficient $c_k$ est born\'e sur le domaine $V_{u_0,k+1}$. On conclut en remarquant que le complémentaire de ces domaines dans n'importe quel domaine  de type puissance, est une partie relativement compact.\qed\par

\vskip 3mm

\proclaim{Lemme VA4} Le germe de $F$ est un \'el\'ement de $QR{\Cal H}^{1,(1,q+1)}$.\endproclaim

\noindent {\bf Preuve.} La d\'emarche est classique ([I-Y], [M], [Ro2] ) pour la partie concernant la structure asymptotique formelle de type Hilbert. Elle utilise les formes prénormales de $\omega_\alpha$. Les propriétés du reste découlent d'une deuxième application des opérateurs intégrals de Dulac.\par

\vskip 3mm

 On pr\'epare analytiquement $\omega_\alpha$ \`a un certain ordre $2N>1$ par un diff\'eomorphisme qui pr\'eserve la coordonn\'ee $x$

$$
a(x,y,\alpha)=\sum_{n\geq 1} a_n(x,\alpha)y^n
\quad\text{avec}\quad
a_n=\left \lbrace 
  \aligned x^n \widetilde{a}_n(\alpha) & \text{ pour}\quad n\leq 2N \\
    x^{2N+1}\widetilde{a}_n (x,\alpha) & \text{ pour}\quad n>2N \endaligned
\right.
$$

\noindent Soit $f_{1,N}$ l'int\'egrale premi\`ere de la forme prénormale

$$
\omega_{\alpha,N}=
xdy+y(1+\mu+{{\sum}}_{n=1}^{2N}a_n(x,\alpha)y^n)dx
\tag 16
$$

\noindent telle que $f_{1,N}(1,y,\alpha)\equiv y$. Il est connu ([I-Y] par exemple) qu'il existe $F_{1,N} \in QR{\Cal H}^{1,(1,q+1)}_{\text{cvg}}$ telle que

$$
 f_{1,N}=x^r y(1+F_{1,N}) \quad \text{ et }\quad F_{1,N}(0)=0
$$

\noindent En effet, dans le morphisme $\Phi_\alpha(x,y)=(x_1,y_1)=(xy,z(x,\mu)y)$, la 1-forme $\omega_{\alpha,N}$ se d\'esingularise en 

$$
 \widetilde{\omega}_{\alpha,N}=
    (\mu+{{\sum}}_{n=1}^N \widetilde{a}_n x_1^n )dy_1 +(1-
             y_1{{\sum}}_{n=1}^N\widetilde{a}_n x_1^{n-1})dx_1
$$

\noindent et cette 1-forme admet une int\'egrale premi\`ere analytique $g$ qui v\'erifie $g(x_1,0,\alpha )=x_1$. Soit $g_0=x_1+\mu y_1$ l'int\'egrale premi\`ere de la partie lin\'eaire. Par un calcul simple 

$$
 dg_0\wedge \widetilde{\omega}_{\alpha,N}=g_0 \sum_{n=1}^N \widetilde{a}_n x_1^{n-1} dx_1\wedge dy_1
$$

\noindent et ceci montre que $g$ est divisible par $g_0$.\par

\vskip 2mm

 Soit $\Psi_\alpha (x,y)=(X,Y)$ le changement de coordonn\'ees ramifi\'e

$$ 
 \left\lbrace
 \aligned
      X &=  x   \\
      Y &=  y(1+F_{1,N}) 
   \endaligned
 \right.
\tag 17
$$

\noindent la fonction $f_{1,N}$ \'etant une int\'egrale premi\`ere de $\omega_{\alpha,N}$, la 1-forme ramifi\'ee $\eta_\alpha =(\Psi_\alpha^{-1})^* \omega_\alpha$ s'\'ecrit

$$
 \eta_\alpha =XdY+Y(r+(XY)^{2N+1} b(X,Y,\alpha))dX 
$$

\noindent avec $b \in QR{\Cal H}^{1,(1,q+1)}_{\text{cvg}}$. Soit $f_{2,N}$ l'int\'egrale premi\`ere de $\eta_\alpha$ telle que $f_{2,N}(1,Y,\alpha)\equiv Y$. Soit $\Psi(x,y,\alpha)=(\Psi_\alpha(x,y),\alpha)$. L'int\'egrale premi\`ere $f$ de $\omega_\alpha$ est donn\'ee par 

$$
 f=f_{2,N} \circ \Psi
$$

\noindent  Or il est facile de voir que $f_{2,N}$ s'\'ecrit

$$
 f_{2,N}=X^rY(1+ X^N H_N) 
\tag 18
$$

\noindent On montre alors que $H_N=\sum_{n>2N} h_n Y^n$ est un \'el\'ement de $SB^{1,q+2}$ (et même de $QA^{1,q+2}$), en appliquant la m\'ethode des lemmes pr\'ec\'edents aux coefficients $h_n$ par l'interm\'ediaire des op\'erateurs ${\Cal L}_{(n-1)r-N}$ sur le domaine non-singulier $\Omega_{2r-1}$. Ceci implique en particulier que $F\in QA^{1,q+2}$ (pour un bon choix de $N$), et donc que $F\in QR{\Cal H}^{1,(1,q+1)}$\qed\par

\vskip 3mm

\noindent {\bf §3. Remarque VA1.}\par

\vskip 3mm

 On peut construire des chemins sur les feuilles de $\omega_\alpha$ qui ne quittent pas un certain voisinage de $0$, par exemple $D(0,1)\times D(0,1)$, et qui rencontrent une seule fois chacune des transversales $\{ x=1 \}$ et $\{ y=1\}$. En effet, soit $Y=y(1+F)$. D'apr\`es l'\'etude pr\'ec\'edente

$$
 c_1|y|\leq ||Y(.,y,.)||_{V\times W_1\times {\bold W}_q}\leq c_2|y|
\tag 19
$$

\noindent Soit $y_0<c_1/c_2$. L'intégration réelle au dessus du segment $[1,y_0]$ est presque une translation réelle dans la coordonnée $w$. Maintenant, en utilisant (19) et l'int\'egrale premi\`ere $x^r Y$ de la partie lin\'eaire de $\omega_\alpha$, on montre que l'int\'egration au dessus des chemins $\gamma_w$ ne quitte pas le voisinage $D(0,1)\times D(0,1)$.\par

\vskip 5mm

\noindent {\bf B. Les théorèmes de division.}\par

\vskip 3mm

 Ils sont bas\'es sur l'algorithme de division d'Hironaka [B-M]. Rappelons quelques r\'esultats de ce travail: soit $a\in {\bold{\RR}}\{\alpha \}$ avec $\alpha=(\alpha_1,\ldots,\alpha_q)$. On note $e_L(a)$ le plus petit indice $m\in {\bold{\NN}}^q$ de coefficient non nul dans la s\'erie

$$
 a=\sum_{m}a_{m} \alpha^{m}
$$

\noindent les \'el\'ements de ${\bold{\NN}}^q$ \'etant ordonn\'es par l'ordre lexicographique $(L(m),m_1,\ldots,m_q)$ o\`u $L(m)=\sum_j \lambda_j m_j$ est une forme lin\'eaire positive. Soit $J$ un id\'eal de ${\bold{\RR}}\{\alpha \}$ et $N(J)=\{ e_L(a);\quad a\in J \}$ son diagramme des exposants initiaux. Il existe une liste minimale $m^1,\ldots,m^l\in {\bold{\NN}}^q$ telle que

$$
 N(J)=\cup_{i=1}^l \{ m^i +{\bold{\NN}}^q \}
$$

\noindent On d\'efinit une partition de ${\bold{\NN}}^q$ par

$$
 \Delta_1=m^1 +{\bold{\NN}}^q
\quad \Delta_i=m^i +{\bold{\NN}}^q -\cup_{k=1}^{i-1} \Delta_k
\quad\text{et}\quad
 \Delta={\bold{\NN}}^q -\cup_{i=1}^l \Delta_i
$$

\noindent alors\par

\roster

\item"$(i)$" toute famille $a_1,\ldots,a_l\in J$ telle que $e_L(a_i)=m^i$ est une base de $J$\par

\item"$(ii)$" tout $f\in {\bold{\RR}}\{ \alpha \}$ se divise de mani\`ere unique dans $J$ sous la forme

\endroster

$$
 f=\sum_{i=1}^l Q_i a_i + R \qquad Q_i,\ R\in {\bold{\RR}}\{ \alpha \}
$$

$$
 m^i + \text{Supp}(Q_i)\subset \Delta_i \quad \text{et} \quad \text{Supp}(R)\subset \Delta
$$

\noindent La partie analytique de cette algorithme fournit des estimations pr\'ecises: soit $L$ une forme linéaire positive, $\sigma>0$ et

$$
 {\bold{\RR}}\{ \alpha \}_{L,\sigma}=\{ f=\sum_{m} f_{m} \alpha^{m};\quad ||f||_{L,\sigma}<\infty \} \quad\text{avec}\quad ||f||_{L,\sigma}=\sum |f_{m}| \sigma^{L(m)}
$$

\noindent alors il existe $L$ et $\varepsilon>0$ tels que si $f\in {\bold{\RR}}\{ \alpha \}_{L,\sigma}$ et $\sigma\leq \varepsilon$

$$
 ||Q_i(f)||_{L,\sigma}\leq \frac{2}{\sigma^{L(m^i)}} ||f||_{L,\sigma}
\quad\text{et}\quad ||R(f)||_{L,\sigma}\leq 2||f||_{L,\sigma}
$$

\vskip 3mm

\noindent Soit $\Omega\in {\Cal E}{\Cal I}$ un domaine quasi-analytique et $QA^{1,q}[\Omega]\subset QA^{1,q}$ l'algèbre des germes $f=\sum f_m(x)\alpha^m$ dont les coefficients $f_m$ admettent un prolongement holomorphe et born\'e sur $\Omega$. Pour l'action de $\chi_0=x\partial/\partial x$ sur ces algèbres, on a le

\vskip 3mm

\proclaim{Théorème VB1 (théorème de de division 1)} Soit ${\Cal B}=QA^{1,q}[\Omega]$ ou $SB^{1,q}$ et $J$ un id\'eal de ${\bold{\RR}}\{\alpha\}$. Soit $a_1,\ldots,a_l$ une base de $J$. Alors pour tout $f\in {\Cal B}$, il existe de mani\`ere unique $Q_i,\ R\in {\Cal B}$ tels que

$$
 f=\sum_i Q_i a_i + R
$$

$$
 m^i + \text{Supp}(Q_i(x,.))\subset \Delta_i\qquad
 \text{Supp}(R(x,.))\subset \Delta
$$

\endproclaim

\vskip 3mm

\noindent La partie formelle de cet algorithme implique le

\vskip 3mm

\proclaim{Lemme VB1} Soit $f=\sum_{m} f_{m} \alpha^{m}\in {\Cal B}$. Si pour tout $m$, on a $f_{m}=o(x^n)$ (dans l'anneau $SB^{1,q}$), alors il en est de m\^eme pour les s\'eries des $Q_i$ et de $R$.\endproclaim

\vskip 3mm

 Soit $s\in\NN^*$ et soit la dérivation 

$$
 \chi=x\frac{\partial}{\partial x}-s\sum_{j=1}^\ell u_j\frac{\partial}{\partial u_j}
$$

\noindent Quitte à effectuer une ramification en $x$, on peut supposer que $s=1$. La dérivation $\chi$ agit sur l'anneau $SB^{1,.}(x,\alpha,u)$. Soit $\pi_\chi:(x,\alpha,u)\in U\mapsto (\alpha,(\lambda_j))=(\alpha,(xu_j))\in W$ son morphisme intégral. Les idéaux $\chi$-transverses le long de $\gamma=\{(\alpha,u)=0\}$ sont des idéaux de l'anneau $\RR\{\alpha,\lambda\}$.

\proclaim{Théorème VB2 (théorème de division 2)} Soit $J$ un idéal de $\RR\{\alpha,\lambda\}$. Il existe un entier $n(J)$ tel que pour tout $f\in SB(x,\alpha,u)$ dont l'idéal $\chi$-transverse est inclus dans $J$, alors $(x^{n(J)})I_{\chi,f}\subset \pi_{\chi}^*(J)$.
\endproclaim

\noindent{\bf Preuve.} Soit $S_\theta$ un secteur dans la coordonnée $x$ et $P_\varepsilon$ un polydisque dans les coordonnées $(\alpha,u)$ tels que la série $f=\sum c_{n,m}(x)\alpha^n u^m$ soit convergente sur $S_\theta\times P_\varepsilon$. Soit $F=\sum x^{-|m|}c_{n,m}\alpha^n\lambda^m$. On a $\pi_{\chi}^*(F)=f$ et pour tout $x\in S_\theta$, la série de $F(x,.)$ est convergente sur le produit de polydisques $P_\varepsilon\times P_{\varepsilon |x|}$. De plus, l'idéal $\chi_0$-transverse de $F$ est inclus dans $J$. Soit $(\varphi_j)$ une base de l'idéal $J$ dans l'anneau $\RR\{\alpha,\lambda\}$. D'après le théorème de division VB1

$$
 F=\sum Q_j\varphi_j
\tag 1
$$

\noindent Ce théorème se généralise facilement aux produits de polydisques: soit $L(n,m)=L_1(n)+L_2(m)$ une forme linéaire positive et soit $\sigma=(\sigma_1,\sigma_2)$ avec $\sigma_i>0$. Notons

$$
 \RR\{\alpha,\lambda\}_{L,\sigma}=\{ g=\sum g_{n,m}\alpha^n\lambda^m,\ ||g||_{L,\sigma}=\sum |g_{n,m}|\sigma_1^{L_1(n)}\sigma_2^{L_2(m)}<\infty\}
$$

\noindent alors il existe $L$, des entiers $\ell_{1,j}$, $\ell_{2,j}$ et $\varepsilon_0>0$ tels que si $\sigma_i<\varepsilon_0$

$$
 ||Q_j(x,.)||_{L,\sigma}<\frac{2}{\sigma_1^{\ell_{1,j}}\sigma_2^{\ell_{2,j}}}||F(x,.)||_{L,\sigma}
$$

\noindent Soit $n(J)=\max\{\ell_{2,j}\}$, on a $||x^{n(J)}Q_j(x,.)||_{L,\sigma}<c_j||F(x,.)||_{L,\sigma}$, par conséquent si $\varepsilon$ est suffisament petit, on a $q_j=\pi_{\chi}^*(x^{n(J)}Q_j)\in SB(x,\alpha,u)$. On obtient le résultat en relevant la relation (1)

$$
 x^{n(J)}f=\sum q_j\pi_{\chi'}^*(\varphi_j)
\qed
$$

\vskip 3mm

 Soit $\nu=(\nu_1,\ldots,\nu_p)$, $\nu'=(\nu'_1,\ldots,\nu'_{p'})$, $\alpha=(\mu,\nu')$ et $\alpha'=(\alpha,\nu)$. L'algèbre $QR{\Cal H}(x,\alpha')$ s'identifie à une sous-algèbre de ${\Cal B}=QR{\Cal H}(x,\alpha)\{\nu\}$. Les séries d'élé\-ments de ${\Cal B}$ sont convergentes sur un voisinage produit. Soit $f=\sum f_m\nu^m\in{\Cal B}$, l'algèbre $QR{\Cal H}(x,\alpha)$ étant locale de topologie de Krull séparée, on peut définir, comme dans [B-M], un ordre sur les monômes de ${\Cal B}$: soit ${\Cal M}$ l'idéal maximal de $QR{\Cal H}(x,\alpha)$, l'ordre $e(f_m)$ de $f_m$ est le plus grand entier $e$ tel que $f_m\in{\Cal M}^{e}$. Soit $L(m)=\sum \lambda_j m_j$ une forme linéaire positive et soit $(L(m),m_1,\ldots,m_p)$ l'ordre lexicographique sur les monômes $\nu^m$. Pour tout entier $\ell$, l'application $L'(f_m\nu^m)=(\ell e(f_m),L(m),m)$ est un ordre sur les monômes de ${\Cal B}$.\par

\vskip 3mm

 Soit $I=\langle g_1,\ldots,g_q\rangle$ un idéal de ${\Cal B}$ et soit $g_j=\sum g_{j,m}\nu^m$ la série de $g_j$. On suppose que pour tout $j$, il existe un coefficient $g_{j,m}$ d'ordre 0. Soit $m^j$ le plus petit de ces entiers et $(\Delta_j,\Delta)$ la partition de $\NN^p$ associée. On définit le support de $f\in{\Cal B}$ par $\text{supp}(f)=\{m;\quad f_m\neq 0\}$. Les algorithmes formel et analytique de [B-M] s'adaptent à $I$ et à ${\Cal B}$ et on obtient

\proclaim{Lemme VB2} Il existe $L$ et $\ell$ tels que tout $f\in{\Cal B}$ se divise de manière unique dans $I$ sous la forme

$$
 f=\sum_{j=1}^q Q_j g_j +R\quad\quad Q_j,\ R\in{\Cal B}
$$

$$
m^j+\text{supp}(Q_j)\subset\Delta_j\quad\text{et}\quad\text{supp}(R)\subset\Delta
$$

\endproclaim

 L'algèbre $QR{\Cal H}(x,\alpha')$ est isomorphe à une sous-algèbre ${\Cal B}_*$ de ${\Cal B}$ définie comme suit: Soit $X(x,\mu)$ les fonctions élémentaires de $QR{\Cal H}$. Alors $f\in{\Cal B}_*$ si et seulement si pour tout $n$, la série $\tau_n(f)=\sum{\bold j}_X^n(f_m)\nu^m$ est un élément de $QR{\Cal H}_{cvg}$. Si $I$ est un idéal de ${\Cal B}_*$ satisfaisant aux mêmes hypothèses que ci-dessus, alors

\proclaim{Théorème VB3 (théorème de division 3)} Tout $f\in{\Cal B}_*$ se divise de manière unique dans $I$ sous la forme

$$
 f=\sum_{j=1}^q Q_j g_j +R\quad\quad Q_j,\ R\in{\Cal B}_*
$$

$$
m^j+\text{supp}(Q_j)\subset\Delta_j\quad\text{et}\quad\text{supp}(R)\subset\Delta
$$

\endproclaim

\noindent{\bf Preuve.} D'après la partie formel de l'algorithme (lemme VB1), si $f_m=o(x^n)$ pour tout $m$, il en est de même des séries des $Q_j$ et de $R$. Soit $Q_j(f)$ et $R(f)$ donnés par le lemme VB2. Par l'unicité de la division, $Q_j(f)=Q_j(\tau_n(f))+Q_j(f-\tau_n(f))$ et $R(f)=R(\tau_n(f))+R(f-\tau_n(f))$. Par conséquent $\tau_n(Q_{j}(f))=\tau_n(Q_{j}(\tau_n(f)))$ et $\tau_n(R(f))=\tau_n(R(\tau_n(f)))$. Soit l'idéal ${I}_n=\langle \tau_n(g_{1}),\ldots,\tau_n(g_{q})\rangle$, la partition de $\NN^p$ qui lui est associée est la même que celle de $I$. Soit $Q_{j,n}(\tau_n(f))$ et $R_n(\tau_n(f))$ donnés par le lemme VB2 appliqué à $I_n$, ce sont des éléments de $QR{\Cal H}_{cvg}$ d'après la partie analytique de l'algorithme. Or chaque coefficient de la série de $\sum Q_{j,n}(\tau_n(f))(g_j-\tau_n(g_{j}))$ est un $o(x^n)$, donc par l'unicité de la division, $\tau_n(Q_{j}(\tau_n(f)))={\bold j}_X^n(Q_{j,n}(\tau_n(f)))$ et $\tau_n(R(\tau_n(f)))={\bold j}_X^n(R_n(\tau_n(f)))$.\qed

\vskip 3mm

 D'après [B-M], une conséquence de ce théorème est le théorème des fonctions implicites pour les applications régulières dans la coordonnée $\nu$. Une autre conséquence est

\proclaim{Théorème VB4 (théorème d'inversion)} Soit $f=x(1+O(x))\in QR{\Cal H}(x,\alpha)$. Il existe un unique $g=y(1+O(y))\in QR{\Cal H}(y,\alpha)$ qui inverse $f$.
\endproclaim

\noindent{\bf Preuve.} Posons $f=x(1+F)$ avec $F\in QR{\Cal H}$. Il existe un unique inverse en classe ${\Cal C}^1$. Cherchons $G\in QR{\Cal H}$ tel que $g=y(1+G)$ soit un inverse. La condition $f\circ g=Id$ donne l'équation $G+(1+G)F(y(1+G),\alpha)=0$ qui est analytique en $G$ et régulière d'ordre 1. Le théorème des fonctions implicites permet de conclure.\qed

\vskip 5mm

\noindent {\bf R\'ef\'erences.}\par

\vskip 3mm

\roster

\item"[B-M]" Bierstone E., Milman P.D., {\sl The local geometry of analytic mappings.} Dottorato di Ricerca in Matematica. [Doctorate in Mathematical Research] ETS Editrice, Pisa, (1988), iv+77pp. ISBN: 88-7741-419.\par

\vskip 2mm

\item"[D]" Dulac H., {\sl Sur les cycles limites.} Bull. Soc. Math. France {\bf 51} (1923), pp. 45-188.\par

\vskip 2mm

\item"[Dr]" Dries L.V.D., {\sl Remarks on Tarski's problem concerning } $(\RR,+,.,exp)$, Logic Colloquium '82, G. Lolli, G. Longo and A. Marcja, eds., North-Holland (1984), pp. 97-121.\par

\vskip 2mm

\item"[D-M-M]" Dries L.V.D., Macintyre A. and Marker D., {\sl The elementary theory of restricted analytic fields with exponentiation}, Annals of Maths., {\bf 140} (1994), pp. 183-205.\par

\vskip 2mm

\item"[D-S]" Dries L.V.D. and P. Speissegger, {\sl The real field with convergent generalized power series}, Trans. Amer. Math. Soc., {\bf 350} (1998), pp.4377-4421.\par

\vskip 2mm

\item"[E]" Ecalle J., {\sl Introduction aux fonctions analysables et preuve constructive de la conjecture de Dulac.} Actualit\'es Math\'ematiques, Hermann (1992).\par

\vskip 2mm

\item"[E-M-R]" El Morsalani M., Mourtada A., Roussarie R., {\sl Quasi-regularity property for unfolding of hyperbolic polycycles.} Ast\'erisque {\bf 220} (1994), pp. 303-326.\par

\vskip 2mm

\item"[H]" Herv\'e M., {\sl Several complex variables, local theory.} Tata Institute of Fundamental Research, Bombay, Oxford University Press (1963).\par

\vskip 2mm

\item"[Hi]" Hilbert D., {\sl Mathematische probleme (lecture), the second international congress of mathematicians, Paris 1900.} Nachr. Ges. Wiss. Gottingen Math.-Phys. Kl. (1900), p. 253-297; {\sl Mathematical developmentsarising from Hilbert's problems.} Proceedings of Symposium in Pure Mathematics, A.M.S., F. Browder editor, {\bf 28} (1976), pp. 50-51.\par

\vskip 2mm

\item"[Hir1]" Hironaka H., {\sl Resolution of th singularities of an algebraic variety over a field of characteristic zero, I, II.} Ann. of Math., {\bf 79} (1964), pp. 109-326.\par

\vskip 2mm

\item"[Hir2]" Hironaka H., {\sl Introduction to real-analytic sets and real-analytic maps.} Instituto Matematico "L. Tonelli", Pisa (1973).\par

\vskip 2mm

\item"[I]" Il'yashenko Y. S., {\sl Limits cycles of polynomial vector fields with non degenerate singular points on the real plane.} Func. Ana. and Appl., {\bf 18.3} (1985), pp. 199-209.\par

\vskip 2mm

\item"[I-Y1]" Il'yashenko Y. S., Yakovenko S., {\sl Finitely smooth normal forms for local diffeomorphisms and vector fields.} Rus. Math. Surveys, {\bf 46} (1991), pp. 1-43.\par

\vskip 2mm

\item"[I-Y2]" Il'yashenko Y. S., Yakovenko S., {\sl Hilbert-Arnold problem for elementary polycycles. In Around Hilbert's 16th problem,} Advances in Sov. Math. (1994).\par                       

\vskip 2mm

\item"[Ka]" Kaloshin V., {\sl The Hilbert 16th problem and an estimate for cylicity of an elementary polycycle}, Preprint.\par

\vskip 2mm

\item"[K1]" Khovanskii A., {\sl Fewnomials.} A.M.S., Providence, RI (1991).\par

\vskip 2mm

\item"[K2]" Khovanskii A., {\sl Real analytic varieties with the finiteness property and complex abelian integrals.} Funct. Anal. and Appl., {\bf 18} (1984), pp. 199-207.\par

\vskip 2mm

\item"[L]" Lojasiewicz S., {\sl Introduction to Complex Analytic Geometry}, Birkhauser (1991).\par

\vskip 2mm

\item"[L-R]" Lion J.M. et Rolin J.P., {\sl Volumes, feuilles de Rolle de feuilletages analytiques et théorème de Wilkie}, Ann. de Toulouse {\bf 7} (1998), pp. 93-112.\par

\vskip 2mm

\item"[L-S]" Lion J.M. and Speissegger P., {\sl Analytic stratification in the Pfaffian closure of an o-minimal structure}, Duke Math. J. {\bf 103}, no 2 (2000), pp. 215-231.\par

\vskip 2mm

\item"[M]" Mourtada A., {\sl Cyclicit\'e finie des polycycles hyperboliques de champs de vecteurs du plan; algorithme de finitude. Ann. Inst. Four. tome 41, fasc. {\bf 3} (1991), pp. 719-753.\par

\vskip 2mm

\item"[M1]" Mourtada A., {\sl Extension aux cycles singuliers du théorème de Khovanski-Varchenko}, Préprint Janvier 2007.\par

\vskip 2mm

\item"[M2]" Mourtada A., {\sl Application de Dulac en dimension quelconque: quasi-analy-\par\noindent cité, monodromie et synthèse}, En cours de rédaction.\par

\vskip 2mm

\item"[M-M]" Mourtada A. et Moussu R., {\sl Applications de Dulac et applications pfaffiennes.} Bull. Soc. Math. France {\bf 125} (1997), pp. 1-13.\par

\vskip 2mm

\item"[Ro1]" Roussarie R., {\sl A note on finite cyclicity and Hilbert's 16th problem.} Springer Lecture Notes in Mathematics {\bf 1331} (1988), pp. 161-168.\par

\vskip 2mm

\item"[Ro2]" Roussarie R., {\sl On the number of limit cycles which appear by perturbation of separatrix loop of planar vector fields.} Bol. Soc. Bras. Mat. {\bf 17} (1986), pp. 67-101.\par

\vskip 2mm

\item"[Ru]" Rudin W., {\sl Real and complex analysis.} McGraw-Hill, (1966).\par

\vskip 2mm

\item"[S]" Speissegger P., {\sl The pfaffian closure of an o-minimal structure}, J. reine angew. Math. {\bf 508} (1999), pp. 189-211.\par

\vskip 2mm

\item"[T]" Tougeron J.-C.,{\sl Alg\`ebres analytiques topologiquement n\oe      th\'eriennes et th\'e\-orie de Khovanski}. Ann. Inst. Four., tome 41, fasc. {\bf 4} (1991), pp. 823-840.\par

\vskip 2mm

\item"[W]" Wilkie  A., {\sl Model completeness results for expansions of the ordered field of real numbers by restricted Pfaffian functions and the exponential function}, J. Amer. Math. Soc. {\bf 9} (1996), pp. 1051-1094.\par

\endroster

\enddocument